Fosco Loregian

# $t$-structures on stable $(\infty, 1)$-categories

a thesis in Mathematics
with 105 diagrams
coded in koDi

# Contents

В родстве со всем, что есть, уверясь
И знаясь с будущим в быту,
Нельзя не впасть к концу, как в ересь,
В неслыханную простоту.

Но мы пощажены не будем,
Когда ее не утаим.
Она всего нужнее людям,
Но сложное понятней им...

*Б. Л. Пастернак*

四拳

波羅蜜

大光明.

# Introduction

The present work re-enacts the classical theory of *t-structures* reducing the classical definition given in [BBD82, KS] to a rather primitive categorical gadget: suitable *reflective factorization systems* (Def. **2.3.1**, **2.3.9**), which we call *normal torsion theories* following [CHK85, RT07]. A relation between these two objects has previously been noticed by other authors [RT07, HPS97, BR07] on the level of homotopy categories. The main achievement of the present thesis is to observe and prove that this relation exists genuinely when the definition is lifted to the higher-dimensional world where the notion of triangulated category comes from, i.e. *stable* $(\infty, 1)$*-categories.*

Stable $(\infty, 1)$-categories provide a far more natural setting to interpret the language of homological algebra: the main conceptual aim of the present work is to give explicit examples of this meta-principle.

To achieve this result, it seemed unavoidable to adopt a preferential model for $(\infty, 1)$-category theory: instead of working in a 'model-free' setting, we choose the ubiquitous dialect of Lurie's *stable quasicategories*; discussing to which extent (if any) the results we prove are affected by this choice, and establishing a meaningful dictionary between the validity of the general statement **3.1.1** in various different flavours on $\infty$-category theory occupies sections **A.5** and **3.3**; despite the fact that this is one of the most important issues from a categorical point of view, a rapid convergence of the present thesis into its final form has to be ensured; hence, we will defer a torough examination of the topic of model (in)dependence to subsequent works.

The first part of the thesis (Ch. **1**–**3**) builds (or rather, 'reinterprets') the calculus of factorization in the setting of $\infty$-categories. The desire to link this calculus with homological algebra and higher algebra deserves further explanation.

The language of factorization systems proved to be ubiquitous inside and outside category theory (among various different applications now established in the mathematical practice, the 'modern view' in algebraic topology revolves around the notion of orthogonality and lifting/extension problem, as it is said in the first pages of [**?**]. The modern 'synthetic' approach to homotopy theory inescapably relies on the notion of a (weak) factorization system ([Qui67, DS95, Rie11]).

In light of this, finding 'concrete' means of application for the calculus of factorization should be a natural step towards a popularization of this pervasive and deep language. And among all the various fields of application, homological algebra, a notable kind of 'abelian' homotopy theory, should be the most natural test bench to measure the validity of this effort. Despite the intrinsic simplicity, almost a triviality, of Thm. **3.1.1**, and despite the fact that the author feels he had failed at such an ambitious task, the pages you're about to read should be interpreted in this spirit.

## Structure of the thesis

The thesis is the results of a re-organization and methodical arrangement of the papers [FL16b, FL15a, FL15b, FL16a] (all written having my advisor as co-author) that have appeared on the arXiv since August 2014; the content is essentially unchanged; some sections and subsections (like e.g. **1.5**, **1.5.2.1**, **3.2**, a renewed proof of **4.3.20**, and Ch. **6**) do not appear anywhere at the moment of writing[(1)], but contain little new material and serve as linking sections making the discussion more complete and streamlined, developing certain natural derivations of the basic theory which would have easily exceeded the average length of a research paper.

Figure (**1**) below depicts the dependencies among the various chapters: a dashed line indicates a feeble logical dependence, whereas a thick line indicates a stronger one, unavoidable at first reading.

The first three chapters outline the main result of the present work, summarized as follows:

> For each stable $\infty$-category $\mathbf{C}$ there is a bijective correspondence between $t$-structures on the triangulated homotopy category $\mathrm{Ho}(\mathbf{C})$ and suitable orthogonal factorization systems on $\mathbf{C}$ called *normal torsion theories.*

This constitutes the backbone and the basic environment in which every subsequent application (the theory of *recollements* in stable $\infty$-categories in Ch. **5**, and Bridgeland's theory of *stability conditions* in Ch. **7**) takes place.

The main original contribution given in the present work is the 'Rosetta stone' theorem proving the quoted remark above; this is the main result of [FL16b], the only preprint that, at the moment of writing, has also been published by a peer-reviewed journal.

There are several minor results following from the 'Rosetta stone', like the fact that constructions one can perform on normal torsion theories are (at least to the categorically-minded) more natural and canonical than the corresponding construction in homological algebra, done on bare $t$-structures.

[(1)] May 24, 2020

## A word on model dependency

Ideally speaking, if there is an equivalence between two models for $\infty$-categories (say, *red* and *blue* $\infty$-categories), these two models both possess a notion of factorization systems and a calculus[(2)] thereof; moreover, these two notions of factorization system correspond to each other under the equivalence of models. Turning this principle of equivalence and correspondence into a genuine theorem is often a subtle matter (apart from being inherently difficult and a delicate issue, this is perhaps due to the fact that the author is ignorant of how to retrieve such a result in the existing literature): it is however possible to recognize at least three different settings having each its own 'calculus of factorization':

- stable model categories, where one can speak about *homotopy factorization systems* following [Bou77, Joy08]; this leads to the definition of a *homotopy $t$-structures* on stable model categories as suitable analogues of normal torsion theories in the set $\text{HFS}(\mathbf{M})$ of homotopy factorization systems on a model category $\mathbf{M}$.
- DG-categories, where we speak about enriched (over $\mathsf{Ch}(\mathbf{k})$) factorization systems (see [DK74]); this leads to the definition of DG-*$t$-structures* as enriched analogues of normal torsion theories in the set of $\text{DG-FS}(\mathbf{D})$ of enriched factorization systems on a DG-category $\mathbf{D}$.
- derivators, where we can define *$t$-derivators* via a (genuinely new) notion of factorization system on a derivator, and recognize the analogue of normal torsion theory in this setting.

At the moment of writing, all these points are being studied, and will hopefully appear as separate results in the near future.

## A word on the state of the art

Drawing equally from homological algebra, algebraic geometry, topology and category theory, the present work has not a single, well-defined flavour. Several sources of inspirations came from classical literature in algebraic topology [HPS97, Tie69, Hel68]; several others belong to the classical and less classical literature on algebraic geometry [Ver96, Bri07, Bri09, BO95]; others belong to pure category theory [RT07, CHK85, JM09, KT93, LW, Zan04], and others (see below) do not even belong to what is canonically recognized as mathematical literature.

The approach to the theory of $\infty$-categories taken here will certainly appear rather unorthodox to some readers: [Lur09, Lur17] have taught the author more about 1-categories than he did about $\infty$-categories. This, again, must be attributed to the ignorance of the author, which is more comfortable with the language of categories rather than with homotopy theory.

[(2)]By a 'calculus' of factorization systems we naïvely mean an analogue of the major results expressed in Ch. **1**, translated from the red to the blue model.

## Notation and Conventions

Categories (in the broad sense of 'categories and $\infty$-categories') are denoted as boldface letters $\mathbf{C}, \mathbf{D}$ and suchlike, opposed to generic, variable simplicial sets which are denoted by capital Latin letters (this creates an extremely rare, harmless conflict with the same notation adopted for objects in a category: the context always allows us to avoid confusion); functors between categories are always denoted as capital Latin letters in a sufficiently large neighbourhood[(3)] of $F, G, H, K$ and suchlike; the category of functors $\mathbf{C} \to \mathbf{D}$ is denoted as $\mathrm{Fun}(\mathbf{C}, \mathbf{D})$, $\mathbf{D}^{\mathbf{C}}$, $[\mathbf{C}, \mathbf{D}]$ (or, at the risk of being pedantic, as $(\mathbf{Q})\mathbf{Cat}(\mathbf{C}, \mathbf{D})$); morphisms in $\mathrm{Fun}(\mathbf{C}, \mathbf{D})$ (i.e. natural transformations between functors) are often written in the Greek alphabet; the simplex category $\Delta$ is the *topologist's delta*, having objects *nonempty* finite ordinals $\Delta[n] := \{0 < 1 \cdots < n\}$ regarded as categories in the obvious way; we adopt [Lur09] as a main reference for $\infty$-category theory, even if we can't help but confess that we profited from every single opportunity to deviate from the aesthetic of that book; in particular, we accept the (alas!) settled abuse to treat 'quasicategory' and '$\infty$-category' as synonyms; any other unexplained choice of notation belongs to folklore, or leans on common sense.

A general working principle of stable $\infty$-category theory is that homological algebra becomes easier and better motivated when looked at from a higher perspective[(4)]. To refer to this more natural environment we will often call *the stable setting* any theory of stable $\infty$-categories.

A not completely standard choice of notation is the following: each time a concept notion appears together with its dual, we write co/notion to denote that we refer to notion and conotion at the same time. So, if we write '$\mathbf{C}$ is a co/complete category' we mean that $\mathbf{C}$ is *both* complete and cocomplete, and if we write 'co/limit' we are speaking about limits *and* colimits at the same time.

## About the *kamon* on the titlepage

The titlepage contains the *kamon* of the Tachibana branch of Yoneda (!) family, traditionally drawn [TM02] as a tea-berry (a $t$-berry!) inside a circle ('丸茶の実', *maru Cha no Mi*):

[(3)]The set $A$ of letters of the English alphabet admits an obvious monotone bijection $A \xrightarrow{\varphi} \Delta[26]$; define a distance on $A$ by putting $d(-,=) \stackrel{\triangle}{=} |\varphi(-) - \varphi(=)|$.

[(4)]This rather operative and meta-linguistic principle is sketched in our Appendix **A**, where a complete proof of how triangulated category axioms follow from the 'pullout axiom' **A.2.3** is worked out in full detail.

## A word on the way I drew diagrams

Basically every existing package to draw commutative diagrams sucks. Starting from this undeniable truth, I spurred P. B. (see the acknowledgements) to write a `tikzlibrary` capable of producing beautiful and readable diagrams on both the coders' and the readers' side. The result is repo-ed here under the name `koDi`.

`koDi` acts via three different kinds of command: a `\lattice` environment, describing where to put the objects of the commutative diagram: each object of a `\lattice` is included in a `\obj #;` environment, and a command `\mor`, which produces a chain of morphisms of variable length, all linked by arrows `->` having different styles (basically those of `TiKZ`).

Each `\obj #;` environment allows the user to label the node with a tag which can be internally referred to inside a `\mor` environment: so, for example, an intelligent way to rename the node $\gamma(\widehat{X}^{\mathsf{s}}, \lambda_0)$ is `\obj (gX-l0):{\gamma(\widehat{X}^\textsf{s},\lambda_0)};` whereas an arrow $\gamma(\widehat{X}^{\mathsf{s}}, \lambda_0) \to Y$ can be written `\mor gX-l0 -> Y;`.

Since an example is worth a thousand words, here is the code producing diagram (**3.6**).

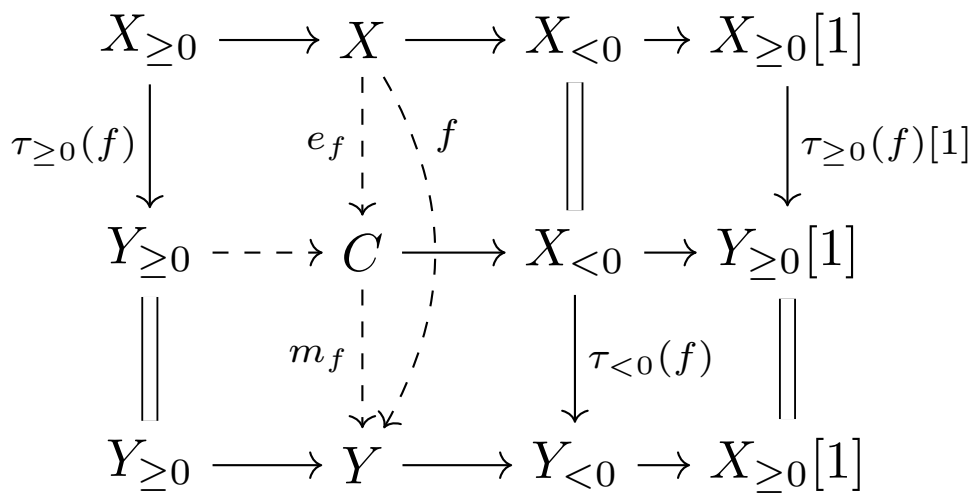

```
\begin{kD}
\lattice[mesh]{
        \obj  (Xge):X_{\ge 0}; & \obj X; & \obj  (Xle):X_{<0};
                      & \obj  (Xge+):X_{\ge 0}[1]; \\
        \obj  (Yge):Y_{\ge 0}; & \obj C; & \obj (Xle'):X_{<0};
                      & \obj  (Yge+):Y_{\ge 0}[1]; \\
        \obj (Yge'):Y_{\ge 0}; & \obj Y; & \obj  (Yle):Y_{<0};
                      & \obj (Yge+'):X_{\ge 0}[1]; \\};
\mor Xge -> X -> Xle -> Xge+
                {\tau_{\ge 0}(f)[1]}:-> Yge+ 2- Yge+';
\mor Xge swap:{\tau_{\ge 0}(f)}:-> Yge
                dashed,-> C {crossing over},-> Xle'
                {\tau_{<0}(f)}:-> Yle -> Yge+';
\mor Yge 2- Yge' -> Y -> Yle; \mor Xle 2- Xle' -> Yge+;
\mor[swap] X e_f:dashed,-> C m_f:dashed,-> Y;
\mor[dashed,near start] X f:r> Y;
\end{kD}
```

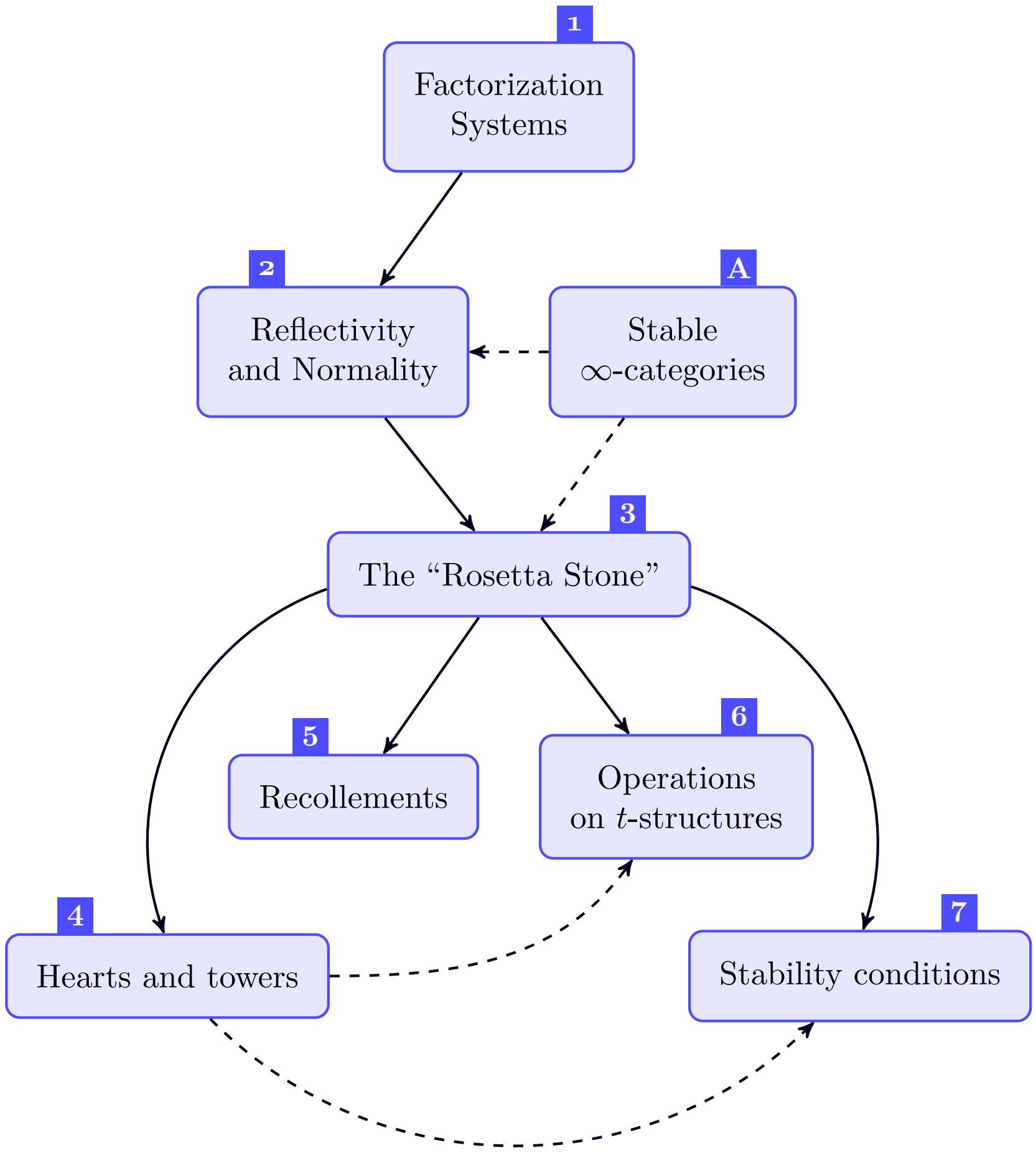

Figure 1: Dependencies between chapters

# Chapter 1

# Factorization Systems

## 1.1 Overview of factorization systems.

It isn't that they can't see the solution.
It is that they can't see the problem.

G.K. Chesterton

This first chapter deals with the general definition of a *factorization system* in the setting of $\infty$-categories, as given by [Lur09] and [Joy08, p. **178**].

We do not claim originality here, aiming only at a balance between the creation of a flexible and natural formalism, to be used along the subsequent chapters, and the necessity of rigor and generality.

The current literature seems to be too poor and too rich at the same time, when dealing with factorization systems; several authors often decide to rebuild the basic theory from scratch when they prove a new result, as there are several slightly different flavours in which one wants to interpret the basic idea behind the definition (that is, *factor every arrow of a category into two distinguished pieces*).

Since simple and extremely pervasive structures are often discovered independently, an accurate overview of the topic is somewhat impossible; we can however try to date back to the pioneering [ML48], published in 1948 (!) a forerunner to the modern notion of factorization system (see in particular axioms BC1–5).

With the passing of time, it became clear that together with unique factorization, the *orthogonality* relation between the arrows of two classes $\mathcal{E}, \mathcal{M} \subseteq \hom(\mathbf{C})$ played an essential rôle in the definition of a factorization system.[(1)] Another forerunner of the modern theory is Isbell's [Isb64]

(1) However, the rôle of uniqueness is much more essential (see Remark **1.4.8**): even if factorization of arrows with respect to a prefactorization is unique, a strictly unique factorization (there, the author doesn't mention orthogonality, but clearly refers to what in the subsequent [FK72] will be called in this way); his work was first popularized in the lucid and methodical presentation of the latter paper, which together with [CHK85] has been a fundamental starting point, and a source of suggestions for the main result exploited throughout this work: *reflective subcategories of a category* $\mathbf{C}$ *originate from the calculus of factorization systems*. This will be the main theme developed in the following chapters, and will culminate in Ch. **3** with the proof of our "Rosetta stone" theorem.

## 1.2 Markings and prefactorizations.

.

Matthew 25:33

DEFINITION 1.2.1. (MARKED SIMPLICIAL SET): Recall that a *marked simplicial set* $\underline{X}$ consists of a pair $(X, \mathcal{S})$, where $X$ is a simplicial set, and $\mathcal{S} \subseteq X_1$ is a class of distinguished 1-simplices on $X$, which contains every degenerate 1-simplex.

REMARK 1.2.2. The class of all marked simplicial sets forms a category $\mathbf{sSet}^{\varsigma}$, where a morphism is a simplicial map $f : (X, \mathcal{S}_X) \to (Y, \mathcal{S}_Y)$ which *respects* the markings, in the sense that $f(\mathcal{S}_X) \subseteq \mathcal{S}_Y$; the obvious forgetful functor

$$U : \mathbf{sSet}^{\varsigma} \to \mathbf{sSet} \tag{1.1}$$

admits both a right adjoint $X \mapsto X^\sharp = (X, X_1)$ and a left adjoint $X \mapsto X^\flat = (X, s_0(X_0))$, given by choosing the maximal and minimal markings, respectively (mnemonic trick: **r**ight adjoint is sha**r**p, **l**eft adjoint is f**l**at).

NOTATION 1.2.3. A *marked* $\infty$*-category* simply consists of a marked simplicial set which, in addition, is an $\infty$-category. From now on, we will consider only marked $\infty$-categories, and occasionally confuse $X$ and $X^\flat$.

DEFINITION 1.2.4. (ORTHOGONALITY): Let $f, g$ be two edges in an $\infty$-category $\mathbf{C}$. We will say that $f$ is *left orthogonal* to $g$ (or equivalently that $g$ is *right orthogonal* to $f$) if in any commutative square $\Delta[1] \times \Delta[1] \to \mathbf{C}$ like the following,

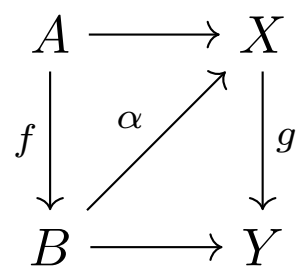

(1.2)

with respect to two classes *implies* orthogonality between the classes; from time to time we will need to exploit this useful remark, first observed by Joyal and used in [Joy] as definition of a factorization system. In any case, assuming the orthogonality relation and the factorization property as primitive and unrelated properties is a common practice.

the space of liftings $a$ rendering the two triangles (homotopy) commutative is contractible.[(2)]

Remark 1.2.5. This is Def. [Lur09, **5.2.8.1**]; compare also the older [JM09, Def. **3.1**].

Remark 1.2.6. "Being orthogonal" defines a binary relation on the set of edges in a marked $\infty$-category $\mathbf{C}$, denoted $f \perp g$.

Notation 1.2.7. We denote ${}^{\perp}(-) \dashv (-)^{\perp}$ the (antitone) Galois connection induced by the relation $\perp$ on subsets $\mathcal{S} \subseteq \mathbf{C}_1$;[(3)] more explicitly, we denote

$$\begin{aligned} \mathcal{S}^{\perp} &= \{f\colon \Delta[1] \to \mathbf{C} \mid s \perp f, \forall s \in \mathcal{S}\} \\ {}^{\perp}\mathcal{S} &= \{f\colon \Delta[1] \to \mathbf{C} \mid f \perp s, \forall s \in \mathcal{S}\}, \end{aligned}$$

and we consider the adjunction ${}^{\perp}(-)\colon \mathcal{P}(\hom(\mathbf{C}) \leftrightarrows \mathcal{P}(\hom(\mathbf{C}))\colon (-)^{\perp}$.

Definition 1.2.8. (Category of markings): If $\mathbf{C}$ is a small $\infty$-category we can define a poset $\mathrm{Mrk}(\mathbf{C})$ whose objects are markings of $\mathbf{C}$ and whose arrows are given by inclusions as subsets of $\mathbf{C}_1$.

Remark 1.2.9. The maximal and the minimal markings are, respectively, the terminal and initial object of $\mathrm{Mrk}(\mathbf{C})$; this category can also be characterized as the fiber over $\mathbf{C}$ of the forgetful functor $U\colon \mathbf{sSet}^{\varsigma} \to \mathbf{sSet}$. Moreover, the Galois connection ${}^{\perp}(-) \dashv (-)^{\perp}$ defined above induces an analogous adjunction on $\mathrm{Mrk}(\mathbf{C})$, via the obvious identification.

This remark leads to a second

Remark 1.2.10. The correspondence ${}^{\perp}(-) \dashv (-)^{\perp}$ forms a Galois connection in the category of markings of $X$; the maximal marking, and the marking Eqv made by all isomorphisms in $\mathbf{C}$ exchange each other under these correspondences. More precisely,

Proposition 1.2.11. The following conditions are equivalent, for $f\colon \Delta[1] \to \mathbf{C}$:

(1) $f$ is an isomorphism;
(2) $f \in \mathbf{C}_1^{\perp}$;
(3) $f \in {}^{\perp}\mathbf{C}_1$;

---

(2) By requiring that the space of liftings $\alpha$ is only *nonempty* one obtains the notion of weak orthogonality. In the following discussion we will only cope with the stronger request, but we rapidly address the issue in Def. **1.3.12** and in the subsequent points.

(3) Recall that if $R \subseteq A \times B$ is a relation, it induces a Galois connection

$${}^{R}(-)\colon \mathcal{P}(A) \rightleftarrows \mathcal{P}(B)\colon (-)^{R};$$

"negative thinking" tells us that this is simply the nerve-realization adjunction generated by $R$ regarded as a $\Omega$-profunctor ($\Omega$ can, but must not, be the Boolean algebra $\{0, 1\}$). This remark has, however, little importance in the ongoing discussion, and only serves the purpose of using the word "profunctor" in the present thesis.

(4) $f \perp f$.

REMARK 1.2.12. The technique applied here (devise suitable lifting problems which, solved, prove the claim) is a standard trick in the calculus of factorization systems: we will often use arguments like the following.

*Proof.* This case is extremely simple and paradigmatic. It is evident that the implications $1 \Rightarrow 2 \Rightarrow 4$ and $1 \Rightarrow 3 \Rightarrow 4$ (the inverse of $f$, composed with the upper horizontal arrow of a lifting problem, does the trick); to close the circle of implications, it is enough to show that $4 \Rightarrow 1$: this is evident, since the solution to the lifting problem

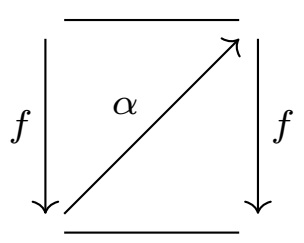

(1.3)

(where horizontal arrows are identity maps) must be the inverse of $f$ (in $\infty$-categories, there is a contractible space of such inverses, agreeing with Def. **1.2.4**). □

PROPOSITION 1.2.13. There exists an adjunction

$$\underline{{}^{\boxslash}(-)} \dashv \underline{(-)^{\boxslash}} \colon \mathbf{QCat}/X^{\Delta[1]} \leftrightarrows \left(\mathbf{QCat}/X^{\Delta[1]}\right)^{\mathrm{op}} \tag{1.4}$$

"lifting" the Galois connection ${}^{\boxslash}(-) \dashv (-)^{\boxslash} \colon \mathrm{Mrk}(X) \leftrightarrows \mathrm{Mrk}(X)$.

This can be seen as an $\infty$-categorical version of [Gar09, Prop. **3.8**].

REMARK 1.2.14. (ORTHOGONALITY AND LOCALITY): There is another notion of orthogonality of an object $X$ with respect to a morphism $f \in \hom(\mathbf{C})$; given these data, we say that $X$ is *right-orthogonal* to $f$ (or that $X$ is an *$f$-local* object) if the hom functor $\hom(-, X)$ inverts $f$.

If $\mathbf{C}$ has a terminal object $*$, this notion is related to Def. **1.2.4**, in the sense that $X$ is right-orthogonal to $f$ if and only if the terminal arrow $X \to *$ is right orthogonal to $f$. For this reason, we always refer to *object-orthogonality* as orthogonality with respect to terminal arrows. (Obviously, there is a dual notion of left-object-orthogonality between $f$ and $B \in \mathbf{C}$, which in the presence of an initial object reduces to left orthogonality with respect to $\varnothing \to B$).

NOTATION 1.2.15. Extending this notation a little bit more, we can speak about orthogonality between two objects, without introducing new definitions: in a category with both a terminal and initial object (which, since our blanket assumption in all the remaining chapters is to work in a stable $\infty$-category, will always be the case) we can say that

- Two objects $B$ and $X$ are orthogonal if $\hom(B, X)$ is contractible; we denote this (non-symmetric) relation as $B \perp X$.

- Two classes of object $\mathcal{H}$ and $\mathcal{K}$ in $\mathbf{C}$ are orthogonal if each object $H \in \mathcal{H}$ is orthogonal to each object $K \in \mathcal{K}$; we denote this situation by $\mathcal{H} \perp \mathcal{K}$.

This notation will greatly help us in Ch. **3** and **4**

The following nomenclature is modeled on the analogous categorical notion of a *prefactorization system* introduced in [FK72].

Definition 1.2.16. A pair of markings $(\mathcal{E}, \mathcal{M})$ in an $\infty$-category $\mathbf{C}$ is said to be a *($\infty$-categorical) prefactorization* when $\mathcal{E} = {}^{\perp}\mathcal{M}$ and $\mathcal{M} = \mathcal{E}^{\perp}$. In the following we will denote a prefactorization on $\mathbf{C}$ as $\mathbb{F} = (\mathcal{E}, \mathcal{M})$.

Remark 1.2.17. The collection of all prefactorizations on a given $\infty$-category $\mathbf{C}$ forms a poset, which we will call $\textsc{pf}(\mathbf{C})$, with respect to the order $\mathbb{F} = (\mathcal{E}, \mathcal{M}) \preceq \mathbb{F}' = (\mathcal{E}', \mathcal{M}')$ iff $\mathcal{M} \subset \mathcal{M}'$ (or equivalently, $\mathcal{E}' \subset \mathcal{E}$).

Remark 1.2.18. It is evident (as an easy consequence of adjunction identities) that any marking $\mathcal{S} \in \mathrm{Mrk}(\mathbf{C})$ induces two *canonical* prefactorizations on $\mathbf{C}$, obtained by sending $\mathcal{S}$ to $({}^{\perp}\mathcal{S}, ({}^{\perp}\mathcal{S})^{\perp})$ and $({}^{\perp}(\mathcal{S}^{\perp}), \mathcal{S}^{\perp})$. These two prefactorizations are denoted $\mathbb{S}_{\perp}$ e ${}_{\perp}\mathbb{S}$, respectively, and termed the *right* and *left prefactorization* associated to $\mathcal{S}$.

Definition 1.2.19. If a prefactorization $\mathbb{F}$ on $\mathbf{C}$ is such that there exists a marking $\mathcal{S} \in \mathrm{Mrk}(\mathbf{C})$ such that $\mathbb{F} = \mathbb{S}_{\perp}$ (resp., $\mathbb{F} = {}_{\perp}\mathbb{S}$) then $\mathbb{F}$ is said to be *right* (resp., *left*) *generated* by $\mathcal{S}$.

Remark 1.2.20. Since in a prefactorization system $\mathbb{F} = (\mathcal{E}, \mathcal{M})$ each class uniquely determines the other, the prefactorization is characterized by any of the two markings $\mathcal{E}, \mathcal{M}$ and the poset $\textsc{pf}(\mathbf{C})$ defined in **1.2.17** can be confused with a sub-poset of $\mathrm{Mrk}(\mathbf{C})$ defined in **1.2.8**; accordingly with **1.2.17** the class of all prefactorizations $\mathbb{F} = (\mathcal{E}, \mathcal{M})$ on an $\infty$-category $\mathbf{C}$ is a poset whose greatest and smallest elements are respectively

$$ {}_{\perp}(\mathbf{C}^{\sharp}) = (\textsc{Eqv}, \mathbf{C}_1) \text{ and } (\mathbf{C}^{\sharp})_{\perp} = (\mathbf{C}_1, \textsc{Eqv}). \tag{1.5} $$

Definition 1.2.21. If $f\colon \Delta[1] \to \mathbf{C}$ is an arrow in $\mathbf{C}$, a *factorization* of $f$ is an element of the simplicial set $\mathrm{Fact}(f)$ defined to be the upper-left corner on the following pullback diagram

$$ \begin{array}{ccc} \mathrm{Fact}(f) & \longrightarrow & \mathbf{C}^{\Delta[1]} \times \mathbf{C}^{\Delta[1]} \\ \downarrow & \lrcorner & \downarrow \\ \Delta[0] & \longrightarrow & \mathbf{C}^{\Delta[1]} \end{array} \tag{1.6} $$

A factorization $\sigma \in \mathrm{Fact}(f)$ will usually be denoted as $(u, v)$ and the sentence "$(v, u)$ is a factorization of $f\colon A \to B$" will be shortened into "$f = u \circ v\colon A \to F \to B$" for an object $F$ determined from time to time, and called the *factor* of $F$.

REMARK 1.2.22. Notice that since $\mathbf{C}$ is an $\infty$-category, so is $\mathrm{Fact}(f)$ for $f\colon A \to B$; its morphisms can be depicted as commutative squares

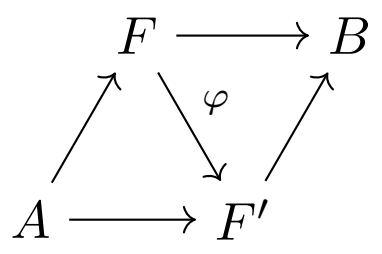

(1.7)

where $\varphi\colon F \to F'$ is a morphism between the factors, such that the two triangles are commutative. Moreover (see the introduction), we will only consider factorizations which are *functorial*, in the obvious sense of being given as functors out of the arrow category $\mathbf{C}^{\Delta[1]}$.

This isn't too restrictive an assumption, since the factorization systems relevant to the present work are all functorial (see, in particular, the proof of our **3.1.1**).

REMARK 1.2.23. A useful characterization of orthogonality available in 1-categorical world is the following: given $f\colon A \to B, g\colon X \to Y$, we have $f \perp g$ if and only if the following square

$$\begin{array}{ccc} \hom(B,X) & \longrightarrow & \hom(B,Y) \\ \downarrow & & \downarrow \\ \hom(A,X) & \longrightarrow & \hom(A,Y) \end{array} \tag{1.8}$$

is a pullback (the proof of this fact is immediate). This characterization can be used to define *enriched* factorization systems (see [DK74, LW] and our discussion in §**3.3**).

This characterization exports, *mutatis mutandis*, to the $\infty$-categorical setting: see [MG14, **A.4.(41)**]

## 1.3 Factorization systems.

A basic transition step from prefactorizations to factorizations is the following result, which is the analogue of the classical results about uniqueness of a factorization with respect to a prefactorization system:

REMARK 1.3.1. If $f\colon \Delta[1] \to \mathbf{C}$ is a morphism, and $\mathbb{F} = (\mathcal{E}, \mathcal{M}) \in \mathrm{PF}(\mathbf{C})$ is a prefactorization system on the $\infty$-category $\mathbf{C}$, then the subspace $\mathrm{Fact}_{\mathbb{F}}(f) \subseteq \mathrm{Fact}(f)$ of factorizations $(e, m)$ such that $e \in \mathcal{E}, m \in \mathcal{M}$, is a contractible simplicial set as soon as it is nonempty.

*Proof.* It all boils down to solving the right lifting problem: if $f\colon A \to B$ can be factored in two ways $(e, m), (e', m') \in \mathrm{Fact}_{\mathbb{F}}(f)$, the first lifting problem gives "comparison" arrows $X \leftrightarrows X'$, and the other two (together with essential uniqueness of the factorization) entail that $u, v$ are mutually inverse.

$$\begin{array}{ccc} A & \xrightarrow{e'} & X' \\ {\scriptstyle e}\downarrow & {\scriptstyle u}\nearrow \swarrow{\scriptstyle v} & \downarrow{\scriptstyle m'} \\ X & \xrightarrow[m]{} & B \end{array} \qquad \begin{array}{ccc} A & \xrightarrow{e} & X \\ {\scriptstyle e}\downarrow & {\scriptstyle vu}\nearrow & \downarrow{\scriptstyle m} \\ X & \xrightarrow[m]{} & B \end{array} \qquad \begin{array}{ccc} A & \xrightarrow{e'} & X' \\ {\scriptstyle e'}\downarrow & {\scriptstyle uv}\nearrow & \downarrow{\scriptstyle m'} \\ X' & \xrightarrow[m]{} & B \end{array} \tag{1.9}$$

□

Definition 1.3.2. ($\mathbb{F}$-crumbled morphisms): Given a prefactorization $\mathbb{F} \in \text{PF}(\mathbf{C})$ we say that an arrow $f: X \to Y$ is $\mathbb{F}$-*crumbled*, (or $(\mathcal{E}, \mathcal{M})$-*crumbled* for $\mathbb{F} = (\mathcal{E}, \mathcal{M})$) when there exists a(n essentially unique, in view of Remark **1.3.1**) factorization for $f$ as a composition $m \circ e$, with $e \in \mathcal{E}$, $m \in \mathcal{M}$; let $\sigma_{\mathbb{F}}$ be the class of all $\mathbb{F}$-crumbled morphisms, and define

$$\text{PF}_{\mathcal{S}}(\mathbf{C}) = \{\mathbb{F} \mid \sigma_{\mathbb{F}} \supset \mathcal{S}\} \subset \text{PF}(\mathbf{C}). \tag{1.10}$$

Definition 1.3.3. A prefactorization system $\mathbb{F} = (\mathcal{E}, \mathcal{M})$ in $\text{PF}(\mathbf{C})$ is said to be a *($\infty$-categorical) factorization system* on $\mathbf{C}$ if $\sigma_{\mathbb{F}} = \hom(\mathbf{C})$; factorization systems, identified with $\text{PF}_{\hom(\mathbf{C})}(\mathbf{C})$, form a sub-poset $\text{FS}(\mathbf{C}) \leq \text{PF}(\mathbf{C})$.

This last definition (factorizations "crumble everything", i.e. split every arrow in two) justifies the form of a more intuitive presentation for a factorization system on $\mathbf{C}$, modeled on the classical, 1-categorical definition:

Definition 1.3.4. ($\infty$-categorical Factorization System): Let $\mathbf{C}$ be an $\infty$-category; a *factorization system* (FS for short) $\mathbb{F}$ on $\mathbf{C}$ consists of a pair of markings $\mathcal{E}, \mathcal{M} \in \text{Mrk}(\mathbf{C})$ such that

(1) For every morphism $h: X \to Z$ in $\mathbf{C}$ we can find a factorization $X \xrightarrow{e} Y \xrightarrow{m} Z$, where $e \in \mathcal{E}$ and $m \in \mathcal{M}$; an evocative notation for this condition, which we sometimes adopt, is $\mathbf{C} = \mathcal{M} \circ \mathcal{E}$;
(2) $\mathcal{E} = {}^{\perp}\mathcal{M}$ and $\mathcal{M} = \mathcal{E}^{\perp}$.

It is useful to introduce the following alternative formalism to express the class of $\mathbb{F}$-crumbled morphisms:

Definition 1.3.5. (The "$\fatsemi$" symbol): Let $\mathbf{C}$ be an $\infty$-category and $\mathcal{A}, \mathcal{B} \subseteq \hom(\mathbf{C})$; we denote $\mathcal{A} \fatsemi \mathcal{B}$ the class of all $f \in \hom(\mathbf{C})$ such that there exists a factorization $f = b \circ a$ with $a \in \mathcal{A}, b \in \mathcal{B}$.

It is obvious that for each prefactorization $\mathbb{F} = (\mathcal{E}, \mathcal{M})$, $\sigma_{\mathbb{F}} = \mathcal{E} \fatsemi \mathcal{M}$ and that a prefactorization $\mathbb{F} = (\mathcal{E}, \mathcal{M})$ is a factorization if and only if $\mathcal{E} \fatsemi \mathcal{M} = \hom(\mathbf{C})$.

The proof of the following lemma is an immediate consequence of Def. **1.3.5**, **1.2.4**:

Lemma 1.3.6. Let $\mathbf{C}$ be an $\infty$-category, and let $\mathcal{A}, \mathcal{B}, \mathcal{C}, \mathcal{D} \in \text{Mrk}(\mathbf{C})$.

- If $\mathcal{A} \perp \mathcal{C}$ and $\mathcal{A} \perp \mathcal{D}$, then $\mathcal{A} \perp (\mathcal{C} \fatsemi \mathcal{D})$;
- If $\mathcal{A} \perp \mathcal{C}$ and $\mathcal{B} \perp \mathcal{C}$, then $(\mathcal{A} \fatsemi \mathcal{B}) \perp \mathcal{C}$.

REMARK 1.3.7. The collection of all factorization systems on an $\infty$-category $\mathbf{C}$ form a poset $\text{FS}(\mathbf{C})$ with respect to the partial order induced by $\text{PF}(\mathbf{C})$.

REMARK 1.3.8. In the presence of condition (1) of Definition **1.3.3**, the second condition may be replaced by

(2a) $\mathcal{E} \perp \mathcal{M}$ (namely $\mathcal{E} \subset {}^{\perp}\mathcal{M}$ and $\mathcal{M} \subset \mathcal{E}^{\perp}$);
(2b) $\mathcal{E}$ and $\mathcal{M}$ are closed under isomorphisms in $\mathbf{C}^{\Delta[1]}$.

Notice that this is precisely [Lur09, Def. **5.2.8.8**].

REMARK 1.3.9. Condition (2) of the previous Definition (or the equivalent pair of conditions (2a), (2b)) entails that each of the two classes $(\mathcal{E}, \mathcal{M})$ in a factorization system on $\mathbf{C}$ uniquely determines the other (compare the analogous statement about prefactorizations): this is [Lur09, Remark **5.2.8.12**].

It is often of great interest to determine whether a given class right-generates or left-generates (Def. **1.2.19**) a factorization system and not only a prefactorization: a general procedure to solve this problem is to invoke the *small object argument.*

THEOREM 1.3.10. (SMALL OBJECT ARGUMENT): Let $\mathbf{C}$ be an $\infty$-category, and $\mathcal{J} \in \text{Mrk}(\mathbf{C})$. If for each $f: I \to J$ in $\mathcal{J}$ the functor $\hom(I, -)$ commutes with filtered colimits, then $\mathcal{J}_\perp$ is a factorization system on $\mathbf{C}$.

REMARK 1.3.11. Let $\mathbf{C}$ be an $\infty$-category with initial object $\varnothing$. If a class $\mathcal{K}$ is generated via the small object argument by a small set then $0/\mathcal{K}$ is generated as the object-orthogonal of a small set.

### 1.3.1 Weak factorization systems.

We will not make use of the content of this subsection, as our main results unavoidably need uniqueness for solutions of lifting problems and factorizations; we only record this definition for the sake of completeness.

There is a more general notion of *weak* factorization system on an $\infty$-category, again modeled on the 1-categorical notion. This "weakness" shows up in two respects:

- Orthogonality is no longer strong: in $\infty$-categories, this means that the space of solution is no longer contractible, but only *nonempty.*
- Factorization is no longer unique up to a unique equivalence: this means that there are possibly many different connected components on the space of factorizations.

A 1-dimensional example of such a structure is the pair of classes (MONO, EPI) in the category **Set** of sets and functions; there seems to be no mention of this condition in the world of $\infty$-categories, in [Lur09] or [Joy08] (this could possibly be related to the fact that the concept of "monomorphism" does not naturally belong to the world of $\infty$-categories, see also **1.5.6** below).

However, a tentative definition of a *model* $\infty$*-category* has been given in [MG14], with applications to Goerss-Hopkins obstruction theory. [MG14, §**2**] contains plenty of examples of weak factorization systems in $\infty$-categories.

Definition 1.3.12. (Weak orthogonality): Let $f, g$ be two edges in an $\infty$-category $\mathbf{C}$. We will say that $f$ is *weakly left orthogonal* to $g$ (or equivalently that $g$ is *weakly right orthogonal* to $f$) if in any lifting problem like (**1.2**) the space of solutions $a$ is nonempty. We denote this binary relation as $f \pitchfork g$, and the resulting Galois connection ${}^{\pitchfork}(-) \dashv (-)^{\pitchfork}$.

Now, a *weak prefactorization system* is a pair of classes $\mathbb{F} = (\mathcal{E}, \mathcal{M}) \subseteq \hom(\mathbf{C}) \times \hom(\mathbf{C})$ such that $\mathcal{E} = {}^{\pitchfork}\mathcal{M}$, $\mathcal{M} = \mathcal{E}^{\pitchfork}$. A *weak factorization system* is a weak prefactorization system such that every arrow is $\mathbb{F}$-crumbled (Def. **1.3.2**).

Weak factorization systems are organized in a poset $\text{WFS}(\mathbf{C})$ which contains as a sub-poset $\text{FS}(\mathbf{C})$ of Def. **1.2.8**.

It is possible to relate weak orthogonality to strong orthogonality, and give conditions ensuring that a weak factorization system is indeed strong: see [RT07, §**1**] for further details, and consult [MG14] for further information about model $\infty$-categories.

Example 1.3.13. A *Quillen model structure* on a small-bicomplete $\infty$-category $\mathbf{C}$ is defined by three markings $(\mathcal{W}\|, \mathcal{F}\rangle\lfloor, \mathcal{C}\wr\{)$ such that

- $\mathcal{W}\|$ is a 3-for-2 class (Def. **1.4.6**) containing all isomorphisms and closed under retracts;
- The markings $(\mathcal{C}\wr\{, \mathcal{W}\| \cap \mathcal{F}\rangle\lfloor)$ and $(\mathcal{W}\| \cap \mathcal{C}\wr\{, \mathcal{F}\rangle\lfloor)$ both form a weak factorization system on $\mathbf{C}$.

## 1.4 Closure properties.

Definition 1.4.1. (Closure operators associated to markings): Let $\mathbf{C}$ be an $\infty$-category. A marking $\mathcal{J} \in \text{Mrk}(\mathbf{C})$ is called w.) **wide** if it contains all the isomorphisms and it is closed under composition;

A wide marking $\mathcal{J}$ (in an $\infty$-category $\mathbf{C}$ which admits in each case the $\infty$-categorical co/limits needed to state the definition) is called p.) **presaturated** if is closed under *co-base change*, i.e. whenever we are given arrows $j \in \mathcal{J}$, and $h$ such that we can form the pushout

$$\begin{array}{ccc} A & \xrightarrow{h} & X \\ {\scriptstyle j}\downarrow & & \downarrow{\scriptstyle j'} \\ B & \longrightarrow & Y \end{array} \tag{1.11}$$

then the arrow $j'$ is in $\mathcal{J}$;

Q.) **almost saturated** if it is presaturated and closed under *retracts* (in the category $\mathbf{C}^{\Delta[1]}$), i.e. whenever we are given a diagram like

$$\begin{array}{ccccc} A & \xrightarrow{i} & A' & \xrightarrow{r} & A \\ \downarrow{\scriptstyle u} & & \downarrow{\scriptstyle v} & & \downarrow{\scriptstyle u} \\ B & \xrightarrow{i'} & B' & \longrightarrow & B \end{array} \tag{1.12}$$

where $ri = \mathrm{id}_A$ and $r'i' = \mathrm{id}_B$, if $v$ lies in $\mathcal{J}$, then the same is true for $u$;

C.) **cellular** if it is presaturated and closed under *transfinite composition*, namely whenever we have a cocontinuous functor $F\colon \alpha \to \mathcal{J}^{(4)}$ defined from a limit ordinal $\alpha$ admits a composite in $\mathcal{J}$, i.e. the canonical arrow

$$F(0) \longrightarrow F(\alpha) = \varinjlim_{i<\alpha} F(i) \tag{1.13}$$

lies in $\mathcal{J}$;

S.) **saturated** if it is almost saturated and cellular.

All these properties induce suitable closure operators, encoded as suitable (idempotent) monads on $\mathrm{Mrk}(\mathbf{C})$, defined for any property $\mathsf{p}$ among $\{\mathrm{W}, \mathrm{P}, \mathrm{Q}, \mathrm{C}, \mathrm{S}\}$ as

$$\mathsf{p}(-)\colon \mathcal{S} \mapsto \mathsf{p}(\mathcal{S}) = \bigcap_{\mathcal{U} \supseteq \mathcal{S}} \Big\{ \mathcal{U} \in \mathrm{Mrk}(\mathbf{C}) \mid \mathcal{U}\ P \text{ property has} \Big\} \tag{1.14}$$

In classical category theory, the *cellularization* $\mathsf{c}(-)$ and the *saturation* $\mathsf{s}(-)$ of a marking $\mathcal{J}$ on $\mathbf{C}$ are of particular interest (especially in homotopical algebra), in view of what we state in Prop. **1.4.3**.

REMARK 1.4.2. A little more generality is gained by supposing that the cardinality of the coproducts or the transfinite compositions in $\mathbf{C}$ is bounded by some (regular) cardinal $\alpha$. In this case we speak of $\alpha$-saturated or $\alpha$-cellular classes, and define the closure operators of *$\alpha$-cellularization* and *$\alpha$-saturation*, etc.

(4) This notation is a shorthand to denote the fact that each edge $F(i \le j)\colon F(i) \to F(j)$ is an element of $\mathcal{J}$; alternatively, we can regard this notation as consistent, via the obvious identification between markings on $\mathbf{C}$ and full subcategories of $\mathbf{C}$.

The following Proposition is a standard result in the theory of factorization systems, which we will often need throughout the discussion; a proof of the 1-categorical version of the statement can be found in any of the references about factorization systems provided in the bibliography.

PROPOSITION 1.4.3. Let $(\mathbf{C}, \mathcal{S})$ be a marking of a cocomplete $\infty$-category $\mathbf{C}$; then the marking ${}^{\perp}\mathcal{S}$ of $\mathbf{C}$ is a saturated class. In particular, the left class of a weak factorization system in a cocomplete $\infty$-category is saturated.

Completely dual definitions give rise to co-p-classes.[(5)] again, suitable monads acting as co-p-closure operators are defined on $\mathrm{Mrk}(\mathbf{C})$, giving the dual of Proposition **1.4.3**:

PROPOSITION 1.4.4. Let $(\mathbf{C}, \mathcal{S})$ be a marking of a cocomplete $\infty$-category $\mathbf{C}$; then the marking $\mathcal{S}^{\perp}$ of $\mathbf{C}$ is a co-saturated class. In particular, the right class of a weak factorization system in a complete category is co-saturated.

PROPOSITION 1.4.5. Let $\mathbf{C}$ be an $\infty$-category and $\mathbb{F} = (\mathcal{E}, \mathcal{M}) \in \mathrm{FS}(\mathbf{C})$; then $\mathcal{E} \cap \mathcal{M}$ equals the class of all equivalences in $\mathbf{C}$.

*Proof.* Again, the proof in a 1-category case can be found in any reference about factorization systems. The idea is extremely simple: if $g \in \mathcal{E} \cap \mathcal{M}$ then it is orthogonal to itself, and we can invoke **1.2.11**. □

DEFINITION 1.4.6. Let $\mathcal{S} \in \mathrm{Mrk}(\mathbf{C})$; then for each 2-simplex in $\mathbf{C}$ representing a composable pair of arrows, whose edges are labeled $f, g$, and $fg$ we say that

- $\mathcal{S}$ is L32 if $f, fg \in \mathcal{S}$ imply $g \in \mathcal{S}$;
- $\mathcal{S}$ is R32 if $fg, g \in \mathcal{S}$ imply $f \in \mathcal{S}$.

A marking $\mathcal{S}$ which is closed under composition and both L32 and R32 is said to *satisfy the 3-for-2 property*, or a *3-for-2 class*.

PROPOSITION 1.4.7. Given a FS $(\mathcal{E}, \mathcal{M})$ in the $\infty$-category $\mathbf{C}$, then

(1) If the $\infty$-category $\mathbf{C}$ has $K$-colimits, for $K$ a given simplicial set, then the full subcategory of $\mathbf{C}^{\Delta[1]}$ spanned by $\mathcal{E}$ has $K$-colimits; dually, if the $\infty$-category $\mathbf{C}$ has $K$-limits, then the full subcategory of $\mathbf{C}^{\Delta[1]}$ spanned by $\mathcal{M}$ has $K$-limits;
(2) The class $\mathcal{E}$ is R32, and the class $\mathcal{M}$ is L32 (see Def. **1.4.6**).

*Proof.* Point (1) is [Lur09, Prop. **5.2.8.6**]; point (2) is easy to prove for 1-categories, and then the translation to the $\infty$-categorical setting is straightforward.[(6)] □

[(5)] Obviously, wideness and closure under retracts are auto-dual properties. The definition of *transfinite op-composition* needed to define co-cellularity may be difficult to guess; see [Joy08] for reference.

[(6)] This translation process being often straightforward, here and everywhere a bibliographic support is needed, we choose to rely on classical sources to prove most of the result involving $\infty$-categorical factorization systems. This should cause no harm to the reader.

It is a remarkable, and rather useful result, that each of the properties (1) and (2) of the above Proposition characterizes factorizations among weak factorizations: see [RT07, Prop. **2.3**] for more details.

We close this section with an useful observation, showing that in Def. **1.3.4** "factorization is all what matters": asking two classes $(\mathcal{E}, \mathcal{M})$ to uniquely crumble every morphism $f \in \hom(\mathbf{C})$ entails that $(\mathcal{E}, \mathcal{M})$ are mutually orthogonal classes.

REMARK 1.4.8. Let $\mathbb{F} = (\mathcal{E}, \mathcal{M})$ be a pair of wide markings such that every $f \in \hom(\mathbf{C})$ has a unique factorization $f = m \circ e$ with $m \in \mathcal{M}$, $e \in \mathcal{E}$; then, $\mathbb{F}$ is a prefactorization, i.e. $\mathcal{E} \perp \mathcal{M}$.

*Proof.* Given a lifting problem

$$\begin{array}{ccc} A & \xrightarrow{u} & X \\ {\scriptstyle e''}\downarrow & & \downarrow{\scriptstyle m''} \\ B & \xrightarrow[v]{} & Y \end{array} \tag{1.15}$$

we can factor $u$ as $m \circ e$ and $v$ as $m' \circ e'$, so that the square becomes

$$\begin{array}{ccccc} A & \xrightarrow{e} & U & \xrightarrow{m} & X \\ {\scriptstyle e'}\downarrow & & & & \downarrow{\scriptstyle m''} \\ B & \xrightarrow[e']{} & V & \xrightarrow[m']{} & Y \end{array} \tag{1.16}$$

Now, the factorizations $(m'' \circ m, e)$ and $(m', e' \circ e'')$ must be isomorphic by the uniqueness assumption, so that there exists an isomorphism $U \cong V$ which composed with $e', m$ gives a solution to the lifting problem. □

## 1.5 A second glance at factorization.

一条の矢は折るべく十条の矢は折るべからず

Japanese proverb

We add here a different presentation of $\infty$-categorical factorization systems, faithfully following [Joy08, pp. 178—].

DEFINITION 1.5.1. (ORTHOGONALITY AND FILLERS): Let $\mathbf{C}$ be an $\infty$-category, and $u: A \to B, f: X \to Y$ two edges of $\mathbf{C}$. We define the space $\mathrm{Sq}(u, f)$ of commutative squares associated to $(u, f)$ to be the space of simplicial maps $s: \Delta[1] \times \Delta[1] \to \mathbf{C}$ such that $s|_{\Delta^0 \times \Delta[1]} = u, s|_{\Delta[1] \times \Delta^0} = f$.

A *diagonal filler* for $s \in \mathrm{Sq}(u, f)$ consists of an extension $\bar{s}: \Delta[1] \star \Delta[1] \to \mathbf{C}$ (where $\star$ denotes the *join* of simplicial sets, see [Joy08, §**3.1** and **3.2**]) of $s$ along the natural inclusion $\Delta[1] \times \Delta[1] \subset \Delta[1] \star \Delta[1]$.

REMARK 1.5.2. Denote by $\mathrm{Fill}(s)$ the top-left corner of the fiber sequence

$$\begin{array}{ccc} \mathrm{Fill}(s) & \longrightarrow & X^{\Delta[1]\star\Delta[1]} \\ \downarrow & \lrcorner & \downarrow q \\ \Delta[0] & \xrightarrow{s} & X^{\Delta[1]\times\Delta[1]}. \end{array} \tag{1.17}$$

The simplicial set $\mathrm{Fill}(s)$ is a Kan complex, since $q$ is a Kan fibration (as a consequence of [Joy08, Prop. **2.18**]).

This leads us to the following

DEFINITION 1.5.3. We say that the edge $u$ is *left orthogonal* to the edge $f$ in the $\infty$-category $\mathbf{C}$ (or $f$ is *right orthogonal* to $u$) if $\mathrm{Fill}(s)$ is a *contractible* Kan complex for any $s \in \mathrm{Sq}(u, f)$. We denote this relation between $u$ and $f$ as $u \boxslash f$.

A first and natural task is to prove that the two relations $\perp$ and $\boxslash$ defined on the set of edges $\mathbf{C}_1$ of an $\infty$-category coincide: this is immediate since $\Delta[1] \star \Delta[1] = \Delta[3]$ ([Joy08, p. 244]) and since “solved commutative squares” can be identified with simplicial maps $\Delta[3] \to \mathbf{C}$ (there is a unique edge outside the image of $\Delta[1] \times \Delta[1] \subset \Delta[1] \star \Delta[1]$; this is the solution to the lifting problem).

Given this, for the rest of the section we will stick to the notation $f \perp g$ to denote orthogonality in this sense ([Joy08] uses the same symbol and takes **1.5.1** as a definition).

PROPOSITION 1.5.4. Factorization systems can be lifted along left or right fibrations: if $p\colon \mathbf{C} \to \mathbf{D}$ is such a simplicial map, and $(\mathcal{E}, \mathcal{M}) \in \mathrm{FS}(\mathbf{D})$, then $(p^{\leftarrow}(\mathcal{E}), p^{\leftarrow}(\mathcal{M}))$ is a factorization system on $\mathbf{C}$.

COROLLARY 1.5.5. As a consequence, since the simplicial maps $\mathbf{C}_{/X} \to \mathbf{C}$ and $\mathbf{C}_{Y/} \to \mathbf{C}$ are left/right fibrations, every factorization system on $\mathbf{C}$ *lifts* to a factorization system on the slice/coslice $\infty$-category. This is the $\infty$-categorical version of the classical statement saying that a factorization system on $\mathbf{C}$ induces factorization systems on all co/slice categories $\mathbf{C}_{/X}$ and $\mathbf{C}_{Y/}$.

### 1.5.1 A factor-y of examples.

[Joy08, pp. 178—] is an invaluable source of examples for factorization systems on $\infty$-categories; a standard technique to provide such examples is to reduce suitable “niceness” properties for categories (like regularity or exactness, or the possibility to find “Postnikov towers” for morphisms) to the presence of suitable factorization systems on it.

This general tenet remains valid in an $\infty$-categories.

REMARK 1.5.6. We must observe, here, that it is rather difficult (i.e. rather more difficult than in 1-categories) to produce intuitive examples of factorization systems in $\infty$-category, as many of the 1-dimensional examples rely on the intuition that (EPI, MONO) is a well-behaved and paradigmatic example of such a structure in many categories (such as sets, toposes, abelian categories...: the factorization systems such that every $\mathcal{E}$ is a epimorphism, and every $\mathcal{M}$ is a monomorphism are called *proper* in [FK72, Kel80] to suggest how this notion is common and familiar).

This cannot be achieved in $\infty$-category theory, as [Lur09, p. 562] conveys the intuition that the notion of monomorphism is not as meaningful in $\infty$-category theory as it is in 1-category theory (compare, however, the statement that every topos has an (EPI, MONO)-factorization system with the existence of a "Postnikov" factorization system on each $\infty$-topos, [Lur09]).

Several construction can be performed inside the category of categories with factorization system: these are classical definitions that can be recovered in every text about factorization systems (especially those with an interest towards model categories).

EXAMPLE 1.5.7. (CO/PRODUCTS, CO/SLICES): Let $\mathbf{C}$ be an $\infty$-category; then every co/slice of $\mathbf{C}$ inherits a factorization system from $\mathbb{F} = (\mathcal{E}, \mathcal{M}) \in \text{FS}(\mathbf{C})$ obtained by putting

$$\begin{aligned} \mathcal{E}_{/X} &= \{(Y, f) \xrightarrow{\varphi} (Z, g) \mid \varphi \in \mathcal{E}\} \\ \mathcal{M}_{/X} &= \{(Y, f) \xrightarrow{\psi} (Z, g) \mid \psi \in \mathcal{M}\} \end{aligned} \tag{1.18}$$

(the definition for coslices $\mathbf{C}_{X/}$ is analogous).

Let $\{\mathbf{C}_i\}$ be any small family of $\infty$-categories; the product $\prod \mathbf{C}_i$ of all the elements of the family inherits a factorization system from a family $\mathbb{F}_i = (\mathcal{E}_i, \mathcal{M}_i) \in \text{FS}(\mathbf{C}_i)$, defined by putting

$$\begin{aligned} \prod \mathcal{E}_i &= \{(f_i)_{i\in I} \mid f_i \in \mathcal{E}_i \ \forall i \in I\} \\ \prod \mathcal{M}_i &= \{(g_i)_{i\in I} \mid g_i \in \mathcal{M}_i \ \forall i \in I\}. \end{aligned} \tag{1.19}$$

A similar definition works for coproducts: the coproduct $\coprod \mathbf{C}_i$ inherits a factorization system defined in a dual fashion (an arrow $f \in \coprod \mathbf{C}_i$ lies in one and only one $\mathbf{C}_{i^*}$; $f \in \coprod \mathcal{E}$ if and only if $f \in \mathcal{E}_{i^*}$).

EXAMPLE 1.5.8. (SURJECTION-MONO FACTORIZATIONS): We say that an arrow $f \colon X \to Y$ in an $\infty$-category is *monic* if the square

$$\begin{array}{ccc} X & = & X \\ \| & & \downarrow{\scriptstyle f} \\ X & \xrightarrow[f]{} & Y \end{array} \tag{1.20}$$

is cartesian. The class of monic arrows in $\mathbf{C}$ is collected in a marking $\text{MONO}(\mathbf{C}) = \text{MONO}$.

The class of *surjective* arrows is defined to be the class ${}^{\perp}\text{MONO}(\mathbf{C})$; we say that the $\infty$-category $\mathbf{C}$ *has a surjection-mono factorization* if the prefactorization $({}^{\perp}\text{MONO}(\mathbf{C}), \text{MONO}(\mathbf{C}))$ is also a factorization system.

DEFINITION 1.5.9. (REGULAR $\infty$-CATEGORY): A finitely complete $\infty$-category is said to be *regular* if it admits a *pullback-stable* surjection-mono factorization system; the coherent nerve of the category of Kan complexes, as a full sub-$\infty$-category of the nerve of the whole **sSet**, is regular. Notice that this is the $\infty$-categorical counterpart of *Barr-regular categories.*

### 1.5.2 Chains of factorization systems.

DEFINITION 1.5.10. ($k$-ARY FACTORIZATION SYSTEM): Let $k \geq 2$ be a natural number. A *$k$-fold factorization system* on a category $\mathbf{C}$ consists of a monotone map $\phi \colon \Delta[k-2] \to \text{FS}(\mathbf{C})$, where the codomain has the partial order of Def. **1.2.17**; denoting $\phi(i) = \mathbb{F}_i$, a $k$-fold factorization system on $\mathbf{C}$ consists of a chain

$$\mathbb{F}_1 \preceq \cdots \preceq \mathbb{F}_{k-1}, \tag{1.21}$$

This means that if we denote $\mathbb{F}_i = (\mathcal{E}_i, \mathcal{M}_i)$ we have two chains –any of which determines the other– in $\hom(\mathbf{C})$:

$$\begin{gathered} \mathcal{E}_1 \supset \mathcal{E}_2 \supset \cdots \supset \mathcal{E}_{k-1}, \\ \mathcal{M}_1 \subset \mathcal{M}_2 \subset \cdots \subset \mathcal{M}_{k-1}. \end{gathered}$$

The definition of a $k$-ary factorization system is motivated by the fact that a chain in $\text{FS}(\mathbf{C})$ results in a way to factor each arrow "coherently" as the composition of $k$ pieces, coherently belonging to the various classes of arrows. This is explained by the following simple result:

LEMMA 1.5.11. Every arrow $f \colon A \to B$ in a category endowed with a $k$-ary factorization system $\mathbb{F}_1 \preceq \cdots \preceq \mathbb{F}_{k-1}$ can be uniquely factored into a composition

$$A \xrightarrow{\mathcal{E}_1} X_1 \xrightarrow{\mathcal{E}_2 \cap \mathcal{M}_1} X_2 \to \cdots \to X_{k-2} \xrightarrow{\mathcal{E}_{k-1} \cap \mathcal{M}_{k-2}} X_{k-1} \xrightarrow{\mathcal{M}_{k-1}} B, \tag{1.22}$$

where each arrow is decorated with the class it belongs to.

*Proof.* For $k = 1$ this is the definition of factorization system: given $f \colon X \to Y$, we have its $\mathbb{F}_{i_1}$-factorization

$$X \xrightarrow{\mathcal{E}_{i_1}} Z_{i_1} \xrightarrow{\mathcal{M}_{i_1}} Y. \tag{1.23}$$

Then we work inductively on $k$. Given an arrow $f \colon X \to Y$ we first consider its $\mathbb{F}_{i_k}$-factorization

$$X \xrightarrow{\mathcal{E}_{i_k}} Z_{i_k} \xrightarrow{\mathcal{M}_{i_k}} Y, \tag{1.24}$$

and then observe that the chain $i_1 \le \cdots \le i_{k-1}$ induces a $(k-1)$-ary factorization system on $\mathbf{C}$, which we can use to decompose $Z_{i_k} \to Y$ as

$$Z_{i_k} \xrightarrow{\mathcal{E}_{i_{k-1}}} Z_{i_{k-1}} \xrightarrow{\mathcal{E}_{i_{k-2}}\cap\mathcal{M}_{i_{k-1}}} Z_{i_{k-2}} \to \cdots \to Z_{i_2} \xrightarrow{\mathcal{E}_{i_1}\cap\mathcal{M}_{i_2}} Z_{i_1} \xrightarrow{\mathcal{M}_{i_1}} Y, \tag{1.25}$$

and we are only left to prove that $Z_{i_k} \to Z_{i_{k-1}}$ is actually in $\mathcal{E}_{i_{k-1}} \cap \mathcal{M}_{i_k}$. This is an immediate consequence of the left cancellation property for the class $\mathcal{M}_{i_1}$. Namely, since $\mathcal{M}_{i_1} \subseteq \mathcal{M}_{i_2} \subseteq \cdots \subseteq \mathcal{M}_{i_k}$, and $\mathcal{M}_{i_k}$ is closed for composition, the morphism $Z_{i_{k-1}} \to Y$ is in $\mathcal{M}_{i_k}$. Then the L32 property applied to

$$Z_{i_k} \to Z_{i_k-1} \xrightarrow{\mathcal{M}_{i_k}} Y \tag{1.26}$$

concludes the proof. ☐

### 1.5.2.1 The transfinite case.

We are now interested to refine the previous theory in order to deal with possibly infinite chains of factorization systems. From a conceptual point of view, it seems natural how to extend the former definition to an infinite ordinal $\alpha$; it must consists on a "suitable" functor $F\colon \alpha \to \mathrm{FS}(\mathbf{C})$.

The problem is that suitable necessary co/continuity assumptions for such a $F$ might be covered by the fact that its domain is finite (and in particular admits an initial and a terminal object): in principle, dealing with infinite quantities could force such $F$ to fulfill some other properties in order to preserve the basic intuition behind factorization.

We start, now, by recalling a number of properties motivating Def. **1.5.10** below.

NOTATION 1.5.12. A factorization system on $\mathbf{C}$ naturally defines a pair of pointed/co-pointed endofunctors on $\mathbf{C}^{\Delta[1]}$, starting from the factorization

$$\begin{array}{ccc} X & \longrightarrow & Y \\ & {\scriptstyle \overleftarrow{F}(f)} \searrow \quad \nearrow {\scriptstyle \overrightarrow{F}(f)} & \\ & F(f). & \end{array} \tag{1.27}$$

(This has also been noticed in [Lur09]). A refinement of this notion (see [GT06, ?, Rie11]) regards this pair of functors as monad/comonad on $\mathbf{C}^{\Delta[1]}$: in this case $F\colon \mathbf{C}^{\Delta[1]} \to \mathbf{C}$ is a functor and $f \mapsto \overleftarrow{F}(f)$ has the structure of a (idempotent) comonad, whose comultiplication is

$$\begin{array}{ccc} Ff & \xrightarrow{\overleftarrow{F}(\overrightarrow{F}(f))} & FFf \\ {\scriptstyle \overrightarrow{F}(f)}\big\downarrow & & \big\downarrow{\scriptstyle \overrightarrow{F}(\overrightarrow{F}(f))} \\ Y & = & Y \end{array} \tag{1.28}$$

and $f \mapsto \overrightarrow{F}(f)$ has the structure of a (idempotent) monad, whose multiplication is

$$\begin{array}{ccc} X & = & Ff \\ {\scriptstyle \overleftarrow{F}(\overleftarrow{F}(f))}\downarrow & & \downarrow{\scriptstyle \overleftarrow{F}(f)} \\ FFf & \xrightarrow[\overrightarrow{F}(\overleftarrow{F}(f))]{} & Ff. \end{array} \tag{1.29}$$

Remark 1.5.13. (On functors to posetal categories): (Small) posets form the category **PCat** of (small) posetal categories (categories where every hom-set is either empty or has one element).

This category is reflective in **Cat**, since we have an adjunction

$$\mathbf{PCat}(p\mathbf{C}, P) \cong \mathbf{Cat}(\mathbf{C}, P). \tag{1.30}$$

(The poset $p(\mathbf{J})$ results as the partially ordered set $\mathrm{Ob}(\mathbf{J})$ where $A \leq B$ iff there is an arrow from $A$ to $B$.) Hence functors $\mathbf{J} \to P$ are uniquely determined by a monotone function $p(\mathbf{J}) = J \to P$, with respect to this order on $J$.

Remark 1.5.14. (On (co)limits in slice categories): Slice and coslice categories $\mathbf{C}_{X/}$ and $\mathbf{C}_{X/}$ are complete and cocomplete whenever $\mathbf{C}$ is: colimits in $\mathbf{C}_{X/}$ and limits $\mathbf{C}_{X/}$ are simply reflected by the natural forgetful functor $U\colon \mathbf{C}_{X/}, \mathbf{C}_{X/} \to \mathbf{C}$, so that the limit of a diagram $j \mapsto \begin{bmatrix} X \\ \downarrow \\ A_j \end{bmatrix}$ is simply the arrow $\begin{bmatrix} X \\ \downarrow \\ \varprojlim_j^{\mathbf{C}} A_j \end{bmatrix}$ (and dually for colimits in $\mathbf{C}_{X/}$); limits in a slice category, and colimits in a coslice category are, generally, more difficult to compute.

The general recipe for (say) colimits of a functor $F\colon \mathbf{J} \to \mathbf{C}_{X/}$ exploits the isomorphism

$$\mathrm{Fun}(\mathbf{J}, \mathbf{C}_{X/}) \cong \mathrm{Fun}_F(\mathbf{J}^{\triangleleft}, \mathbf{C}) \tag{1.31}$$

where $\mathbf{J}^{\triangleleft}$ is the category $[0] \star \mathbf{J}$ obtained freely adding an initial object, and $\mathrm{Fun}_F(\mathbf{J}^{\triangleleft}, \mathbf{C})$ is the category of functors $\mathbf{J}^{\triangleleft} \to \mathbf{C}$ which coincide with $F$ when restricted to $\mathbf{J} \subset [0] \star \mathbf{J}$.

Now fortunately, whenever the indexing category is *connected*, limits in slice categories, and colimits in coslice categories are again reflected along the natural forgetful $U$: a particular application of this result, when $J$ is an ordinal regarded as a category, serves to state the following definition.

Definition 1.5.15. Let $\alpha$ be an ordinal. A *$\alpha$-ary factorization system*, or *factorization system in $\alpha$-stages*, on $\mathbf{C}$ consists of a monotone function $\alpha \to \mathrm{FS}(\mathbf{C})\colon i \mapsto \mathbb{F}_i$ such that, if we denote by

$$\begin{array}{ccc} X & \longrightarrow & Y \\ & {\scriptstyle \overleftarrow{F}_i(f)} \searrow \quad \nearrow {\scriptstyle \overrightarrow{F}_i(f)} & \\ & F_i(f) & \end{array} \tag{1.32}$$

the $\mathbb{F}_i$-factorization of $f\colon X \to Y$, we have the following two "tame convergence" conditions:

$$\varprojlim_{i\in\alpha} \overrightarrow{F}_i(f) = \varprojlim_{i\in\alpha} \begin{bmatrix} F_i(f) \\ \downarrow \\ Y \end{bmatrix} = \begin{bmatrix} X \\ \downarrow \\ Y \end{bmatrix}; \qquad \varinjlim_{i\in\alpha} \overleftarrow{F}_i(f) = \varinjlim_{i\in\alpha} \begin{bmatrix} X \\ \downarrow \\ F_i(f) \end{bmatrix} = \begin{bmatrix} X \\ \downarrow \\ Y \end{bmatrix}$$

$$\varinjlim_{i\in\alpha} \overrightarrow{F}_i(f) = \varinjlim_{i\in\alpha} \begin{bmatrix} F_i(f) \\ \downarrow \\ Y \end{bmatrix} = 1_Y; \qquad \varprojlim_{i\in\alpha} \overleftarrow{F}_i(f) = \varprojlim_{i\in\alpha} \begin{bmatrix} X \\ \downarrow \\ F_i(f) \end{bmatrix} = 1_X$$

(all the diagrams have to be considered defined in suitable slice and coslice categories) which can be summarized in the presence of "extremal" factorizations

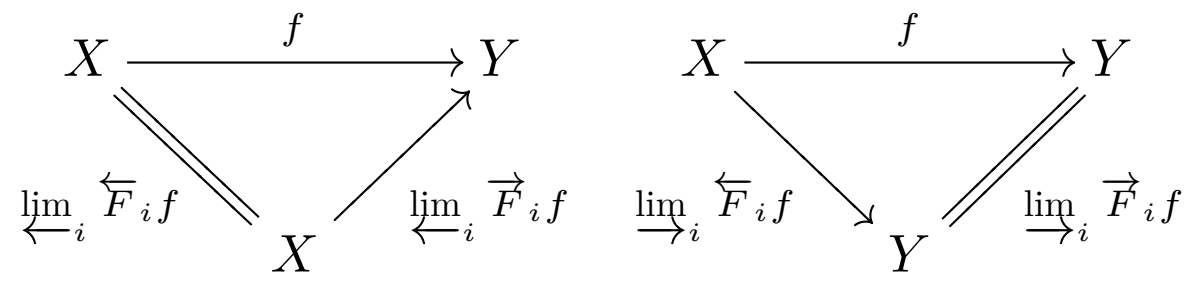

(1.33)

THEOREM 1.5.16. (THE MULTIPLE SMALL OBJECT ARGUMENT): Let $\mathcal{J}_1 \subseteq \cdots \subseteq \mathcal{J}_n$ be a chain of markings on $\mathbf{C}$; if each class $\mathcal{J}_\alpha$ has small domains then applying $n$ times the small object argument, the extensivity of the $_\perp((-)^\perp)$ and $(^\perp(-))^\perp$ closure operators entails that there exists a chain of factorization systems

$$\big(^\perp(J_n^\perp), J_n^\perp\big) \preceq \cdots \preceq \big(^\perp(J_1^\perp), J_1^\perp\big) \tag{1.34}$$

(the order relation is that of Def. **1.2.17**).

# Chapter 2

# Reflectivity and Normality

We now translate in the setting of $\infty$-categories the main definitions and results outlined in [CHK85], with a special attention to the setting of Lurie's stable $\infty$-categories. We take for granted all the basic definition of stable $\infty$-category, $t$-structure and properties thereof, outlined in appendix **A**.

The paper [CHK85] extensively describes various types of reflective subcategories of a given category $\mathbf{C}$ obtained by means of factorization systems on $\mathbf{C}$. A number of results are discussed and applied to additive and abelian categories, pointed categories, etc. leading to the notion of a *normal torsion theory*.

Among these results, one the most interesting for the present purposes is the antitone bijection established between localizations of $\mathbf{C}$, collected in the poset $\mathrm{Rex}(\mathbf{C})$[(1)], and factorization system $\mathbb{F} = (\mathcal{E}, \mathcal{M})$ such that both classes are 3-for-2 in the sense of our **1.4.6**: this analysis paves the way to the foundations for a "theory of torsion and torsion-free classes" in non-additive categories, and it is a starting point to motivate the $\infty$-categorical translation of the theory.

The present chapter profits from the blanket assumption of stability for the $\infty$-category $\mathbf{C}$; here a triangulated structure on the homotopy category $\mathrm{Ho}(\mathbf{C})$ is induced by easy and categorically natural axioms, verified at the "higher" level, and universal properties utterly simplify the proof of the analogy "co/reflective pairs" = "$t$-structures". From this we deduce a rather primitive statement (hinted at in [AHHK07, BR07], and others, but never extracted from the land of *folklore*: for a discussion on this point, see **2.3.8**): this result is called "Rosetta stone" theorem in **3.1.1**, and constitutes the backbone of the thesis.

We now sketch the content of [CHK85, §**6**], and offer an $\infty$-categorical

(1) This notation may appear deceiving: "Rex" stands here for **re**fle**x**ions, and not for *right exact*.

counterpart thereof: given a *reflective factorization system* $\mathbb{F} = (\mathcal{E}, \mathcal{M})$ (Prop. **2.1.6**) on an $\infty$-category with initial and terminal objects, the classes

$$\begin{aligned} \varnothing/\mathcal{E} &\stackrel{\triangle}{=} \{X \in \mathbf{C} \mid \varnothing \to X \in \mathcal{E}\}, \\ \mathcal{M}/1 &\stackrel{\triangle}{=} \{Y \in \mathbf{C} \mid Y \to 1 \in \mathcal{M}\} \end{aligned} \tag{2.1}$$

are respectively a coreflective and reflective subcategory of $\mathbf{C}$. A number of additional requests on $\mathbb{F}$ ensure that these two subcategories behave well under several other constructions, or enjoy additional properties of mutual interaction (e.g. determining each other up to equivalence, via the object-orthogonality relation).

This is, again, a chapter devoted to purely categorical results; we can nevertheless outline a couple of interesting points, even at this level of abstraction.

In their review of [CHK85], the authors of [RT07] outline a sequence of implications between the properties of *(semi-left/right-)exactness*, *simplicity* and *normality* of a torsion theory $\mathbb{F}$, and confess a certain difficulty in exhibiting a non-artificial example of a *non-normal* torsion theory; they conclude, then (with a certain coherence in the choice of notation), that the notion of non-normality is somewhat pathological, and suggest ([RT07, Remark **4.11**]) that there are few (if any) examples of non-normal torsion theories.

In **2.3.16** we prove that, in the setting of stable $\infty$-categories, the three notion of exactness, simplicity and normality collapse into a single notion (simply called *normality*); this result deserves further investigation in light of the use of reflective factorization systems in [BJ01] and in view of the fact that any category $\mathbf{A}$ has a (canonically constructed) stabilization $\mathrm{S{\scriptsize P}}(\mathbf{A})$, where the asymmetry between normality, semi-exactness and simplicity stated in [CHK85, §**4.4**] disappears.

Notation 2.0.1. A blanket assumption throughout all this chapter is that $\mathbf{C}$ is an $\infty$-category with an initial and terminal object, respectively denoted $\varnothing$ and 1: subsequently we will specialize this assumption by asking that $\mathbf{C}$ is stable (so in particular it is pointed and finitely co/complete). Other specializations (like in Def. **2.2.1** or **2.2.3**) will always be notified to the reader; here we do not strive for a particular sharpness in statements and proofs: several results are still valid outside our main case of interest (i.e. when $\mathbf{C}$ is not stable, but still has finite limits or is at least pointed).

We denote by $\boldsymbol{\tau}_{\mathbf{C}}$ the class of the terminal morphisms $\{t_X : X \to 1 \mid X \in \mathbf{C}\}$, and $\mathrm{Rex}(\mathbf{C})$ be the poset of reflective subcategories $(\mathbf{B}, R)$ of $\mathbf{C}$ (where $R: \mathbf{C} \to \mathbf{B}$ is the reflection functor, left adjoint to the inclusion).

## 2.1 The fundamental connection.

菩提本無樹，明鏡亦非臺。
本來無一物，何處惹塵埃。

Huìnéng

The aim of the present section is to re-enact [CHK85, Prop. **2.2**], where the authors build a correspondence between $\mathrm{PF}_{\mathcal{T}}(\mathbf{C})$ (see Def. **1.3.2**) and $\mathrm{Rex}(\mathbf{C})$.

Proposition 2.1.1. There exists a(n antitone) Galois connection $\Phi \dashv \Psi$ between the posets $\mathrm{Rex}(\mathbf{C})$ and $\mathrm{PF}_{\mathcal{T}}(\mathbf{C})$, where $\Psi$ sends $\mathbb{F} = (\mathcal{E}, \mathcal{M})$ to the subcategory $\mathcal{M}/1 = \{B \in \mathbf{C} \mid (B \to 1) \in \mathcal{M}\}$, and $\Phi$ is defined by sending $(\mathbf{B}, R) \in \mathrm{Rex}(\mathbf{C})$ to the prefactorization *right generated* (see Definition **1.2.19**) by $\hom(\mathbf{B})$.

*Proof.* A complete proof can be found in [CHK85]; we prefer to give only a sketch of such argument. The definition of the two functions $\Phi, \Psi$ turns the verification that the two form a Galois connection into a straightforward check, and all the other main steps of the proofs are resumed in the following remarks. $\square$

Remark 2.1.2. The action of the functor $R\colon \mathbf{C} \to \mathcal{M}/1$ is induced on objects by a choice of $\mathbb{F}$-factorizations of terminal morphisms: $X \xrightarrow{e} RX \xrightarrow{m} 1$. On arrows it is obtained from a choice of solutions to lifting problems

$$\begin{array}{ccc} A & \xrightarrow{ef} & RB \\ \downarrow & \overset{Rf}{\nearrow} & \downarrow{\scriptstyle m} \\ RA & \xrightarrow[m]{} & 1 \end{array} \tag{2.2}$$

Remark 2.1.3. Showing that there is an adjunction $R\colon \mathbf{C} \rightleftarrows \mathcal{M}/1\colon i$ boils down to showing that $\mathbf{C}(-, X)$ inverts each reflection $A \to RA$; this is an easy consequence of the arrow-orthogonality between $\left[\begin{smallmatrix} A \\ \downarrow \\ RA \end{smallmatrix}\right]$ and $\left[\begin{smallmatrix} X \\ \downarrow \\ 1 \end{smallmatrix}\right]$, equivalent to the object-orthogonality on $\left[\begin{smallmatrix} A \\ \downarrow \\ RA \end{smallmatrix}\right]$ and $X \in \mathcal{M}/1$.

Remark 2.1.4. The unit $\mathrm{id}_{\mathrm{Rex}(\mathbf{C})} \Rightarrow \Psi\Phi$ of this adjunction is an isomorphism. The comonad $\Phi\Psi \Rightarrow \mathrm{id}_{\mathrm{PF}_{\mathcal{T}}(\mathbf{C})}$ is much more interesting, as it acts like an *interior operator* on the poset $\mathrm{PF}_{\mathcal{T}}(\mathbf{C})$, sending $\mathbb{F}$ to a new prefactorization $\mathbb{F}^\circ = (\mathcal{E}^\circ, \mathcal{M}^\circ)$ which is by construction *reflective*, i.e. satisfies $\mathbb{F}^\circ = \mathbb{F}$ (whereas in general we have only a proper inclusion $\mathbb{F}^\circ \preceq \mathbb{F}$ deduced from $\mathcal{M}^\circ \subseteq \mathcal{M}$).

What we said so far entails that

PROPOSITION 2.1.5. The adjunction $\Phi \dashv \Psi$ restricts to an equivalence (a bijection between posets) between the reflective prefactorizations in $\mathbb{F} \in \text{PF}_{\boldsymbol{\tau}}(\mathbf{C})$ and the poset $\text{Rex}(\mathbf{C})$.

PROPOSITION 2.1.6. $\mathbb{F} \in \text{PF}_{\boldsymbol{\tau}}(\mathbf{C})$ is reflective if and only if $\mathcal{E}$ is a 3-for-2 class (see Definition **1.4.6**), or equivalently (since each $\mathcal{E}$-class of a factorization system is R32) if and only if $\mathcal{E}$ is L32.

*Proof.* It is an immediate consequence of [CHK85, Thm. **2.3**], where it is stated that $g \in \mathcal{E}^\circ$ iff $fg \in \mathcal{E}$ for a suitable $f \in \mathcal{E}$. □

REMARK 2.1.7. We can also state a completely dual antitone bijection between the poset of *co*reflective subcategories, $\text{CoRex}(\mathbf{C})$, and the poset of (pre)factorization systems $\text{PF}_{\boldsymbol{\iota}}(\mathbf{C})$ factoring *initial* arrows $\boldsymbol{\iota} = \{\varnothing \to X \mid X \in \mathbf{C}\}$; this is defined via the correspondence $\mathbb{F} \mapsto \varnothing/\mathcal{E} = \{Y \in \mathbf{C} \mid (\varnothing \to Y) \in \mathcal{E}\}$; the coreflection of $\mathbf{C}$ along $\varnothing/\mathcal{E}$ is given by a functor $S$ defined by a choice of $\mathbb{F}$-factorization $\varnothing \xrightarrow{e} SX \xrightarrow{m} X$.

REMARK 2.1.8. We can also define *co*reflective factorization systems, and prove that $\mathbb{F}$ is coreflective iff $\mathcal{M}$ is R32, and *bi*reflective factorization systems as those which are reflective *and* coreflective at the same time: as these will consistute the main object of study of the present and subsequent chapters, we gather these remarks into a precise definition.

DEFINITION 2.1.9. (REFLECTIVE FACTORIZATION SYSTEM): A *bireflective (pre)factorization system* $\mathbb{F} = (\mathcal{E}, \mathcal{M}) \in \text{PF}(\mathbf{C})$ is a (pre)factorization system such that both classes $\mathcal{E}, \mathcal{M}$ are 3-for-2 classes.

## 2.2 Semiexactness and simplicity.

> The guiding motto in the life of every natural philosopher should be, seek simplicity and distrust it.
>
> A.N. Whitehead

A fairly general theory, subsumed in [CHK85], stems from the above construction, and several notable subclasses of (co)reflective factorization systems become of interest. We now concentrate on *semi-exact* and *simple* factorizations:

DEFINITION 2.2.1. A *semi-left-exact* factorization system on a finitely complete $\mathbf{C}$ consists of a reflective $\mathbb{F} = (\mathcal{E}, \mathcal{M}) \in \text{FS}(\mathbf{C})$ such that the left class $\mathcal{E}$ is closed under pulling back by $\mathcal{M}$ arrows; more explicitly, in the pullback

$$\begin{array}{ccc} A & \longrightarrow & B \\ {\scriptstyle e'}\downarrow & \lrcorner & \downarrow{\scriptstyle e} \\ C & \xrightarrow[m]{} & D \end{array} \tag{2.3}$$

the arrow $e'$ lies in $\mathcal{E}$.

Equivalent conditions for $\mathbb{F}$ to be semi-left-exact are given in [CHK85, Thm. **4.3**]. There is a dual definition of a semi-*right*-exact factorization system.

Notation 2.2.2. We call *semiexact* a factorization system which is both left and right exact.

Another important class of factorization systems is made by *simple* ones in categories with finite limits and colimits, where $\mathbb{F}$ gives "a simple rule to factor morphisms". More precisely, if $\mathbf{C}$ has pullbacks, we can define

Definition 2.2.3. A *left simple* factorization system on $\mathbf{C}$ is a reflective $\mathbb{F} \in \text{FS}(\mathbf{C})$ such that, if we denote by $R$ the reflection $\mathbf{C} \to \mathcal{M}/1$, with unit $\eta\colon 1_{\mathbf{C}} \Rightarrow iR$ (often denoted simply as $\eta\colon 1_{\mathbf{C}} \Rightarrow R$ with a harmless abuse of notation), associated to $\mathbb{F}$, then the $\mathbb{F}$-factorization of $f\colon X \to Y$ can be obtained as $X \to RX \times_{RY} Y \to Y$ in the diagram

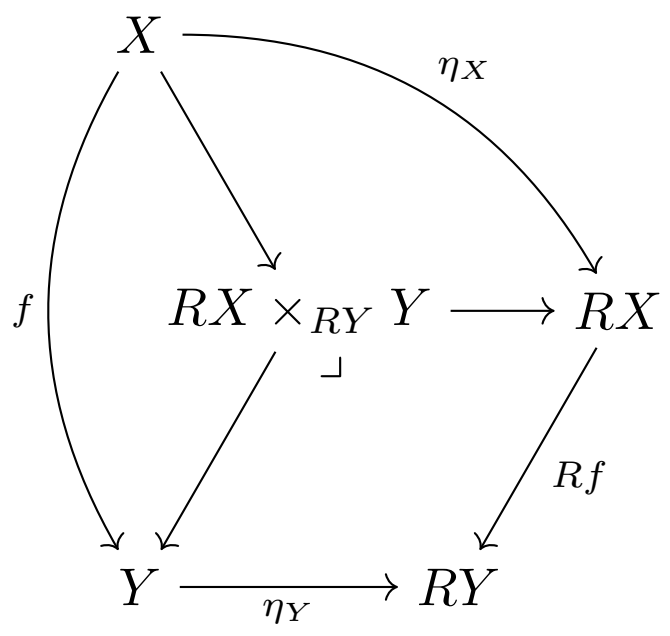

(2.4)

obtained from the naturality square for $f$.

Simple factorization systems are, in other words, those such that the canonical arrow $X \to RX \times_{RY} Y$ lies in $\mathcal{E}$ (the pullback arrow $RX \times_{RY} Y \to Y$ always lies in $\mathcal{M}$, by the 3-for-2 closure property of $\mathcal{M}$).

Remark 2.2.4. Every semi-left-exact factorization system is left simple, as proved in [CHK85, Thm. **4.3**]. In the 1-categorical setting, the converse doesn't hold in general (see [CHK85, Example **4.4**]), whereas our Prop. **2.3.15** shows that in the stable $\infty$-categorical world the two notions coincide. This is a first evidence of the notable and really symmetric "internal behaviour" of a stable $\infty$-category (the proof of our **2.3.15** makes essential use of the *pullout axiom*, which is only valid and nontrivial in an $\infty$-categorical setting).

Remark 2.2.5. There is an analogous notion of *right simple* factorization system: it is enough to dualize Def. **2.2.3**; dualizing also [CHK85, Thm. **4.3**], it is possible to prove that semi-right-exact factorization systems are right simple.

A useful result follows from the semi-exactness of a factorization system $\mathbb{F}$ both of whose classes are 3-for-2 (these last are called *torsion theories* in [RT07]; see our Def. **2.3.1** for an extensive discussion).

PROPOSITION 2.2.6. Let $\mathbb{F}$ be a semiexact (Def. **2.2.1**; its domain of definition is, in particular, finitely co/complete) torsion theory with reflection functor $R\colon \mathbf{C} \to \mathcal{M}/1$ and whose coreflection is $S$; then we have that

$$SY \amalg_{SX} X \cong RX \times_{RY} Y \qquad (2.5)$$

for any $f\colon X \to Y$.

*Proof.* The claim holds simply because semiexactness gives the $\mathbb{F}$-factorization of $f\colon X \to Y$ as $X \to RX \times_{RY} Y \to Y$ (on the left), and $X \to SY \amalg_{SX} X \to Y$ (on the right).

There is a more explicit argument which makes explicit use of the orthogonality and 3-for-2 closure property: consider the diagram

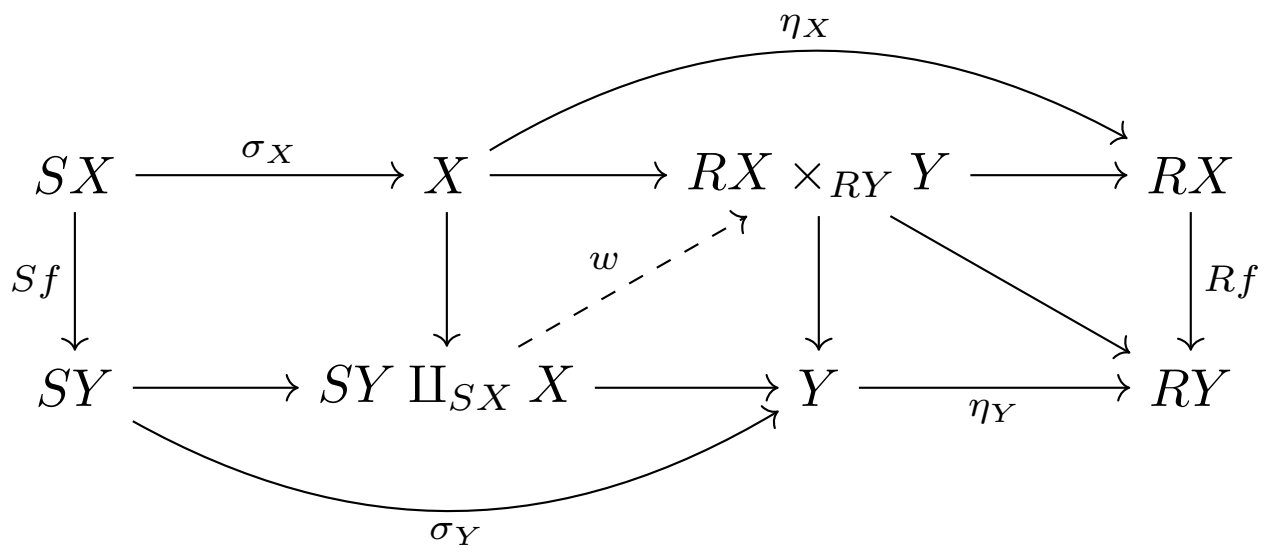

(2.6)

where $\eta$ is the unit of the reflection $R$, $\sigma$ is the counit of the coreflection $S$, and the diagonal of the central square is filled by $f\colon X \to Y$. Now, denote $P = RX \times_{RY} Y$ and $Q = SY \amalg_{SX} X$ the arrow $\begin{bmatrix} X \\ \downarrow \\ Q \end{bmatrix}$ is in $\mathcal{E}$, and the arrow $\begin{bmatrix} P \\ \downarrow \\ Y \end{bmatrix}$ is in $\mathcal{M}$, as a consequence of stability under cobase and base change (see Prop. **1.4.3**); this entails that there is a unique $w\colon Q \to P$ making the central square commute. Now, semiexactness entails that $X \to P \to Y$ and $X \to Q \to Y$ are both $\mathbb{F}$-factorizations of $f\colon X \to Y$, and since both classes $\mathcal{E}, \mathcal{M}$ are 3-for-2, we can now conclude that $w\colon Q \to P$ lies in $\mathcal{E} \cap \mathcal{M}$, and hence is an equivalence (see Prop. **1.4.5**). $\square$

## 2.3 Normal torsion theories.

> We have normality. I repeat: we have normality. Anything you still can't cope with is therefore your own problem.
>
> D. Adams

Refining the blanket assumption of the initial section, we now assume that $\mathbf{C}$ is a stable $\infty$-category, with zero object $0 = \varnothing = *$. Following (and slightly adapting to our particular case) [RT07, §4] we give the following definitions

DEFINITION 2.3.1. (TORSION THEORY, TORSION CLASSES): A *torsion theory* in **C** consists of a factorization system $\mathbb{F} = (\mathcal{E}, \mathcal{M})$ (see Remark **2.1.8** and Def. **2.1.9**), where both classes are 3-for-2 (in the sense of Definition **1.4.6**). We define $\mathcal{T}(\mathbb{F}) = 0/\mathcal{E}$ and $\mathcal{F}(\mathbb{F}) = \mathcal{M}/0$ (see Prop. **2.1.1**, and Remark **2.1.7**) to be respectively the *torsion* and *torsion-free* classes associated to the torsion theory.

REMARK 2.3.2. [RT07, **3.1**] Let **C** be an $\infty$-category with terminal object $*$; then the class $\mathcal{F}(\mathbb{F})$ is *firmly $\mathcal{E}$-reflective*, meaning that any morphism $A \to F$ with $F \in \mathcal{F}(\mathbb{F})$ is isomorphic to the reflection $A \to RA$. This directly follows from the uniqueness of the $\mathbb{F}$-factorization.

REMARK 2.3.3. In view of Prop. **2.1.6** and its dual, the torsion and torsion-free classes of a torsion theory $\mathbb{F} \in \text{FS}(\mathbf{C})$ are respectively a coreflective and a reflective subcategory of **C**.

If we $\mathbb{F}$-factor the terminal and initial morphisms of any object $X \in \mathbf{C}$, we obtain the reflection $R\colon \mathbf{C} \to \mathcal{M}/0$ and coreflection $S\colon \mathbf{C} \to 0/\mathcal{E}$, and a "complex"

$$SX \to X \to RX \tag{2.7}$$

(in the sense of pointed categories), i.e., a homotopy commutative diagram

$$\begin{array}{ccc} SX & \longrightarrow & X \\ \downarrow & & \downarrow \\ 0 & \longrightarrow & RX \end{array} \tag{2.8}$$

We deduce this commutativity from the orthogonality condition: the lifting problem

$$\begin{array}{ccc} 0 & \longrightarrow & RX \\ \downarrow & & \downarrow \\ SX & \longrightarrow & 0 \end{array} \tag{2.9}$$

has unique solution the zero arrow $SX \to RX$, so that the space $\mathbf{C}(SX, RX)$ is contractible: since there cannot be nonzero arrows $SX \to RX$, the claim is proved.

PROPOSITION 2.3.4. Let **C** be a stable $\infty$-category with a normal torsion theory $\mathbb{F} = (\mathcal{E}, \mathcal{M})$, having coreflection $S\colon \mathbf{C} \to 0/\mathcal{E}$. Then the following conditions are equivalent for an object $X \in \mathbf{C}$:

(1) $X$ is an $S$-coalgebra, i.e. there exists an arrow $c\colon X \to SX$ such that $SX \xrightarrow{\sigma_X} X \xrightarrow{c} SX$ is the identity of $SX$;

(2) $X \in \mathcal{T} = 0/\mathcal{E}$;
(3) $X \cong SX$;
(4) $X \in {}^{\perp}\{SA \to A\}$, i.e. $X$ is left-object-orthogonal (Def. **1.2.15**) to each coreflection arrow $SA \to A$.

The present statement results from a mixture of [RT07] and [Kel80, Prop. 5.2].

Obviously, a dual result can be stated and proved with basically no effort:

PROPOSITION 2.3.5. Let $\mathbf{C}$ be a stable $\infty$-category with a normal torsion theory $\mathbb{F} = (\mathcal{E}, \mathcal{M})$. Then the following conditions are equivalent for an object $X \in \mathbf{C}$

(1) $X$ is an $R$-algebra;
(2) $X \in \mathcal{M}/0$;
(3) $X \cong RX$;
(4) $X \in \{A \to RA\}^{\perp}$, i.e. $X$ is object-orthogonal (Def. **1.2.15**) to each reflection arrow $A \to RA$.

REMARK 2.3.6. Given the closure properties of the classes $\mathcal{E}, \mathcal{M}$, we can define natural functors $F\colon \mathbf{C} \to \mathcal{F}$ and $T\colon \mathbf{C} \to \mathcal{T}$ taking $FX$ as the homotopy pullback, and $TX$ as the homotopy pushout in the diagrams below

$$\begin{array}{ccc} FX & \longrightarrow & SX \\ \downarrow & \lrcorner & \downarrow \\ 0 & \longrightarrow & X \end{array} \qquad \begin{array}{ccc} X & \longrightarrow & 0 \\ \downarrow & \ulcorner & \downarrow \\ RX & \longrightarrow & TX. \end{array} \tag{2.10}$$

We now come to the gist of the present chapter, i.e. the definition of a *normal* torsion theory and its relation with the notion of $t$-structures, which will occupy entirely Chapter **3** with the proof of the *Rosetta stone* theorem.

An initial step to motivate the quest for a class of factorization system describing $t$-structures (identified with the pair of subcategories called *aisle* and *coaisle* in the literature, see [KV88]) in stable $\infty$-categories starts precisely from the observation that suitable additional properties of a co/-reflective subcategory $\mathbf{B} \subseteq \mathbf{C}$ translate into properties of the associated co/reflective factorization system $\Phi(\mathbf{B}) = \mathbb{F}$.

Torsion theories in a stable $\infty$-category, in the form of bireflective factorization systems, produce such pairs of well-behaved coreflective/reflective subcategories via the correspondence $(\mathcal{E}, \mathcal{M}) \mapsto (0/\mathcal{E}, \mathcal{M}/0)$; so we are only one step away from characterizing $t$-structures: we only lack axiom (iii) of Def. **A.3.2**.

It turns out that the possibility of putting every object $X$ into a distin-guished triangle (or, better to say in our setting, a fiber sequence)

$$\begin{array}{ccc} X_{\geq 0} & \longrightarrow & X \\ \downarrow & \lrcorner \quad \ulcorner & \downarrow \\ 0 & \longrightarrow & X_{<0} \end{array} \tag{2.11}$$

is equivalent to the request that $(\mathcal{E}, \mathcal{M})$ be a *normal* factorization system on $\mathbf{C}$; in a nutshell, the idea is the following.

General torsion theories generate a sequence $SX \to X \to RX$ whose composition is the zero morphism; the factorization systems rendering this composition also an *exact sequence* are called *normal* (the term is borrowed from [CHK85] who first studied the notion, reprised in [RT07]).

> A normal torsion theory is a factorization system $\mathbb{F} = (\mathcal{E}, \mathcal{M})$ such that the diagram
>
> $$\begin{array}{ccc} SX & \longrightarrow & X \\ \downarrow & \lrcorner \quad \ulcorner & \downarrow \\ 0 & \longrightarrow & RX \end{array} \tag{2.12}$$
>
> is a pullout.

As discussed above, it is fairly natural to define functors $F$ and $T$ taking respectively the *fiber of the coreflection* and the *cofiber of the reflection* morphism. Normality involves the alternate procedure, considering the fiber $KX$ of the reflection $X \to RX$ and the cofiber $QX$ of the coreflection $SX \to X$. A priori, there is no way to control the subcategory where the functors $K, Q$ take value: the idea behind a normal torsion theory is that in certain situation this is possible, as the two functors $K$ and $Q$ do not introduce new information, as they are respectively isomorphic to $S$ and $R$.

REMARK 2.3.7. This terse characterization of normality, and especially our Remark **2.3.16** which states that left, right and two-sided normality all coincide in a stable $\infty$-category, seems to shed a light on [CHK85, Remark **7.8**] and [RT07, Remark **4.11**], where the non-existence of a non-artificial example of a non-normal torsion theory is conjectured.

REMARK 2.3.8. The present analysis owes to [RT07, CHK85, BR07] an infinite debt; it may appear strange, hence, that such many different sources ignore the possibility of turning this suggestion into a precise statement.

Indeed, somehow mysteriously, [RT07, §**4**] seems to ignore application of the formalism of torsion theories to the triangulated world, even if its authors point out clearly (see [RT07, Remark **4.11.(2)**]) that

> It [our definition of torsion theory, *auth.*] applies, for example, to a triangulated category $\mathbf{C}$. Such a category has only weak kernels and weak cokernels and our definition precisely corresponds to torsion theories considered there as pairs $\mathcal{F}$ and $\mathcal{T}$ of colocalizing and localizing subcategories (see [HPS97]).

Even more mysteriously, another encyclopedic source for a "calculus of torsion theories" in triangulated categories, [BR07], explicitly says (p. 17) that

> Torsion pairs in triangulated categories are used in the literature mainly in the form of $t$-structures.

and yet it avoids, in a certain sense, to offer a more primitive characterization for $t$-structures than the one given *ibi*, Thm **2.13**.

This situation indicates well a general tenet according to which working in the stable setting gives more symmetric and better motivated results.

The "Rosetta stone" theorem casts a shadow on the homotopy category $\mathbf{T} = \mathrm{Ho}(\mathbf{C})$, giving a similar but insufficient characterization of $t$-structures as those factorization systems in $\mathbf{T}$ which are closed under homotopy pullback and pushouts in $\mathbf{T}$.[(2)]

DEFINITION 2.3.9. We call *left normal* a torsion theory $\mathbb{F} = (\mathcal{E}, \mathcal{M})$ on $\mathbf{C}$ such that the fiber $KX \to 0$ of a reflection morphism $X \to RX$ lies in $\mathcal{E}$, as in the diagram

$$\begin{array}{ccc} KX & \longrightarrow & X \\ \big\downarrow & \lrcorner & \big\downarrow \\ 0 & \longrightarrow & RX \end{array} \tag{2.13}$$

In other words, the $\mathcal{E}$-morphisms arising as components of the unit $\eta\colon 1 \Rightarrow R$ are stable under pullback along the initial $\mathcal{M}$-morphism $0 \to RX$.

REMARK 2.3.10. This last sentence deserves a deeper analysis: by the very definition of $RX$ it is clear that $RX \to 0$ lies in $\mathcal{M}$; but more is true (and this seemingly innocuous result is a key step of most of the proofs we are going to present): since $\mathcal{M}$ enjoys the 3-for-2 property, and it contains all isomorphisms of $\mathbf{C}$, it follows immediately that an initial arrow $0 \to A$ lies in $\mathcal{M}$ *if and only if* the terminal arrow $A \to 0$ on the same object lies in $\mathcal{M}$. The same reasoning applied to $\mathcal{E}$ gives a rather notable "specularity" property for both classes $\mathcal{E}, \mathcal{M}$:

LEMMA 2.3.11. (SATOR LEMMA): In a pointed $\infty$-category $\mathbf{C}$, an initial arrow $0 \to A$ lies in a class $\mathcal{E}$ or $\mathcal{M}$ of a bireflective (see Remark **2.1.8**)

(2) This result is part of a work in progress [?] and will hopefully introduce a subsequent joint work exploring the shape of the "Rosetta stone" in the setting of *stable derivators*.

factorization system $\mathbb{F}$ if and only if the terminal arrow $A \to 0$ lies in the same class.(3)

Notation 2.3.12. This motivates a little abuse of notation: we can say that an object $A$ of $\mathbf{C}$ *lies in* a 3-for-2 class $\mathcal{K}$ if its initial or terminal arrow lies in $\mathcal{K}$: in this sense, a left normal factorization system is an $\mathbb{F}$ such that the fiber $KX$ of $X \to RX$ lies in $\mathcal{E}$, for every $X$ in $\mathbf{C}$.

Equivalent conditions for $\mathbb{F}$ to be left normal are given in [RT07, Thm. **4.10**] and [CHK85, **7.3**].

Remark 2.3.13. There is, obviously, a notion of *right* normal factorization system: it is an $\mathbb{F}$ such that the cofiber $QX$ of $SX \to X$ lies in $\mathcal{M}$, for every $X$ in $\mathbf{C}$. In the following we call simply *normal*, or *two-sided normal* a factorization system $\mathbb{F} \in \text{FS}(\mathbf{C})$ which is both left and right normal.

Now we come to an interesting point: in a stable $\infty$-category the three notions of simple, semiexact and normal torsion theory collapse to be three equivalent conditions.

To see this, we have to prove a preliminary result:

Proposition 2.3.14. For every object $X$, consider the following diagram in $\mathbf{C}$, where every square is a pullout.

$$\begin{array}{ccccc} SX \oplus RX[-1] & \longrightarrow & SX & \longrightarrow & 0 \\ {\scriptstyle m''}\downarrow & & \downarrow{\scriptstyle \sigma_X} & & \downarrow \\ KX & \longrightarrow & X & \longrightarrow & QX \\ \downarrow & & \downarrow{\scriptstyle \rho_X} & & \downarrow{\scriptstyle e''} \\ 0 & \longrightarrow & RX & \longrightarrow & SX[1] \oplus RX \end{array} \tag{2.15}$$

Then the following conditions are equivalent for a bireflective factorization system $\mathbb{F} = (\mathcal{E}, \mathcal{M})$ on $\mathbf{C}$:

(1) $\mathbb{F}$ is left normal;
(2) $\mathbb{F}$ is right normal;
(3) $\mathbb{F}$ is normal;
(4) $RX \simeq QX$;

(3)The so-called *Sator square*, first found in the ruins of Pompeii, consists of the $5 \times 5$ matrix

$$\begin{array}{ccccc} \text{S} & \text{A} & \text{T} & \text{O} & \text{R} \\ \text{A} & \text{R} & \text{E} & \text{P} & \text{O} \\ \text{T} & \text{E} & \text{N} & \text{E} & \text{T} \\ \text{O} & \text{P} & \text{E} & \text{R} & \text{A} \\ \text{R} & \text{O} & \text{T} & \text{A} & \text{S} \end{array} \tag{2.14}$$

where the letters are arranged in such a way that the same phrase ("sator arepo tenet opera rotas", approximately "Arepo, the farmer, drives carefully the plough") appears when it is read top-to-bottom, bottom-to-top, left-to-right, and right-to-left.

(5) $SX = KX$;
(6) $SX \to X \to RX$ is a fiber sequence.

*Proof.* We start by proving that the first three conditions are equivalent. If we assume left normality, then the arrow $\begin{bmatrix} QX \\ \downarrow \\ SX[1]\oplus RX \end{bmatrix}$ lies in $\mathcal{E}$, since it results as a pushout of an arrow in $\mathcal{E}$. So we can consider

$$\begin{array}{ccc} QX & \xrightarrow{\quad e' \quad} & RQX \\ \Big\downarrow{\scriptstyle e''} & & \Big\downarrow{\scriptstyle m'} \\ SX[1]\oplus RX \xrightarrow[e]{} R(SX[1]\oplus RX) & \xrightarrow[m]{} & 0 \end{array} \tag{2.16}$$

$\mathbb{F}$-factoring the morphisms involved (notice that $R(SX[1] \oplus RX) \cong RX$): $R(SX[1] \oplus RX) = RRX = RX$ since $RS = 0$. Thus $RQX \cong RX$, which entails $\begin{bmatrix} 0 \\ \downarrow \\ QX \end{bmatrix} \in \mathcal{M}$, which entails right normality. A dual proof gives that (2) $\Rightarrow$ (1), thus right normality equals left normality and hence two-sided normality. Now it is obvious that (6) is equivalent to (4) and (5) together; the non-trivial part of the proof consists of the implications (1) $\Rightarrow$ (4), and dually (2) $\Rightarrow$ (5).

Once this is noticed, start with the diagram

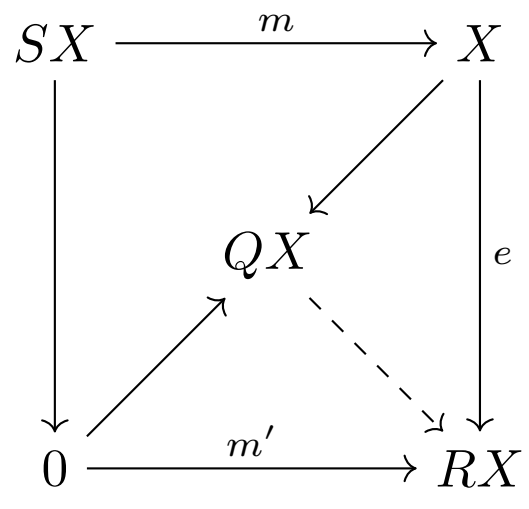

(2.17)

and consider the canonical arrow $QX \to RX$ obtained by universal property: the arrow $\begin{bmatrix} 0 \\ \downarrow \\ RX \end{bmatrix}$ lies in $\mathcal{M}$ (this is a general fact); left normality now entails that $\begin{bmatrix} 0 \\ \downarrow \\ QX \end{bmatrix} \in \mathcal{M}$, so that $\begin{bmatrix} QX \\ \downarrow \\ RX \end{bmatrix}$ lies in $\mathcal{M}$ too by reflectivity.

A similar argument shows that since both $\begin{bmatrix} X \\ \downarrow \\ QX \end{bmatrix}$, $\begin{bmatrix} X \\ \downarrow \\ RX \end{bmatrix}$ lie in $\mathcal{E}$, $\begin{bmatrix} QX \\ \downarrow \\ RX \end{bmatrix}$ lies in $\mathcal{E}$ too by reflectivity. This entails that $\begin{bmatrix} QX \\ \downarrow \\ RX \end{bmatrix}$ is an equivalence. Conversely, if we start supposing that $QX \cong RX$, then we have (left) normality. This concludes the proof, since in the end we are left with the equality (4) $\iff$ (5). □

As previewed before, the three notions of simplicity, semiexactness and normality collapse in a single notion in the stable setting:

Proposition 2.3.15. A torsion theory $\mathbb{F}$ is left normal if and only it is semi-left-exact in the sense of [CHK85, **4.3.i**], namely if and only if in the pullout square

$$\begin{array}{ccc} E & \longrightarrow & X \\ {\scriptstyle e'}\downarrow & & \downarrow{\scriptstyle \rho_X} \\ Q & \xrightarrow[m]{} & RX \end{array} \tag{2.18}$$

the arrow $e'$ lies in $\mathcal{E}$. Dually, a factorization system $\mathbb{F}$ is right normal if and only it is semi-right-exact in the sense of (the dual of) [CHK85, **4.3.i**].

*Proof.* Consider the diagram

$$\begin{array}{ccccc} KX & \longrightarrow & E & \longrightarrow & X \\ \downarrow & & \downarrow{\scriptstyle e'} & & \downarrow{\scriptstyle e} \\ 0 & \longrightarrow & Q & \xrightarrow{m} & RX \end{array} \tag{2.19}$$

where the arrow $Q \to RX$ belongs to $\mathcal{M}$. On the one hand it is obvious that if $\mathbb{F}$ is semi-left-exact, then it is normal (just pull back two times $e$ along $\mathcal{M}$-arrows). On the other hand, the converse implication relies on the pullout axiom: if $\mathbb{F}$ is normal, then $KX$ lies in $\mathcal{E}$; but now since the left square is a pullout, the arrow $\left[\begin{smallmatrix} E \\ \downarrow \\ Q \end{smallmatrix}\right]$ belongs to $\mathcal{E}$ too, giving semi-left-exactness. □

Remark 2.3.16. The three notions coincide since "classically" we have

$$\text{SLEX} \to \text{SIMPLE} \to \text{NORMAL}, \tag{2.20}$$

whereas in our setting the chain of implication proceeds one step further and closes the circle:

$$\text{SLEX} \to \text{SIMPLE} \to \text{NORMAL} \xrightarrow{\star} \text{SLEX}. \tag{2.21}$$

This gives a pleasant consequence:

Remark 2.3.17. In a stable $\infty$-category the $\mathbb{F}$-factorization of $f\colon A \to B$ with respect to a normal torsion theory is always

$$A \to RA \times_{RB} B \to B, \tag{2.22}$$

or equivalently (see Prop. **2.2.6**)

$$A \to SB \amalg_{SA} A \to B. \tag{2.23}$$

A useful remark appearing in [RT07, §**4.6**, **5**] (here adapted to the stable case) is the following: torsion and torsionfree classes of a torsion theory in a stable $\infty$-category are closed under extensions.

DEFINITION 2.3.18. Let $\mathcal{K} \subseteq \mathrm{Ob}(\mathbf{C})$ be a class of objects in a stable $\infty$-category; $\mathcal{K}$ is said to be *closed under extensions* if for each pullout square

$$\begin{array}{ccc} A & \longrightarrow & B \\ \downarrow & & \downarrow \\ 0 & \longrightarrow & C \end{array} \tag{2.24}$$

such that $A, C \in \mathcal{K}$, then also $B \in \mathcal{K}$.

PROPOSITION 2.3.19. Let $\mathbb{F} = (\mathcal{E}, \mathcal{M})$ be a torsion theory in a stable $\infty$-category $\mathbf{C}$; then the classes $0/\mathcal{E}$, $\mathcal{M}/0$ of Def. **2.3.1** are closed under extension.

*Proof.* We only prove the statement if $A, C$ of diagram (**2.24**) lie in $0/\mathcal{E}$; the proof for $\mathcal{M}/0$ is identical. Now, it is enough to consider the diagram

$$\begin{array}{ccccc} 0 & \longrightarrow & A & \longrightarrow & B \\ & & \downarrow & & \downarrow \\ & & 0 & \longrightarrow & C \\ & & & & \downarrow \\ & & & & 0 \end{array} \tag{2.25}$$

where we have $A, C \in 0/\mathcal{E}$, i.e. $A \to 0, C \to 0$ lie in $\mathcal{E}$; the arrow $B \to C$ is in $\mathcal{E}$ since $\mathcal{E}$ is closed under pushout; so $B \to C \to 0$ is in $\mathcal{E}$. $\square$

# Chapter 3

# The "Rosetta stone"

## 3.1 $t$-structures are factorization systems.

> Acaso un arquetipo no revelado aún a los hombres, un objeto eterno (para usar la nomenclatura de Whitehead), esté ingresando paulatinamente en el mundo; su primera manifestación fue el palacio; la segunda el poema. Quien los hubiera comparado habría visto que eran esencialmente iguales.
>
> [Bor97], *El sueño de Coleridge*

This is (both form a conceptual and order-theoretical point of view) the central chapter of the thesis, where we prove our main result: we gathered enough material and mastery of the *iaidō* of factorization to embark on a complete, exhaustive proof of our "Rosetta stone" theorem **3.1.1**, i.e. to prove that in a stable quasicategory, normal torsion theories correspond to $t$-structures, via the following dictionary.

| **Normal torsion theories** | **$t$-structures** |
|---|---|
| $\mathbb{F} = (\mathcal{E}, \mathcal{M})$ | $\mathfrak{t}$ |
| $(\mathcal{T}, \mathcal{F})$ | $(\mathbf{C}_{\geq 0}, \mathbf{C}_{<0})$ |
| $\hom(\mathcal{T}, \mathcal{F}) \simeq *$ | $\hom(\mathbf{C}_{\geq 0}, \mathbf{C}_{<0}) = 0$ |
| factorization of initial/terminal | reflection/coreflection functors |

We provide an introduction to $t$-structures in **A.3**; the interested reader can also consult classical references as [KS, BBD82] and the unique (at the moment of writing) reference for $t$-structures in stable $\infty$-categories, [Lur17].

In some sense, the present result, which turned out to be the main conceptual achievement of the present work, arose from the innocuous desire to better understand [Lur17, **1.2.1.4**], which defines $t$-structures on a stable $\infty$-category $\mathbf{C}$ as classical $t$-structures on the homotopy category $\mathrm{Ho}(\mathbf{C})$.

Albeit true, this result seems to hide part of the story. A deeper analysis of it, motivated by the desire for a more intrinsic characterization of $t$-structures, motivated the following statement:

THEOREM 3.1.1. (THE ROSETTA STONE): Let $\mathbf{C}$ be a stable $\infty$-category. There is a bijective correspondence between the class of normal torsion theories $\mathbb{F} = (\mathcal{E}, \mathcal{M})$ on $\mathbf{C}$ (in the sense of Definition **2.3.9**) and the class of $t$-structures on $\mathbf{C}$ (in the sense of Definition **A.3.2**).

The proof of this result will occupy the entire chapter, and will be followed by examples coming from homological algebra and algebraic topology, showing how to reinterpret classical constructions in light of this result.

To simplify the discussion we will deduce **3.1.1** as a consequence of a number of separate statements.

The strategy is simple: we first construct the pair of correspondences

$$\text{normal torsion theories} \underset{\mathbb{F}(-)}{\overset{\mathfrak{t}(-)}{\rightleftarrows}} t\text{-structures} \tag{3.1}$$

We are obviously led to exploit the fundamental connection outlined in §**2.1**:

- Given a normal, bireflective factorization system $\mathbb{F} = (\mathcal{E}, \mathcal{M})$ on $\mathbf{C}$ we define the two classes $(\mathbf{C}_{\geq 0}(\mathbb{F}), \mathbf{C}_{<0}(\mathbb{F}))$ of the $t$-structure $\mathfrak{t}(\mathbb{F})$ to be the torsion and torsion-free classes $(0/\mathcal{E}, \mathcal{M}/0)$ associated to $\mathbb{F}$, in the sense of Definition **2.3.1**.
- On the other hand, given a $t$-structure $\mathfrak{t} = (\mathbf{C}_{\geq 0}, \mathbf{C}_{<0})$ in the sense of Definition **A.3.2**, we have to define classes $\mathbb{F}(\mathfrak{t}) = (\mathcal{E}(\mathfrak{t}), \mathcal{M}(\mathfrak{t}))$ which form a factorization system. If $\tau_{\geq 0}, \tau_{<0}$ denote, respectively, the co/-truncation of the $t$-structures (Remark **A.3.4**), we set:

$$\begin{aligned} \mathcal{E}(\mathfrak{t}) &= \{f \in \mathbf{C}^{\Delta[1]} \mid \tau_{<0}(f)\text{equivalence an is } \}; \\ \mathcal{M}(\mathfrak{t}) &= \{f \in \mathbf{C}^{\Delta[1]} \mid \tau_{\geq 0}(f)\text{equivalence an is } \}. \end{aligned} \tag{3.2}$$

The language developed throughout the previous chapter will give a manageable (in fact, several) characterizations of these two classes of morphisms.

Half of the proof for Thm. **3.1.1** consists in a mere recasting of the definition of normal torsion theory, to check that the pair $(\mathbf{C}_{\geq 0}(\mathbb{F}), \mathbf{C}_{<0}(\mathbb{F}))$ really is a $t$-structure:

PROPOSITION 3.1.2. The pair $\mathfrak{t}(\mathbb{F})$ is a $t$-structure on $\mathbf{C}$ in the sense of Definition **A.3.2**.

*Proof.* The orthogonality condition is immediate by definition of the two classes (see Remark **2.3.3**). As for the closure under positive/negative shifts, $(A \to B) \in \mathcal{E}$ entails that $(A[1] \to B[1]) \in \mathcal{E}$ since left classes in factorization systems are closed under (homotopy) colimits in the arrow category (see Prop. **1.4.7**) and in particular under the homotopy pushout defining the shift $A \mapsto A[1]$ on $\mathbf{C}$. This justifies the chain of implications

$$X \in \mathbf{C}_{\geq 0}(\mathbb{F}) \iff \begin{bmatrix} 0 \\ \downarrow \\ X \end{bmatrix} \in \mathcal{E} \implies \begin{bmatrix} 0 \\ \downarrow \\ X[1] \end{bmatrix} \in \mathcal{E} \iff X[1] \in \mathbf{C}_{\geq 0}(\mathbb{F}). \tag{3.3}$$

The case of $\mathbf{C}_{<0}$ is completely dual: since $\mathcal{M}$ admits any limit, $\begin{bmatrix} X \\ \downarrow \\ 0 \end{bmatrix} \in \mathcal{M}$ implies that $\begin{bmatrix} X[-1] \\ \downarrow \\ 0 \end{bmatrix} \in \mathcal{M}$, so that $\mathbf{C}_{<0}(\mathbb{F})[-1] \subset \mathbf{C}_{<0}(\mathbb{F})$.

To see that any object $X \in \mathbf{C}$ fits into a fiber sequence $X_{\geq 0} \to X \to X_{<0}$, with $X_{\geq 0}$ in $\mathbf{C}_{\geq 0}(\mathbb{F})$ and $X_{<0}$ in $\mathbf{C}_{<0}(\mathbb{F})$, it suffices to $\mathbb{F}$-factor the terminal morphism of $X$ obtaining a diagram like

$$X \xrightarrow{e} RX \xrightarrow{m} 0 \tag{3.4}$$

and then to take the fiber of $e$,

$$\begin{array}{ccc} KX & \longrightarrow & X \\ \downarrow & \lrcorner & \downarrow \\ 0 & \ulcorner\longrightarrow & RX \end{array} \tag{3.5}$$

Set $X_{\geq 0} = KX$ and $X_{<0} = RX$. Then $X_{<0} \in \mathbf{C}_{<0}(\mathbb{F})$ by construction and $SX \cong X_{\geq 0} \in \mathbf{C}_{\geq 0}(\mathbb{F})$ by normality. □

In order to prove, now, that the pair of markings $\mathbb{F}(\mathfrak{t})$ is a factorization system on the stable $\infty$-category $\mathbf{C}$, we use the data of the $t$-structure to produce a functorial factorization of morphisms, and we recall ([Lur17, Def. **1.2.1.4**] and our Remark **A.3.7**) that a $t$-structure on $\mathbf{C}$ corresponds to a classical $t$-structure on the triangulated homotopy category of $\mathbf{C}$; this gives us a certain freedom in moving between data living in $\mathbf{C}$ and their "shadow" living in Ho($\mathbf{C}$), at least as soon as these data involve only homotopy invariant information associated to the $t$-structure.

Finally, we use the characterization outlined in Remark **2.3.17** of the factorization functor in terms of its pair reflection/coreflection.

Recall that by Def. **A.3.2.(iii)** every object $X \in \mathbf{C}$ fits into a distinguished triangle $X_{\geq 0} \to X \to X_{<0} \to X_{\geq 0}[1]$. This triangle in Ho($\mathbf{C}$) is the image of a fiber sequence (denoted with the same symbols) in $\mathbf{C}$ via the homotopy-category functor, and can be lifted to such a sequence: this entails that given $f\colon X \to Y$ we can build the diagram[(1)]

$$\begin{array}{ccccccc}
X_{\geq 0} & \longrightarrow & X & \longrightarrow & X_{<0} & \to & X_{\geq 0}[1] \\
{\scriptstyle \tau_{\geq 0}(f)}\downarrow & & \downarrow{\scriptstyle e_f} \searrow{\scriptstyle f} & & \| & & \downarrow{\scriptstyle \tau_{\geq 0}(f)[1]} \\
Y_{\geq 0} & \dashrightarrow & C & \longrightarrow & X_{<0} & \to & Y_{\geq 0}[1] \\
\| & & \downarrow{\scriptstyle m_f} & & \downarrow{\scriptstyle \tau_{<0}(f)} & & \| \\
Y_{\geq 0} & \longrightarrow & Y & \longrightarrow & Y_{<0} & \to & X_{\geq 0}[1]
\end{array} \tag{3.6}$$

where the decorated square is a pullout (so $C \cong X_{<0} \times_{Y_{<0}} Y$, a characterization reminiscent of simplicity for the would-be factorization of $f$), and hence the dotted arrows are determined by the obvious universal property. Now, mapping $f$ to the pair $(e_f, m_f)$ is a candidate factorization functor (a tedious but easy check) in the sense of [KT93].

Now, we have to summon a rather easy but subtle result, [KT93, Thm. **A**], which in a nutshell says that a factorization system on a category $\mathbf{C}$ is determined by a functorial factorization $F$ such that the arrows $m_{e_f}$, $e_{m_f}$ are invertible (the meaning of this notation is self-evident). Functors satisfying this property are called *Eilenberg-Moore factorization functors* in [KT93].[(2)] More precisely, if one defines

$$\begin{aligned}
\mathcal{E}_F &= \{h \in \mathbf{C}^{\Delta[1]} \mid m_h \text{is invertible}\} \\
\mathcal{M}_F &= \{h \in \mathbf{C}^{\Delta[1]} \mid e_h \text{is invertible}\},
\end{aligned} \tag{3.7}$$

then $(\mathcal{E}_F, \mathcal{M}_F)$ is a factorization system as soon as $e_f \in \mathcal{E}_F$ and $m_f \in \mathcal{M}_F$ for any morphism $f$ in $\mathbf{C}$.

REMARK 3.1.3. Before we go on with the proof notice that by the very definition of the factorization functor $F$ in (**3.6**) associated with a $t$-structure above, we have that $\mathcal{M}_F$ coincides with the class of arrows $f$ such that the naturality square of $f$ with respect to the “truncation” functor $\tau_{<0}$ of the $t$-structure is cartesian: we denote this marking of $\mathbf{C}$ as CART$(\tau_{<0})$ adopting the same notation as [RT07, §**4**]. This is reminiscent of our characterization of simplicity via the pullbacks given in Def. **2.2.3**.

LEMMA 3.1.4. The homotopy commutative sub-diagram

$$\begin{array}{ccc}
X_{\geq 0} & \longrightarrow & X \\
\downarrow & & \downarrow \\
Y_{\geq 0} & \longrightarrow & C
\end{array} \tag{3.8}$$

---

(1) This construction, and the link with Remark **2.3.17** was suggested to us by E. Wofsey in a public discussion on `Mathoverflow` [Wof].

(2) These are not the weakest assumptions to ensure that $\mathbb{F}(F) = (\mathcal{E}_F, \mathcal{M}_F) \in \text{FS}(\mathbf{C})$: see the final remark in [KT93] and [JT99, **1.3**].

in the diagram (**3.6**) is a pullout.

*Proof.* Consider the diagram

$$\begin{array}{ccccc} X_{\geq 0} & \longrightarrow & X & & \\ {\scriptstyle \tau_{\geq 0}(f)}\big\downarrow & & \big\downarrow{\scriptstyle e_f} & & \\ Y_{\geq 0} & \longrightarrow & C & \xrightarrow{\ m_f\ } & Y \\ \big\downarrow & & \big\downarrow & & \big\downarrow \\ 0 & \longrightarrow & X_{<0} & \xrightarrow[\tau_{<0}(f)]{} & Y_{<0} \end{array} \tag{3.9}$$

where all the squares are homotopy commutative and apply twice the 3-for-2 law for pullouts. □

LEMMA 3.1.5. Let $F : f \mapsto (e_f, m_f)$ be the factorization functor associated to a $t$-structure by the diagram (**3.6**). Then $\tau_{<0}(e_f)$ and $\tau_{\geq 0}(m_f)$ are equivalences.

*Proof.* Since $\tau_{<0}\tau_{\geq 0} = 0$, by applying $\tau_{<0}$ to the pullout diagram in $\mathbf{C}$ given by lemma **3.1.4**, we get the pushout diagram

$$\begin{array}{ccc} 0 & \longrightarrow & X_{<0} \\ \big\downarrow & & \big\downarrow{\scriptstyle \tau_{<0}(e_f)} \\ 0 & \longrightarrow & C_{<0} \end{array} \tag{3.10}$$

in $\mathbf{C}_{<0}$ which tells us that $\tau_{<0}(e_f)$ is an equivalence. The proof that $\tau_{\geq 0}(m_f)$ is a equivalence is perfectly dual and is obtained by applying $\tau_{\geq 0}$ to the marked pullout diagram in (**3.6**). □

It is now rather obvious that a proof of the equations

$$\mathcal{E}_F = \tau_{<0}^{-1}(\text{EQV}); \qquad \mathcal{M}_F = \tau_{\geq 0}^{-1}(\text{EQV}) \tag{3.11}$$

will imply that $F$ is an Eilenberg-Moore factorization functor. Once this is proved, it is obvious that the preimage of a 3-for-2 class along a functor is again a 3-for-2 class in $\mathbf{C}$, and this entails that both classes in $\mathbb{F}(\mathfrak{t})$ are 3-for-2. We are now ready to prove

PROPOSITION 3.1.6. The pair of markings $\mathbb{F}(\mathfrak{t})$ is a factorization system on the quasicategory $\mathbf{C}$, in the sense of Definition **1.3.3**.

*Proof.* By the very definition of the factorization procedure, and invoking the pullout axiom, we can deduce that the arrow $f$ lies in $\mathcal{E}_F$ if and only if it is inverted by $\tau_{<0}$; this entails that $\mathcal{E}_F = \tau_{<0}^{-1}(\text{EQV})$. So it remains to show that $\mathcal{M}_F = \tau_{\geq 0}^{-1}(\text{Eqv})$. We have already remarked that $\mathcal{M}_F = \text{Cart}(\tau_{<0})$, so we are reduced to showing that $\tau_{\geq 0}^{-1}(\text{Eqv}) = \text{Cart}(\tau_{<0})$. But again, this is easy because on the one side, if $f \in \text{Cart}(\tau_{<0})$ then the square

$$\begin{array}{ccc} \cdot & \longrightarrow & \cdot \\ {\scriptstyle \tau_{\geq 0}(f)}\downarrow & \lrcorner \quad \ulcorner & \downarrow \\ \cdot & \longrightarrow & \cdot \end{array} \tag{3.12}$$

is a pullout since $\tau_{\geq 0}$ preserves pullouts, and yet $\tau_{\geq 0}\tau_{<}0(f)$ is the identity of the zero object. So $\tau_{\geq 0}(f)$ must be an equivalence. On the other hand, the stable $\infty$-categorical analogue of the triangulated five lemma (see [Nee01, Prop. **1.1.20**]), applied to the diagram (**3.6**) shows that if $\tau_{\geq 0}(f)$ is an equivalence then $e_f$ is an equivalence and so $C \cong X$, i.e., $f \in \text{Cart}(\tau_{<0})$. □

Remark 3.1.7. As a side remark, we notice that a completely dual proof would have arisen using $D = Y_{\geq 0} \amalg_{X_{\geq 0}} X$ (see Lemma **3.1.4**) and then showing first that $\mathbb{F}(\mathfrak{t})$ is the factorization system $(\text{Cocart}(\tau_{\geq 0}), \tau_{\geq 0}^{-1}(\text{Eqv}))$ and that $\text{Cocart}(\tau_{\geq 0}) = \tau_{<0}^{-1}(\text{Eqv})$.

This is in line with remark **2.3.17**.

To check that $\mathbb{F}(\mathfrak{t})$ is normal, it only remains to verify that any of the equivalent conditions for normality given in Proposition **2.3.14** holds, which is immediate. This concludes the proof that there is a correspondence between normal torsion theories and $t$-structures: it remains to show that this correspondence is bijective, i.e., that the following proposition holds.

Proposition 3.1.8. In the notations above, we have $\mathbb{F}(\mathfrak{t}(\mathbb{F})) = \mathbb{F}$ and $\mathfrak{t}(\mathbb{F}(\mathfrak{t})) = \mathfrak{t}$.

*Proof.* On the one side, consider the factorization system

$$\mathbb{F}(\mathfrak{t}(\mathbb{F})) = (\tau_{<0}^{-1}(\text{Eqv}), \tau_{\geq 0}^{-1}(\text{Eqv})), \tag{3.13}$$

where the functor $\tau_{<0}$ is the reflection $R$ obtained from the $\mathbb{F}$-factorization of each $X \to 0$, as in the fundamental connection of §**2.1**: $X \xrightarrow{e} X_{<0} \xrightarrow{m} 0$. Recall (Remark **2.1.2**) that the action of $\tau_{<0}\colon \mathbf{C} \to \mathcal{M}/0$ on arrows is obtained from a choice of solutions to lifting problems

$$\begin{array}{ccc} A & \xrightarrow{e'f} & \tau_{<0}B \\ {\scriptstyle e}\downarrow & \nearrow & \downarrow{\scriptstyle m'} \\ \tau_{<0}A & \xrightarrow[m]{} & 0. \end{array} \tag{3.14}$$

It is now evident that $\tau_{<0}^{-1}(\text{Eqv}) = \mathcal{E}$. Indeed:

- If $f \in \tau_{<0}^{-1}(\text{Eqv})$, then in the above square $e'f = \tau_{<0}(f)\,e$, which is in $\mathcal{E}$ since $\mathcal{E}$ contains equivalences and is closed for composition. But $e'$ lies in $\mathcal{E}$, so that $f \in \mathcal{E}$ by the 3-for-2 property of $\mathcal{E}$;
- If $f \in \mathcal{E}$, then $e'f$ is in $\mathcal{E}$ and so in the same square we read two lifting problems with unique solutions, which implies that $\tau_{<0}(f)$ is invertible.

On the other side, we have to compare the $t$-structures $\mathfrak{t} = (\mathbf{C}_{\geq 0}, \mathbf{C}_{<0})$ and $\mathfrak{t}(\mathbb{F}(\mathfrak{t}))$. We have $X \in \mathbf{C}_{\geq 0}(\mathbb{F}(\mathfrak{t}))$ if and only if $\begin{bmatrix} 0 \\ \downarrow \\ X \end{bmatrix} \in \mathcal{E}(\mathfrak{t})$. Since $\mathcal{E}(\mathfrak{t}) = \tau_{<0}^{-1}(\text{Eqv})$, we see that $X \in \mathbf{C}_{\geq 0}(\mathbb{F}(\mathfrak{t}))$ if and only if $X_{<0} \cong 0$. But it is a direct consequence of Lemma **2.3.4** that $X_{<0} \cong 0$ if and only if $X \in \mathbf{C}_{\geq 0}$. Dually, one can prove that $\mathbf{C}_{<0}(\mathbb{F}(\mathfrak{t})) = \mathbf{C}_{<0}$ (but this, in view of Remark **A.3.8**, is superfluous). □

## 3.2 Examples.

> Stand firm in your refusal to remain conscious during algebra. In real life, I assure you, there is no such thing as algebra.
>
> F. Leibowitz

We gather here a series of classical and less classical examples (more will be given in the subsequent chapters), heavily relying on existing literature. As a consequence, this section is more sketchy and gives several (even non trivial) statements without proof.

Example 3.2.1. (Bousfield localization of spectra): The category **Sp** of spectra furnishes the most natural example of a stable $\infty$-category; a classical construction in [Bou79] endows **Sp** with a $t$-structure for each object $E$, whose right class (and whose reflection functor) is called *$E$-localization*; we define the subcategories

$$\mathcal{T}_E = \{X \in \mathbf{Sp} \mid X \wedge E \simeq *\} \tag{3.15}$$

$$\mathcal{F}_E = \{Y \in \mathbf{Sp} \mid [X, Y] \simeq * \; \forall X \in \mathcal{T}_E\} = \mathcal{T}_E^{\perp} \tag{3.16}$$

These two classes form a stable $t$-structure $\mathfrak{t}_E$ in the sense of **4.4.3** (the notation is chosen to inspire the correspondence between $\mathcal{T}_E$ and torsion objects, and between $\mathcal{F}_E$ and free objects.

We now want to characterize the factorization system corresponding to this (stable) $t$-structure under the Rosetta stone theorem. We start by recalling that [Bou79, Lemma **1.13**] ensures that $\mathcal{T}_E$ is generated under homotopy colimits by a single element $G_E$, and that $\mathcal{F}_E$ is precisely the right object-orthogonal to this single object; now let $g\colon * \to G_E$ be the initial morphism in **Sp**, and let

$$\mathbb{F}_E = \left({}^{\perp}(\{g\}^{\perp}), \{g\}^{\perp}\right) \in \text{PF}(\mathbf{Sp}). \tag{3.17}$$

THEOREM 3.2.2. The pair of markings $\mathbb{F}_E$ is a normal torsion theory, and corresponds to the $E$-localization of **Sp** under Thm. **3.1.1**.

*Proof.* It is basically a way to rewrite [Bou79, **1.13**, **1.14**] replacing object-orthogonality and generation with arrow-orthogonality and generation (this can be done in view of **1.3.11**), and subsequently to check that the prefactorization left generated by $g$ coincides with $\mathbb{F}(\mathfrak{t}_E)$ of Def. **3.2**. □

The above example survives to the category of chain complexes of abelian groups, giving the $p$-localization of the category $\mathbf{Ch}(\mathbb{Z})$; the two contexts are linked by [Bou79, §**2**] (see in particular [Bou79, **2.4**, **2.5**]). For another glance to $p$-localization see Example **3.2.6** below.

EXAMPLE 3.2.3. (THE $p$-ACYCLIC $t$-STRUCTURE ON $\mathbf{Ch}(\mathbb{Z})$): Let $p \in \mathbb{Z}$ be a prime, and let $A \in \mathbf{Ch}(\mathbb{Z})$ be a chain complex of abelian groups. We say that $A$ is *$p$-acyclic* if (i) $A$ is projective and (ii) the tensor product $A \otimes_{\mathbb{Z}} \mathbb{Z}/p\mathbb{Z}$ is nullhomotopic; the class of $p$-acyclic complexes is denoted $p{-}\mathbf{Ac}$. We call *$p$-local* complexes the elements of $(p{-}\mathbf{Ac})^{\perp}$.

The pair $\big(p{-}\mathbf{Ac}, (p{-}\mathbf{Ac})^{\perp}\big)$ induces a $t$-structure on the category of chain complexes; the reflection with respect to this $t$-structure is called *$p$-localization*, and it is defined by

$$A \mapsto \widehat{A} := \varprojlim_n \big(A \otimes \mathbb{Z}/p^n\mathbb{Z}\big) \tag{3.18}$$

Since it is a homotopy limit of $p$-local chain complexes, we conclude that $\widehat{A}$ is again $p$-local.

EXAMPLE 3.2.4. (THE STANDARD $t$-STRUCTURE ON CHAIN COMPLEXES): **A.3.1** defines the canonical $t$-structure on the derived category $\mathbf{D}(R)$ of a ring $R$ as the pair of subcategories

$$\begin{aligned} \mathbf{D}_{\geq 0}(R) &= \{A_* \in \mathbf{D}(R) \mid H^n(A_*) = 0;\ n \leq 0\} \\ \mathbf{D}_{\leq 0}(R) &= \{B_* \in \mathbf{D}(R) \mid H^n(B_*) = 0;\ n \geq 0\}. \end{aligned}$$

The construction of $\mathbb{F}(\mathfrak{t})$ provided by (**3.6**) gives the following definition for the two classes of chain maps in $\mathbf{Ch}(R)$: $\mathcal{E}(\mathfrak{t})$ (resp. $\mathcal{M}(\mathfrak{t})$) is the class of arrows such that the negative (resp. positive) truncation is Working out the details, this means that the factorization of $f \colon X_* \to Y_*$ is defined via the pullout

$$X \longrightarrow X_{<0} \oplus_{Y_{<0}} Y \longrightarrow Y \tag{3.19}$$

where the object $X_{<0} \oplus_{Y_{<0}} Y$ is defined to be the mapping cone of the map $(f_{<0}, \rho_Y) \colon X_{<0} \oplus Y \to Y_{<0}$.

EXAMPLE 3.2.5. (THE STANDARD $t$-STRUCTURE ON SPECTRA): The stable $\infty$-category of spectra carries another $t$-structure, whose left class is determined by those objects whose homotopy groups vanish in negative dimension (recall that a spectrum has homotopy groups in each, possibly negative, degree).

We can reproduce the above argument to find the corresponding factorization system.

Example 3.2.6. (The $p$-local/$p$-complete arithmetic square): Let $p \in \mathbb{Z}$ be a prime number; a spectrum $E \in \mathbf{Sp}$ is called *$p$-torsion* if for every $x \in \pi_*(E)$ there exists a $n = n_x$ such that $p^n x = 0$. The full sub-$\infty$-category of $p$-torsion spectra is coreflective in $\mathbf{Sp}$, via a coreflection $G_p(-) \to (-)$; this means that every spectrum $X$ has a *$p$-torsion approximation* fitting into a fiber sequence

$$\tau_p X \to X \to X\left[\tfrac{1}{p}\right] \tag{3.20}$$

the rightmost object of which is called the *$p$-localization* of $X$. The class of $p$-torsion and $p$-local spectra form mutually (object-)orthogonal subcategories of $\mathbf{Sp}$, and together they form a $t$-structure called the *$p$-local* $t$-structure.

Let again $p \in \mathbb{Z}$ be a prime number; a spectrum $E \in \mathbf{Sp}$ is called *$p$-complete* if the homotopy limit of the tower

$$E \xrightarrow{p} E \xrightarrow{p} E \xrightarrow{p} \cdots \tag{3.21}$$

vanishes. The full sub-$\infty$-category of $p$-complete spectra is reflective in $\mathbf{Sp}$, via a reflection $X \to \widehat{X}_p$; this means that every spectrum has a *$p$-completion* fitting into a fiber sequence $G_p X \to X \to \widehat{X}_p$, the leftmost object of which is called *$p$-torsion approximation.* These data determine another $t$-structure on $\mathbf{Sp}$, called the *$p$-complete* $t$-structure.

These two $t$-structures can be arranged into a so-called *arithmetic square* or *fracturing square*, i.e. in the following diagram

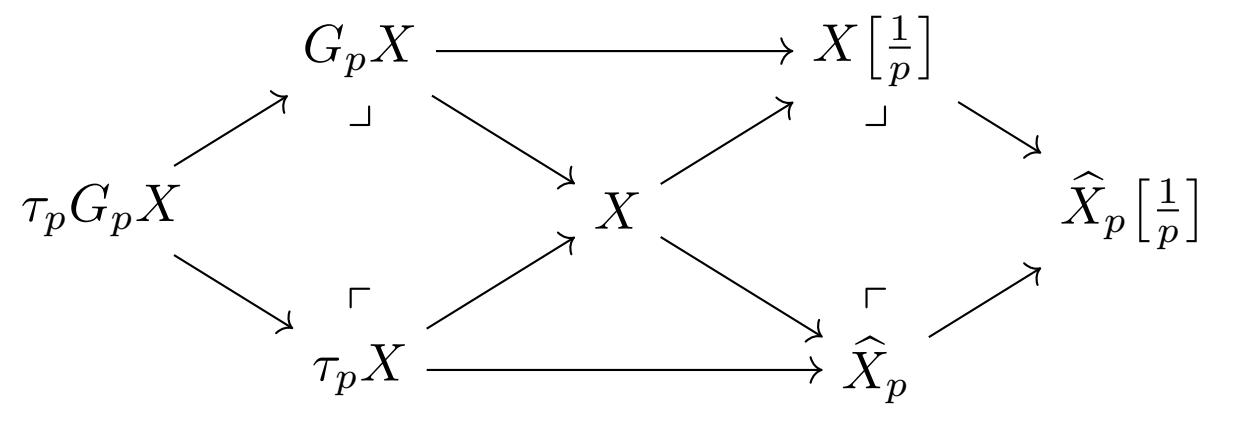

(3.22)

Such a diagram, canonically built from the prime number $p$ alone and the spectra $E$ (and functorial in this argument), contains an impressive amount of informations that we now attempt to characterize more explicitly:

1. the two squares are pullout;
2. the two sequences $\tau_p G_p X \to G_p X \to X\left[\tfrac{1}{p}\right] \to \widehat{X}_p\left[\tfrac{1}{p}\right]$ and $\tau_p G_p X \to \tau_p X \to \widehat{X}_p \to \widehat{X}_p\left[\tfrac{1}{p}\right]$ are long exact fiber sequences (this means that $\widehat{X}_p\left[\tfrac{1}{p}\right] \cong \tau_p G_p X[1]$);
3. the diagonals are fiber sequences by construction.

Motivated by this example, we give the following

Definition 3.2.7. (Crimson $t$-structures): Let $\mathfrak{t}_1, \mathfrak{t}_2 \in \text{TS}(\mathbf{C})$ be two $t$-structures; the two are called *crimson*, or *fracturing*, if the two fiber sequences $S_1X \to X \to R_1X$ and $S_2X \to X \to R_2X$ arrange into an hexagonal diagram

$$\begin{array}{ccccccc} & & S_1X & \longrightarrow & R_2X & & \\ & \nearrow & \lrcorner & \searrow \quad \nearrow & \lrcorner & \searrow & \\ S_2S_1X & & & X & & & R_2R_1X \\ & \searrow & \ulcorner & \nearrow \quad \searrow & \ulcorner & \nearrow & \\ & & S_2X & \longrightarrow & R_1X & & \end{array} \tag{3.23}$$

natural in the object $X$, such that properties (1)–(3) above hold.

## 3.3 Model dependency

One might wonder, at this point, to which extent the "Rosetta stone" theorem is true in other models for $(\infty, 1)$-category theory. Apart from stable $\infty$-categories, extensively treated in the present work, we know (see , **A.5**) there are many, well suited to the description of homological algebra:

(1) (stable) model categories;
(2) (DG-)enriched categories;
(3) (stable) derivators.

It is really tempting to think that a "generic object" $\mathsf{C}$ of any of these higher categories is a "model-free" (stable) $(\infty, 1)$-category, and possesses a natural notion of $t$-structure; with the possible exception of stable derivators(3), each of these models is rich enough to interpret a notion of "factorization system on $\mathsf{C}$", and then the fundamental connection between reflective (pre)factorization systems on $\mathsf{C}$ and reflective sub-$\infty$-categories of $\mathsf{C}$; each of these models is powerful enough to interpret the notion of normal torsion theory, and *subsequently* of $t$-structure, taking the former as the definition of the latter.

A major achievement of our Rosetta stone **3.1.1** is, hence, the possibility to give the notion of $t$-structure a meaning in several different categorical contexts, like enriched categories and model categories.

The scope of the present section is to pave the way to speculations in this respect, and will hopefully be a starting point for future investigations. We start recalling the various flavours of factorization systems we have to deal with, in studying (stable) $(\infty, 1)$-categories.

(3) As mentioned elsewhere, at the moment of writing there is a work in progress in this direction, [?].

### 3.3.1 Enriched factorization systems.

Intuitively, an *enriched* factorization system in an enriched category $\mathbf{C} \in \mathcal{V}\text{-}\mathbf{Cat}$ consists, according to [DK74, LW] of a pair $\mathbb{F} = (\mathcal{E}, \mathcal{M})$ of classes of morphisms in $\mathbf{C}$ such that $\mathcal{E} = {}^{\perp}\mathcal{M}$, and $\mathcal{M} = \mathcal{E}^{\perp}$, where the orthogonality relation is defined in $\mathcal{V}$-**Cat** by an enriched analogue of Remark **1.2.23**, and such that every arrow in $\mathbf{C}$ is $\mathbb{F}$-crumbled in the obvious sense. More explicitly, if $\mathcal{V}$ is an enriched symmetric monoidal category with finite limits, then $f \perp g$ in $\mathcal{V}$-**Cat** if and only if the square in (**1.8**) is a pullback in $\mathcal{V}$:

$$\begin{array}{ccc} \mathbf{C}(B,X) & \longrightarrow & \mathbf{C}(B,Y) \\ \downarrow & \lrcorner & \downarrow \\ \mathbf{C}(A,X) & \longrightarrow & \mathbf{C}(A,Y) \end{array} \tag{3.24}$$

This formalism applies well to simplicial(ly enriched) categories, and more precisely in the stable setting, to DG-categories, which can be regarded as particular examples of simplicial categories via the Dold-Kan correspondence.

REMARK 3.3.1. In the case of simplicially enriched categories the above definition admits an equivalent reformulation relying on the adjunction

$$\mathfrak{C} \colon \mathbf{sSet} \leftrightarrows \mathbf{sSet\text{-}Cat} \colon N_{\mathbf{sSet}} \tag{3.25}$$

In particular, for each $\mathbf{C} \in \mathbf{sSet\text{-}Cat}$ we define:

- a "lifting problem" as a map $\mathfrak{C}(\Delta[1] \times \Delta[1]) \to \mathbf{C}$;
- a "solution" to the lifting problem is presented by an extension over $\mathfrak{C}(\Delta[3]) = \mathfrak{C}(\Delta[1] \star \Delta[1])$ (which nevertheless is only **sSet**-equivalent, and not isomorphic, to $\mathfrak{C}(\Delta[1]) \star \mathfrak{C}(\Delta[1])$).

(These definitions work well only when $\mathbf{C}$ is Bergner-cofibrant [Ber10])

Mild assumptions on $\mathbf{C}$ (see [**?**]) ensure that enriched factorization systems on $\mathbf{C}$ and 1-dimensional factorization systems on $|\mathbf{C}|$ (the **Set**-category naturally associated to $\mathbf{C}$) are in bijection. This paves the way to the following definition of $t$-structure in a DG-category:

DEFINITION 3.3.2. A $t$-structure on a DG-category $\mathbf{A}$ is an enriched factorization system $(\mathcal{E}, \mathcal{M})$ such that

(1) the two classes of morphisms $\mathcal{E}, \mathcal{M}$ are 3-for-2;
(2) the coreflective/reflective pair $0/\mathcal{E}$, $\mathcal{M}/0$ have co/reflection functors $S, R$ respectively, and each object $X \in \mathbf{A}$ fits into a pullback and pushout square

$$\begin{array}{ccc} SX & \longrightarrow & X \\ \downarrow & \lrcorner \quad \ulcorner & \downarrow \\ 0 & \longrightarrow & RX \end{array} \tag{3.26}$$

### 3.3.2 Homotopy factorization systems.

Model categories $\mathbf{M}$ possess a notion of “homotopy” factorization system, which induces a 1-dimensional factorization system on the homotopy category Ho($\mathbf{M}$); the following definition is taken from [Joy08, Def. **F.1.3**]:

DEFINITION 3.3.3. Let $\mathbf{M}$ be a model category with model structure $(\mathcal{C}\wr\{, \mathcal{W}\|, \mathcal{F}\rangle\lfloor)$. A pair $(\mathcal{E}, \mathcal{M})$ of classes of maps in $\mathbf{M}$ is a *homotopy factorisation system* if

(HFS1) the classes $\mathcal{E}, \mathcal{M}$ are *homotopy replete*;
(HFS2) the pair $(\mathcal{E} \cap \mathcal{C}\wr\{_{\text{cf}}, \mathcal{M} \cap \mathcal{F}\rangle\lfloor_{\text{cf}})$ is a weak factorisation system in $\mathbf{M}_{\text{cf}}$, where for $\mathcal{K} \subseteq \hom(\mathbf{M})$ we denote $\mathcal{K}_{\text{cf}}$ the morphisms in $\mathcal{K}$ having co/fibrant co/domain;
(HFS3) the class $\mathcal{E}$ is R32, and the class $\mathcal{M}$ is L32.

It can be shown ([Joy08, Prop. **F.2.6**]) that a homotopy factorization system determines a unique factorization system on the homotopy category Ho($\mathbf{M}$); also, several theorems of the calculus of factorization survive to this setting, and most notably the closure properties of §**1.4** taking care to replace every co/limit appearing there with the appropriate homotopy version: so, in particular we have ([Joy08, Prop. **F.4.8**])

PROPOSITION 3.3.4. The right class of a homotopy factorisation system is closed under homotopy base change. Dually, the left class is closed under homotopy cobase change.

This paves the way to the following definition of a $t$-structure in a stable model category:

DEFINITION 3.3.5. Let $\mathbf{M}$ be a stable model category; a *homotopy normal torsion theory* on $\mathbf{M}$ is a homotopy factorization system $(\mathcal{E}, \mathcal{M})$ on $\mathbf{M}$ such that

(1) both $\mathcal{E}, \mathcal{M}$ are 3-for-2 classes;
(2) the subcategories $0/\mathcal{E}$, $\mathcal{M}/0$ (defined in the same fashion as (**2.1**)) are respectively coreflective and reflective, and the co/reflection fit into the homotopy-pullback-and-pushout diagram

$$\begin{array}{ccc} SX & \longrightarrow & X \\ \downarrow & \lrcorner \quad \ulcorner & \downarrow \\ 0 & \longrightarrow & RX \end{array} \tag{3.27}$$

# Chapter 4

# Hearts and towers

In the present section we exploit the description of $t$-structures as normal torsion theories of Ch. **3** to discuss two apparently separated constructions in the theory of triangulated categories: the characterization of *bounded $t$-structures* in terms of their hearts, and *semiorthogonal decompositions* on triangulated categories. In the stable setting both notions stem as particular cases of a single construction.

In analogy with the example of the Postnikov decomposition of a morphism $f: X \to Y$ of spaces (or spectra, or objects of an $\infty$-topos), we construct (Def. **4.2.6**) the *tower* $♖_{\{i_j\}}(f)$[(1)] of a morphism induced by a $\mathbb{Z}$-equivariant *$J$-family* of normal torsion theories $\{\mathbb{F}_i\}_{i\in J}$, i.e. a monotone function $J \to \text{FS}(\mathbf{C})$ "taking normal values", which is *equivariant* with respect to an action of the group $\mathbb{Z}$ on both sets.

As we will see along the chapter, a natural way to encompass these structures is to vary the action on the domain of the $J$-family (choosing diffferent $J$s and different actions on $J$ will result in different kinds of $t$-structures for the values $J(\lambda)$. We will concentrate on the following two "extremal" examples:

- For $J = \mathbb{Z}$ with its obvious self-action, we recover the classical notion of Postnikov towers in a triangulated category endowed with a $t$-structure (and a fortiori, the notion of Postnikov tower in the category **Sp** of spectra), and subsequently we give a neat, conceptual proof of the the abelianity of the *heart* of a $t$-structure in the stable setting, basically relying on the uniqueness of a suitable factorization.
- For $J$ a finite totally ordered set, or more generally any set $J$ with trivial $\mathbb{Z}$-action, we recover the theory of *semiorthogonal decompositions* [BO95, Kuz11], showing in Thm. **4.4.9** that such a $J \to \text{FS}_\nu(\mathbf{C})$ consists of a family $\{\mathbb{F}_i\}_{i\in J}$ of *stable* $t$-structures. This is a classical result.

---

(1) Pron. *rook*; it is the same rook of the game of chess.

## 4.1 Posets with $\mathbb{Z}$-actions.

為無為。事無事。味無味。

Laozi LXIII

This section has an introductory purpose, aiming to introduce the terminology about partially ordered groups and their actions, and then specialize the discussion to $\mathbb{Z}$-actions on partially ordered sets.

We do not aim at reaching a complete generality, but instead at gathering a number of useful results and nomenclature we can refer to along the present chapter. Among various possible choices, we mention specialized references as [Bly05, Gla99, Fuc63] for an extended discussion of the theory of actions on ordered groups.

Definition 4.1.1. A *partially ordered group* ("po-group" for short) consists of a group $\mathbf{G} = (G, \cdot, 1)$ endowed with a relation $\preceq$ which is a partial order and a (two-sided) congruence on $G$, namely for any $g \preceq h$ and $a, b \in \mathbf{G}$ we have

(i) $a \cdot g \preceq a \cdot h$ and
(ii) $g \cdot b \preceq h \cdot b$.

Remark 4.1.2. We should draw a distinction between a *left* po-group (satisfying only property *i* above) and a *right* po-group (satisfying only *ii*). At the level of generality we need ignoring this subtlety is absolutely harmless.

A supplementary motivation to choose this slightly looser definition is that it seems more natural for a group to be ordered by a two-sided congruence, since in this case inversion $(-)^{-1} \colon \mathbf{G} \to \mathbf{G}$ is an antitone antiautomorphism of groups, i.e. we have that

- $g \preceq h \iff h^{-1} \preceq g^{-1}$;
- The set $G^+$ of *positive* elements, i.e. the set $\{g \in G \mid 1 \preceq g\}$ is closed under conjugation.

Definition 4.1.3. A *homomorphism* of po-groups consists of a group morphism $f \colon G \to H$ which is also a monotone mapping. This, with the obvious choices of identities and composition, defines a category **POGrp** of partially ordered groups and their morphisms.

Definition 4.1.4. Let $G$ be any group. A $G$-poset is a partially ordered set $(P, \leq)$ endowed with a group homomorphism $G \to \mathrm{Aut}_{\leq}(P)$ to the group of order isomorphisms of $P$.

Remark 4.1.5. Obviously, the former definition of $G$-poset is equivalent to the following one: a $G$-poset consists of a poset $(P, \leq)$ with a map $a \colon G \times P \to P$ satisfying the well-known properties of a group action, and furthemore such that for each $g \in G, p \leq q \in P$ one has $a(g, p) \leq a(g, q)$.

Lemma 4.1.6. The category **Pos** of partially ordered sets and monotone maps is cartesian closed.

*Proof.* This is a classical result; there is only one way to endow the underlying set of a product $P \times Q$ of posets with a partial order in such a way that the universal property of the product is satisfied, and there is only one way to endow the set of all monotone functions between two posets with a partial order relation to obtain the adjunction

$$\mathbf{Pos}(P \times Q, R) \cong \mathbf{Pos}(P, R^Q). \tag{4.1}$$

□

Proposition 4.1.7. When $G$ is a po-group $(G, \preceq)$, the action map $a\colon G \times P \to P$ defining a $G$-poset is a monotone mapping if we endow $G \times P$ with the product order; equivalently, the map $G \to \mathrm{Aut}_{\leq}(P)$ is monotone if we endow the codomain with the order inherithed from the inclusion $\mathrm{Aut}_{\leq}(P) \subseteq P^P$.

*Proof.* Straightforward, unwinding the definitions: Def. **4.1.4** can be reinterpreted in light of this viewing $G$ endowed with the trivial partial order where $x \preceq y$ if and only if $x = y$. □

Definition 4.1.8. A $\mathbb{Z}$-*poset* is a partially ordered set $(P, \leq)$ together with a group action

$$+_P\colon P \times \mathbb{Z} \to P\colon (x, n) \mapsto x +_P n \tag{4.2}$$

which is a morphism of partially ordered sets, when $\mathbb{Z}$ is regarded with its usual total order.

Remark 4.1.9. It is immediate to see that a $\mathbb{Z}$-poset is equivalently the datum of a poset $(P, \leq)$ together with a monotone bijection $\rho\colon P \to P$ such that $x \leq \rho(x)$ for any $x$ in $P$. The function $\rho$ and the action are related by the identity $\rho(x) = x +_P 1$.

Notation 4.1.10. To avoid a cumbersome accumulation of indices, the action $+_P$ will be often denoted as a simple "+". This is meant to evoke in the reader the two most natural examples of a $\mathbb{Z}$-poset, described below:

Example 4.1.11. The poset $(\mathbb{Z}, \leq)$ of integers with their usual order is a $\mathbb{Z}$-poset with the action given by the usual sum of integers. The poset $(\mathbb{R}, \leq)$ of real numbers with their usual order is a $\mathbb{Z}$-poset for the action given by the sum of real numbers with integers (seen as a subring of real numbers).

Remark 4.1.12. If $(P, \leq)$ is a finite poset, then the only $\mathbb{Z}$-action it carries is the trivial one. Indeed, if $\rho\colon P \to P$ is the monotone bijection associated with the $\mathbb{Z}$-action, one sees that $\rho$ is of finite order by the finiteness of $P$.

Therefore there exists an $n \geq 1$ such that $\rho^n = \mathrm{id}_P$. It follows that, for any $x$ in $P$,

$$x \leq x+1 \leq \cdots \leq x+n = x \tag{4.3}$$

and so $x = x+1$.

NOTATION 4.1.13. An obvious terminology: a *$G$-fixed point* for a $G$-poset $P$ is an element $p \in P$ kept fixed by all the elements of $G$ under the action $+_P$. An important observation is that an element $p$ of a $\mathbb{Z}$-poset $P$ is a $\mathbb{Z}$-fixed point if and only if $p +_P 1 = p$.

LEMMA 4.1.14. If $k \in P$ is a $\leq$-maximal or $\leq$-minimal element in the $\mathbb{Z}$-poset $(P, \leq)$, then it is a $\mathbb{Z}$-fixed point.

REMARK 4.1.15. Given a poset $P$ we can always define a partial order on the set $P \cup \{-\infty, +\infty\}$ which extends the partial order on $P$ by the rule $-\infty \leq x \leq +\infty$ for any $x \in P$.

LEMMA 4.1.16. If $(P, \leq)$ is a $\mathbb{Z}$-poset, then $(P \cup \{\pm\infty\}, \leq)$ carries a natural $\mathbb{Z}$-action extending the $\mathbb{Z}$-action on $P$, by declaring both $-\infty$ and $+\infty$ to be $\mathbb{Z}$-fixed points.

*Proof.* Adding a fixed point always gives an extension of an action, so we only need to check that the extended action is compatible with the partial order. This is equivalent to checking that also on $P \cup \{\pm\infty\}$ the map $x \to x+1$ is a monotone bijection such that $x \leq x+1$, which is immediate. □

Posets with $\mathbb{Z}$-actions naturally form a category, whose morphisms are $\mathbb{Z}$-*equivariant* morphisms of posets. More explicitly, if $P$ and $Q$ are $\mathbb{Z}$-posets with actions $+_P$ and $+_Q$, then a morphism of $\mathbb{Z}$-posets between them is a morphism of posets $\varphi\colon P \to Q$ such that

$$\varphi(x +_P n) = \varphi(x) +_Q n, \tag{4.4}$$

for any $x \in P$ and any $n \in \mathbb{Z}$.

LEMMA 4.1.17. The choice of an element $x$ in a $\mathbb{Z}$-poset $P$ is equivalent to the datum of a $\mathbb{Z}$-equivariant morphism $\varphi\colon (\mathbb{Z}, \leq) \to (P, \leq)$. Moreover $x$ is a $\mathbb{Z}$-fixed point if and only if the corresponding morphism $\varphi$ factors $\mathbb{Z}$-equivariantly through $(*, \leq)$, where $*$ denotes the terminal object of **Pos**.

*Proof.* To the element $x$ one associates the $\mathbb{Z}$-equivariant morphism $\varphi_x$ defined by $\varphi_x(n) = x + n$. To the $\mathbb{Z}$-equivariant morphism $\varphi$ one associates the element $x_\varphi = \varphi(0)$. It is immediate to check that the two constructions are mutually inverse. The proof of the second part of the statement is straightforward. □

LEMMA 4.1.18. Let $\varphi\colon (\mathbb{Z}, \leq) \to (P, \leq)$ be a $\mathbb{Z}$-equivariant morphism of $\mathbb{Z}$-posets. Then $\varphi$ is either injective or constant.

*Proof.* Assume $\varphi$ is not injective. then there exist two integers $n$ and $m$ with $n > m$ such that $\varphi(n) = \varphi(m)$. By $\mathbb{Z}$-equivariancy we therefore have

$$x_\varphi + (n - m) = x_\varphi, \tag{4.5}$$

with $n - m \geq 1$ and $x_\varphi = \varphi(0)$. The conclusion then follows by the same argument used in Remark **4.1.12**. □

Lemma 4.1.19. Let $\varphi\colon (P, \leq) \to (Q, \leq)$ be a morphism of $\mathbb{Z}$-posets. Assume $Q$ has a minimum and a maximum. Then $\varphi$ extends to a morphism of $\mathbb{Z}$-posets $(P \cup \{\pm\infty\}, \leq) \to (Q, \leq)$ by setting $\varphi(-\infty) = \min(Q)$ and $\varphi(+\infty) = \max(Q)$.

*Proof.* Since $\min(Q)$ and $\max(Q)$ are $\mathbb{Z}$-fixed points by Lemma **4.1.14**, the extended $\varphi$ is a morphism of $\mathbb{Z}$-posets. Moreover, since $\min(Q)$ and $\max(Q)$ are the minimum and the maximum of $Q$, respectively, the extended $\varphi$ is indeed a morphism of posets, and so it is a morphism of $\mathbb{Z}$-posets. □

### 4.1.1 $J$-families of $t$-structures.

The main reason why we are interested in the theory of $\mathbb{Z}$-poset is the following result we already used in Ch. **2**, Ch. **3** (and recalled also in **A.3.11**):

Remark 4.1.20. Let $\mathbf{C}$ be a stable $\infty$-category. Then, the collection $\textsc{ts}(\mathbf{C})$ of all $t$-structures on $\mathbf{C}$ is a poset with respect to following order relation: given two $t$-structures $\mathfrak{t}_a = (\mathbf{C}_{\geq_a 0}, \mathbf{C}_{<_a 0})$[(2)] and $\mathfrak{t}_b = (\mathbf{C}_{\geq_b 0}, \mathbf{C}_{<_b 0})$, one has $\mathfrak{t}_a \preceq \mathfrak{t}_b$ iff $\mathbf{C}_{<_a 0} \subseteq \mathbf{C}_{<_b 0}$.

The ordered group $\mathbb{Z}$ acts on $\textsc{ts}(\mathbf{C})$ in a way that is fixed (Remark **4.1.9**) by the action of the generator $+1$; this maps a $t$-structure $\mathfrak{t} = (\mathbf{C}_{\geq 0}, \mathbf{C}_{<0})$ to the *shifted* $t$-structure $\mathfrak{t}[1] = (\mathbf{C}_{\geq 0}[1], \mathbf{C}_{<0}[1])$.

Since $\mathfrak{t} \preceq \mathfrak{t}[1]$ one sees that $\textsc{ts}(\mathbf{C})$ is naturally a $\mathbb{Z}$-poset (this follows from **A.3.11**).

Notation 4.1.21. If $\mathfrak{t} = (\mathbf{C}_{\geq 0}, \mathbf{C}_{<0})$ is a $t$-structure on $\mathbf{C}$, it is customary to write $\mathbf{C}_{\geq 1}$ for $\mathbf{C}_{\geq 0}[1]$ and $\mathbf{C}_{<1}$ for $\mathbf{C}_{<0}[1]$, so that $\mathfrak{t}[1] = (\mathbf{C}_{\geq 1}, \mathbf{C}_{<1})$, and more generally $\mathbf{C}_{\geq n} := \mathbf{C}_{\geq 0}[n]$, $\mathbf{C}_{<n} := \mathbf{C}_{<0}[n]$ for each $n \in \mathbb{Z}$, so that $\mathfrak{t}[n] = (\mathbf{C}_{\geq n}, \mathbf{C}_{<n})$.

We now have the natural desire to consider families of $t$-structures on $\mathbf{C}$ indexed by an *arbitrary* $\mathbb{Z}$-poset $J$, as in the following

Definition 4.1.22. Let $(J, \leq)$ be a $\mathbb{Z}$-poset. A *$J$-family* of $t$-structures on a stable $\infty$-category $\mathbf{C}$ is a $\mathbb{Z}$-equivariant morphism of posets $\mathfrak{t}\colon J \to \textsc{ts}(\mathbf{C})$.

More explicitly, a $J$-family is a family $\{\mathfrak{t}_j\}_{j \in J}$ of $t$-structures on $\mathbf{C}$ such that

(1) $\mathfrak{t}_i \preceq \mathfrak{t}_j$ if $i \leq j$ in $J$;

[(2)]The baffled reader is invited to look at Notation **4.1.24**.

(2) $\mathfrak{t}_{i+1} = \mathfrak{t}_i[1]$ for any $i \in J$.

REMARK 4.1.23. A natural choice of notation, motivated by the "Rosetta stone" **3.1.1**, is the following: a $J$-family of $t$-structures is the same as a $J$-family of normal torsion theories on $\mathbf{C}$ (or, more formally, the maps $\mathbb{F}(-)$ and $\mathfrak{t}(-)$ defined in the proof of the Rosetta stone become isomorphisms *in the category* $\mathbb{Z}$-**Pos** for a suitable choice of partial order and $\mathbb{Z}$-action on normal torsion theories).

Motivated by this remark, we feel free to call "$J$-family of normal torsion theories" any monotone function $J \to \text{NTT}(\mathbf{C})$ which is also $\mathbb{Z}$-equivariant.

NOTATION 4.1.24. For $i \in J$, we will write $\mathbf{C}_{\leq i}$ and $\mathbf{C}_{>i}$ for $\mathbf{C}_{\leq_i 0}$ and $\mathbf{C}_{<_i 0}$, respectively. With this notation we have that $\mathfrak{t}_i = (\mathbf{C}_{\geq i}, \mathbf{C}_{<i})$. Note that, by $\mathbb{Z}$-equivariancy, this notation is consistent. Namely $\mathfrak{t}_{i+1} = \mathfrak{t}_i[1]$ implies $\mathbf{C}_{\geq_{i+1} 0} = \mathbf{C}_{\geq_i 0}[1]$ and so

$$\mathbf{C}_{\geq i+1} = \mathbf{C}_{\geq i}[1]. \tag{4.6}$$

Similarly, one has

$$\mathbf{C}_{<i+1} = \mathbf{C}_{<i}[1]. \tag{4.7}$$

We underline how in this choice of notation the condition $\mathfrak{t}_i \preceq \mathfrak{t}_j$ for $i \leq j$ translates to the very natural condition $\mathbf{C}_{<i} \subseteq \mathbf{C}_{<j}$ for $i \leq j$. Notice that this is basically [GKR04, Def. **3.1**].

EXAMPLE 4.1.25. A $\mathbb{Z}$-family of $t$-structures is, by Lemma **4.1.9**, equivalent to the datum of a $t$-structure $\mathfrak{t}_0 = (\mathbf{C}_{\geq 0}, \mathbf{C}_{<0})$. One has $\mathfrak{t}_1 = (\mathbf{C}_{\geq 1}, \mathbf{C}_{<1})$ consistently with the notations in Remark **4.1.20**. Notice that by our Remark **4.1.12**, as soon as $\mathbf{C}_{\geq 0}[1] \subset \mathbf{C}_{\geq 0}$ (proper inclusion), then this proper inclusion is valid for all $n \in \mathbb{Z}$, i.e. the orbit $\mathfrak{t} + \mathbb{Z}$ is an infinite set.

EXAMPLE 4.1.26. An $\mathbb{R}$-family of $t$-structures is the datum of a $t$-structure $\mathfrak{t}_\lambda = (\mathbf{C}_{\geq \lambda}, \mathbf{C}_{<\lambda})$ on $\mathbf{C}$ for any $\lambda \in \mathbb{R}$ in such a way that $\mathfrak{t}_{\lambda+1} = \mathfrak{t}_\lambda[1]$. Such a structure is called a *slicing* of $\mathbf{C}$ in [Bri07].[(3)]

EXAMPLE 4.1.27. (A TAUTOLOGICAL EXAMPLE): By taking $J = \text{TS}(\mathbf{C})$ and $\mathfrak{t}$ to be the identity of $\text{TS}(\mathbf{C})$ one sees that the whole $\text{TS}(\mathbf{C})$ can be looked at as a particular $J$-family of $t$-structures on $\mathbf{C}$.

REMARK 4.1.28. The poset $\text{TS}(\mathbf{C})$ has a minimum and a maximum given by

$$\min(\text{TS}(\mathbf{C})) = (\mathbf{C}, \mathbf{0}); \qquad \max(\text{TS}(\mathbf{C})) = (\mathbf{0}, \mathbf{C}). \tag{4.8}$$

which correspond under the bijection of Thm. **3.1.1** to the maximal and minimal factorizations on $\mathbf{C}$ respectively, and will be called the *trivial* factorizations/$t$-structures.

(3) This is not entirely true, as will appear clear in Ch. **7**, but it's a good approximation of the definition given there. [Bri07] imposes more restrictive conditions to ensure "compactness" of the factorization. Compare also [GKR04].

Hence, by Lemma **4.1.19**, any $J$-family of $t$-structures $\mathfrak{t}\colon J \to \mathrm{TS}(\mathbf{C})$ extends to a $(J\cup\{\pm\infty\})$-family by setting $\mathfrak{t}_{-\infty} = (\mathbf{C}, \mathbf{0})$ and $\mathfrak{t}_{+\infty} = (\mathbf{0}, \mathbf{C})$.

DEFINITION 4.1.29. Let $\mathfrak{t}$ be a $J$-family of $t$-structures. For $i$ and $j$ in $J$ we set

$$\mathbf{C}_{[i,j)} = \mathbf{C}_{\geq i} \cap \mathbf{C}_{<j}. \tag{4.9}$$

Consistently with Remark **4.1.28** and Notation **4.1.24**, we also set

$$\mathbf{C}_{[i,+\infty)} = \mathbf{C}_{\geq i}; \qquad \mathbf{C}_{[-\infty,i)} = \mathbf{C}_{<i} \tag{4.10}$$

for any $i$ in $J$. We say that $\mathbf{C}$ is *$J$-bounded* if

$$\mathbf{C} = \bigcup_{i,j\in J} \mathbf{C}_{[i,j)}. \tag{4.11}$$

Similarly, we say that $\mathbf{C}$ is *$J$-left-bounded* if $\mathbf{C} = \bigcup_{i\in J} \mathbf{C}_{[i,+\infty)}$ and *$J$-right-bounded* if $\mathbf{C} = \bigcup_{i\in J} \mathbf{C}_{[-\infty,i)}$. This notion is well known in the classical as well as in the quasicategorical setting: see [BBD82, Lur17].

REMARK 4.1.30. Since $\mathbf{C}_{[i,j)} = \mathbf{C}_{[i,+\infty)} \cap \mathbf{C}_{[-\infty,j)}$ one immediately sees that $\mathbf{C}$ is $J$-bounded if and only if $\mathbf{C}$ is both $J$-left- and $J$-right-bounded.

REMARK 4.1.31. As it is natural to expect, if $i \geq j$, then $\mathbf{C}_{[i,j)}$ is contractible. Namely, since $j \leq i$ one has $\mathbf{C}_{<j} \subseteq \mathbf{C}_{<i}$ and so

$$\mathbf{C}_{[i,j)} = \mathbf{C}_{\geq i} \cap \mathbf{C}_{<j} \subseteq \mathbf{C}_{\geq i} \cap \mathbf{C}_{<i} = \mathbf{C}_{\geq_i 0} \cap \mathbf{C}_{<_i 0} \tag{4.12}$$

which corresponds to the contractible subcategory of zero objects in $\mathbf{C}$ (this is immediate, in view of the definition of the two classes).

REMARK 4.1.32. Let $\mathfrak{t}$ be a $\mathbb{Z}$-family of $t$-structures on $\mathbf{C}$. Then $\mathbf{C}$ is $\mathbb{Z}$-bounded (resp., $\mathbb{Z}$-left-bounded, $\mathbb{Z}$-right-bounded) if and only if $\mathbf{C}$ is bounded (resp., left-bounded, right-bounded)with respect to the $t$-structure $\mathfrak{t}_0$, agreeing with the classical definition of boundedness as given, e.g., in [BBD82].

REMARK 4.1.33. If $\mathfrak{t}$ is an $\mathbb{R}$-family of $t$-structures on $\mathbf{C}$, then one can define

$$\mathbf{C}_\lambda = \bigcap_{\epsilon>0} \mathbf{C}_{[\lambda,\lambda+\epsilon)}. \tag{4.13}$$

These subcategories $\mathbf{C}_\lambda$ are the *slices* of $\mathbf{C}$ in the terminology of [Bri07].

REMARK 4.1.34. For any $i, j, h, k$ in $J$ with $j \leq h$ one has

$$\mathbf{C}_{[i,j)} \subseteq \mathbf{C}^{\perp}_{[h,k)}, \tag{4.14}$$

i.e., $\mathbf{C}(X, Y)$ is contractible whenever $X \in \mathbf{C}_{[h,k)}$ and $Y \in \mathbf{C}_{[i,j)}$ (one says that $\mathbf{C}_{[i,j)}$ is *right-orthogonal* to $\mathbf{C}_{[h,k)}$, see Notation **1.2.15**). Indeed, since $\mathbf{C}_{<j} = \mathbf{C}_{<_j 0} = \mathbf{C}^{\perp}_{\geq_j 0} = \mathbf{C}^{\perp}_{\geq j}$, and passing to the orthogonal reverses the inclusions, we have

$$\mathbf{C}_{[i,j)} \subseteq \mathbf{C}_{<j} = \mathbf{C}^{\perp}_{\geq j} \subseteq \mathbf{C}^{\perp}_{\geq h} \subseteq \mathbf{C}^{\perp}_{[h,k)}. \tag{4.15}$$

DEFINITION 4.1.35. Let $(\mathbf{C}, \mathfrak{t})$ be a stable $\infty$-category endowed with a $t$-structure, arising from the normal torsion theory $\mathbb{F} = (\mathcal{E}, \mathcal{M})$. For each $n \in \mathbb{Z}$, let $\mathbf{C}_{\geq n}$ and $\mathbf{C}_{<n}$ be the reflective and coreflective subcategories of $\mathbf{C}$ determined by the $t$-structure $\mathfrak{t}$.

Then $\mathfrak{t}$ is said to be

- *bounded* if $\bigcup \mathbf{C}_{\geq n} = \mathbf{C}$;
- *limited* if every $f\colon X \to Y$ fits into a fiber sequence

$$\begin{array}{ccccc} F & \longrightarrow & X & \longrightarrow & 0 \\ {\scriptstyle m[a]}\downarrow & & \downarrow{\scriptstyle f} & & \downarrow{\scriptstyle e[b]} \\ 0 & \longrightarrow & Y & \longrightarrow & C \end{array} \tag{4.16}$$

  where $F = \mathrm{fib}(f), C = \mathrm{cofib}(f)$, and $m[a] \in \mathcal{M}[a], e[b] \in \mathcal{E}[b]$ for suitable integers $a, b \in \mathbb{Z}$;
- *narrow* if $\mathbf{C} = \bigcup_{a \leq b} \mathbf{C}_{[a,b)}$, where $\mathbf{C}_{[a,b)} = \mathbf{C}_{\geq a} \cap \mathbf{C}_{<b}$.

PROPOSITION 4.1.36. Let $(\mathbf{C}, \mathfrak{t})$ be a stable $\infty$-category endowed with a $t$-structure. Then $\mathfrak{t}$ is narrow if and only if it is bounded, if and only if it is limited.

REMARK 4.1.37. We say that an $f\colon X \to Y$ in $(\mathbf{C}, \mathfrak{t})$ is *limited between* $a, b \in \mathbb{Z}$ if there exists a diagram like (**4.16**) for $f$; we say that $f$ is *limited* if it is limited between $a, b$ for some $a, b \in \mathbb{Z}$. In this terminology, a $t$-structure $\mathfrak{t}$ is limited if and only if every $f\colon X \to Y$ is limited with respect to $\mathfrak{t}$.

*Proof.* It is rather obvious that $\mathfrak{t}$ is narrow if and only if it is limited, so we can reduce ourselves to prove that bounded and limited $t$-structures coincide.

This is a consequence of the application of the following

LEMMA 4.1.38. Let $f\colon X \to Y$ be limited between $a, b$; then $f$ belongs to $\mathcal{M}[a+1] \cap \mathcal{E}[b-1]$.

*Proof.* We can reduce the result to an easy consequence of the Sator Lemma **2.3.11**. Moreover, we only prove that $f \in \mathcal{M}[a+1]$, the proof that $f \in \mathcal{E}[b-1]$ being dual.

By the abovementioned Sator Lemma, $\left[\begin{smallmatrix} F \\ \downarrow \\ 0 \end{smallmatrix}\right] \in \mathcal{M}[a]$ if and only if $\left[\begin{smallmatrix} 0 \\ \downarrow \\ F \end{smallmatrix}\right] \in \mathcal{M}[a]$; but now $F \simeq C[-1]$ in diagram (**4.16**), and $\left[\begin{smallmatrix} 0 \\ \downarrow \\ C[-1] \end{smallmatrix}\right] \in \mathcal{M}[a]$ implies that $\left[\begin{smallmatrix} 0 \\ \downarrow \\ C \end{smallmatrix}\right] \in \mathcal{M}[a+1]$. $\square$

Now we can return to the proof of the initial result, implicitly invoking Lemma **4.1.38** when needed: if $\mathfrak{t}$ is a limited $t$-structure, then every $\left[\begin{smallmatrix} X \\ \downarrow \\ 0 \end{smallmatrix}\right]$ is limited between $a_X, b_X$, hence $\left[\begin{smallmatrix} X \\ \downarrow \\ 0 \end{smallmatrix}\right] \in \mathcal{M}[a_X + 1]$, so that $X \in \mathbf{C}_{<a_X}$; in the same way $\left[\begin{smallmatrix}0\\ \downarrow\\ X\end{smallmatrix}\right] \in \mathcal{E}[b_X - 1]$, so that $X \in \mathbf{C}_{\geq b_X - 1}$ and $X \in \bigcup_{u,v} \mathbf{C}_{[u,v)}$. The other inclusion is obvious.

Conversely, if $\mathfrak{t}$ is bounded, we have that each object $X$ lies in $\mathcal{E}[u_X] \cap \mathcal{M}[v_X]$; so if we consider the following diagram of pullout squares

$$\begin{array}{ccccc} Y[-1] & \longrightarrow & 0 & & \\ \downarrow & & \downarrow{\scriptstyle m[v_X]} & & \\ F & \longrightarrow & X & \longrightarrow & 0 \\ \downarrow & & \downarrow & {\scriptstyle m[v_Y]} & \downarrow \\ 0 & \longrightarrow & Y & \longrightarrow & C \end{array} \tag{4.17}$$

we deduce that the arrow $\left[\begin{smallmatrix}F\\ \downarrow\\ 0\end{smallmatrix}\right]$ belongs to $\mathcal{M}[v]$, where $v = \max\{v_X, v_Y\}$, as a consequence of the stability under pullbacks and the 3-for-2 closure property of each class $\mathcal{M}[n]$.

Reasoning in a perfectly dual fashion, we deduce that $\left[\begin{smallmatrix}0\\ \downarrow\\ C\end{smallmatrix}\right] \in \mathcal{E}[u]$, where $u = \min\{u_X, u_Y\}$, so that each $f\colon X \to Y$ is limited between $u, v$. □

## 4.2 Towers of morphisms.

הָבָה, נֵרְדָה, וְנָבְלָה שָׁם, שְׂפָתָם--אֲשֶׁר לֹא יִשְׁמְעוּ, אִישׁ שְׂפַת רֵעֵהוּ.

[ER77], Genesis 11:7

In the remainder of this section, $J$ will be a fixed $\mathbb{Z}$-poset and $\mathfrak{t}_i$ will be the $i^{\text{th}}$ element of a $J$-family of $t$-structures on $\mathbf{C}$; $\mathbb{F}_i$ will denote the corresponding $J$-family of factorization systems.

We recall Def. **1.5.10**, and in particular that

Lemma 4.2.1. *The chain $i_1 \leq i_2 \leq \cdots \leq i_k$ determines a $k$-fold factorization system in $\mathbf{C}$ in the sense of Def. **1.5.10**. Namely, every arrow $f\colon X \to Y$ in $\mathbf{C}$ can be uniquely factored into a composition*

$$X \xrightarrow{\mathcal{E}_{i_k}} Z_{i_k} \xrightarrow{\mathcal{E}_{i_{k-1}} \cap \mathcal{M}_{i_k}} Z_{i_{k-1}} \to \cdots \to Z_{i_2} \xrightarrow{\mathcal{E}_{i_1} \cap \mathcal{M}_{i_2}} Z_{i_1} \xrightarrow{\mathcal{M}_{i_1}} Y. \tag{4.18}$$

*Proof.* This was proved in Lemma **1.5.11**. □

Lemma 4.2.2. *Let $i, j$ be elements in $J$ and let $X$ be an object in $\mathbf{C}_{\geq j}$ (see Definition **4.1.29**). If a morphism $f\colon X \to Y$ is in $\mathcal{E}_i \cap \mathcal{M}_j$, then $\operatorname{cofib}(f)$ is in $\mathbf{C}_{[i,j)}$.*

*Proof.* Since $X$ is in $\mathbf{C}_{\geq j}$, $0 \to X \xrightarrow{f} Y$ is the $(\mathcal{E}_j, \mathcal{M}_j)$-factorization of $0 \to Y$ (in particular, $X \cong S_j Y$ if $S_j$ denotes the coreflection of $\mathbf{C}$ on $\mathbf{C}_{\geq j}$; see our Def. **2.3.2** of "firm reflectivity" and [RT07, Prop **3.2**]). Since the factorization system $\mathbb{F}_j$ is normal, hence semi-right-exact, we have the following pullout diagram:

$$\begin{array}{ccc} X & \xrightarrow{\mathcal{E}_j} & Y \\ {\scriptstyle \mathcal{E}_j}\big\downarrow & & \big\downarrow{\scriptstyle \mathcal{M}_j} \\ 0 & \underset{\mathcal{M}_j}{\longrightarrow} & \operatorname{cofib}(f) \end{array} \tag{4.19}$$

Hence $\operatorname{cofib}(f)$ is in $\mathbf{C}_{<j}$. On the other hand, $f$ is in $\mathcal{E}_i$, which is closed under pushouts, and so $0 \to \operatorname{cofib}(f)$ is in $\mathcal{E}_i$, i.e., $\operatorname{cofib}(f)$ is in $\mathbf{C}_{\geq i}$. □

An immediate corollary of **4.2.2** is that the cofibers of each $f_j \colon Y_{i_j} \to Y_{i_{j-1}}$ in the $k$-fold factorization obtained via **4.2.1** belong to the subcategories $\mathbf{C}_{[i_{j-1}, i_j)}$. This remark is the basic building block of the *tower* of $f$.

COROLLARY 4.2.3. Let $i_1 \leq i_2 \leq \cdots \leq i_k$ an ascending chain in $J$. Then for any object $Y$ in $\mathbf{C}$, the arrows $f_j \colon Y_{i_j} \to Y_{i_{j-1}}$ in the $k$-fold factorization of the initial morphism $0 \to Y$ are such that $\operatorname{cofib}(f_j) \in \mathbf{C}_{[i_{j-1}, i_j)}$, where we have set $i_{k+1} = +\infty$ and $Y_{+\infty} = 0$ (and, similarly, $i_0 = -\infty$ and $Y_{-\infty} = Y$) consistently with Remark **4.1.28** (and its dual).

*Proof.* From the $k$-fold factorization

$$0 \xrightarrow{\mathcal{E}_{i_k}} Y_{i_k} \xrightarrow{\mathcal{E}_{i_{k-1}} \cap \mathcal{M}_{i_k}} Y_{i_{k-1}} \to \cdots \to Y_{i_2} \xrightarrow{\mathcal{E}_{i_1} \cap \mathcal{M}_{i_2}} Y_{i_1} \xrightarrow{\mathcal{M}_{i_1}} Y, \tag{4.20}$$

and from the fact that $\mathcal{E}_{i_1} \supseteq \mathcal{E}_{i_2} \supseteq \cdots \supseteq \mathcal{E}_{i_k}$ and each class $\mathcal{E}_{i_j}$ is closed for composition, we see that $Y_{i_j}$ is in $\mathbf{C}_{i_j}$ and the previous lemma applies. □

Firm reflectivity implies the converse of **4.2.2**:

LEMMA 4.2.4. Let $i \leq j$ be elements in $J$ and let $f \colon X \to Y$ be a morphism in $\mathbf{C}$. If $X$ is in $\mathbf{C}_{[j,+\infty)}$ and $\operatorname{cofib}(f)$ is in $\mathbf{C}_{[i,j)}$ then $0 \to X \xrightarrow{f} Y$ is the $(\mathcal{E}_j, \mathcal{M}_j)$-factorization of the initial morphism $0 \to Y$ and $Y$ is in $\mathbf{C}_{[i,+\infty)}$. In particular $f$ is in $\mathcal{E}_i \cap \mathcal{M}_j$.

*Proof.* Since $X$ is in $\mathbf{C}_{\geq j}$, the morphism $0 \to X$ is in $\mathcal{E}_j$, and so (reasoning up to equivalence) to show that $0 \to X \to Y$ is the $(\mathcal{E}_j, \mathcal{M}_j)$-factorization of $0 \to Y$ we are reduced to showing that $f \colon X \to Y$ is in $\mathcal{M}_j$. Since $\operatorname{cofib}(f)$ is in $\mathbf{C}_{[i,j)}$, we have in particular that $\operatorname{cofib}(f) \to 0$ is in $\mathcal{M}_j$ and so $0 \to \operatorname{cofib}(f)$ is in $\mathcal{M}_j$ by the Sator lemma. Then we have a homotopy

pullback diagram

$$\begin{array}{ccc} X & \xrightarrow[f]{\quad} & 0 \\ \big\downarrow & & \big\downarrow{\scriptstyle \mathcal{M}_j} \\ Y & \longrightarrow & \mathrm{cofib}(f) \end{array} \tag{4.21}$$

and so $f$ is in $\mathcal{M}_j$ by the fact that $\mathcal{M}_j$ is closed under pullbacks.

To show that also $f \in \mathcal{E}_i$ let $0 \to X \to T \to Y$ be the ternary factorization of $f$. We can consider the diagram

$$\begin{array}{ccccccc} 0 & & & & & & \\ {\scriptstyle \mathcal{E}_j}\big\downarrow & & & & & & \\ X & \xrightarrow{\mathcal{E}_j} & 0 & & & & \\ {\scriptstyle \mathcal{E}_i\cap\mathcal{M}_j}\big\downarrow & & \big\downarrow{\scriptstyle \mathcal{E}_i\cap\mathcal{M}_j} & & & & \\ T & \xrightarrow{\mathcal{E}_j} & U & \xrightarrow{\mathcal{E}_i\cap\mathcal{M}_j} & 0 & & \\ {\scriptstyle \mathcal{M}_i}\big\downarrow & & \big\downarrow{\scriptstyle \mathcal{M}_i} & & \big\downarrow{\scriptstyle \mathcal{M}_i} & & \\ Y & \xrightarrow[\mathcal{E}_j]{} & \mathrm{cofib}(f) & \xrightarrow[\mathcal{E}_i\cap\mathcal{M}_j]{} & V & \xrightarrow[\mathcal{M}_i]{} & 0 \end{array} \tag{4.22}$$

where all the squares are pullouts, and where we have used the Sator lemma, the fact that the classes $\mathcal{E}$ are closed for pushouts while the classes $\mathcal{M}$ are closed for pullbacks, and the 3-for-2 property for both classes.

□

Lemma 4.2.5. Let $Y$ an object in $\mathbf{C}$ and let $i_1 \leq i_2 \leq \cdots \leq i_k$ be an ascending chain in $J$. If a factorization

$$0 \xrightarrow{f_{k+1}} Y_{i_k} \xrightarrow{f_k} Y_{i_{k-1}} \to \cdots \to Y_{i_2} \xrightarrow{f_2} Y_{i_1} \xrightarrow{f_1} Y, \tag{4.23}$$

of the initial morphism $0 \to Y$ is such that $\mathrm{cofib}(f_j)$ is in $\mathbf{C}_{[i_{j-1}, i_j)}$ (with $i_{k+1} = +\infty$ and $i_0 = -\infty$) then this factorization is the $k$-fold factorization of $0 \to Y$ associated with the chain $i_1 \leq \cdots \leq i_k$.

*Proof.* By uniqueness of the $k$-fold factorization we only need to prove that $f_j \in \mathcal{E}_{i_{k-1}} \cap \mathcal{M}_{i_k}$, which is immediate by repeated application of Lemma **4.2.4**. □

This paves the way to the definition of the tower of $f$: the basic idea is to “pull back” the factorization of the initial morphism $0 \to \mathrm{cofib}(f)$ using Lemma **4.2.5**.

Definition 4.2.6. (Tower of a morphism): Let $f\colon X \to Y$ be a morphism in $\mathbf{C}$ and let $i_1 \le i_2 \le \cdots \le i_k$ be an ascending chain in $J$. We say that a factorization

$$X \xrightarrow{f_{k+1}} Z_{i_k} \xrightarrow{f_k} Z_{i_{k-1}} \to \cdots \to Z_{i_2} \xrightarrow{f_2} Z_{i_1} \xrightarrow{f_1} Y, \tag{4.24}$$

of $f$ is a *tower* of $f$ relative to the chain $\{i_j\} = \{i_1 \le i_2 \le \cdots \le i_k\}$ if for any $j = 1, \dots, k+1$ one has $\operatorname{cofib}(f_j) \in \mathbf{C}_{[i_{j-1}, i_j)}$ (with $i_{k+1} = +\infty$ and $i_0 = -\infty$).

Proposition 4.2.7. Let $f\colon X \to Y$ be a morphism in $\mathbf{C}$ and let $i_1 \le i_2 \le \cdots \le i_k$ be an ascending chain in $J$. Then a tower for $f$ relative to $\{i_j\}$, denoted $\text{♖}_{\{i_j\}}(f)$, exists and it is unique up to isomorphisms.

*Proof.* We split the proof in two parts: existence and uniqueness of the tower;

(1) Consider the pullout diagram

$$\begin{array}{ccc} X & \xrightarrow[f]{} & 0 \\ \big\downarrow & & \big\downarrow \mathcal{M}_j \\ Y & \longrightarrow & \operatorname{cofib}(f) \end{array} \tag{4.25}$$

By Corollary **4.2.3**, the $k$-fold factorization

$$0 \xrightarrow{\varphi_{k+1}} A_{i_k} \xrightarrow{\varphi_k} A_{i_{k-1}} \to \cdots \to A_{i_2} \xrightarrow{\varphi_2} A_{i_1} \xrightarrow{\varphi_1} \operatorname{cofib}(f) \tag{4.26}$$

of the initial morphism $0 \to \operatorname{cofib}(f)$ is such that $\operatorname{cofib}(\varphi_{i_j}) \in \mathbf{C}_{[i_{j-1}, i_j)}$. Pulling back this factorization along $Y \to \operatorname{cofib}(f)$ we obtain a factorization

$$\begin{array}{ccc} X & \longrightarrow & 0 \\ {\scriptstyle f_{k+1}}\big\downarrow & & \big\downarrow{\scriptstyle \phi_{k+1}} \\ Z_{i_k} & \longrightarrow & A_{i_k} \\ {\scriptstyle f_k}\big\downarrow & & \big\downarrow{\scriptstyle \phi_k} \\ \vdots & & \vdots \\ {\scriptstyle f_2}\big\downarrow & & \big\downarrow{\scriptstyle \phi_2} \\ Z_{i_1} & \longrightarrow & A_{i_1} \\ {\scriptstyle f_1}\big\downarrow & & \big\downarrow{\scriptstyle \phi_1} \\ Y & \longrightarrow & \operatorname{cofib}(f) \end{array} \tag{4.27}$$

of $f$, and the pasting of pullout diagrams

$$\begin{array}{ccccc} Z_{i_j} & \longrightarrow & A_{i_j} & \longrightarrow & 0 \\ {\scriptstyle f_j}\big\downarrow & & \big\downarrow{\scriptstyle \phi_j} & & \big\downarrow \\ Z_{i_{j-1}} & \to & A_{i_{j-1}} & \rightarrowtail & \operatorname{cofib}(\varphi_j) \end{array} \tag{4.28}$$

shows that $\operatorname{cofib}(f_j) = \operatorname{cofib}(\varphi_j)$ and so $\operatorname{cofib}(f_j) \in \mathbf{C}_{[i_{j-1}, i_j)}$. This proves the existence of the tower.

(2) To prove uniqueness, start with a tower $♖_{\{i_j\}}(f)$ for $f$ and push it out along $Y \to \operatorname{cofib}(f)$ to obtain a tower for the initial morphism $0 \to \operatorname{cofib}(f)$. By Lemma **4.2.5**, this is the $k$-fold factorization of $0 \to \operatorname{cofib}(f)$ associated with the chain $\{i_j\}$ and so $♖_{\{i_j\}}(f)$ is precisely the tower constructed in the first part of the proof. Note how the pullout axiom of stable $\infty$-categories plays a crucial role. □

Remark 4.2.8. A tower for $f$ relative to an ascending chain $\{i_j\}$ can be equivalently defined as a factorization of $f$ such that $\operatorname{fib}(f) \in \mathbf{C}_{[i_{j-1}-1, i_j-1)}$, for any $j = 0, \dots, k+1$.

Remark 4.2.9. It's an unavoidable temptation to think of the tower $♖_{\{i_j\}}(f)$ relative to an ascending chain $\{i_j\}$ as the $k$-fold factorization of $f$ associated with the chain $\{i_j\}$.

As the following counterexample shows, when $f$ is not an initial morphism this is in general not true.(4) Let $J = \mathbb{Z}$ and take an ascending chain consisting of solely the element 0. Now take a morphism $f\colon X \to Y$ between two elements in $\mathbf{C}_{[-1,0)}$. The object $\operatorname{cofib}(f)$ will lie in $\mathbf{C}_{[-1,+\infty)}$, since $\mathcal{E}_{-1}$ is closed for pushouts, but in general it will not be an element in $\mathbf{C}_{[0,+\infty)}$. In other words, we will have, in general, a nontrivial $(\mathcal{E}_0, \mathcal{M}_0)$-factorization of the initial morphism $0 \to \operatorname{cofib}(f)$. Pulling this back along $Y \to \operatorname{cofib}(f)$ we obtain the tower $X \xrightarrow{f_2} Z \xrightarrow{f_1} Y$ of $f$, and this factorization will be nontrivial since its pushout is nontrivial. It follows that $(f_2, f_1)$, cannot be the $(\mathcal{E}_0, \mathcal{M}_0)$-factorization of $f$. Indeed, by the 3-for-2 property of $\mathcal{M}_0$, the morphism $f$ is in $\mathcal{M}_0$, so its $(\mathcal{E}_0, \mathcal{M}_0)$-factorization is trivial.

(4) When $f\colon A \to 0$ is the terminal morphism, our notation and construction is in line with the classical [Lur17], where the "Postnikov tower" of $A$ is the sequence

$$A \to \cdots \to R_2 A \to R_1 A \to R_0 A \to 0 \tag{4.29}$$

of factorizations obtained from the (stable image of) the $n$-connected factorization system of [Joy08].

## 4.3 Hearts of $t$-structures.

> I watched a snail crawl along the edge of a straight razor. That's my dream. That's my nightmare. Crawling, slithering, along the edge of a straight razor... and surviving.
>
> Col. Walter E. Kurtz

We now focus in the case $J = \mathbb{Z}$. As indicated in remark **4.1.9** this is equivalent to a single distinguished $t$-structure $\mathfrak{t} = \mathfrak{t}_0$ on the stable $\infty$-category $\mathbf{C}$, together with its orbit $\{\mathfrak{t}_j = \mathfrak{t}_0[j]\}_{j\in\mathbb{Z}}$. As the set of indices for our family of $t$-structures is the ordered set of integers, we will always consider complete ascending chains of the form

$$n < n+1 < n+2 < \cdots n+k-1 \tag{4.30}$$

in what follows. In particular Proposition **4.2.7** becomes

PROPOSITION 4.3.1. Let $f\colon X \to Y$ be a morphism in $\mathbf{C}$. Then for any integer $n$ and any positive integer $k$ there exists a unique tower for $f$ associated with the ascending chain $n < n+1 < \cdots < n+k-1$. Denoting this tower by

$$X \xrightarrow{f_{n+k}} Z_{n+k-1} \xrightarrow{f_{n+k-1}} Z_{n+k-2} \to \cdots \to Z_{n+1} \xrightarrow{f_n} Z_n \xrightarrow{f_{n-1}} Y, \tag{4.31}$$

one has $\mathrm{cofib}(f_j) \in \mathbf{C}_{[j,j+1)}$ for any $j = n, \ldots, n+k-1$, $\mathrm{cofib}(f_{n-1}) \in \mathbf{C}_{<n}$ and $\mathrm{cofib}(f_{n+k}) \in \mathbf{C}_{\geq n+k}$.

Since $\mathbf{C}_{[j,j+1)} = \mathbf{C}_{[0,1)}[j]$ for any $j \in \mathbb{Z}$, the above Proposition suggests to focus on the subcategory $\mathbf{C}_{[0,1)}$ of $\mathbf{C}$. This subcategory has a special name and special properties (it is an *abelian* subcategory).

DEFINITION 4.3.2. Let $\mathbf{C}$ be a stable $\infty$-category equipped with a $t$-structure $\mathfrak{t} = (\mathbf{C}_{\geq 0}, \mathbf{C}_{<0})$; the *heart* $\mathbf{C}^{\heartsuit}$ of $\mathfrak{t}$ is the subcategory $\mathbf{C}_{[0,1)}$ defined following Def. **4.1.29**.

REMARK 4.3.3. There is a rather evocative pictorial representation of the heart of a $t$-structure, manifestly inspired by [Bri07]: if we depict $\mathbf{C}_{<0}$ and $\mathbf{C}_{\geq 0}$ as contiguous half-planes, like in the following picture, then the action of the shift functor can be represented as an horizontal shift, and the closure properties of the two classes $\mathbf{C}_{\geq 0}, \mathbf{C}_{<0}$ under positive and negative shifts are a direct consequence of the shape of these areas. With these notations, an object $Z$ is in the heart of $\mathfrak{t}$ if it lies in a "boundary region", i.e. if it lies in $\mathbf{C}_{\geq 0}$, but $Z[-1]$ lies in $\mathbf{C}_{<0}$.

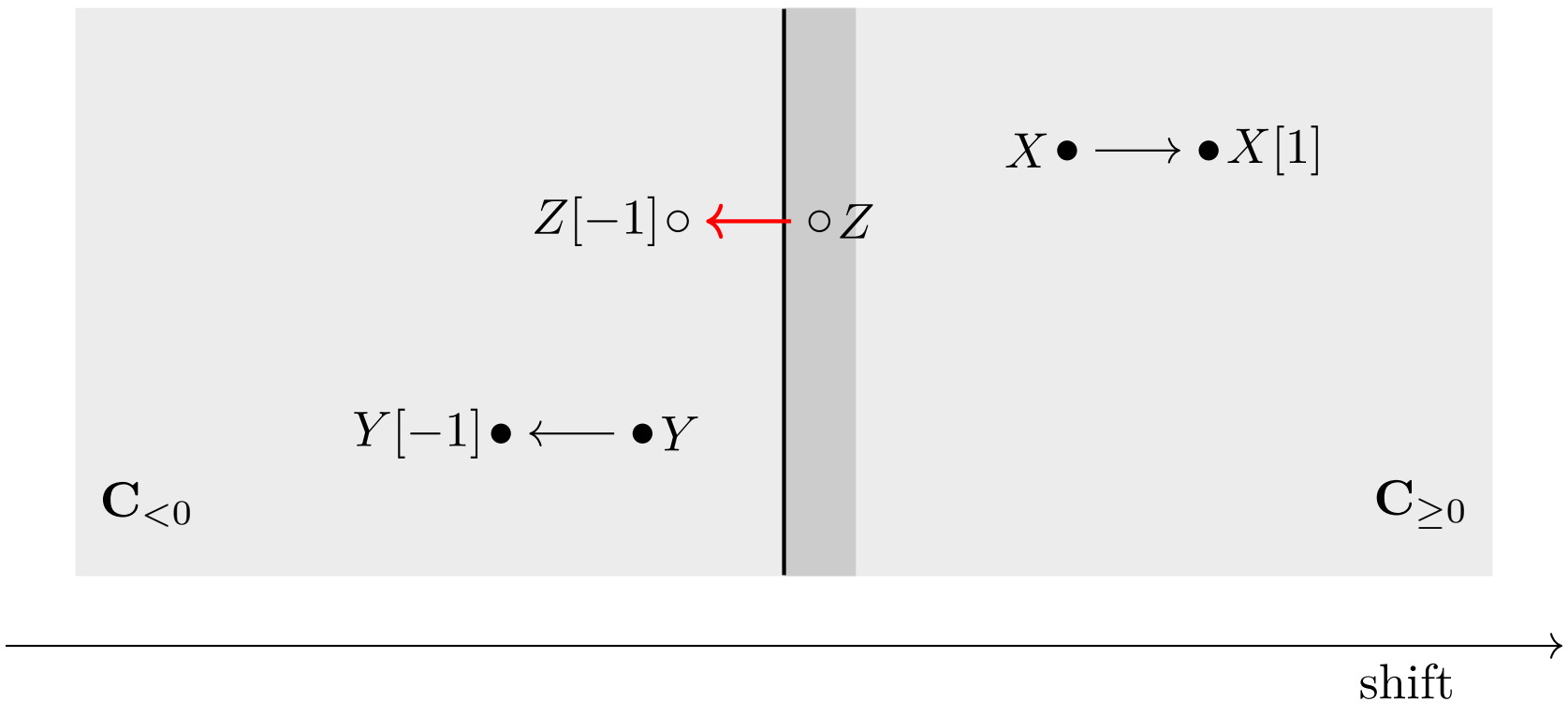

Figure 4.1: Heart of a $t$-structure

Having introduced this notation, we can rephrase the existence of the tower for $f$ as follows: given a morphism $f\colon X \to Y$ in $\mathbf{C}$, for any integer $n$ and any positive integer $k$ there exists a unique factorization of $f$

$$X \xrightarrow{f_{n+k}} Z_{n+k-1} \xrightarrow{f_{n+k-1}} Z_{n+k-2} \to \cdots \to Z_{n+1} \xrightarrow{f_n} Z_n \xrightarrow{f_{n-1}} Y, \quad (4.32)$$

such that $\mathrm{cofib}(f_j) \in \mathbf{C}^{\heartsuit}[j]$ for any $j = n, \dots, n+k-1$, $\mathrm{cofib}(f_{n-1}) \in \mathbf{C}_{<n}$ and $\mathrm{cofib}(f_{n+k}) \in \mathbf{C}_{\geq n+k}$.

The content of this statement becomes more interesting when $\mathbf{C}$ is *bounded* with respect to the $t$-structure $\mathfrak{t}$ (see Definition **4.1.29**). If $\mathbf{C}$ is bounded, then the $(\mathcal{E}_n, \mathcal{M}_n)$-factorizations of an initial morphism $0 \to Y$ are trivial (see Definition **4.1.29** and the subsequent Remark) for $|n| \gg 0$.

As an immediate consequence, the morphisms $X \xrightarrow{f_{n+k}} Z_{n+k-1}$ and $Z_n \xrightarrow{f_{n-1}} Y$ in the tower of $f$ associated with the chain $n < n+1 < \cdots < n+k-1$ are isomorphisms for $n \ll 0$ and $k \gg 0$. One notices, as it is obvious, that the class of isomorphisms in $\mathbf{C}$ is closed under transfinite composition this leads to the following

Proposition 4.3.4. Let $\mathbf{C}$ be a stable $\infty$-category which is bounded with respect to a given $t$-structure $\mathfrak{t}$. Then for any morphism $f\colon X \to Y$ in $\mathbf{C}$ there exists an integer $n_0$ and a positive integer $k_0$ such that for any integer $n \leq n_0$ and any positive integer $k$ with $k \geq n_0 - n + k_0$ there exists a unique factorization of $f$

$$X \xrightarrow{\sim} Z_{n+k-1} \xrightarrow{f_{n+k-1}} Z_{n+k-2} \to \cdots \to Z_{n+1} \xrightarrow{f_n} Z_n \xrightarrow{\sim} Y \quad (4.33)$$

such that $\text{cofib}(f_j) \in \mathbf{C}^\heartsuit[j]$ for any $j = n, \dots, n+k-1$.

REMARK 4.3.5. By uniqueness in Proposition **4.3.4**, one has a well defined $\mathbb{Z}$-factorization

$$X = \lim(Z_j) \to \cdots \to Z_{j+1} \xrightarrow{f_j} Z_j \xrightarrow{f_{j-1}} Z_{j-1} \to \cdots \to \text{colim}(Z_j) = Y \tag{4.34}$$

with with $j$ ranging over the integers, $\text{cofib}(f_j) \in \mathbf{C}^\heartsuit[j]$ for any $j \in \mathbb{Z}$ and with $f_m$ being an isomorphism for $|j| \gg 0$. We will refer to this factorization as the $\mathbb{Z}$-*tower* of $f$. Notice how the boundedness of $\mathbf{C}$ has played an essential role: when $\mathbf{C}$ is not bounded, one still has towers for any finite ascending chain, but in general they do not stabilize.

REMARK 4.3.6. Since we know that the tower of an initial morphism is its $k$-fold $(\mathcal{E}_j, \mathcal{M}_j)$-factorization, we see that in a stable $\infty$-category $\mathbf{C}$ which is bounded with respect to a $t$-structure $\mathfrak{t} = (\mathbf{C}_{\geq 0}, \mathbf{C}_{<0})$ the $\mathbb{Z}$-tower of $0 \to Y$,

$$0 = \lim(Y_j) \to \cdots \to Y_{j+1} \xrightarrow{f_j} Y_j \xrightarrow{f_{j-1}} Y_{j-1} \to \cdots \to \text{colim}(Y_j) = Y \tag{4.35}$$

is such that $f_j \in \mathcal{E}_j \cap \mathcal{M}_{j+1}$ for any $j \in \mathbb{Z}$. It follows that an object $Y$ is in $\mathbf{C}_{\geq 0}$ if and only if the $\mathbb{Z}$-tower of $0 \to Y$ satisfies $\text{cofib}(f_j) = 0$ for any $j < 0$, while $Y$ is in $\mathbf{C}_{<0}$ if and only if $\text{cofib}(f_j) = 0$ for any $j \geq 0$.

### 4.3.1 Abelianity of the heart.

In the following section we present a complete proof, in the stable setting, of the fact that the heart of a $t$-structure, as defined in [Lur17, Def. **1.2.1.11**], is an abelian $\infty$-category.

In other words, $\mathbf{C}^\heartsuit$ is homotopy equivalent to its homotopy category $h\mathbf{C}^\heartsuit$, which is an abelian category; this is the higher-categorical counterpart of a classical result, first proved in [BBD82, Thm. **1.3.6**], which only relies on properties stated in terms of normal torsion theories in a stable $\infty$-category. We begin with the following

DEFINITION 4.3.7. (ABELIAN $\infty$-CATEGORY): An *abelian $\infty$-category* is a quasicategory $\mathbf{A}$ such that

(1) the hom space $\mathbf{A}(X, Y)$ is a homotopically discrete infinite loop space for any $X, Y$, i.e., there exists an infinite sequence of $\infty$-groupoids $Z_0, Z_1, Z_2, \dots$, with $Z_0 \cong \mathbf{C}(X, Y)$ and homotopy equivalences $Z_i \cong \Omega Z_{i+1}$ for any $i \geq 0$, such that $\pi_n Z_0 = 0$ for any $n \geq 1$;
(2) $\mathbf{A}$ has a zero object, (homotopy) kernels, cokernels and biproducts;
(3) for any morphism $f$ in $\mathbf{A}$, the natural morphism from the *coimage* of $f$ to the *image* (see Definition **4.3.15**) of $f$ is an equivalence.

REMARK 4.3.8. Axiom (i) is the homotopically-correct version of $\mathbf{A}(X, Y)$ being an abelian group. For instance, if the abelian group is $\mathbb{Z}$, then the corresponding homotopy discrete space is the Eilenberg-Mac Lane spectrum $\mathbb{Z}, K(\mathbb{Z},1), K(\mathbb{Z},2), \dots$. The homotopy category of such an $\mathbf{A}$ is an abelian category in the classical sense (note that $\mathbf{A}(X,Y)$ being homotopically discrete is necessary in order that kernels and cokernels in $\mathbf{A}$ induce kernels and cokernels in $h\mathbf{A}$). Moreover, since the hom spaces $\mathbf{A}(X,Y)$ are homotopically discrete, the natural morphism $\mathbf{A} \to h\mathbf{A}$ is actually an equivalence.

The rest of the section is devoted to the proof of the following result:

Theorem 4.3.9. The heart $\mathbf{C}^{\heartsuit}$ of a $t$-structure $\mathfrak{t}$ on a stable $\infty$-category $\mathbf{C}$ is an abelian $\infty$-category; its homotopy category $h\mathbf{C}^{\heartsuit}$ is the abelian category arising as the heart of the $t$-structure $h(\mathfrak{t})$ on the triangulated category $h\mathbf{C}$.

Lemma 4.3.10. For any $X$ and $Y$ in $\mathbf{C}^{\heartsuit}$, the hom space $\mathbf{C}^{\heartsuit}(X,Y)$ is a homotopically discrete infinite loop space.

*Proof.* Since $\mathbf{C}^{\heartsuit}$ is a full subcategory of $\mathbf{C}$, we have $\mathbf{C}^{\heartsuit}(X,Y) = \mathbf{C}(X,Y)$, which is an infinite loop space since $\mathbf{C}$ is a stable $\infty$-category.

So we are left to prove that $\pi_n\mathbf{C}(X,Y) = 0$ for $n \geq 1$. Since $\pi_n\mathbf{C}(X,Y) = \pi_{n-1}\Omega\mathbf{C}(X,Y) = \pi_{n-1}\mathbf{C}(X,Y[-1])$, this is equivalent to showing that $\mathbf{C}(X,Y[-1])$ is contractible. Since $X$ and $Y$ are objects in $\mathbf{C}^{\heartsuit}$, we have $X \in \mathbf{C}_{[0,1)}$ and $Y[-1] \in \mathbf{C}_{[-1,0)}$. But $\mathbf{C}_{[-1,0)}$ is right object-orthogonal to $\mathbf{C}_{[0,1)}$ (see Remark **4.1.34**), therefore $\mathbf{C}(X,Y[-1])$ is contractible. □

The subcategory $\mathbf{C}^{\heartsuit}$ inherits the 0 object and biproducts (in fact, all finite limits) from $\mathbf{C}$, so in order to prove it is is abelian we are left to prove that it has kernels and cokernels, and that the canonical morphism from the coimage to the image is an equivalence.

Lemma 4.3.11. Let $f\colon X \to Y$ be a morphism in $\mathbf{C}^{\heartsuit}$. Then $\mathrm{fib}(f)$ is in $\mathbf{C}_{<1}$ and $\mathrm{cofib}(f)$ is in $\mathbf{C}_{\geq 0}$.

*Proof.* Since both $X \to 0$ and $Y \to 0$ are in $\mathcal{M}[1]$, by the 3-for-2 property also $f$ is in $\mathcal{M}[1]$. Since $\mathcal{M}[1]$ is closed for pullbacks, $\mathrm{fib}(f) \to 0$ is in $\mathcal{M}[1]$ and so $\mathrm{fib}(f)$ is in $\mathbf{C}_{<1}$. The proof for $\mathrm{cofib}(f)$ is completely dual. □

Definition 4.3.12. Denote by

$$0 \xrightarrow{\mathcal{E}} \ker(f) \xrightarrow{\mathcal{M}} \mathrm{fib}(f) \tag{4.36}$$

the $(\mathcal{E},\mathcal{M})$-factorization of the morphism $0 \to \mathrm{fib}(f)$ and by

$$\mathrm{cofib}(f) \xrightarrow{\mathcal{E}[1]} \mathrm{coker}(f) \xrightarrow{\mathcal{M}[1]} 0 \tag{4.37}$$

the $(\mathcal{E}[1],\mathcal{M}[1])$-factorization of the morphism $\mathrm{cofib}(f) \to 0$. We call $S\mathrm{fib}(f) = \ker(f)$ and $R_{[1]}\mathrm{cofib}(f) = \mathrm{coker}(f)$ respectively the *kernel* and the *cokernel* of $f$ in $\mathbf{C}^{\heartsuit}$.

REMARK 4.3.13. Since $\mathrm{cofib}(f)[-1] \cong \mathrm{fib} f$, one can equivalently define $\mathrm{coker}(f)$ by declaring the $(\mathcal{E}, \mathcal{M})$-factorization of $\mathrm{fib}(f) \to 0$ to be $\mathrm{fib}(f) \xrightarrow{\mathcal{E}} \mathrm{coker}(f)[-1] \xrightarrow{\mathcal{M}} 0$. Similarly, one can define $\ker(f)$ by declaring the $(\mathcal{E}[1], \mathcal{M}[1])$-factorization of $0 \to \mathrm{cofib}(f)$ to be $0 \xrightarrow{\mathcal{E}[1]} \ker(f)[1] \xrightarrow{\mathcal{M}[1]} \mathrm{cofib}(f)$. By normality of the factorization system we therefore have the homotopy commutative diagram

$$
\begin{array}{ccccc}
0 & \xrightarrow{\mathcal{E}} & \ker(f) & \xrightarrow{\mathcal{M}} & \mathrm{fib}(f) \\
 & & \downarrow{\scriptstyle \mathcal{E}} & & \downarrow{\scriptstyle \mathcal{E}} \\
 & & 0 & \xrightarrow{\mathcal{M}} & \mathrm{coker}(f)[-1] \\
 & & & & \downarrow{\scriptstyle \mathcal{M}} \\
 & & & & 0
\end{array}
\tag{4.38}
$$

whose square sub-diagram is a homotopy pullout.

LEMMA 4.3.14. Both $\ker(f)$ and $\mathrm{coker}(f)$ are in $\mathbf{C}^{\heartsuit}$.

*Proof.* By construction $\ker(f)$ is in $\mathbf{C}_{\geq 0}$, so we only need to show that $\ker(f)$ is in $\mathbf{C}_{<1}$. By definition of $\ker(f)$, we have that $\ker(f) \to \mathrm{fib}(f)$ is in $\mathcal{M}$. Since $\mathcal{M}[-1] \subseteq \mathcal{M}$, we have that also $\ker(f)[-1] \to \mathrm{fib}(f)[-1]$ is in $\mathcal{M}$. By Lemma **4.3.11**, $\mathrm{fib}(f)[-1] \to 0$ is in $\mathcal{M}$ and so we find that also $\ker(f)[-1] \to 0$ is in $\mathcal{M}$. The proof for $\mathrm{coker}(f)$ is perfectly dual. □

By definition of $\ker(f)$ and $\mathrm{coker}(f)$, the defining diagram of $\mathrm{fib}(f)$ and $\mathrm{cofib}(f)$ can be enlarged as

$$
\begin{array}{ccccccccc}
0 & \longrightarrow & \ker(f) & \rightarrow & \mathrm{fib}(f) & \longrightarrow & X & \longrightarrow & 0 \\
 & & & & \downarrow & & \downarrow{\scriptstyle f} & & \downarrow \\
 & & & & 0 & \longrightarrow & Y & \longrightarrow & \mathrm{cofib}(f) & \rightarrow & \mathrm{coker}(f) & \longrightarrow & 0
\end{array}
\qquad \ker(f) \xrightarrow{k_f} X, \quad Y \xrightarrow{c_f} \mathrm{coker}(f)
\tag{4.39}
$$

where $k_f$ and $c_f$ are morphisms in $\mathbf{C}^{\heartsuit}$.

DEFINITION 4.3.15. Let $f \colon X \to Y$ be a morphism in $\mathbf{C}^{\heartsuit}$. The *image* $\mathrm{im}(f)$ and the *coimage* $\mathrm{coim}(f)$ of $f$ are defined as $\mathrm{im}(f) = \ker(c_f)$ and $\mathrm{coim}(f) = \mathrm{coker}(k_f)$.

The following lemma shows that $\ker(f)$ does indeed have the defining property of a kernel:

Lemma 4.3.16. The homotopy commutative diagram

$$\begin{array}{ccc} \ker(f) & \xrightarrow{k_f} & X \\ \downarrow & & \downarrow{\scriptstyle f} \\ 0 & \longrightarrow & Y \end{array} \tag{4.40}$$

is a pullback diagram in $\mathbf{C}^\heartsuit$.

*Proof.* A homotopy commutative diagram

$$\begin{array}{ccc} K & \longrightarrow & X \\ \downarrow & & \downarrow{\scriptstyle f} \\ 0 & \longrightarrow & Y \end{array} \tag{4.41}$$

between objects in the heart is in particular a homotopy commutative diagram in $\mathbf{C}$ so it is equivalent to the datum of a morphism $k' \colon K \to \mathrm{fib}(f)$ in $\mathbf{C}$, with $K$ an object in $\mathbf{C}^\heartsuit$. By the orthogonality of $(\mathcal{E}, \mathcal{M})$, this is equivalent to a morphism $\tilde{k} \colon K \to \ker(f)$:

$$\begin{array}{ccc} 0 & \longrightarrow & \ker(f) \\ {\scriptstyle \mathcal{E}}\downarrow & \overset{\tilde{k}}{\nearrow} & \downarrow{\scriptstyle \mathcal{M}} \\ K & \xrightarrow[k']{} & \mathrm{fib}(f) \end{array} . \tag{4.42}$$

□

There is, obviously, a dual result showing that $\mathrm{coker}(f)$ is indeed a cokernel.

Lemma 4.3.17. The homotopy commutative diagram

$$\begin{array}{ccc} X & \longrightarrow & 0 \\ {\scriptstyle f}\downarrow & & \downarrow \\ Y & \xrightarrow{c_f} & \mathrm{coker}(f) \end{array} \tag{4.43}$$

is a pushout diagram in $\mathbf{C}^\heartsuit$.

Lemma 4.3.18. For $f \colon X \to Y$ a morphism in $\mathbf{C}$, there is a homotopy commutative diagram where all squares are homotopy pullouts:

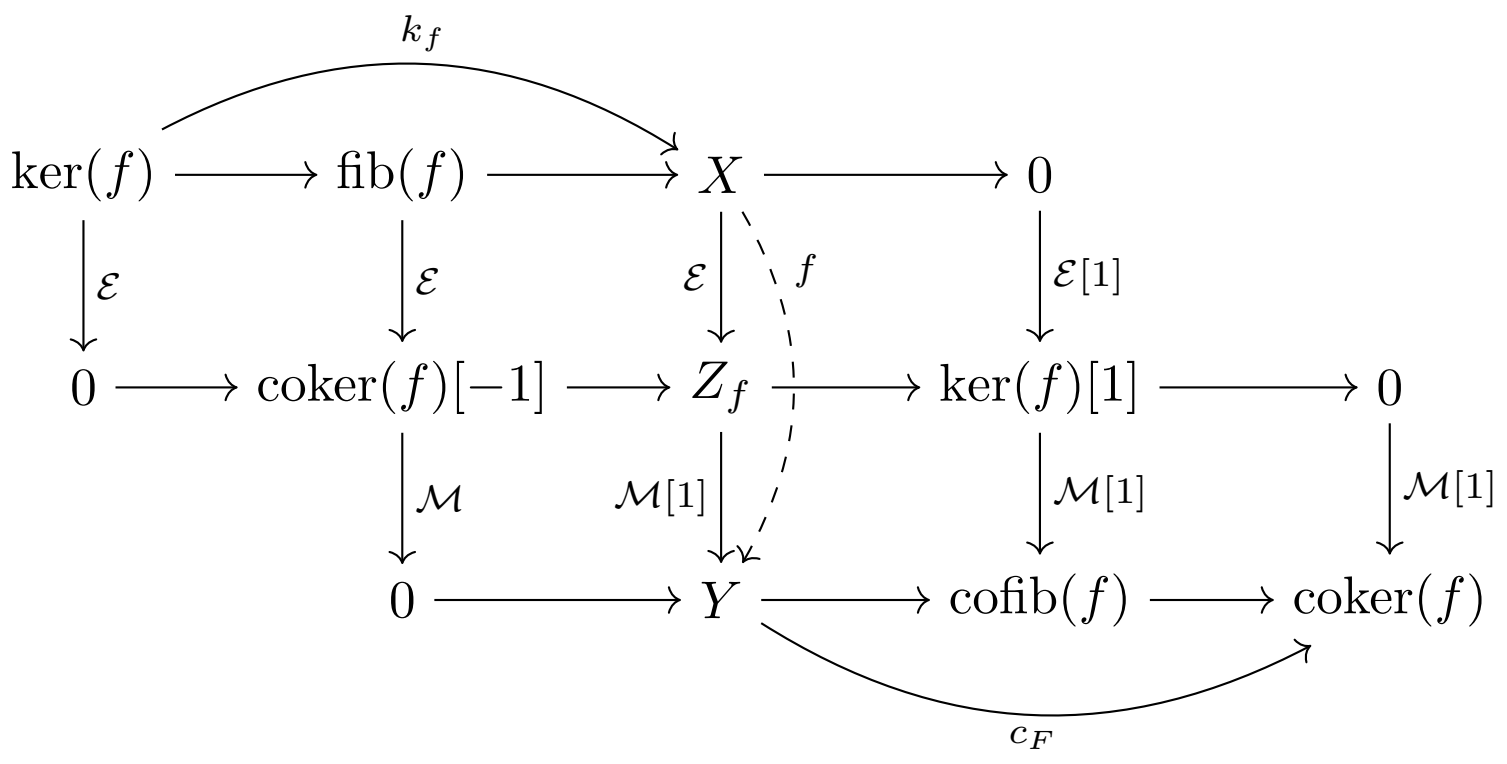

(4.44)

uniquely determining an object $Z_f \in \mathbf{C}^{\heartsuit}$.

*Proof.* Define $Z_f$ as the homotopy pullout

$$\begin{array}{ccc} \mathrm{fib}(f) & \longrightarrow & X \\ {\scriptstyle\mathcal{E}}\downarrow & \mathrm{pp} & \downarrow{\scriptstyle\mathcal{E}} \\ \mathrm{coker}(f)[-1] & \longrightarrow & Z_f \end{array} \tag{4.45}$$

Here the vertical arrow on the right is in $\mathcal{E}$ since the vertical arrow on the left is in $\mathcal{E}$ by definition of $\mathrm{coker}(f)$ (see Remark **4.3.13**) and $\mathcal{E}$ is preserved by pushouts. Next, paste on the left of this diagram the pullout given by Remark **4.3.13** and build the rest of the diagram by taking pullbacks or pushouts. Use again Remark **4.3.13** and the fact that $\mathcal{M}[1]$ is closed under pullbacks to see that $Z_f \to Y$ is in $\mathcal{M}[1]$. Finally, we have

$$0 \xrightarrow{\mathcal{E}} X \xrightarrow{\mathcal{E}} Z_f \xrightarrow{\mathcal{M}[1]} Y \xrightarrow{\mathcal{M}[1]} 0, \tag{4.46}$$

and so $Z_f$ is in $\mathbf{C}^{\heartsuit}$. □

PROPOSITION 4.3.19. There is an isomorphism $\mathrm{im}(f) \cong \mathrm{coim}(f)$.

*Proof.* By definition, $\mathrm{im}(f)$ and $\mathrm{coim}(f)$ are defined by the factorizations

$$0 \xrightarrow{\mathcal{E}} \mathrm{im}(f) \xrightarrow{\mathcal{M}} \mathrm{fib}(c_f) \tag{4.47}$$

and

$$\mathrm{cofib}(k_f) \xrightarrow{\mathcal{E}[1]} \mathrm{coim}(f) \xrightarrow{\mathcal{M}[1]} 0 \tag{4.48}$$

The diagram in Lemma **4.3.18** shows that we have $\mathrm{fib}(c_f) = Z_f = \mathrm{cofib}(k_f)$. Therefore, what we need to exhibit are the $(\mathcal{E}, \mathcal{M})$ factorizations of $0 \to Z_f$ and the $(\mathcal{E}[1], \mathcal{M}[1])$ factorization of $Z_f \to 0$. Since $Z_f$ is an object in $\mathbf{C}^\heartsuit$, these are

$$0 \xrightarrow{\mathcal{E}} Z_f \xrightarrow{\mathrm{id}_{Z_f}} Z_f \tag{4.49}$$

and

$$Z_f \xrightarrow{\mathrm{id}_{Z_f}} Z_f \xrightarrow{\mathcal{M}[1]} 0, \tag{4.50}$$

respectively, thus giving $\mathrm{im}(f) \cong Z_f \cong \mathrm{coim}(f)$. □

### 4.3.2 Abelian subcategories as hearts.

PROPOSITION 4.3.20. Let $\mathbf{A}$ be an abelian full subcategory of a stable $\infty$-category $\mathbf{C}$, such that any morphism $f \colon X \to Y$ in $\mathbf{C}$ has a unique $\mathbf{A}$-weaved $\mathbb{Z}$-Postnikov tower. Let $\mathbf{C}_{\mathbf{A},\geq 0}$ be the full subcategory of $\mathbf{C}$ on those objects $Y$ such that the $\mathbf{A}$-weaved $\mathbb{Z}$-Postnikov tower

$$0 = \lim(Y_j) \to \cdots \to Y_{j+1} \xrightarrow{f_j} Y_j \xrightarrow{f_{j-1}} Y_{j-1} \to \cdots \to \mathrm{colim}(Y_j) = Y \tag{4.51}$$

of the initial morphism $0 \to Y$ is such that $\mathrm{cofib}(f_j) = 0$ for any $j < 0$, and let $\mathbf{C}_{\mathbf{A},<0}$ be the full subcategory of $\mathbf{C}$ on those objects $Y$ such that $\mathrm{cofib}(f_j) = 0$ for any $j \geq 0$. Then $\mathfrak{t}_{\mathbf{A}} = (\mathbf{C}_{\mathbf{A},\geq 0}, \mathbf{C}_{\mathbf{A},<0})$ is a $t$-structure on $\mathbf{C}$, the stable $\infty$-category $\mathbf{C}$ is bounded with respect to $\mathfrak{t}_{\mathbf{A}}$, and the heart of $\mathfrak{t}_{\mathbf{A}}$ is (equivalent to) $\mathbf{A}$.

The proof is split in several Lemmas. We begin introducing the following

NOTATION 4.3.21. For $\mathbf{S}$ a subcategory of $\mathbf{C}$, we write $\langle \mathbf{S} \rangle$ for the smallest extension closed full subcategory of $\mathbf{C}$ containing $S$.

REMARK 4.3.22. Set $\langle \mathbf{S} \rangle_0 = \mathbf{0}$, define $\langle \mathbf{S} \rangle_1$ as the full subcategory of $\mathbf{C}$ generated by $\mathbf{S}$ and $\mathbf{0}$, and define inductively $\langle \mathbf{S} \rangle_n$ as the full subcategory of $\mathbf{C}$ on those objects $X$ which fall into a homotopy fiber sequence

$$\begin{array}{ccc} X_h & \longrightarrow & X \\ \downarrow & & \downarrow \\ 0 & \longrightarrow & X_k \end{array} \tag{4.52}$$

with $h, k \geq 1$, $X_h$ in $\langle \mathbf{S} \rangle_h$, $X_k$ in $\langle \mathbf{S} \rangle_k$ and $h + k = n$. One clearly has

$$\langle \mathbf{S} \rangle_0 \subseteq \langle \mathbf{S} \rangle_1 \subseteq \langle \mathbf{S} \rangle_2 \subseteq \cdots \subseteq \langle \mathbf{S} \rangle. \tag{4.53}$$

Moreover $\bigcup_n \langle \mathbf{S} \rangle_n$ is clearly extension closed, so that

$$\langle \mathbf{S} \rangle = \bigcup_n \langle \mathbf{S} \rangle_n. \tag{4.54}$$

LEMMA 4.3.23. Let $\mathbf{S}_1, \mathbf{S}_2$ be two subcategories of $\mathbf{C}$ with $\mathbf{S}_1 \perp \mathbf{S}_2$. Then $\mathbf{S}_1 \perp \langle \mathbf{S}_2 \rangle$ and $\langle \mathbf{S}_1 \rangle \perp \mathbf{S}_2$, and so $\langle \mathbf{S}_1 \rangle \perp \langle \mathbf{S}_2 \rangle$

*Proof.* By Remark **4.3.22**, to prove the first statement we are reduced to show that, if $Y \in \mathbf{S}_1$ and $X \in \langle \mathbf{S}_2 \rangle_n$ then $\mathbf{C}(Y, X)$ is contractible. We prove this by induction on $n$. For $n = 0, 1$ there is nothing to prove by the assumption $\mathbf{S}_1 \perp \mathbf{S}_2$. For $n \geq 2$, consider a fiber sequence $X_h \to X \to X_k$ with $1 \leq h, k$ and $h + k = n$ as in Remark **4.3.22**. Since $\mathbf{C}(Y, -)$ preserves homotopy fiber sequences, we get a homotopy fiber sequence of $\infty$-groupoids

$$\begin{array}{ccc} \mathbf{C}(Y, X_h) & \longrightarrow & \mathbf{C}(Y, X) \\ \downarrow & & \downarrow \\ * & \longrightarrow & \mathbf{C}(Y, X_k) \end{array} \quad . \tag{4.55}$$

By the inductive hypothesis both $\mathbf{C}(Y, X_h)$ and $\mathbf{C}(Y, X_k)$ are contractible, so $\mathbf{C}(Y, X)$ also is. The proof of the second statement is perfectly dual, due to the fact that in $\mathbf{C}$ every fiber sequence is also a cofiber sequence, and $\mathbf{C}(-, Y)$ transforms a cofiber sequence into a fiber sequence. □

Lemma 4.3.24. Let $\mathbf{A}$ be an abelian full subcategory of $\mathbf{C}$. Then $\langle \{\mathbf{A}[s]\}_{s \geq 0} \rangle \perp \langle \{\mathbf{A}[s]\}_{s<0} \rangle$. In particular, in the hypothesis of Proposition **4.3.20** we have $\mathbf{C}_{\mathbf{A}, \geq 0} \perp \mathbf{C}_{\mathbf{A}, <0}$

*Proof.* By Lemma **4.3.23**, we only need to show that $\mathbf{A}[s_1] \perp \mathbf{A}[s_2]$ whenever $s_1 \geq 0 > s_2$. Let $X \in \mathbf{A}[s_1]$ and $Y \in \mathbf{A}[s_2]$. Then $X = Z_1[s_1]$ and $Y = Z_2[s_2]$ for suitable $Z_1, Z_2 \in \mathbf{A}$ and so

$$\begin{aligned} \mathbf{C}(X, Y) = \mathbf{C}(Z_1[s_1], Z_2[s_2]) &\cong \mathbf{C}(Z_1, Z_2[s_2 - s_1]) \\ &\cong \Omega^{s_1 - s_2} \mathbf{C}(Z_1, Z_2) = \Omega^{s_1 - s_2} \mathbf{A}(Z_1, Z_2), \end{aligned}$$

where in the last equality we used the fact that $\mathbf{A}$ is full. Since $s_1 - s_2 > 0$, the space $\Omega^{s_1 - s_2} \mathbf{A}(Z_1, Z_2)$ is contractible by definition of abelian $\infty$-category. Finally, in the hypothesis of Proposition **4.3.20** one clearly has $\mathbf{C}_{\mathbf{A}, <0} \subseteq \langle \{\mathbf{A}[s]\}_{s<0} \rangle$ and $\mathbf{C}_{\mathbf{A}, \geq 0} \subseteq \langle \{\mathbf{A}[s]\}_{s \geq 0} \rangle$. □

Lemma 4.3.25. In the hypothesis of Proposition **4.3.20** every object $Y$ of $\mathbf{C}$ sits into a homotopy fiber sequence $Y_{\geq 0} \to Y \to Y_{<0}$ with $Y_{\geq 0} \in \mathbf{C}_{\mathbf{A}, \geq 0}$ and $Y_{<0} \in \mathbf{C}_{\mathbf{A}, <0}$.

*Proof.* Let

$$0 = \lim(Y_j) \to \cdots \to Y_1 \xrightarrow{f_0} Y_0 \xrightarrow{f_{-1}} Y_{-1} \to \cdots \to \operatorname{colim}(Y_j) = Y \tag{4.56}$$

be teh $\mathbf{A}$-weaved Postnikov tower of $0 \to Y$ and consider the pullout diagram

$$\begin{array}{ccc} Y_0 & \longrightarrow & 0 \\ {\scriptstyle f_{<0}} \downarrow & & \downarrow \\ Y & \longrightarrow & \operatorname{cofib}(f_{<0}) \end{array} \tag{4.57}$$

together with the $\mathbf{A}$-weaved $\mathbb{Z}$-Postnikov towers

$$0 = \lim(Y_j) \to \cdots \to Y_1 \xrightarrow{f_0} Y_0 \tag{4.58}$$

and

$$Y_0 \xrightarrow{f_{-1}} Y_{-1} \to \cdots \to \operatorname{colim}(Y_j) = Y. \tag{4.59}$$

The first Postnikov tower shows that $Y_0 \in \mathbf{C}_{\mathbf{A},\geq 0}$ while the second Postnikov tower shows that $\operatorname{cofib}(f_{<0}) \in \mathbf{C}_{\mathbf{A},<0}$. □

Lemma 4.3.26. In the hypothesis of Proposition **4.3.20**, for any $\lambda \in \mathbb{R}$ let $\mathbf{C}_{\mathbf{A},\geq\lambda}$ be the full subcategory of $\mathbf{C}$ on those objects $Y$ such that the $\mathbf{A}$-weaved $\mathbb{Z}$-Postnikov tower

$$0 = \lim(Y_j) \to \cdots \to Y_{j+1} \xrightarrow{f_j} Y_j \xrightarrow{f_{j-1}} Y_{j-1} \to \cdots \to \operatorname{colim}(Y_j) = Y \tag{4.60}$$

of the initial morphism $0 \to Y$ is such that $\operatorname{cofib}(f_j) = 0$ for any $j < \lambda$, and let $\mathbf{C}_{\mathbf{A},<\lambda}$ be the full subcategory of $\mathbf{C}$ on those objects $Y$ such that $\operatorname{cofib}(f_j) = 0$ for any $j \geq \lambda$. Then, for any $n \in \mathbb{Z}$, one has $\mathbf{C}_{\mathbf{A},<\lambda}[n] = \mathbf{C}_{\mathbf{A},<\lambda+n}$ and $\mathbf{C}_{\mathbf{A},\geq\lambda}[n] = \mathbf{C}_{\mathbf{A},\geq\lambda+n}$. In particular, $\mathbf{C}_{\mathbf{A},<0}[-1] \subseteq \mathbf{C}_{\mathbf{A},<0}$ and $\mathbf{C}_{\mathbf{A},\geq 0}[1] \subseteq \mathbf{C}_{\mathbf{A},\geq 0}$.

*Proof.* Since the shift functor commutes with the formation of $\mathbf{A}$-weaved $\mathbb{Z}$-Postnikov towers, an object $Y$ lies in $\mathbf{C}_{\mathbf{A},<\lambda}[n]$ if and only if $\operatorname{cofib}(f_{j+n}[-n]) = 0$ for any $j \geq \lambda$, i.e., if and only if $\operatorname{cofib}(f_j) = 0$ for any $j \geq \lambda + n$. The proof for $\mathbf{C}_{\mathbf{A},\geq\lambda}[n]$ is identical. □

*Proof of Proposition 4.3.20.* Lemmas **4.3.24**, **4.3.25** and **4.3.26** together show that $\mathfrak{t}_{\mathbf{A}} = (\mathbf{C}_{\mathbf{A},\geq 0}, \mathbf{C}_{\mathbf{A},<0})$ is a bounded $t$-structure on $\mathbf{C}$. To see that the heart of $\mathfrak{t}_{\mathbf{A}}$ is $\mathbf{A}$ notice that an object $Y$ lies in $\mathbf{C}_{\mathbf{A},[0,1)}$ if and only if the $\mathbf{A}$-weaved $\mathbb{Z}$-Postnikov tower

$$0 = \lim(Y_j) \to \cdots \to Y_{j+1} \xrightarrow{f_j} Y_j \xrightarrow{f_{j-1}} Y_{j-1} \to \cdots \to \operatorname{colim}(Y_j) = Y \tag{4.61}$$

of its initial morphism has $\operatorname{cofib}(f_j) = 0$ for every $j \neq 0$, and so it is of the form

$$\cdots 0 \to 0 \to \cdots \to 0 \xrightarrow{f_0} Y \xrightarrow{\mathrm{id}_Y} Y \xrightarrow{\mathrm{id}_Y} \cdots \xrightarrow{\mathrm{id}_Y} Y \xrightarrow{\mathrm{id}_Y} \cdots, \tag{4.62}$$

with $Y = \operatorname{cofib}(f_0) \in \mathbf{A}$. □

Remark 4.3.27. The same reasoning used in the proof of Proposition **4.3.20**, shows that $(\mathbf{C}_{\mathbf{A},\geq\lambda}, \mathbf{C}_{\mathbf{A},<\lambda})$ is a bounded $t$-structure on $\mathbf{C}$ for every $\lambda \in \mathbb{R}$, and that the assignment $\lambda \mapsto (\mathbf{C}_{\mathbf{A},\geq\lambda}, \mathbf{C}_{\mathbf{A},<\lambda})$ is a $\mathbb{Z}$-equivariant morphisms of posets $\mathbb{R} \to \mathrm{TS}(\mathbf{C})$, so it is a slicing of $\mathbf{C}$. The heart of $(\mathbf{C}_{\mathbf{A},\geq\lambda}, \mathbf{C}_{\mathbf{A},<\lambda})$ is $\mathbf{A}[\lceil\lambda\rceil]$, where $\lceil\lambda\rceil = \min\{n \in \mathbb{Z} \mid n \geq \lambda\}$.

## 4.4 Semiorthogonal decompositions.

La vie c'est ce qui se décompose à tout moment; c'est une perte monotone de lumière, une dissolution insipide dans la nuit, sans sceptres, sans auréoles, sans nimbes.

E. Cioran, *Précis de décomposition.*

At the opposite end of the transitive case studied in the previous section, there is the *finite case*, where $J$ is a finite totally ordered set. As we are going to show, this is another well investigated case in the literature: $J$-familes of $t$-structures with a finite $J$ capture (and slightly generalize) the notion of *semiorthogonal decompositions* for the stable $\infty$-category $\mathbf{C}$ (see [BO95, Kuz11] for the notion of semiorthogonal decomposition in the classical triangulated context).

To fix notations for this section, let $J = \Delta[k-1]$ be the totally ordered set on $k$ elements seen as a poset, i.e., $J = \{i_1, i_2, \dots, i_k\}$ with $i_1 \leq i_2 \leq \dots \leq i_k$, and let $\mathfrak{t}\colon \Delta[k-1] \to \text{TS}(\mathbf{C})$ be a $\mathbb{Z}$-equivariant $\Delta[k-1]$-family of $t$-structures on $\mathbf{C}$. We also set, for any $j = 1, \dots, k+1$,

$$\mathbf{A}_j = \mathbf{C}_{[i_{j-i}, i_j)} \tag{4.63}$$

where, as usual, $i_0 = -\infty$ and $i_{k+1} = +\infty$. We have that any morphism $f\colon X \to Y$ in $\mathbf{C}$ has a unique factorization

$$X \xrightarrow{f_{k+1}} Z_{i_k} \xrightarrow{f_k} Z_{i_{k-1}} \to \dots \to Z_{i_2} \xrightarrow{f_2} Z_{i_1} \xrightarrow{f_1} Y, \tag{4.64}$$

with $\text{cofib}(f_j) \in \mathbf{A}_j$, and $\mathbf{A}_j \subseteq \mathbf{A}_h^{\perp}$, for any $1 \leq j < h \leq k+1$.

What we are left to investigate are therefore the special features of the $t$-structures $\mathfrak{t}_{i_j} = (\mathbf{C}_{\geq i_j}, \mathbf{C}_{<i_j})$ coming from the finiteness assumption on $J$. As we noticed in Remark **4.1.12**, a $\mathbb{Z}$-action on a finite poset $J$ is necessarily trivial. By $\mathbb{Z}$-equivariancy of the map $\Delta[k-1] \to \text{TS}(\mathbf{C})$ we have therefore that all the $t$-structures $\mathfrak{t}_{i_j}$ are $\mathbb{Z}$-fixed points for the natural $\mathbb{Z}$-action on $\text{TS}(\mathbf{C})$.

Now, a rather pleasant fact is that fixed points of the $\mathbb{Z}$-action on $\text{TS}(\mathbf{C})$ are precisely those $t$-structures $\mathfrak{t} = (\mathbf{C}_{\geq 0}, \mathbf{C}_{<0})$ for which $\mathbf{C}_{\geq 0}$ is a stable sub-$\infty$-category of $\mathbf{C}$. We will make use of the following

Lemma 4.4.1. Let $\mathbf{B}$ be a full sub-$\infty$-category of the stable $\infty$-category $\mathbf{C}$; then, $\mathbf{B}$ is a stable sub-$\infty$-category of $\mathbf{C}$ if and only if $\mathbf{B}$ is closed under shifts in both directions and under pushouts in $\mathbf{C}$.

*Proof.* The "only if" part is trivial, so let us prove the "if" part.

First of all let us see that under these assumptions $\mathbf{B}$ is closed under fibers. This is immediate: if $f\colon X \to Y$ is an arrow in $\mathbf{B}$ (i.e. an arrow of $\mathbf{C}$ between objects of $\mathbf{B}$, by fullness), then $f[-1]$ is again in $\mathbf{B}$ since $\mathbf{B}$ is closed with respect to the left shift. Since $\mathbf{B}$ is closed under pushouts in $\mathbf{C}$, also $\mathrm{fib}(f) = \mathrm{cofib}(f[-1])$ is in $\mathbf{B}$. It remains to show how this implies that $\mathbf{B}$ is actually stable, i.e. it is closed under all finite limits and satisfies the pullout axiom. Unwinding the assumptions on $\mathbf{B}$, this boils down to showing that in the square

$$\begin{array}{ccc} B & \longrightarrow & X \\ \big\downarrow & \text{pb} & \big\downarrow f \\ Y & \xrightarrow[g]{} & Z \end{array} \tag{4.65}$$

the pullback $B$ of $f, g \in \hom(\mathbf{B})$ done in $\mathbf{C}$ is actually an object of $\mathbf{B}$; indeed, once this is shown, the square above will satisfy the pullout axiom in $\mathbf{C}$, so *a fortiori* it will have the universal property of a pushout in $\mathbf{B}$. To this aim, let us consider the enlarged diagram of pullout squares in $\mathbf{C}$

$$\begin{array}{ccccc} Z[-1] & \longrightarrow & \mathrm{fib}(g) & \longrightarrow & 0 \\ \big\downarrow & \star & \big\downarrow & & \big\downarrow \\ \mathrm{fib}(f) & \longrightarrow & B & \longrightarrow & X \\ \big\downarrow & & \big\downarrow & & \big\downarrow f \\ 0 & \longrightarrow & Y & \xrightarrow[g]{} & Z. \end{array} \tag{4.66}$$

The objects $Z[-1], \mathrm{fib}(f)$ and $\mathrm{fib}(g)$ lie in $\mathbf{B}$ by the first part of the proof, so the square $(\star)$ is in particular a pushout of morphism in $\mathbf{B}$; by assumption, this entails that $B \in \mathbf{B}$. □

REMARK 4.4.2. Obviously, a completely dual statement can be proved in a completely dual fashion: a full sub-$\infty$-category $\mathbf{B}$ of a stable $\infty$-category $\mathbf{C}$ is a stable sub-$\infty$-category if and only if it is closed under shifts in both directions and under pullbacks in $\mathbf{C}$.

PROPOSITION 4.4.3. Let $\mathfrak{t} = (\mathbf{C}_{\geq 0}, \mathbf{C}_{<0})$ be a $t$-structure on a stable $\infty$-category $\mathbf{C}$; then the following conditions are equivalent:

(1) $\mathfrak{t}$ is a fixed point for the $\mathbb{Z}$-action on $\mathrm{TS}(\mathbf{C})$, i.e., $\mathfrak{t}[1] = \mathfrak{t}$ (or equivalently in view of remark **A.3.8**, $\mathbf{C}_{\geq 1} = \mathbf{C}_{\geq 0}$);
(2) $\mathbf{C}_{\geq 0}$ is a stable sub-$\infty$-category of $\mathbf{C}$.

*Proof.* '(2) implies (1)' is obvious. Namely, if $\mathbf{C}_{\geq 0}$ is a stable sub-$\infty$-category of $\mathbf{C}$, then it is closed under shifts in both directions. Therefore $\mathbf{C}_{\geq 1} = \mathbf{C}_{\geq 0}[1] \subseteq \mathbf{C}_{\geq 0}$. Since, by definition of $t$-structure, $\mathbf{C}_{\geq 1} \subseteq \mathbf{C}_{\geq 0}$, we have $\mathbf{C}_{\geq 1} = \mathbf{C}_{\geq 0}$. To prove that '(1) implies (2)', assume $\mathbf{C}_{\geq 1} = \mathbf{C}_{\geq 0}$. This means that not only $\mathbf{C}_{\geq 0}[1] \subseteq \mathbf{C}_{\geq 0}$ as for any $t$-structure, but also $\mathbf{C}_{\geq 0} \subseteq \mathbf{C}_{\geq 0}[1]$, which implies that $\mathbf{C}_{\geq 0}[-1] \subseteq \mathbf{C}_{\geq 0}$. Therefore $\mathbf{C}_{\geq 0}$ is closed under shifts in both directions. By Lemma **4.4.1**, we then have only to show that $\mathbf{C}_{\geq 0}$ is closed under pushouts in $\mathbf{C}$ to conclude that $\mathbf{C}_{\geq 0}$ is a stable $\infty$-subcategory of $\mathbf{C}$. Consider a pushout diagram

$$\begin{array}{ccc} A & \longrightarrow & B \\ {\scriptstyle h}\downarrow & \text{po} & \downarrow{\scriptstyle k} \\ C & \longrightarrow & P \end{array} \tag{4.67}$$

in $\mathbf{C}$ with $A$, $B$ and $C$ in $\mathbf{C}_{\geq 0}$, and let $\mathbb{F} = (\mathcal{E}, \mathcal{M})$ be the normal torsion theory associated to $\mathfrak{t}$. Since $A$ and $C$ are in $\mathbf{C}_{\geq 0} = 0/\mathcal{E}$ we have that both $0 \to A$ and $0 \to C$ are in $\mathcal{E}$. But $\mathcal{E}$ has the 3-for-2 property, so also $A \to C$ is $\mathcal{E}$. Since $\mathcal{E}$ is closed for pushouts, this implies that also $B \to P$ is in $\mathcal{E}$. But $0 \to B$ in in $\mathcal{E}$ since $B$ is in $\mathbf{C}_{\geq 0}$, and therefore also $0 \to P$ is in $\mathcal{E}$, i.e., $P$ is in $\mathbf{C}_{\geq 0}$. □

REMARK 4.4.4. The statement of Prop. **4.4.3** can easily be dualized: $\mathbb{Z}$-fixed points in $\text{TS}(\mathbf{C})$ as those $t$-structures $(\mathbf{C}_{\geq 0}, \mathbf{C}_{<0})$ for which $\mathbf{C}_{<0}$ is a stable sub-$\infty$-category of $\mathbf{C}$, as well as those such that $\mathbf{C}_{<0} = \mathbf{C}_{<1}$.

Proposition **4.4.3** and remark **4.4.4** characterize $\mathbb{Z}$-fixed points on $\text{TS}(\mathbf{C})$ as the $t$-structures with stable classes $\mathbf{C}_{\geq 0}$ and $\mathbf{C}_{<0}$. By the correspondence between $t$-structures and normal factorization systems, one should expect that these should be equally characterized as the normal factorization systems $\mathbb{F} = (\mathcal{E}, \mathcal{M})$ for which the classes $\mathcal{E}$ and $\mathcal{M}$ are "stable on both sides", i.e., are closed both for pullbacks and for pushouts.

THEOREM 4.4.5. Let $\mathfrak{t}$ be a $t$-structure on a stable $\infty$-category $\mathbf{C}$ and let $\mathbb{F} = (\mathcal{E}, \mathcal{M})$ be the corresponding normal factorization system; then the following conditions are equivalent:

(1) $\mathfrak{t}[1] = \mathfrak{t}$;
(2) $\mathbf{C}_{\geq 0}$ is a stable $\infty$-category;
(3) $\mathbf{C}_{<0}$ is a stable $\infty$-category;
(4) $\mathcal{E}$ is closed under pullback;
(5) $\mathcal{M}$ is closed under pushout.

*Proof.* In view of the previous results, the only implication we need to prove is that '(1) is equivalent to (4)'. Assume $\mathcal{E}$ is closed under pullbacks. Then for any $X$ in $\mathbf{C}_{\geq 0}$ we have that $0 \to X$ is in $\mathcal{E}$, and so $X[-1] \to 0$ is in $\mathcal{E}$. By the Sator lemma this implies that $0 \to X[-1]$ is in $\mathcal{E}$, i.e., that $X[-1]$ is in $\mathbf{C}_{\geq 0}$. This shows that $\mathbf{C}_{\geq 0}[-1] \subseteq \mathbf{C}_{\geq 0}$ and therefore that $\mathfrak{t}[1] = \mathfrak{t}$.

Conversely, assume $\mathfrak{t}[1] = \mathfrak{t}$, and consider a morphism $f \colon X \to Y$ in $\mathcal{E}$. For any morphism $B \to Y$ in $\mathbf{C}$ consider the diagram

$$\begin{array}{ccccccc} \text{fib}(f) & \longrightarrow & A & \longrightarrow & X & \longrightarrow & 0 \\ \downarrow & & \downarrow & & \downarrow{\scriptstyle f} & & \downarrow \\ 0 & \longrightarrow & B & \longrightarrow & Y & \longrightarrow & \text{cofib}(f) \end{array} \tag{4.68}$$

where all the squares are pullouts in $\mathbf{C}$. Since $f$ is in $\mathcal{E}$ and $\mathcal{E}$ is closed for pushouts, also $0 \to \mathrm{cofib}(f)$ is in $\mathcal{E}$. This means that $\mathrm{cofib}(f)$ is in $\mathbf{C}_{\geq 0}$ and so, since we are assuming that $\mathbf{C}_{\geq 0} = \mathbf{C}_{\geq 0}[-1]$, also $\mathrm{fib}(f) = \mathrm{cofib}(f)[-1]$ is in $\mathbf{C}_{\geq 0}$, i.e., $0 \to \mathrm{fib}(f)$ is in $\mathcal{E}$. By the Sator lemma, $\mathrm{fib}(f) \to 0$ is in $\mathcal{E}$, which is closed for pushouts, and so $A \to B$ is in $\mathcal{E}$. □

REMARK 4.4.6. In the literature, a factorization system $(\mathcal{E}, \mathcal{M})$ for which the class $\mathcal{E}$ is closed for pullbacks is sometimes called an *exact reflective* factorization, see, e.g., [CHK85]. This is equivalent to saying that the associated reflection functor is left exact (this is called a *localization* in the jargon of [CHK85]). Dually, one characterizes *colocalizations* of a category $\mathbf{C}$ with an initial object as *coexact coreflective* factorizations where the right class $\mathcal{M}$ of $\mathbb{F}$ is closed under pushouts. Therefore, in the stable $\infty$-case, we see that a $\mathbb{Z}$-fixed point in $\mathrm{TS}(\mathbf{C})$ is a $t$-structure $(\mathbf{C}_{\geq 0}, \mathbf{C}_{<0})$ such that the truncation functors $\tau_{\geq 0} \colon \mathbf{C} \to \mathbf{C}_{\geq 0}$ and $\tau_{<0} \colon \mathbf{C} \to \mathbf{C}_{<0}$ respectively form a colocalizations and a localization of $\mathbf{C}$. In the terminology of [BR07] we therefore find that in the stable $\infty$-case $\mathbb{Z}$-fixed point in $\mathrm{TS}(\mathbf{C})$ correspond to *hereditary torsion pairs* on $\mathbf{C}$. Since we have seen that for a $\mathbb{Z}$-fixed point in $\mathrm{TS}(\mathbf{C})$ both $\mathbf{C}_{\geq 0}$ and $\mathbf{C}_{<0}$ are stable $\infty$-categories, this result could be deduced also from [Lur17, Prop. **1.1.4.1**]: a left (resp., right) exact functor between stable $\infty$-categories is also right (resp., left) exact.

We can now precisely relate semiorthogonal decompositions in a stable $\infty$-category $\mathbf{C}$ to $\Delta[k-1]$-families of $t$-structures on $\mathbf{C}$. The only thing we still need is the following definition, which is an immediate adaptation to the stable setting of the classical definition given for triangulated categories (see, e.g., [BO95, Kuz11] ).

DEFINITION 4.4.7. Let $\mathbf{C}$ be a stable $\infty$-category. A *semiorthogonal decomposition* with $k$ classes on $\mathbf{C}$ is the datum of $k+1$ stable $\infty$-subcategories $\mathbf{A}_1$, $\mathbf{A}_2$,..., $\mathbf{A}_{k+1}$ of $\mathbf{C}$ such that

(1) one has $\mathbf{A}_i \subseteq \mathbf{A}_j^{\perp}$ for $i < j$ (semiorthogonality);
(2) for any object $Y$ in $\mathbf{C}$ there exists a unique $\{\mathbf{A}_i\}$-weaved tower, i.e., a factorization of the initial morphism $0 \to Y$ as

$$0 = Y_0 \to \cdots \to Y_{j+1} \xrightarrow{f_j} Y_j \xrightarrow{f_{j-1}} Y_{j-1} \to \cdots \to Y_{k+1} = Y \tag{4.69}$$

with $\mathrm{cofib}(f_j) \in \mathbf{A}_j$ for any $j = 1, \ldots, k+1$.

REMARK 4.4.8. Since $\{\mathbf{A}_i\}$-weaved Postnikow towers are preserved by pullouts, one can equivalently require that any morphism $f \colon X \to Y$ in $\mathbf{C}$ has a unique factorization of the form

$$X = Z_0 \to \cdots \to Z_{j+1} \xrightarrow{f_j} Z_j \xrightarrow{f_{j-1}} Z_{j-1} \to \cdots \to Z_{k+1} = Y \tag{4.70}$$

with $\mathrm{cofib}(f_j) \in \mathbf{A}_j$ for any $j = 1, \ldots, k+1$.

Theorem 4.4.9. Let $\mathbf{C}$ be a stable $\infty$-category. Then the datum of a semiorthogonal decompositions with $k$ classes on $\mathbf{C}$ is equivalent to the datum of a $\mathbb{Z}$-equivariant $\Delta[k-1]$-family of $t$-structures on $\mathbf{C}$

*Proof.* Let us start with a $\mathbb{Z}$-equivariant $\Delta[k-1]$-family of $t$-structures $\mathfrak{t}$, and write $i_1 < i_2 < \cdots < i_k$ for the elements of $\Delta[k-1]$ and $\mathfrak{t}_{i_j} = (\mathbf{C}_{\geq i_j}, \mathbf{C}_{<i_j})$ for the corresponding $t$-structures on $\mathbf{C}$. Then, setting $\mathbf{A}_j = \mathbf{C}_{[i_{j-1}, i_j)}$ we have semiorthogonality between the $\mathbf{A}_j$'s and the existence of $\{\mathbf{A}_j\}$-weaved Postnikow towers by the general argument recalled at the beginning of this section. So we are only left to prove that the subcategories $\mathbf{A}_j$ are stable. This is immediate: by Theorem **4.4.5** both the sub-$\infty$-categories $\mathbf{C}_{\geq i_{j-1}}$ and $\mathbf{C}_{<i_j}$ are stable, and so also their intersection is stable (see, [Lur17]). Vice versa, if we start with a semiorthogonal decomposition, then repeating verbatim the argument in the proof of Proposition **4.3.20** one defines a $\mathbb{Z}$-equivariant $\Delta[k-1]$-family of $t$-structures on $\mathbf{C}$. □

Remark 4.4.10. By Remark **4.4.6**, we recover in the stable $\infty$-setting the well known fact (see [BR07, **IV.4**]) that semiorthogonal decompositions with a single class correspond to *hereditary torsion pairs* on the category.

# Chapter 5

# Recollements

The present chapter develops the theory of *recollements* in a stable $\infty$-categorical setting. In the axiomatization of Beĭlinson, Bernstein and Deligne, recollement situations provide a generalization of Grothendieck's "six functors" between derived categories.

If a recollement is depicted as a diagram $\mathbf{D}^0 \leftrightarrows \mathbf{D} \leftrightarrows \mathbf{D}^1$, given $t$-structures $\mathfrak{t}_0, \mathfrak{t}_1$ on $\mathbf{D}^0, \mathbf{D}^1$ it is possible to construct a "recollée" $t$-structure $\mathfrak{t}_0 \Yup \mathfrak{t}_1$ (see Def. **5.2.1**) on $\mathbf{D}$, exploiting the adjointness relations between these six functors.

Such a classical result, well-known in the setting of triangulated categories, acquires a new taste when $t$-structures are described as normal torsion theories: outlining the construction of the factorization system related to $\mathfrak{t}_0 \Yup \mathfrak{t}_1$ by the "Rosetta stone" reveals a number of interesting formal properties of the construction, and clarifies its origin.

In the geometric case of what [BBD82] calls a *stratified space*, various recollements arise, and they "interact well" with the combinatorics of the intersections of strata to give a well-defined, associative $\Yup$ operation. From this we deduce a generalized associative property for $n$-fold gluing $\mathfrak{t}_0 \Yup \cdots \Yup \mathfrak{t}_n$, valid in any stable $\infty$-category, provided that a sufficient number of recollement data organize into a diagram (see Def. **5.4.7**) which ensures the possibility of parenthesizing the string $\mathfrak{t}_0 \Yup \cdots \Yup \mathfrak{t}_n$ in different ways and to compare these ways.

## 5.1 Introduction.

Recollements in triangulated categories were introduced by A. Beĭlinson, J. Bernstein and P. Deligne in [BBD82], searching an axiomatization of the Grothendieck's "six functors" formalism for derived categories of sheaves on (the strata of a) stratified topological space. [BBD82] will be our main source of inspiration, and reference for classical results and computations; among other recent but standard references, we mention [KS, Ban07]. Later,

"recollement data" were noticed to appear quite naturally in the context of intersection homology [Pfl01, GM80, GM83] and Representation Theory [PS88, KW01]. In more recent years Beligiannis and Reiten [BR07], adapting to the triangulated setting an old idea of Jans [Jan65], linked recollement data to so-called TTF-triples (i.e. triples $(\mathcal{X}, \mathcal{Y}, \mathcal{Z})$ such that both $(\mathcal{X}, \mathcal{Y})$ and $(\mathcal{Y}, \mathcal{Z})$ are $t$-structures): recollement data, in the form of TTF-triples, appear quite naturally studying derived categories of representations of algebras, see [BR07, Ch. **4**].

We now attempt to translate the basic theory of recollements in the stable setting; adopting this viewpoint clarifies the classical theory and offers a number of interesting results.

Focusing on normal torsion theories as the higher-categorical entities inducing $t$-structures in the triangulated world is, categorically speaking, extremely natural. This "torsio-centric" perspective appears to be very well suited to the description of recollements.

In the presence of a stratification $U_0 \subset U_1 \cdots \subset U_n \subset X$ of a space, with "pure strata" $E_i = U_i \smallsetminus U_{i-1}$ an extended version of the Rorschach lemma entails that an object $\mathcal{F}$ of the derived category of $X$ lies in a class of the glued $t$-structure if and only if $l_i(\mathcal{F})$ lies in the homonym class in $\mathbf{D}(E_i)$, where $l_i$ is a suitable choice of a functor from $\mathbf{D}(X)$ to the category of $i^{\text{th}}$ pure stratum (this appears in **5.4.1**). This associativity result, together with the "compatibility condition" necessary to ensure that two different parenthesizations of $\mathfrak{t}_0 \uplus \cdots \uplus \mathfrak{t}_n$ coincide,

$$(\mathfrak{t}_0 \uplus \cdots \uplus \mathfrak{t}_n)_{\mathfrak{P}} = (\mathfrak{t}_0 \uplus \cdots \uplus \mathfrak{t}_n)_{\mathfrak{Q}} \tag{5.1}$$

is not spelled out explicitly in [BBD82].

## 5.2 Classical Recollements.

> Sitzt ihr nur immer! leimt zusammen,
> Braut ein Ragout von andrer Schmaus,
> und blas't die kümmerlichen Flammen
> aus eurem Aschenhäufchen 'raus!
>
> Faust, I 538-541.

The aim of this subsection is to present the basic features of "classical" recollements in the setting of stable $\infty$-categories ignoring, for the moment, the translation in terms of normal torsion theories which will follow.

Definition 5.2.1. A (*donnée de*) *recollement* consists of the following arrangement of stable $\infty$-categories and functors between them:

$$\mathbf{D}^0 \underset{i_L}{\overset{i_R}{\leftleftarrows\!\!\xrightarrow{\;i\;}}} \mathbf{D} \underset{q_L}{\overset{q_R}{\leftleftarrows\!\!\xrightarrow{\;q\;}}} \mathbf{D}^1 \tag{5.2}$$

satisfying the following axioms:

(1) There are adjunctions $i_L \dashv i \dashv i_R$ and $q_L \dashv q \dashv q_R$;
(2) The counit $\epsilon_{(i_L \dashv i)} \colon i_L i \to 1$ and the unit $\eta_{(i \dashv i_R)} \colon 1 \to i_r i$ are natural isomorphisms; also, the unit $1 \to q q_R$ and counit $q q_L \to 1$ are natural isomorphisms;[(1)]
(3) The (essential) image of $i$ equals the *essential kernel* of $q$, namely the full subcategory of $\mathbf{D}$ such that $qX \cong 0$ in $\mathbf{D}^1$;
(4) The natural homotopy commutative diagrams

$$\begin{array}{ccc} q_L q & \xrightarrow{\epsilon_{(q_L \dashv q)}} & \mathrm{id}_{\mathbf{D}} \\ \downarrow & & \downarrow {\scriptstyle \eta_{(i_L \dashv i)}} \\ 0 & \longrightarrow & i i_L \end{array} \qquad \begin{array}{ccc} i i_R & \xrightarrow{\epsilon_{(i \dashv i_R)}} & \mathrm{id}_{\mathbf{D}} \\ \downarrow & & \downarrow {\scriptstyle \eta_{(q \dashv q_R)}} \\ 0 & \longrightarrow & q_R q \end{array} \tag{5.3}$$

induced by axioms (1), (2) and (3) are pullouts[(2)].

Remark 5.2.2. As an immediate consequence of the axioms, a recollement gives rise to various reflections and coreflections of $\mathbf{D}$: since by axiom (2) the functors $i, q_L, q_R$ are all fully faithful, $q_R q, i i_L$ are reflections and $q_L q, i i_R$ are coreflections. Moreover, axioms (3) and (4) entail that the compositions $i_R q_R, qi, i_L q_L$ are all "exactly" zero, i.e. not only the kernel of $q$ is the essential image of $i$, but also the kernel of $i_{L/R}$ is the essential image of $q_{L/R}$.

Remark 5.2.3. Axioms (2) and (4) together imply that there exists a canonical natural transformation $i_R \to i_L$, obtained as $i_R(\eta_{(i_L \dashv i)})$ (or equivalently, as $i_L(\epsilon_{(i \dashv i_R)})$: it's easy to see that these two arrows coincide). Axiom (4) entails that there is a fiber sequence of natural transformations

$$\begin{array}{ccccc} i_R q_L q & \longrightarrow & i_R & \longrightarrow & 0 \\ \downarrow & & \downarrow & & \downarrow \\ 0 & \longrightarrow & i_L & \longrightarrow & i_L q_R q \end{array}$$

Notation 5.2.4. We will generally use a compact form like

$$(i, q) \colon \mathbf{D}^0 \leftrightarrows \mathbf{D} \leftrightarrows \mathbf{D}^1 \tag{5.4}$$

to denote a recollement (**5.2**), especially in inline formulas. Variations on this are possible, either to avoid ambiguities or to avoid becoming stodgy.

We will for example say that "$(i, q)$ is a recollement on $\mathbf{D}$" or that "$\mathbf{D}$ is the *décollement* of $\mathbf{D}^0, \mathbf{D}^1$" to denote that there exists a diagram like (**5.2**)

(1) With a little abuse of notation we will write $i_L i = \mathrm{id}_{\mathbf{D}^0} = i_R i$, and similarly for $q q_L = \mathrm{id}_{\mathbf{D}} = q q_R$.

(2) Here and everywhere else the category of functors to a stable $\infty$-category becomes a stable $\infty$-category in the obvious way (see [Lur17, Prop. **1.1.3.1**]).

having **D** as a central object. In other situations we adopt an extremely compact notation, referring to a (donné de) recollement with the symbol რ of (the letter *rae* of the Georgian alphabet, in the მხედრული script, see [Hew95]).

**A geometric example.** The most natural example of a recollement comes from the theory of *stratified spaces* [Wei94, Ban07]:

EXAMPLE 5.2.5. Let $X$ be a topological space, $F \subseteq X$ a closed subspace, and $U = X \smallsetminus F$ its open complement.

From the two inclusions $j\colon F \hookrightarrow X$, and $i\colon U \hookrightarrow X$ we obtain the adjunctions $j^* \dashv j_* \dashv j^!$, $i_! \dashv i^* \dashv i_*$ between the categories $\mathbf{Coh}(U), \mathbf{Coh}(X)$ and $\mathbf{Coh}(F)$ of coherent sheaves on the strata. Passing to their (bounded below-)derived versions we obtain functors[3]

$$\mathbf{D}(F) \xrightarrow{j_*} \mathbf{D}(X) \xrightarrow{i^*} \mathbf{D}(U) \tag{5.5}$$

giving rise to reflections and coreflections

$$\mathbf{D}(F) \underset{j^*}{\overset{j_*}{\underset{\longleftarrow}{\overset{\longrightarrow}{\top}}}} \mathbf{D}(X) \underset{i_!}{\overset{i^*}{\underset{\longleftarrow}{\overset{\longrightarrow}{\top}}}} \mathbf{D}(U) \qquad \mathbf{D}(F) \underset{j^!}{\overset{j_*}{\underset{\longleftarrow}{\overset{\longrightarrow}{\bot}}}} \mathbf{D}(X) \underset{i_*}{\overset{i^*}{\underset{\longleftarrow}{\overset{\longrightarrow}{\bot}}}} \mathbf{D}(U). \tag{5.6}$$

These functors are easily seen to satisfy axioms (1)-(4) above: see [BBD82, **1.4.3.1-5**] and [Ban07, **7.2.1**] for details.

REMARK 5.2.6. The above example, first discussed in [BBD82], is in some sense paradigmatic, and it can be seen as a motivation for the abstract definition of recollement: a generalization of Grothendieck's "six functors" formalism. Several sources [Han14, BP13, AHKL11, C+14] convey the intuition that a recollement რ is some sort of "exact sequence" of triangulated categories, thinking **D** as decomposed into two parts, an "open" and a "closed" one. This also motivates the intuition that a donnée de recollement is not symmetric.

**An algebraic example.** The algebraic counterpart of the above example involves derived categories of algebras: we borrow the following discussion from [Han14].

EXAMPLE 5.2.7. Let $A$ be an algebra, and $e \in A$ be an idempotent element; let $J = eAe$ be the ideal generated by $e$, and suppose that

- $Ae \otimes_J eA \cong J$ under the map $(xe, ey) \mapsto xey$;

[3] For a topological space $A$ we denote $\mathbf{D}(A)$ the *derived* $\infty$*-category* of coherent sheaves on $A$ defined in [Lur17, §**1.3.2**]; we also invariably denote as $j^* \dashv j_* \dashv j^!$, $i_! \dashv i^* \dashv i_*$ the functors between stable $\infty$-categories induced by the homonym functors between abelian categories.

- $\mathrm{Tor}_n^J(Ae, eA) \cong 0$ for every $n > 0$.

Then there exists a recollement

$$\mathbf{D}(A/J) \underset{\underset{i_R=\hom(A/J,-)}{\longleftarrow}}{\overset{\overset{i_L=-\otimes_A A/J}{\longleftarrow}}{\xrightarrow{i=-\otimes_{A/J} A/J}}} \mathbf{D}(A) \underset{\underset{q_R=\hom_J(Ae,-)}{\longleftarrow}}{\overset{\overset{q_L=-\otimes_J eA}{\longleftarrow}}{\xrightarrow{q=-\otimes_A Ae}}} \mathbf{D}(eAe) \tag{5.7}$$

between the derived categories of modules on the rings $A/J, A, eAe$.

Interestingly enough, also this example is paradigmatic in some sense; more precisely, *every* recollement $\mathfrak{R}\colon \mathbf{D}(A_1) \leftrightarrows \mathbf{D}(A) \leftrightarrows \mathbf{D}(A_2)$ is equivalent, in a suitable sense, to a "standard" recollement where $i_L$ and $q_L$ act by tensoring with distinguished objects $Y \in \mathbf{D}(A), Y_2 \in \mathbf{D}(A_2)$.

Definition 5.2.8. (Standard recollement): Let $\mathfrak{b} : \mathbf{D}(A_1) \leftrightarrows \mathbf{D}(A) \leftrightarrows \mathbf{D}(A_2)$ be a recollement between algebras; it is called a *standard* recollement generated by a pair $(Y, Y_2)$ if $i_L \cong - \otimes_A Y$, and $q_L \cong - \otimes_{A_2} Y_2$.

Proposition 5.2.9. Let $\mathfrak{R}\colon \mathbf{D}(A_1) \leftrightarrows \mathbf{D}(A) \leftrightarrows \mathbf{D}(A_2)$ be a recollement between algebras; then $\mathfrak{R}$ is equivalent (in the sense of Remark **5.2.14**) to a standard recollement $\mathfrak{b}$ generated by the pair $(Y, Y_2)$.

The proof relies on the following

Lemma 5.2.10. Let $A_1, A, A_2$ be algebras. The derived categories on these algebras are part of a recollement : $\mathbf{D}(A_1) \leftrightarrows \mathbf{D}(A) \leftrightarrows \mathbf{D}(A_2)$ if and only if there exist two objects $X_1, X_2 \in \mathbf{D}(A)$ such that

- $\hom(X_i, X_i) \cong A_i$ for $i = 1, 2$;
- $X_2$ is an exceptional and compact object, and $X_1$ is exceptional and self-compact;
- $X_1 \in \{X_2\}^\perp$;
- $\{X_1\}^\perp \cap \{X_2\}^\perp = (0)$.

See [Han14, §**2**] for details.

**A homotopical example.** Let $\mathrm{Ho}({}_{\mathcal{G}}\mathbf{Sp})$ be the *global stable homotopy category* of [Sch15]; this is defined as the localization of the category of globally equivariant orthogonal spectra at the homotopical class of *global equivalences* ([Sch15, Def. **1.2**]: the homotopical category ${}_{\mathcal{G}}\mathbf{Sp}$ admits a natural forgetful functor $u\colon {}_{\mathcal{G}}\mathbf{Sp} \to \mathbf{Sp}$ which "forgets the equivariancy" (it is the identity on objects, and includes the class of global equivalences in the bigger class of weak equivalences of plain spectra), which has both a left and a right adjoint $u_L, u_R$, and plays the rôle of a $q$-functor in a recollement

$$\mathbf{Sp}_+ \underset{\longleftarrow}{\overset{\longleftarrow}{\xrightarrow[i]{}}} {}_{\mathcal{G}}\mathbf{Sp} \underset{\underset{u_R}{\longleftarrow}}{\overset{\overset{u_L}{\longleftarrow}}{\xrightarrow[u]{}}} \mathbf{Sp} \tag{5.8}$$

where the functor $i\colon \mathbf{Sp}_+ \to {}_{\mathcal{G}}\mathbf{Sp}$ embeds the subcategory of orthogonal spectra that are *stably contractible* in the traditional, non-equivariant sense.

REMARK 5.2.11. Since in a stable $\infty$-category every pullback is a pushout and vice versa, any functor between stable $\infty$-categories preserving either limits or colimits preserves in particular pullout diagrams. Since left adjoints and right adjoints have this property, we find

PROPOSITION 5.2.12. (EXACTNESS OF RECOLLEMENT FUNCTORS): Each of the functors $i, i_L, i_R, q, q_L, q_R$ in a recollement situation preserves pullout diagrams.

This simple remark will be extremely useful in view of the "standard procedure" for proving results in recollement theory outlined in **5.2.25**.

DEFINITION 5.2.13. (THE ($\infty$-)CATEGORY **Recol**): A morphism between two recollements $\mathfrak{r}$ and $\mathfrak{r}'$ consists of a triple of functors $(F_0, F, F_1)$ such that the following square commutes in every part (choosing from time to time homonymous left or right adjoints):

$$\begin{array}{ccccc} \mathbf{D}^0 & \underset{\longleftarrow}{\overset{\longleftarrow}{\xrightarrow{i}}} & \mathbf{D} & \underset{\longleftarrow}{\overset{\longleftarrow}{\xrightarrow{q}}} & \mathbf{D}^1 \\ {\scriptstyle F_0}\downarrow & & {\scriptstyle F}\downarrow & & \downarrow{\scriptstyle F_1} \\ {}'\mathbf{D}^0 & \underset{\longleftarrow}{\overset{\longleftarrow}{\xrightarrow{i'}}} & {}'\mathbf{D} & \underset{\longleftarrow}{\overset{\longleftarrow}{\xrightarrow{q'}}} & {}'\mathbf{D}^1 \end{array} \tag{5.9}$$

This definition turns the collection of all recollement data into an $\infty$-category denoted **Recol** and called the ($\infty$-)category of recollements.

REMARK 5.2.14. The natural definition of equivalence between two recollement data (all three functors $(F_0, F_{01}, F_1)$ are equivalences) has an alternative reformulation (see [PS88, Thm. **2.5**]) asking that only two out of three functors are equivalences; nevertheless (*loc. cit.*) this must not be interpreted as a full 3-for-2 condition.

Equivalently we can define this notion (see [AHKL11, §**1.7**]), asking that the essential images of the fully faithful functors $(i, q_L, q_R)$ are pairwise equivalent with those of $(i', q'_L, q'_R)$.

We now concentrate on other equivalent ways to specify a recollement on a stable $\infty$-category, slightly rephrasing Definition **5.2.1**: first of all, [HJ10, Prop. **4.13.1**] shows that the localization functor $q_R q$, which is an exact localization with reflective kernel, uniquely determines the recollement datum up to equivalence; albeit of great significance as a general result, we are not interested in this perspective, and we address the interested readers to [HJ10] for a thorough discussion.

Another equivalent description of a recollement, nearer to our "torsiocentric" approach, is via a pair of $t$-structures on $\mathbf{D}$ [Nic08]:

DEFINITION 5.2.15. (STABLE TTF TRIPLE): Let $\mathbf{D}$ be a stable $\infty$-category. A *stable* TTF *triple* (short for *torsion-torsionfree triple*) on $\mathbf{D}$ is a triple of full subcategories $(\mathcal{X}, \mathcal{Y}, \mathcal{Z})$ of $\mathbf{D}$ such that both $(\mathcal{X}, \mathcal{Y})$ and $(\mathcal{Y}, \mathcal{Z})$ are $t$-structures on $\mathbf{D}$.

Notice in particular that $\mathbf{D}$ is reflected on $\mathcal{Y}$ via a functor $R^{\mathcal{Y}}$ and coreflected via a functor $S^{\mathcal{Y}}$ . The whole arrangement of categories and functors is summarized in the following diagram

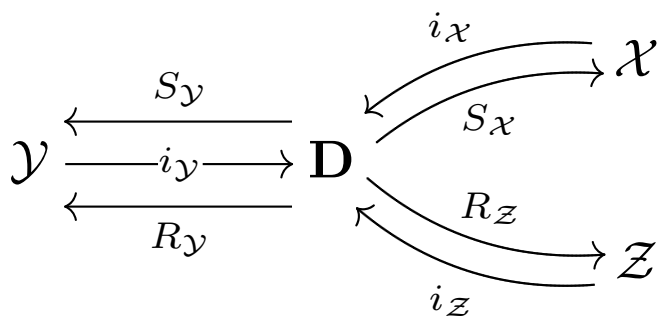

(5.10)

where $S_{\mathcal{Y}} \dashv i_{\mathcal{Y}} \dashv R_{\mathcal{Y}}$, $i_{\mathcal{Z}} \dashv R_{\mathcal{Z}}$ and $S_{\mathcal{X}} \dashv i_{\mathcal{X}}$.

Stable TTF triples are in bijection with equivalence classes of recollements, as it is recalled in [Nic08, Prop. **4.2.4**]; the same bijection holds in the stable setting, *mutatis mutandis*.

We conclude this introductory section with the following Lemma, which will be of capital importance all along §**5.3**: functors in a recollement jointly reflect isomorphisms.

Lemma 5.2.16. (Joint conservativity of recollement data): Let $\mathbf{D}$ be a stable $\infty$-category, and let

$$(i, q) \colon \mathbf{D}^0 \leftrightarrows \mathbf{D} \leftrightarrows \mathbf{D}^1$$

be a recollement on $\mathbf{D}$. Then the following conditions are equivalent for an arrow $f \in \hom(\mathbf{D})$:

- $f$ is an isomorphism in $\mathbf{D}$;
- $q(f)$ is an isomorphism in $\mathbf{D}^1$ and $i_R(f)$ is an isomorphism in $\mathbf{D}^0$;
- $q(f)$ is an isomorphism in $\mathbf{D}^1$ and $i_L(f)$ is an isomorphism in $\mathbf{D}^0$.

In other words, the pairs of functors $\{q, i_R\}$ and $\{q, i_L\}$ *jointly reflect isomorphisms.*

*Proof.* We only prove that if $q(f)$ and $i_L(f)$ are isomorphisms in the respective codomains, then $f$ is an isomorphism in $\mathbf{D}$. We need a preparatory sub-lemma, namely that the pair $\{q, i_L\}$ reflects zero objects; the only non trivial part of this statement is that if $qD \cong 0$ in $\mathbf{D}^1$ and $i_L D \cong 0$ in $\mathbf{D}^0$, then $D \cong 0$ in $\mathbf{D}$, an obvious statement in view of axiom (3) of Def. **5.2.1**, since $qD \cong 0$ entails $D \cong i(D')$, and now $0 \cong i_L(D) = i_L i D' \cong D$.

With this preliminary result, we recall that $f \colon X \to Y$ is an isomorphism if and only if $\mathrm{fib}(f) \cong 0$, and apply the previous result, together with the fact that recollement functors preserve pullouts.

Replacing $i_L$ with $i_R$, the proof shows a similar statement about the joint reflectivity of $\{q, i_R\}$. ☐

Notation 5.2.17. We will often use a rather intuitive shorthand, writing $\{q, i_L\}(f)$, or $\{q, i_R\}(f)$ to both functors applied to the same arrow. For example:

- Given (the left classes of) a pair of $t$-structures $\mathbf{D}^0_{\geq 0}, \mathbf{D}^1_{\geq 0}$ we write "$\{q, i_L\}(D) \in \mathbf{D}_{\geq 0}$" (see Thm. **5.2.19**) to denote that the object $qD \in \mathbf{D}^1_{\geq 0}$ and $i_L(D) \in \mathbf{D}^0_{\geq 0}$; similarly for $\{q, i_R\}(D) \in \mathbf{D}_{<0}$ and other combinations.
- Given (the left classes of) a pair of normal torsion theories $\mathcal{E}_0, \mathcal{E}_1$, we write "$\{q, i_{L/R}\}(f) \in \mathcal{E}$" (see Thm. **5.3.4**) to denote that the arrow $f \in \hom(\mathbf{D})$ is such that $qf \in \mathcal{E}_1$ and $i_{L/R}(f) \in \mathcal{E}_0$; similarly for $\{q, i_{L/R}\}(g) \in \mathcal{M}$ and other combinations.

Remark 5.2.18. The joint reflectivity of the recollement functors $\{q, i_L\}$ or $\{q, i_R\}$ can be seen as an analogue, in the setting of an abstract recollement, of the fact that in the geometric case of the recollement induced by a stratification $\varnothing \subset U \subset X$ one has ([PS88, **2.3**]) that a morphism of sheaves $\varphi\colon \mathcal{F} \to \mathcal{F}'$ on $X$ is uniquely determined by its restrictions $\varphi\big|_U$ and $\varphi\big|_{X\smallsetminus U}$.

### 5.2.1 The classical gluing of $t$-structures.

The main result in the classical theory of recollements is the so-called *gluing theorem*, which tells us how to obtain a $t$-structure $\mathfrak{t} = \mathfrak{t}_0 \between \mathfrak{t}_1$[4] on $\mathbf{D}$ starting from two $t$-structures $\mathfrak{t}_i$ on the categories $\mathbf{D}^i$ of a recollement $\mathfrak{R}$.

Theorem 5.2.19. (Gluing Theorem): Consider a recollement

$$\mathfrak{R} = (i, q)\colon \mathbf{D}^0 \leftrightarrows \mathbf{D} \leftrightarrows \mathbf{D}^1,$$

and let $\mathfrak{t}_i$ be $t$-structures on $\mathbf{D}^i$ for $i = 0, 1$; then there exists a $t$-structure on $\mathbf{D}$, called the *gluing* of the $\mathfrak{t}_i$ (along the recollement $\mathfrak{R}$, but this specification is almost always omitted) and denoted $\mathfrak{t}_0 \between \mathfrak{t}_1$, whose classes $\big((\mathbf{D}^0 \between \mathbf{D}^1)_{\geq 0}, (\mathbf{D}^0 \between \mathbf{D}^1)_{<0}\big)$ are given by

$$\begin{aligned}(\mathbf{D}^0 \between \mathbf{D}^1)_{\geq 0} &= \Big\{X \in \mathbf{D} \mid (qX \in \mathbf{D}^1_{\geq 0}) \wedge (i_L X \in \mathbf{D}^0_{\geq 0})\Big\};\\ (\mathbf{D}^0 \between \mathbf{D}^1)_{<0} &= \Big\{X \in \mathbf{D} \mid (qX \in \mathbf{D}^1_{<0}) \wedge (i_R X \in \mathbf{D}^0_{<0})\Big\}.\end{aligned} \tag{5.11}$$

Remark 5.2.20. Following Notation **5.2.17** we have that $X \in \mathbf{D}_{\geq 0}$ iff $\{q, i_L\}(X) \in \mathbf{D}_{\geq 0}$ and $Y \in \mathbf{D}_{<0}$ iff $\{q, i_R\}(X) \in \mathbf{D}_{<0}$, which is a rather evocative statement: the left/right class of $\mathfrak{t}_0 \between \mathfrak{t}_1$ is determined by the left/right adjoint to $i$.

Remark 5.2.21. The "wrong way" classes

$$\begin{aligned}(\mathbf{D}^0 \between \mathbf{D}^1)^{\bigstar}_{\geq 0} &= \Big\{X \in \mathbf{D} \mid (\{q, i_R\}X \in \mathbf{D}_{\geq 0}\Big\};\\ (\mathbf{D}^0 \between \mathbf{D}^1)^{\bigstar}_{<0} &= \Big\{X \in \mathbf{D} \mid (\{q, i_L\}X \in \mathbf{D}_{<0}\Big\}.\end{aligned} \tag{5.12}$$

[4] The symbol $\between$ (pron. *glue*) recalls the alchemical token describing the process of *amalgamation* between two or more elements (one of which is often mercury): although amalgamation is not recognized as a proper stage of the *Magnum Opus*, several sources testify that it belongs to the alchemical tradition (see [RS76, pp. **409-498**]).

do not define a $t$-structure in general. However they do in the case the recollement situation $\mathfrak{R}$ is the lower part of a *2-recollement*, i.e. there exists a diagram of the form

$$\mathbf{C}^0 \begin{array}{c} \xleftarrow{\;i_1\;} \\ \xrightarrow{\;i_2\;} \\ \xleftarrow{\;i_3\;} \\ \xrightarrow[\;i_4\;]{} \end{array} \mathbf{C} \begin{array}{c} \xleftarrow{\;q_1\;} \\ \xrightarrow{\;q_2\;} \\ \xleftarrow{\;q_3\;} \\ \xrightarrow[\;q_4\;]{} \end{array} \mathbf{C}^1 \tag{5.13}$$

where both

$$\mathfrak{R}_2 = \mathbf{C}^0 \begin{array}{c} \xleftarrow{\;i_1\;} \\ \xrightarrow{\;i_2\;} \\ \xleftarrow[\;i_3\;]{} \end{array} \mathbf{C} \begin{array}{c} \xleftarrow{\;q_1\;} \\ \xrightarrow{\;q_2\;} \\ \xleftarrow[\;q_3\;]{} \end{array} \mathbf{C}^1 \tag{5.14}$$

and

$$\mathfrak{R}_3 = \mathbf{C}^1 \begin{array}{c} \xleftarrow{\;q_2\;} \\ \xrightarrow{\;q_3\;} \\ \xleftarrow[\;q_4\;]{} \end{array} \mathbf{C} \begin{array}{c} \xleftarrow{\;i_2\;} \\ \xrightarrow{\;i_3\;} \\ \xleftarrow[\;i_4\;]{} \end{array} \mathbf{C}^0 \tag{5.15}$$

are recollements, with $\mathfrak{R} = \mathfrak{R}_3$. Indeed, in this situation one has

$$\begin{aligned} (\mathbf{D}^0 \veebar \mathbf{D}^1)^{\bigstar}_{\geq 0} &= \Big\{ X \in \mathbf{D} \mid (\{q, i_R\} X \in \mathbf{D}_{\geq 0} \Big\} \\ &= \Big\{ X \in \mathbf{C} \mid (\{i_3, q_2\} X \in \mathbf{C}_{\geq 0} \Big\} \\ &= (\mathbf{C}^0 \veebar^{\mathfrak{R}_2} \mathbf{C}^1)_{\geq 0}. \end{aligned}$$

More generally, an $n$-recollement is defined as the datum of three stable $\infty$-categories $\mathbf{C}^0, \mathbf{C}, \mathbf{C}^1$ organized in a diagram

$$\mathbf{C}^0 \begin{array}{c} \xleftarrow{\;i_1\;} \\ \xrightarrow{\;i_2\;} \\ \xleftarrow{\;i_3\;} \\ \vdots \\ \xleftarrow[\;i_{n+2}\;]{} \end{array} \mathbf{C} \begin{array}{c} \xleftarrow{\;q_1\;} \\ \xrightarrow{\;q_2\;} \\ \xleftarrow{\;q_3\;} \\ \vdots \\ \xleftarrow[\;q_{n+2}\;]{} \end{array} \mathbf{C}^1 \tag{5.16}$$

with $n+2$ functors on each edge, such that every consecutive three functors form recollements $\mathfrak{R}_{2k} = (i_{2k}, q_{2k})$, $\mathfrak{R}_{2h+1} = (q_{2h+1}, i_{2h+1})$, for $k = 1, \dots, n-1$, $h = 1, \dots, n-2$, see [HQ14, Def. **2**]. Applications of this formalism to derived categories of algebras, investigating the relationships between the recollements of derived categories and the Gorenstein properties of these algebras, can be found in [HQ14, Qin15].

Notation 5.2.22. It is worth noting that $\mathbf{D}^0 \veebar \mathbf{D}^1$ has no real meaning as a category; this is only an intuitive shorthand to denote the pair $(\mathbf{D}, \mathfrak{t}_0 \veebar \mathfrak{t}_1)$; more explicitly, it is a shorthand to denote the following situation:

> The stable $\infty$-category $\mathbf{D}$ fits into a recollement $(i, q)\colon \mathbf{D}^0 \leftrightarrows \mathbf{D} \leftrightarrows \mathbf{D}^1$, $t$-structures on $\mathbf{D}^0$ and $\mathbf{D}^1$ have been chosen, and $\mathbf{D}$ is endowed with the glued $t$-structure $\mathfrak{t}_0 \veebar \mathfrak{t}_1$.

A proof of the gluing theorem in the classical setting of triangulated categories can be found in [Ban07, Thm. **7.2.2**] or in the standard reference [BBD82]. We briefly sketch the argument given in [Ban07] as we will need it in the torsio-centric reformulation of the gluing theorem.

*Proof of Thm. 5.2.19.* We begin showing the way in which every $X \in \mathbf{D}$ fits into a fiber sequence $SX \to X \to RX$ where $SX \in (\mathbf{D}^0 \uplus \mathbf{D}^1)_{\geq 0}, RX \in (\mathbf{D}^0 \uplus \mathbf{D}^1)_{<0}$. Let $\mathbb{F}_i$ denote the normal torsion theory on $\mathbf{D}^i$, inducing the $t$-structure $\mathfrak{t}_i$; let $\eta_1 \colon qX \to R_1 qX$ be the arrow in the fiber sequence

$$S_1 qX \xrightarrow{\epsilon_1} qX \xrightarrow{\eta_1} R_1 qX \tag{5.17}$$

obtained thanks to $\mathbb{F}_1$; let $\hat{\eta}$ be its *mate* $X \to q_R R_1 qX$ in $\mathbf{D}$ under the adjunction $q \dashv q_R$, and let $W = \mathrm{fib}(\hat{\eta})$.

Now, consider $i_L W$ in the fiber sequence

$$S_0 i_L W \xrightarrow{\sigma_0} i_L W \xrightarrow{\theta_0} R_0 i_L W \tag{5.18}$$

induced by $\mathbb{F}_0$ on $\mathbf{D}_0$, and its mate $\hat{\theta} \colon W \to iR_0 i_L W$; take its fiber $SX$, and the object $RX$ defined as the pushout of $iR_0 i_L W \xleftarrow{\hat{\theta}} W \to X$.

To prove that these two objects are the candidate co/truncation we consider the diagram

$$\begin{array}{ccccc} SX & \longrightarrow & W & \longrightarrow & X \\ \downarrow & & \downarrow{\scriptstyle\hat{\theta}} & & \downarrow \\ 0 & \longrightarrow & iR_0 i_L W & \longrightarrow & RX \\ & & \downarrow & & \downarrow \\ & & 0 & \longrightarrow & q_R R_1 qX \end{array} \qquad (\text{dotted } \hat{\eta} \colon X \to q_R R_1 qX) \tag{5.19}$$

where all the mentioned objects fit, and where every square is a pullout. We have to prove that $SX \in (\mathbf{D}^0 \uplus \mathbf{D}^1)_{\geq 0}$ and $RX \in (\mathbf{D}^0 \uplus \mathbf{D}^1)_{<0}$. To do this, apply the functors $q, i_L, i_R$ to (**5.19**), obtaining the following diagram of pullout squares (recall the exactness properties of the recollement functors, stated in Prop. **5.2.12**):

$$\begin{array}{ccccc} qSX & \xrightarrow{\sim} & qW & \longrightarrow & qX \\ \downarrow & & \downarrow & & \downarrow \\ 0 & = & 0 & \longrightarrow & qRX \\ & & \| & & \downarrow \\ & & 0 & \longrightarrow & R_1 qX \end{array} \quad \begin{array}{ccccc} i_L SX & \longrightarrow & i_L W & \longrightarrow & i_L X \\ \downarrow & ① & \downarrow & & \downarrow \\ 0 & \longrightarrow & R_0 i_L W & \longrightarrow & i_L RX \\ & & \downarrow & & \downarrow \\ & & 0 & \longrightarrow & i_L q_R R_1 qX \end{array} \quad \begin{array}{ccccc} i_R SX & \longrightarrow & i_R W & \longrightarrow & i_R X \\ \downarrow & & \downarrow & & \downarrow \\ 0 & \longrightarrow & R_0 i_L W & \xrightarrow{\sim} & i_R RX \\ & & \downarrow & & \downarrow \\ & & 0 & = & 0 \end{array}$$

where we took into account the relations $qi = 0, i_R q_R = 0 = i_L q_L$, we find that

- $qSX \cong qW \cong S_1 qX \in \mathbf{D}^1_{\geq 0}$, since $0 \to S_1 qX$ lies in $\mathcal{M}_1$, and $qRX \cong R_1 qX \in \mathbf{D}^1_{<0}$;
- $i_L SX \cong S_0 i_L W \in \mathbf{D}^0_{\geq 0}$, looking square ①;
- $i_R RX \cong R_0 i_L W \in \mathbf{D}^0_{<0}$.

It remains to show that the two classes $\mathbf{D}_{\geq 0}, \mathbf{D}_{<0}$ are orthogonal; to see this, suppose that $X \in \mathbf{D}_{\geq 0}$ and $Y \in \mathbf{D}_{<0}$. We consider the fiber sequence $ii_R Y \to Y \to q_R qY$ of axiom (4) in Def. **5.2.1**, to obtain (applying the homological functor $\mathbf{D}(X, -)$)

$$\begin{array}{ccccc} \mathbf{D}(X, ii_R Y) & \longrightarrow & \mathbf{D}(X, Y) & \longrightarrow & \mathbf{D}(X, q_R qY) \\ \| & & & & \| \\ \mathbf{D}(i_L X, i_R Y) & & & & \mathbf{D}(qX, qY) \\ \| & & & & \| \\ 0 & & & & 0 \end{array} \tag{5.20}$$

and we conclude, thanks to the exactness of this sequence. □

REMARK 5.2.23. The definition of $\mathfrak{t}_0 \veebar \mathfrak{t}_1$ entails that all the recollement functors $(i, q)\colon (\mathbf{D}^0, \mathfrak{t}_0) \leftrightarrows (\mathbf{D}, \mathfrak{t}_0 \veebar \mathfrak{t}_1) \leftrightarrows (\mathbf{D}^1, \mathfrak{t}_1)$ become $t$-exact in the sense of [Lur17, Def. **1.3.3.1**].

REMARK 5.2.24. Strictly speaking, the domain of definition of the gluing operation $\veebar$ is the set of triples $(\mathfrak{t}_0, \mathfrak{t}_1, \mathfrak{G})$ where $(\mathfrak{t}_0, \mathfrak{t}_1) \in \text{TS}(\mathbf{D}^0) \times \text{TS}(\mathbf{D}^1)$ and $\mathfrak{G} = (i, q)$ is a recollement $\mathbf{D}^0 \leftrightarrows \mathbf{D} \leftrightarrows \mathbf{D}^1$, but unless this (rather stodgy) distinction is strictly necessary we will adopt an obvious abuse of notation.

REMARK 5.2.25. (A STANDARD TECHNIQUE): The procedure outlined above is in some sense paradigmatic, and it's worth to trace it out as an abstract way to deduce properties about objects and arrows fitting in a diagram like (**5.19**). This algorithm will be our primary technique to prove statements in the "torsio-centric" formulation of recollements:

- We start with a particular diagram, like for example (**5.19**) or (**5.22**) below; our aim is to prove that a property (being invertible, being the zero map, lying in a distinguished class of arrows, etc.) is true for an arrow $h$ in this diagram.
- We apply (possibly only some of) the recollement functors to the diagram, and we deduce that $h$ has the above property from
  - The recollement relations between the functors (Def. **5.2.1**);
  - The exactness of the recollement functors (Prop. **5.2.12**);
  - The joint reflectivity of the pairs $\{q, i_L\}$ and $\{q, i_R\}$ (Lemma **5.2.16**);

## 5.3 Stable Recollements.

וַיַּחֲלֹם וְהִנֵּה סֻלָּם מֻצָּב אַרְצָה וְרֹאשׁוֹ מַגִּיעַ הַשָּׁמָיְמָה וְהִנֵּה
מַלְאֲכֵי אֱלֹהִים עֹלִים וְיֹרְדִים בּוֹ׃

[ER77], Genesis 28:12

### 5.3.1 The Jacob's ladder: building co/reflections.

The above procedure to build the functors $R, S$ depends on several choices (we forget half of the fiber sequence $S_1qX \to qX \to R_1qX$) and it doesn't seem independent from these choices, at least at first sight.

The scope of this first subsection is to show that this apparent asymmetry arises only because we are hiding half of the construction, taking into account only half of the fiber sequence (**5.17**). Given an object $X \in \mathbf{D}$ a dual argument yields *another* way to construct a fiber sequence

$$S'X \to X \to R'X \tag{5.21}$$

out of the recollement data, which is naturally isomorphic to the former $SX \to X \to RX$.

We briefly sketch how this dualization process goes: starting from the coreflection arrow $\epsilon_1 \colon S_1qX \to qX$, taking its mate $q_LS_1qX \to X$ under the adjunction $q_L \dashv q$, and reasoning about its cofiber we can build a diagram which is dual to the former one, and where every square is a pullout:

$$\begin{array}{ccccc} q_LS_1qX & \longrightarrow & S'X & \longrightarrow & X \\ \downarrow & & \downarrow & & \downarrow \\ 0 & \longrightarrow & iS_0i_RK & \longrightarrow & K \\ & & \downarrow & & \downarrow \\ & & 0 & \longrightarrow & R'X \end{array} \tag{5.22}$$

Proposition 5.3.1. (The Jacob's ladder): The two squares of the previous constructions fit into a "ladder" induced by canonical isomorphisms $SX \cong S'X, RX \cong R'X$; the construction is functorial in $X$. The "Jacob's ladder" is the following diagram:

$$\begin{array}{ccccccc}
q_L S_1 qX & \longrightarrow & SX & \longrightarrow & W & \longrightarrow & X \\
\downarrow & & \downarrow & & \downarrow & & \downarrow \\
0 & \longrightarrow & iS_0 i_R K & \longrightarrow & C & \longrightarrow & K \\
 & & \downarrow & & \downarrow & & \downarrow \\
 & & 0 & \longrightarrow & iR_0 i_L W & \longrightarrow & RX \\
 & & & & \downarrow & & \downarrow \\
 & & & & 0 & \longrightarrow & q_R R_1 qX
\end{array} \tag{5.23}$$

*Proof.* It suffices to prove that both $SX, S'X$ lie in $\mathbf{D}_{\geq 0}$ and both $RX, R'X$ lie in $\mathbf{D}_{\leq 0}$; given this, we can appeal (a suitable stable $\infty$-categorical version of) [BBD82, Prop. **1.1.9**] which asserts the functoriality of the truncation functors, i.e. that when the same object $X$ fits into *two* fiber sequences arising from the same normal torsion theory, then there exist the desired isomorphisms.[5]

The procedure showing this is actually the same remarked in **5.2.25**: we apply $q, i_L, i_R$ to the diagram (**5.22**) and we exploit exactness of the recollement functors to find pullout diagrams showing that $L \in \mathbf{D}_{<0}$ and $J \in \mathbf{D}_{\geq 0}$.

Once these isomorphisms have been found, it remains only to glue the two sub-diagrams

$$\begin{array}{ccccccc}
q_L S_1 qX & \dashrightarrow & SX & \longrightarrow & W & \longrightarrow & X \\
\downarrow & & | & & | & & | \\
0 & \dashrightarrow & iS_0 i_R K & \dashrightarrow & C & \dashrightarrow & K \\
 & & \downarrow & & \downarrow & & \downarrow \\
 & & 0 & \longrightarrow & iR_0 i_L W & \longrightarrow & RX \\
 & & & & \downarrow & & \downarrow \\
 & & & & 0 & \longrightarrow & q_R R_1 qX
\end{array}
\qquad
\begin{array}{ccccccc}
q_L S_1 qX & \longrightarrow & S'X & — & W & \longrightarrow & X \\
\downarrow & & \downarrow & & \vdots & & \downarrow \\
0 & \longrightarrow & iS_0 i_R K & — & C & \longrightarrow & K \\
 & & \downarrow & & \vdots & & \downarrow \\
 & & 0 & — & iR_0 i_L W & \longrightarrow & R'X \\
 & & & & \vdots & & \vdots \\
 & & & & 0 & \dashrightarrow & q_R R_1 qX
\end{array}$$

to obtain the ladder. Now, this construction is obtained by taking into account the fiber sequence $S_1 qX \to qX \to R_1 qX$ as a whole, and since this latter object is uniquely determined up to isomorphism, we obtain a

[5] In a torsio-centric perspective, this follows from the uniqueness of the factorization of a morphism with respect to the normal torsion theory having reflection $R$ and coreflection $S$; see **1.3.1**.

diagram of endofunctors

$$\begin{array}{ccccccc}
q_L S_1 q & \longrightarrow & S & \longrightarrow & \omega & \longrightarrow & 1 \\
\downarrow & & \downarrow & & \downarrow & & \downarrow \\
0 & \longrightarrow & iS_0 i_R \kappa & \longrightarrow & \gamma & \longrightarrow & \kappa \\
 & & \downarrow & & \downarrow & & \downarrow \\
 & & 0 & \longrightarrow & iR_0 i_L \omega & \longrightarrow & R \\
 & & & & \downarrow & & \downarrow \\
 & & & & 0 & \longrightarrow & q_R R_1 q
\end{array} \tag{5.24}$$

where every square is a pullout (again giving to a category of functors the obvious stable structure [Lur17, Prop. **1.1.3.1**]), and where we commuted in Greek alphabet the functorial dependence of objects $\kappa(X) = K$, $\gamma(X) = C$, $\omega(X) = W$ from the above procedure. Notice also that this latter diagram of functors uses homogeneously all the recollement functors, and that it is "symmetric" with respect to the antidiagonal (it switches left and right adjoints, as well as reflections and coreflections). □

The functors $S, R$ are the co/truncations for the recollée $t$-structure, and the normality of the torsion theory is witnessed by the pullout subdiagram

$$\begin{array}{ccccc}
SX & \longrightarrow & W & \longrightarrow & X \\
\downarrow & & \downarrow & & \downarrow \\
iS_0 i_R K & \longrightarrow & C & \longrightarrow & K \\
\downarrow & & \downarrow & & \downarrow \\
0 & \longrightarrow & iR_0 i_L W & \longrightarrow & RX.
\end{array} \tag{5.25}$$

Notation 5.3.2. From now on, we will always refer to the diagram above as "the Jacob ladder" of an object $X \in \mathbf{D}$, and/or to the diagram induced by a morphism $f\colon X \to Y$ between the ladder of the domain and the codomain, i.e. to three-dimensional diagrams like

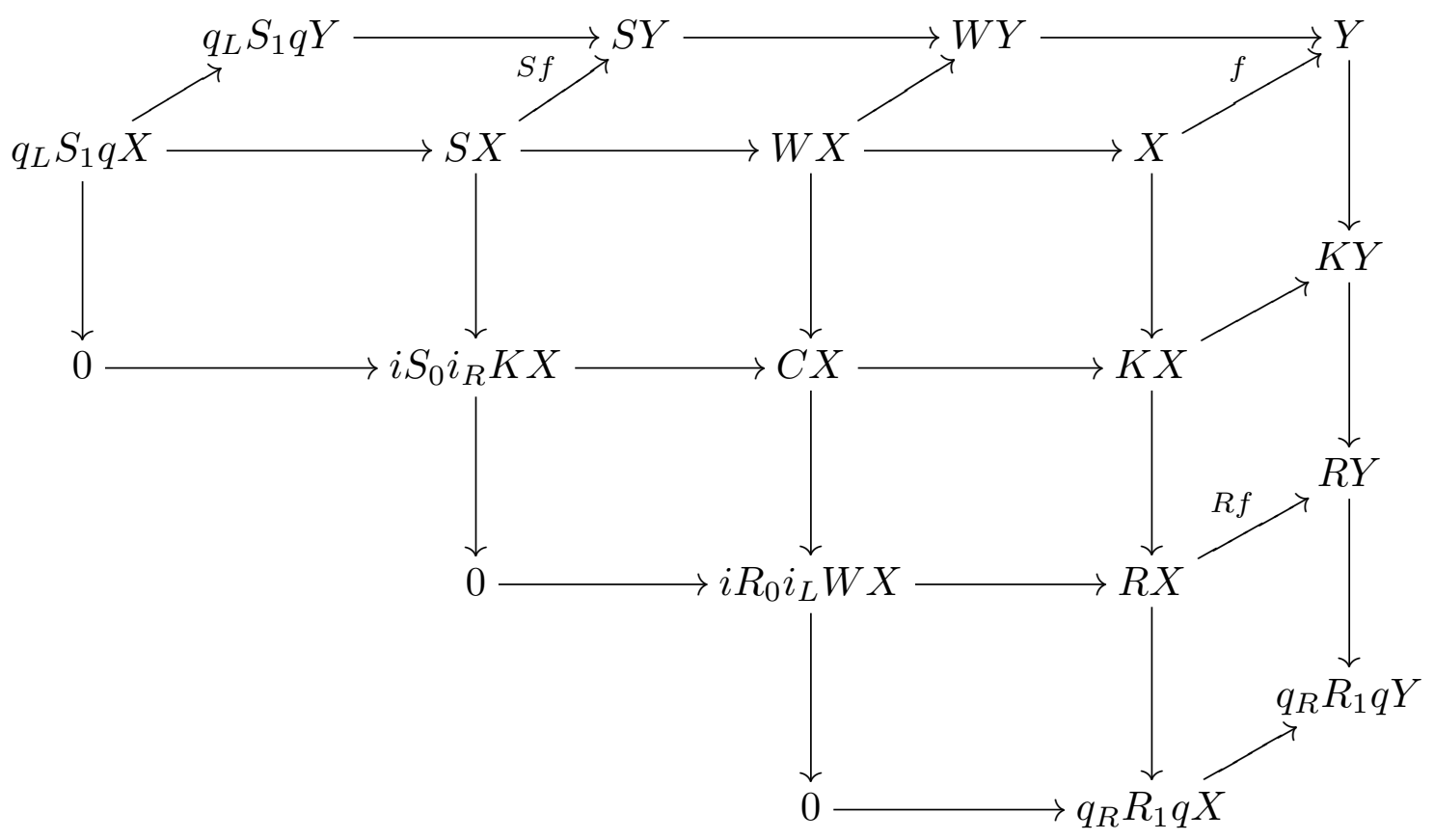

(5.26)

### 5.3.2 The NTT of a recollement.

Throughout this subsection we outline the torsio-centric translation of the classical results recalled above. In particular we give an explicit definition of the $\curlyvee$ operation when it has been "transported" to the set of normal torsion theories, independent from its characterization in terms of the pairs aisle-coaisle of the two $t$-structures. From now on we assume given a recollement

$$\mathbf{D}^0 \underset{\overset{\longleftarrow}{i_L}}{\overset{\overset{i_R}{\longleftarrow}}{\xrightarrow{\;i\;}}} \mathbf{D} \underset{\overset{\longleftarrow}{q_L}}{\overset{\overset{q_R}{\longleftarrow}}{\xrightarrow{\;q\;}}} \mathbf{D}^1.$$

Given $t$-structures $\mathfrak{t}_i \in \textsc{ts}(\mathbf{D}^i)$, in view of our characterization theorem **3.1.1**, there exist normal torsion theories $\mathbb{F}_i = (\mathcal{E}_i, \mathcal{M}_i)$ on $\mathbf{D}^i$ such that $(\mathbf{D}^i_{\geq 0}, \mathbf{D}^i_{<0})$ are the classes $(0/\mathcal{E}_i, \mathcal{M}_i/0)$ of torsion and torsion-free objects of $\mathbf{D}^i$, for $i = 0, 1$; an object $X$ lies in $(\mathbf{D}^0 \curlyvee \mathbf{D}^1)_{\geq 0}$ if and only if $qX \in \mathcal{E}_1$ and $i_LX \in \mathcal{E}_0$[(6)], and similarly an object $Y$ lies in $\mathbf{D}_{\leq 0}$ if and only if $qY \in \mathcal{M}_1$ and $i_RY \in \mathcal{M}_0$.

Remark 5.3.3. The $t$-structure $\mathfrak{t} = \mathfrak{t}_0 \curlyvee \mathfrak{t}_1$ on $\mathbf{D}$ must itself come from a normal torsion theory which we denote $\mathbb{F}_0 \curlyvee \mathbb{F}_1$ on $\mathbf{D}$, so that $\big((\mathbf{D}^0 \curlyvee \mathbf{D}^1)_{\geq 0}, (\mathbf{D}^0 \curlyvee \mathbf{D}^1)_{<0}\big) = \big(0/(\mathcal{E}_0 \curlyvee \mathcal{E}_1), (\mathcal{M}_0 \curlyvee \mathcal{M}_1)/0\big)$; in other words the following three conditions are equivalent for an object $X \in \mathbf{D}$:

- $X$ lies in $(\mathbf{D}^0 \curlyvee \mathbf{D}^1)_{\geq 0}$;
- $X$ lies in $\mathcal{E}_0 \curlyvee \mathcal{E}_1$, i.e. $RX \cong 0$ in the notation of (**5.25**);
- $\{q, i_L\}(X) \in \mathcal{E}$, following Notation **5.2.17**.

(6) Thanks to the Sator lemma we are allowed to use "$X \in \mathcal{K}$" as a shorthand to denote that either the initial arrow $\left[\begin{smallmatrix} 0 \\ \downarrow \\ X \end{smallmatrix}\right]$ or the terminal arrow $\left[\begin{smallmatrix} X \\ \downarrow \\ 0 \end{smallmatrix}\right]$ lie in a 3-for-2 class $\mathcal{K} \subset \hom(\mathbf{C})$. From now on we will adopt this notation.

We now aim to a torsio-centric characterization of the classes $(\mathcal{E}_0 \uplus \mathcal{E}_1, \mathcal{M}_0 \uplus \mathcal{M}_1)$, relying on the factorization properties of $(\mathcal{E}_i, \mathcal{M}_i)$ alone: since we proved Thm. **5.2.19** above, there must be a normal torsion theory $\mathbb{F}_0 \uplus \mathbb{F}_1 = \big(\mathcal{E}_0 \uplus \mathcal{E}_1, \mathcal{M}_0 \uplus \mathcal{M}_1\big)$ inducing $\mathfrak{t}_0 \uplus \mathfrak{t}_1$ as $\big(0/(\mathcal{E}_0 \uplus \mathcal{E}_1), (\mathcal{M}_0 \uplus \mathcal{M}_1)/0\big)$: in other words,

> $\mathbb{F}_0 \uplus \mathbb{F}_1$ is the (unique) normal torsion theory whose torsion/torsionfree classes are $\big((\mathbf{D}^0 \uplus \mathbf{D}^1)_{\geq 0}, (\mathbf{D}^0 \uplus \mathbf{D}^1)_{<0}\big)$ of Thm. **5.2.19**,

Clearly this is only an application of our "Rosetta stone" theorem **3.1.1**, so in some sense this result is "tautological". But there are at least two reasons to concentrate in "proving again" Thm. **5.2.19** from a torsio-centric perspective:

- The construction offered by the Rosetta stone is rather indirect, and only appropriate to show formal statements about the factorization system $\mathbb{F}(\mathfrak{t})$ induced by a $t$-structure;
- In a stable setting, the torsio-centric point of view, using factorization systems, is more primitive and more natural than the classical one using 1-categorical arguments (i.e., $t$-structures $\mathfrak{t}$ on the homotopy category of a stable $\mathbf{D}$ are induced by normal torsion theories in $\mathbf{D}$; in the quotient process one loses important informations about $\mathfrak{t}$).

Both these reasons lead us to adopt a "constructive" point of view, giving an explicit characterization of $\mathbb{F}_0 \uplus \mathbb{F}_1$ which relies on properties of the factorization systems $\mathbb{F}_0$, $\mathbb{F}_1$ alone, independent from triangulated categorical arguments.

In the following section we will discuss the structure and properties of the factorization system $\mathbb{F}_0 \uplus \mathbb{F}_1$, concentrating on a self-contained and categorically well motivated construction of the classes $\mathcal{E}_0 \uplus \mathcal{E}_1$ and $\mathcal{M}_0 \uplus \mathcal{M}_1$ starting from an obvious *ansatz* which follows Remark **5.3.3**.

The discussion above, and in particular the fact that an initial/terminal arrow $0 \leftrightarrows X$ lies in $\mathcal{E}_0 \uplus \mathcal{E}_1$ if and only if $\{q, i_L\}(X) \in \mathcal{E}$, suggests that we define $\mathcal{E}_0 \uplus \mathcal{E}_1 = \big\{f \in \hom(\mathbf{D}) \mid \{q, i_L\}(f) \in \mathcal{E}\big\}$ and $\mathcal{M}_0 \uplus \mathcal{M}_1 = \big\{g \in \hom(\mathbf{D}) \mid \{q, i_R\}(g) \in \mathcal{M}\big\}$. Actually it turns out that this guess is not far from being correct: the correct classes are indeed given by the following:

THEOREM 5.3.4. Let $\mathbf{D}$ be a stable $\infty$-category, in a recollement

$$(i, q)\colon \mathbf{D}^0 \leftrightarrows \mathbf{D} \leftrightarrows \mathbf{D}^1,$$

and let $\mathfrak{t}_i$ be a $t$-structure on $\mathbf{D}^i$. Then the recollée $t$-structure $\mathfrak{t}_0 \uplus \mathfrak{t}_1$ is induced by the normal torsion theory $(\mathcal{E}_0 \uplus \mathcal{E}_1, \mathcal{M}_0 \uplus \mathcal{M}_1)$ with classes

$$\mathcal{E}_0 \uplus \mathcal{E}_1 = \big\{f \in \hom(\mathbf{D}) \mid \{q, i_L W\}(f) \in \mathcal{E}\big\}; \tag{5.27}$$

$$\mathcal{M}_0 \uplus \mathcal{M}_1 = \big\{g \in \hom(\mathbf{D}) \mid \{q, i_R K\}(g) \in \mathcal{M}\big\}. \tag{5.28}$$

*Proof.* We only need to prove the statement for $\mathcal{E}_0 \uplus \mathcal{E}_1$, since the statement for $\mathcal{M}_0 \uplus \mathcal{M}_1$ is completely specular. Thanks to the discussion in section

§**5.2**, an arrow $f \in \hom(\mathbf{D})$ lies in $\mathcal{E}_0 \uplus \mathcal{E}_1$ if and only if $Rf$ (as constructed in the Jacob ladder (**5.26**)) is an isomorphism in $\mathbf{D}$, so we are left to prove that, given $f \in \hom(\mathbf{D})$:

$$\mathbf{D} \text{ in isomorphism } Rf \quad \Leftrightarrow \quad \{q, i_L W\}(f) \in \mathcal{E}. \tag{5.29}$$

Equivalently, we have to prove that

$$Rf \text{ isomorphism} \quad \Leftrightarrow \quad . \text{ isomorphisms are } \{R_1 q, R_0 i_L W\}(f) \tag{5.30}$$

We begin by showing that if $\{R_1 q, R_0 i_L W\}(f)$ are isomorphisms, then also $Rf$ is an isomorphism. By the joint conservativity of the recollement data (Lemma 5.2.16) we need to prove that if $\{R_1 q, R_0 i_L W\}(f)$ are isomorphisms, then both $qRf$ and $i_L Rf$ are isomorphisms. Apply the functor $q$ to the Jacob ladder (**5.26**), to obtain

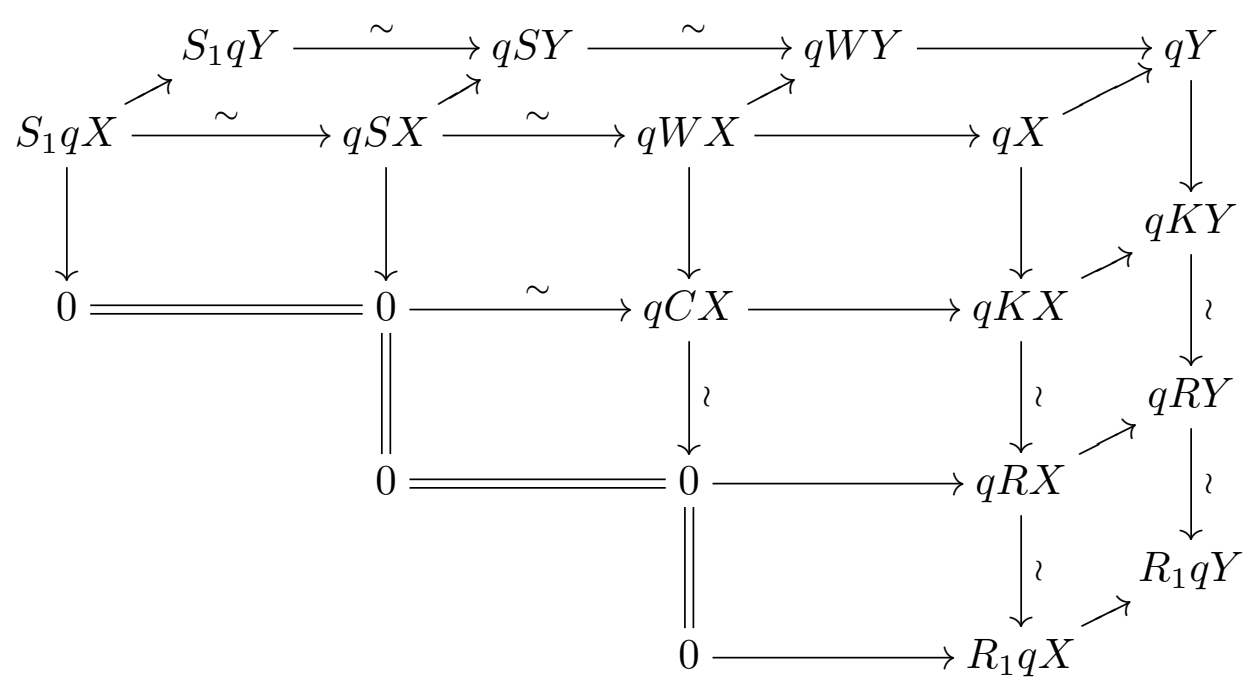

(5.31)

Hence $qRf$ is an isomorphism, since it fits into the square

$$\begin{array}{ccc} qRX & \longrightarrow & qRY \\ \wr\downarrow & & \downarrow\wr \\ R_1 qX & \xrightarrow[\sim]{} & R_1 qY. \end{array} \tag{5.32}$$

Now apply the functor $i_L$ to the Jacob ladder, obtaining

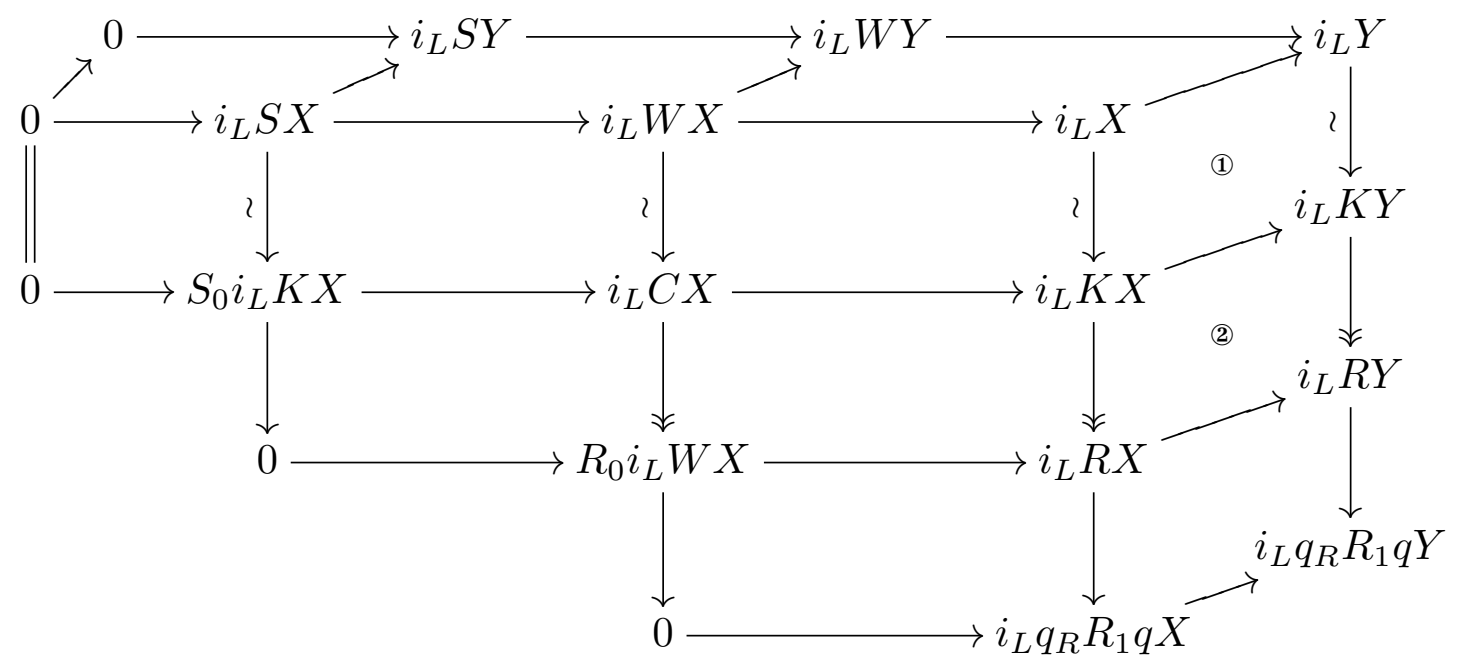

(5.33)

As noticed above, $R_1qf$ is an isomorphism, so also $i_Lq_RR_1qf$ is an isomorphism. Then $i_LRf$ is an isomorphism by the 5-lemma applied to the morphism of fiber sequences

$$\begin{array}{ccccc} & & R_0i_LWY & \longrightarrow & i_LRY \\ & \nearrow & & \nearrow & \downarrow \\ R_0i_LWX & \longrightarrow & i_LRX & & \\ \downarrow & & \downarrow & & i_Lq_RR_1qY \\ & & & \nearrow & \\ 0 & \longrightarrow & i_Lq_RR_1qX & & \end{array} \tag{5.34}$$

Vice versa: assuming $Rf$ is an isomorphism in $\mathbf{D}$, we want to prove that $\{R_1q, R_0i_LW\}(f)$ are isomorphisms. Diagram (**5.31**) gives directly that $R_1qf$ is an isomorphism, since the square

$$\begin{array}{ccc} qRX & \xrightarrow{\sim} & qRY \\ \wr\downarrow & & \downarrow\wr \\ R_1qX & \longrightarrow & R_1qY \end{array} \tag{5.35}$$

is commutative. Then, from diagram (5.34) we see that, since both $i_Lq_RR_1qf$ and $i_LRf$ are isomorphisms, so is also $R_0i_LWf$. $\square$

REMARK 5.3.5. From the sub-diagram

$$\begin{array}{ccc} i_LX & \xrightarrow{i_Lf} & i_LY \\ \wr\downarrow & ① & \downarrow\wr \\ i_LKX & \twoheadrightarrow & i_LKY \\ \downarrow & ② & \downarrow \\ i_LRX & \xrightarrow[\sim]{} & i_LRY \end{array} \tag{5.36}$$

of diagram (**5.33**) one deduces that if $Rf$ is an isomorphism, then $i_Lf \in \mathcal{E}_0$, by the 3-for-2 closure property of $\mathcal{E}_0$. This mean that $\{q, i_LW\}(f) \in \mathcal{E}$ implies that $\{q, i_L\}(f) \in \mathcal{E}$. The converse implication has no reason to be true in general. However it is true for terminal (or initial) morphisms. Namely, from the Rosetta stone one has that $X \in \mathcal{E}_0 \uplus \mathcal{E}_1$ if and only if $X \in (\mathbf{D}^0 \uplus \mathbf{D}^1)_{\geq 0}$, and so if and only if $\{q, i_L\}(X) \in \mathcal{E}$. On the other hand, $X \in \mathcal{E}_0 \uplus \mathcal{E}_1$ if and only if $\{q, i_LW\}(X) \in \mathcal{E}$. The fact that the condition $\{q, i_L\}(X) \in \mathbf{D}_{\geq 0}$ is equivalent to the condition $\{q, i_LW\}(X) \in \mathbf{D}_{\geq 0}$ can actually be easily checked directly. Namely, if $qX \in \mathbf{D}^1_{\geq 0}$, then $q_RR_1qX = 0$ and so $X = WX$ in this case. Specular considerations apply to the right class $\mathcal{M}_0 \uplus \mathcal{M}_1$.

## 5.4 Properties of recollements.

"Do what thou wilt" shall be the whole of the Law. The study of this Book is forbidden. It is wise to destroy this copy after the first reading. Whosoever disregards this does so at his own risk and peril.

Ankh-ef-en-Khonsu i

In this section we address associativity issues for the $\uplus$ operation: it is a somewhat subtle topic, offering examples of several non-trivial constructions even in the classical geometric case: it is our opinion that in a stable setting the discussion can be clarified by simple, well-known categorical properties.

We start by proving a generalization of [Ban07, BBD82] where it is stated that the gluing operation can be iterated in a preferential way determined by a *stratification* of an ambient space $X$. This result hides in fact an associativity property for the gluing operation, in a sense which our Thm. **5.4.2** below makes precise.

Suitably abstracted to a stable setting, a similar result holds true, once we are given a *Urizen compass* (a certain shape of diagram like in Def. **5.4.7**, implying certain relations and compatibilities between different recollements, which taken together ensure associativity).

### 5.4.1 Geometric associativity of the gluing.

An exhaustive account for the theory of stratified spaces can be found in [Pfl01, Ban07, Wei94]. Here, since we do not aim at a comprehensive treatment, we restrict to a sketchy recap of the basic definitions.

A *stratified space* of length $n$ consists of a pair $(X, \mathsf{s})$ where

$$\mathsf{s} \quad : \quad \varnothing = U_{-1} \subset U_0 \subset \cdots \subset U_n \subset X = U_{n+1} \tag{5.37}$$

is a chain of closed subspaces of a space $X$, subject to various technical assumptions which ensure that the homology theory we want to attach to $(X, \mathsf{s})$ is "well-behaved" in some sense.

All along the following section, we will call a *pure stratum* of a stratified space $(X, \mathsf{s})$ the set-theoretical difference $E_i = U_i \smallsetminus U_{i-1}$.

Remark 5.4.1. The definition is intentionally kept somewhat vague in various respects, first of all about the notion of "space": the definition of stratification can obviously be given in different contexts (topological spaces, topological manifolds, PL-manifolds, ...) according to the needs of the specific theory we want to build; when the stratification $\mathsf{s}$ is clear from the context, we indulge to harmless, obvious abuses of notation.

The associativity properties of $\uplus$ are deeply linked with the presence of a stratification on a space $X$, in the sense that a stratification $\mathsf{s}$ is what we need to induce additional recollements "fitting nicely" in the diagram of inclusions determined by $\mathsf{s}$. These recollements define a unique $t$-structure $\mathfrak{t}_0 \curlyvee \cdots \curlyvee \mathfrak{t}_n$, given $\mathfrak{t}_i$ on the derived categories of the pure strata.

To motivate the shape and the strength of the abstract conditions ensuring associativity of $\curlyvee$, outlined in §**5.4.2**, and in particular the definition of a Urizen compass **5.4.7**, we have to dig deep into the argument sketched in the geometric case in [BBD82, **2.1.2-3**]: we start by recalling

THEOREM 5.4.2. [Ban07, p. **158**] Let $(X, \mathsf{s})$ be a stratified space, $\{E_0, \dots, E_n\}$ the set of its pure strata, and $\mathfrak{t}_i$ be a set of $t$-structures, one on each $\mathbf{D}(E_i)$, for $i = 0, \dots, n$.

Then there exists a uniquely determined $t$-structure $\mathfrak{t}_0 \curlyvee \cdots \curlyvee \mathfrak{t}_n$ on $\mathbf{D}(X)$, obtained by an iterated gluing operation as the parenthesization $(\cdots((\mathfrak{t}_0 \curlyvee \mathfrak{t}_1) \curlyvee \mathfrak{t}_2) \curlyvee \cdots \curlyvee \mathfrak{t}_{n-1}) \curlyvee \mathfrak{t}_n$. Following Notation **5.2.22** we will refer to the pair $(\mathbf{D}(X), \mathfrak{t}_0 \curlyvee \cdots \curlyvee \mathfrak{t}_n)$ as $\mathbf{D}(E_0) \curlyvee \cdots \curlyvee \mathbf{D}(E_n)$.

*Proof.* A stratification of $X$ as in (**5.37**) induces a certain triangular diagram $\mathsf{G}_n$ of the following form, where all maps $i_k$ are inclusions of the closed subspaces $U_k$ of $\mathsf{s}$, and all $j_k$ are inclusions of the pure strata $E_k$: in the notation above we obtain

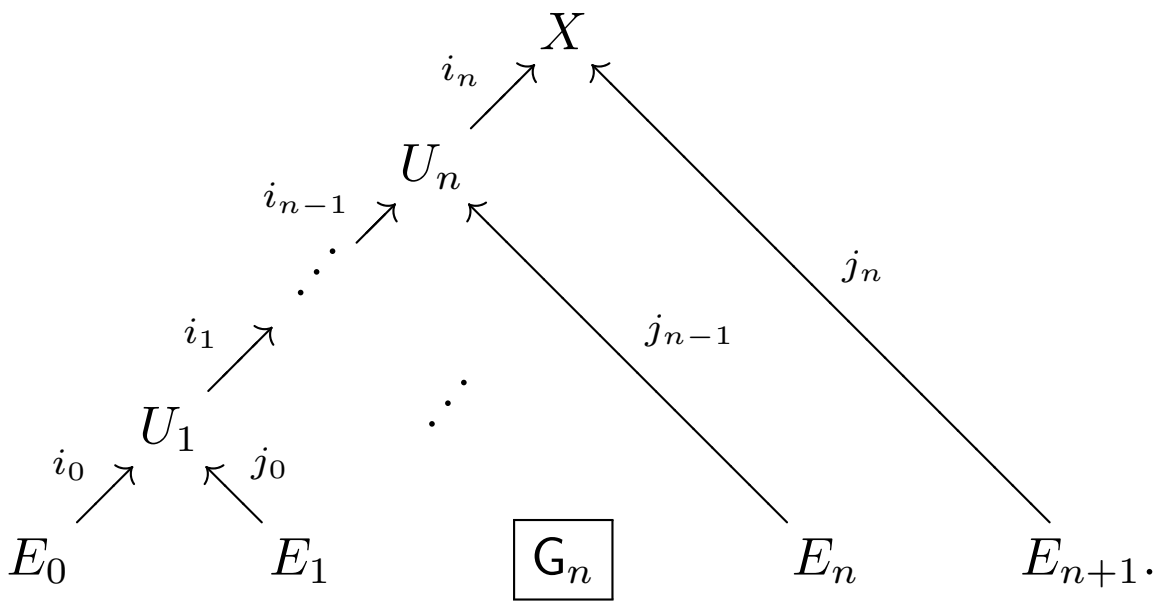

(5.38)

This diagram can clearly be defined inductively starting from $n = 1$ (the diagram of inclusions as in Example **5.2.5**). Given this evident recursive nature, it is sufficient to examine the case $n = 2$ of a stratification $U_0 \subset U_1 \subset X$, depicted as[7]

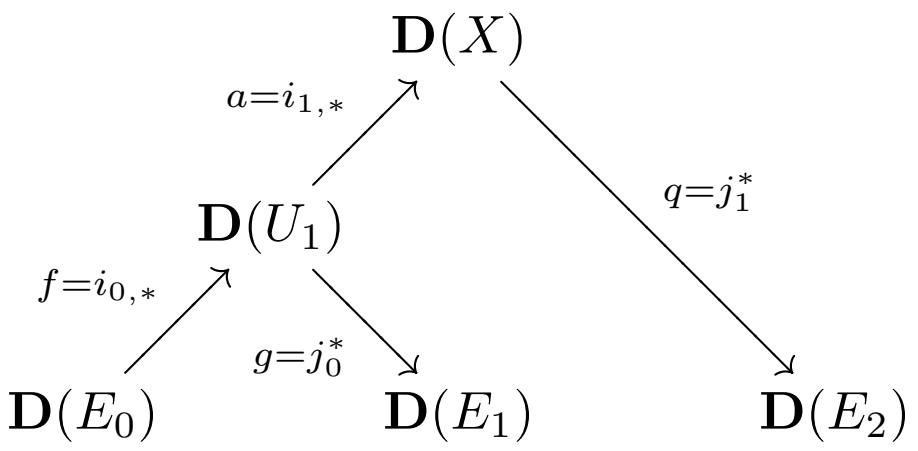

(5.39)

---

[7] Here and for the rest of the section, drawing large diagrams of stable categories, we adopt the following shorthand: every edge $h\colon \mathbf{E} \to \mathbf{F}$ is decorated with an adjoint triple $h_L \dashv h \dashv h_R\colon \mathbf{E} \leftrightarrows \mathbf{F}$.

to notice that the $t$-structure $(\mathfrak{t}_0 \veebar \mathfrak{t}_1) \veebar \mathfrak{t}_2$ obtained by iterated gluing construction is

$$\begin{aligned}
[(\mathbf{D}(E_0) \veebar \mathbf{D}(E_1)) \veebar \mathbf{D}(E_2)]_{\geq 0} &= \left\{ G \in \mathbf{D}(X) \,\middle|\, \begin{array}{l} qG \in \mathbf{D}(E_2)_{\geq 0}, \\ a_L G \in [\mathbf{D}(E_0) \veebar \mathbf{D}(E_1)]_{\geq 0} \end{array} \right\} \\
&= \left\{ G \in \mathbf{D}(X) \,\middle|\, \begin{array}{l} qG \in \mathbf{D}(E_2)_{\geq 0}, \\ g(a_L G) \in \mathbf{D}(E_1)_{\geq 0}, \\ f_L(a_L G) \in \mathbf{D}(E_0)_{\geq 0} \end{array} \right\} \\
(\mathbf{5.2.17}) &= \left\{ G \in \mathbf{D}(X) \mid \{q, g a_L, f_L a_L\}(G) \in \mathbf{D}_{\geq 0} \right\}
\end{aligned}$$

The inductive step simply adds another inclusion (and the obvious maps between derived categories) to these data. □

Remark 5.4.3. In the previous proof, in the case $n = 2$, we could have noticed that two "hidden" recollement data, given by the inclusions

$$(E_1 \hookrightarrow X \smallsetminus U_0, E_2 \hookrightarrow X \smallsetminus U_0) \text{ and } (E_0 \hookrightarrow X, X \smallsetminus U_0 \hookrightarrow X)$$

come into play: the refinement of the inclusions in the diagram above induces an analogous refinement which passes to the derived $\infty$-categories,

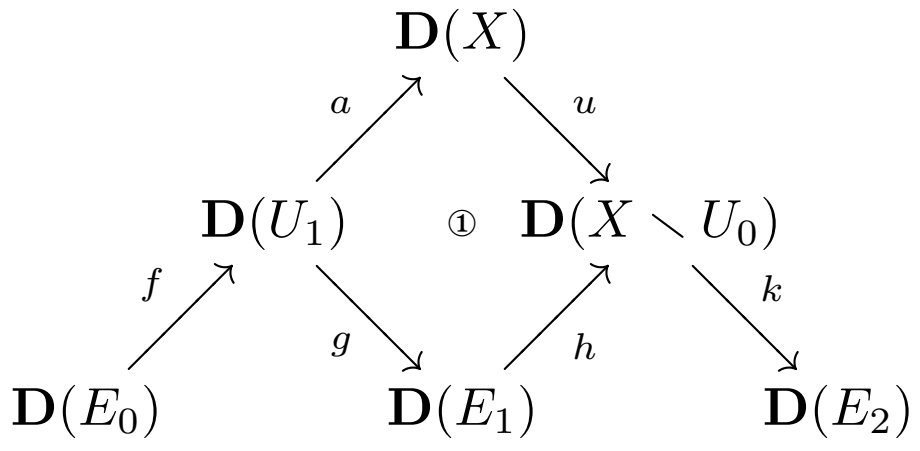

(5.40)

of functors between derived $\infty$-categories on the pure strata. These data induce two additional recollements, $(k, h)$ and $(u, a \circ f)$ which we can use to define a different parenthesization $\mathfrak{t}_0 \veebar (\mathfrak{t}_1 \veebar \mathfrak{t}_2)$.

Remark 5.4.4. When all the recollements data in **5.40** are taken into account, we obtain a graph called the *left-winged* diagram associated with (**5.40**), and defined by taking the left-most adjoint in the string $(-)_L \dashv (-) \dashv (-)_R$, when descending each left "leaf" of the tree represented in diagram (**5.40**). In a completely similar fashion we can define the *right-winged* diagram of (**5.40**). We refer to these diagrams as (**l-5.40**) and (**r-5.40**) respectively.

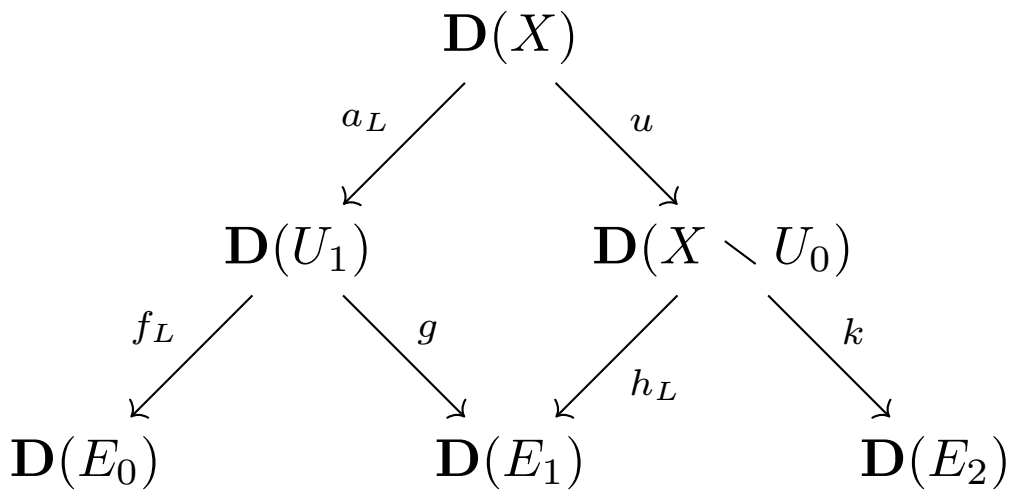

It is now quite natural to speculate about some sort of *comparison* between the two recollements $(\mathfrak{t}_0 \curlyvee \mathfrak{t}_1) \curlyvee \mathfrak{t}_2$ and $\mathfrak{t}_0 \curlyvee (\mathfrak{t}_1 \curlyvee \mathfrak{t}_2)$: in fact we can prove with little effort (once the phenomenon in study has been properly clarified) that the two $t$-structures are equal, since the square

$$\begin{array}{ccc} E_1 & \longrightarrow & X \setminus U_0 \\ \downarrow & & \downarrow \\ U_1 & \longrightarrow & X \end{array} \tag{5.41}$$

is a pullback (in a suitable category of spaces), and so there is a "change of base" morphism $u \circ a \cong h \circ g$ which induces an invertible 2-cell $g \circ a_L \cong h_L \circ u$ filling the square ①in diagram (**5.40**): this is a particular instance of the so-called *Beck-Chevalley condition* for a commutative square, which we now adapt to the $\infty$-categorical setting.

Definition 5.4.5. (Beck-Chevalley condition): Consider the square

$$\begin{array}{ccc} \mathbf{A} & \overset{a_L}{\underset{a_R}{\overset{\longleftarrow}{\underset{\longleftarrow}{\xrightarrow{a}}}}} & \mathbf{B} \\ {\scriptstyle g}\downarrow & & \downarrow{\scriptstyle u} \\ \mathbf{C} & \overset{h_L}{\underset{h_R}{\overset{\longleftarrow}{\underset{\longleftarrow}{\xrightarrow{h}}}}} & \mathbf{D} \end{array} \tag{5.42}$$

in a $(\infty, 2)$-category, filled by an invertible 2-cell $\theta \colon u \circ a \cong h \circ g$ and such that $a_L \dashv a, h_L \dashv h$; then the square is said to satisfy the *left Beck-Chevalley property* (lbc for short) if the canonical 2-cell

$$\hat{\theta} \quad : \quad h_L \circ u \overset{h_L u * \eta}{\Longrightarrow} h_L \circ u \circ a \circ a_L \overset{h_L * \theta * a_L}{\Longrightarrow} h_L \circ h \circ g \circ a_L \overset{\epsilon * g a_L}{\Longrightarrow} g \circ a_L \tag{5.43}$$

is invertible as well. Similarly, when $a \dashv a_R, h \dashv h_R$ we define the 2-cell

$$\tilde{\theta} \quad : \quad g \circ a_R \overset{\eta * g a_R}{\Longrightarrow} h_R \circ h \circ g \circ a_R \overset{h_R * \theta * a_R}{\Longrightarrow} h_R \circ u \circ a \circ a_R \overset{h_R u * \epsilon}{\Longrightarrow} h_R \circ u \tag{5.44}$$

and we say that the square above is *right Beck-Chevalley* (rbc for short) when it is invertible.

In light of this property enjoyed by diagram ①in (**5.40**) it's rather easy to show that the two left classes

$$\big[\big(\mathbf{D}(E_0) \curlyvee \mathbf{D}(E_1)\big) \curlyvee \mathbf{D}(E_2)\big]_{\geq 0} = \big\{G \in \mathbf{D}(X) \mid \{ku, ga_L, f_L a_L\}(G) \in \mathbf{D}_{\geq 0}\big\}$$
$$\big[\mathbf{D}(E_0) \curlyvee \big(\mathbf{D}(E_1) \curlyvee \mathbf{D}(E_2)\big)\big]_{\geq 0} = \big\{G \in \mathbf{D}(X) \mid \{ku, h_L u, f_L a_L\}(G) \in \mathbf{D}_{\geq 0}\big\}$$

coincide up to a canonical isomorphism determined by the Beck-Chevalley 2-cell in ①of diagram (**5.40**).

As a result, both $[(\mathbf{D}(E_0) \uplus \mathbf{D}(E_1)) \uplus \mathbf{D}(E_2)]_{\geq 0}$ and $[\mathbf{D}(E_0) \uplus (\mathbf{D}(E_1) \uplus \mathbf{D}(E_2))]_{\geq 0}$ define the torsion class of the same $t$-structure $(\mathbf{D}^{012}_{\geq 0}, \mathbf{D}^{012}_{<0})$ on $\mathbf{D}(X)$. The previous analysis gives that

Scholium 5.4.6. An object $G \in \mathbf{D}(X)$ lies in $\mathbf{D}^{012}_{\geq 0}$ if and only if $p_0 G \in \mathbf{D}(E_0)_{\geq 0}, p_1 G \in \mathbf{D}(E_1)_{\geq 0}, p_2 G \in \mathbf{D}(E_2)_{\geq 0}$ where $l_i$ is any choice of a functor $\mathbf{D}(X) \to \mathbf{D}(E_i)$ in the left-winged diagram of (**5.40**).

It is now rather easy to repeat the same reasoning with arbitrarily long chains of strata: given a stratified space $(X, \mathsf{s})$ we can induce the diagram

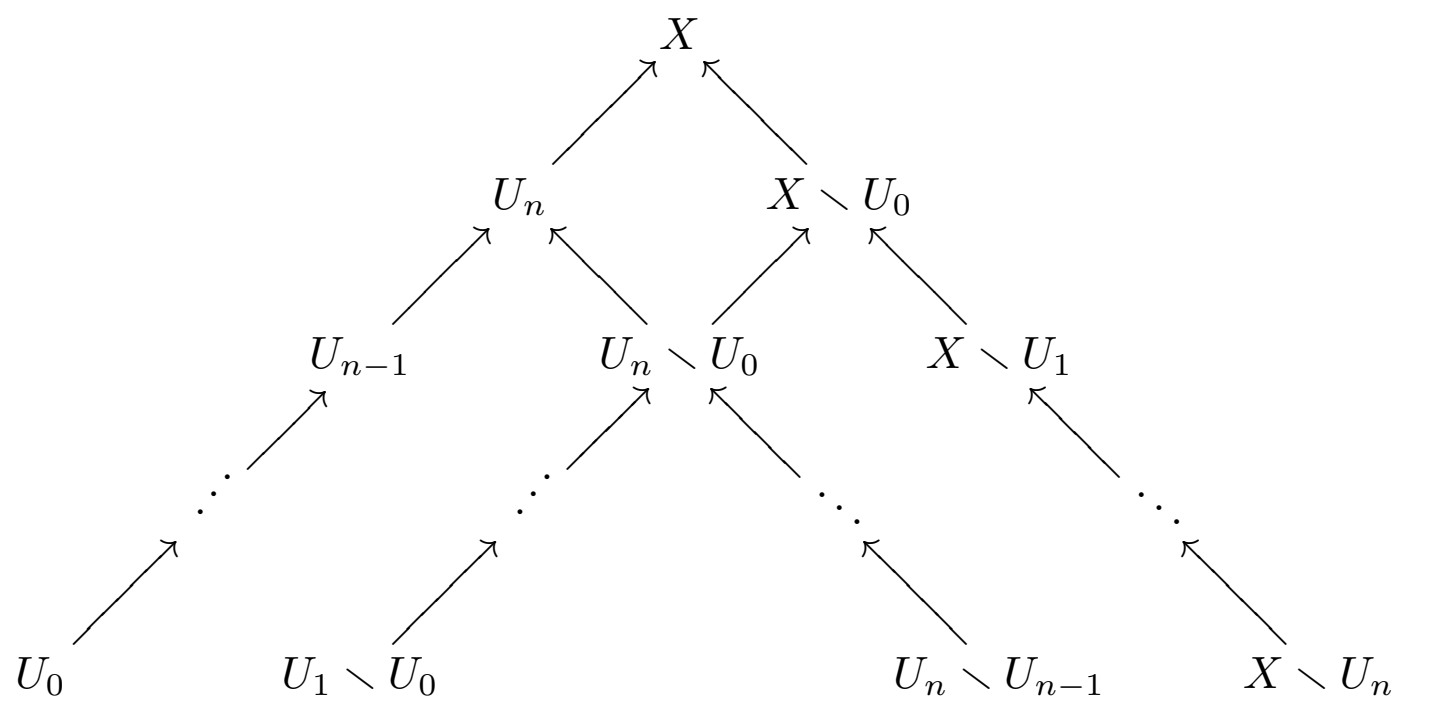

(5.45)

where leaves correspond to pure strata of the stratification of $X$, and every square is a pullback of a proper map along an open embedding, so that the Beck-Chevalley condition is automatically satisfied (inclusions of closed subspaces are proper maps).

Obviously, diagram (**5.45**) induces a diagram $\mathbf{D}$(**5.45**) between the derived categories of the various nodes, and recollement data between some of these nodes; we can again define the left-winged and right-winged version of $\mathbf{D}$(**5.45**), which we will refer as **l-D**(**5.45**) and **r-D**(**5.45**).

Grouping all these considerations we obtain that

(1) There exist "compatible" recollements to give associativity of all the parenthesizations

$$(\mathsf{t}_0 \uplus \cdots \uplus \mathsf{t}_n)_{\mathfrak{P}} = (\mathsf{t}_0 \uplus \cdots \uplus \mathsf{t}_n)_{\mathfrak{Q}} \tag{5.46}$$

for each $\mathfrak{P}, \mathfrak{Q}$ in the set of all possible parenthesizations of $n$ symbols. This is precisely the sense in which, as hinted above, geometric stratifications and recollement data "interact nicely" to give canonical isomorphisms between $(\mathsf{t}_0 \uplus \cdots \uplus \mathsf{t}_n)_{\mathfrak{P}}$ and $(\mathsf{t}_0 \uplus \cdots \uplus \mathsf{t}_n)_{\mathfrak{Q}}$, i.e. a canonical choice for associativity constraints on the $\uplus$ operation.

(2) The following characterization for the class $\big(\mathbf{D}(E_0) \uplus \cdots \uplus \mathbf{D}(E_n)\big)_{\geq 0}$ holds:

$$\big(\mathbf{D}(E_0) \uplus \cdots \uplus \mathbf{D}(E_n)\big)_{\geq 0} = \Big\{ G \mid l_i(G) \in \mathbf{D}(E_i)_{\geq 0},\ \forall i = 0, \dots, n \Big\} \tag{5.47}$$

where $l_i$ is any choice of a functor $\mathbf{D}(X) \to \mathbf{D}(E_i)$ in the left-winged diagram **l-D**(**5.45**).

Similarly, the right class $\big(\mathbf{D}(E_0) \uplus \cdots \uplus \mathbf{D}(E_n)\big)_{<0}$ can be characterized as the class of objects $G$ such that $r_i(G) \in \mathbf{D}(E_i)_{<0}$, where $r_i$ is any choice of a functor $\mathbf{D}(X) \to \mathbf{D}(E_i)$ in the right-winged diagram **r-D**(**5.45**).

### 5.4.2 Abstract associativity of the gluing.

The geometric case studied above gives us enough information to make an ansatz for a general definition, telling us what we have to generalize, and in which way.

In an abstract, stable setting we have the following definition, which also generalizes, in some sense, **5.2.1**.

Let $n \geq 2$ be an integer, and let us denote as $[\![i,j]\!]$ the *interval* between $i, j \in [n]$, i.e., set $\{k \mid i \leq k \leq j\} \subset [n] = \{0, 1, \dots, n\}$ (we implicitly assume $i \leq j$ and we denote $[\![i,i]\!] = \{i\}$ simply as $i$).

Definition 5.4.7. (Urizen compass[(8)]): A *Urizen compass* of length $n$ is an arrangement of stable $\infty$-categories, labeled by intervals $I \subseteq [n]$, and functors in a diagram $\mathsf{G}_n$ of the form

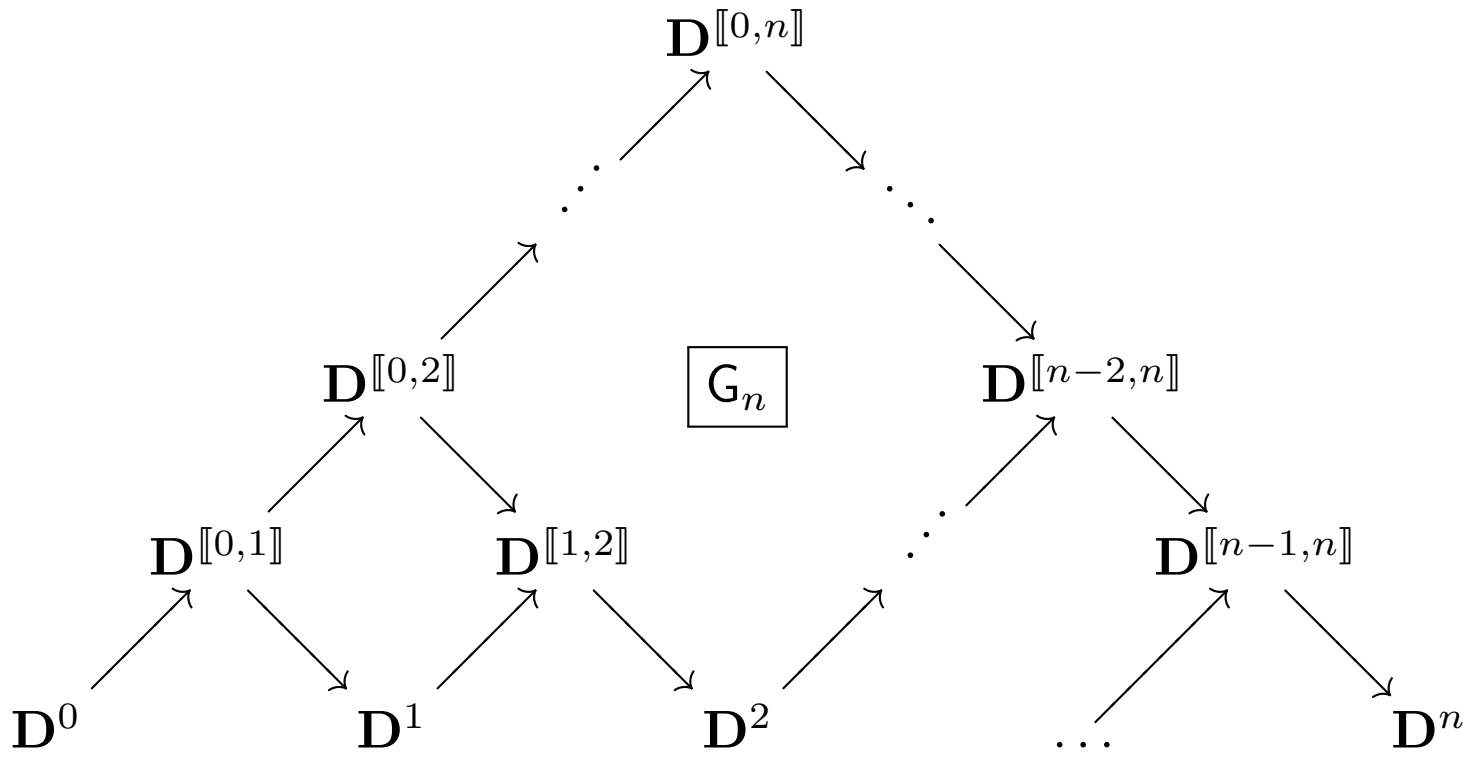

(5.48)

such that the following conditions hold:

- All the triples $\{\mathbf{D}^I, \mathbf{D}^{I \cup J}, \mathbf{D}^J\}$, where $I, J$ are contiguous intervals,[(9)] form different recollements $\mathbf{D}^I \leftrightarrows \mathbf{D}^{I \cup J} \leftrightarrows \mathbf{D}^J$.

(8) In the complicated cosmogony of W. Blake, *Urizen* represents conventional reason and law; it is often represented bearing the same compass of the Great Architect of the Universe postulated by speculative Freemasonry; see for example the painting *The Ancient of Days*, appearing on the frontispiece of the prophetic book "Europe a Prophecy".

(9) Two intervals $I, J \subseteq [n]$ are called *contiguous* if they are disjoint and their union $I \cup J$ is again an interval denoted $I \uplus J$.

- Every square

$$\begin{array}{ccc} \mathbf{D}^{[\![i,j]\!]} & \longrightarrow & \mathbf{D}^{[\![i,j+1]\!]} \\ \downarrow & & \downarrow \\ \mathbf{D}^{[\![i+1,j]\!]} & \longrightarrow & \mathbf{D}^{[\![i+1,j+1]\!]} \end{array} \tag{5.49}$$

  is LBC and RBC in the sense of Definition **5.4.5**.

Note that each row, starting from the base of the diagram, displays all possible intervals of length $k$. We can think of a Urizen compass as a special kind of directed graph (more precisely, a special kind of rooted oriented tree –a *multitree* if we stipulate that each edge shortens a triple of adjunctions); the root of the tree is the category $\mathbf{D}^{[\![0,\dots,n]\!]}$; the leaves are the categories $\{\mathbf{D}^0, \dots, \mathbf{D}^n\}$ (the "generalized pure strata").

Theorem 5.4.8. (The northern hemisphere theorem[(10)]): A Urizen compass of length $n$ induces canonical isomorphisms between the various parenthesizations of $\mathfrak{t}_0 \Yup \cdots \Yup \mathfrak{t}_n$, giving associativity of the gluing operation between $t$-structures.

Rephrasing the above result in a more operative perspective, whenever we have an $n$-tuple $\{(\mathbf{D}^i, \mathfrak{t}_i)\}_{i=0,\dots,n}$ of stable $\infty$-categories with $t$-structure, such that $\{\mathbf{D}^0, \dots, \mathbf{D}^n\}$ are the leaves of a Urizen compass of length $n$, then the gluing operation between $t$-structures gives a unique (up to canonical isomorphism) "glued" $t$-structure on the root $\mathbf{D}^{[\![0,n]\!]}$ of the scheme, resulting as

$$\begin{aligned} \left(\mathbf{D}^0 \Yup \cdots \Yup \mathbf{D}^n\right)_{\geq 0} &= \left\{ X \in \mathbf{D}^{[\![0,n]\!]} \mid l_i(X) \in \mathbf{D}^i_{\geq 0},\ \forall i = 0, \dots, n \right\} \\ \left(\mathbf{D}^0 \Yup \cdots \Yup \mathbf{D}^n\right)_{<0} &= \left\{ X \in \mathbf{D}^{[\![0,n]\!]} \mid r_i(X) \in \mathbf{D}^i_{<0},\ \forall i = 0, \dots, n \right\} \end{aligned} \tag{5.50}$$

where $l_i$ is any choice of a path from the root $\mathbf{D}^{[\![0,n]\!]}$ to the $i^{\text{th}}$ leaf in the left-winged diagram of $\mathsf{G}_n$, and $r_i$ is any choice of a path from the root $\mathbf{D}^{[\![0,n]\!]}$ to the $i^{\text{th}}$ leaf in the right-winged diagram of $\mathsf{G}_n$.

### 5.4.3 Gluing $J$-families.

Our **4.1.20** above shows that the set $\text{TS}(\mathbf{D})$ of $t$-structures on a stable $\infty$-category $\mathbf{D}$ carries a natural action of the ordered group of integers.

(10) In the languages spoken in the northern hemisphere of Tlön, "la célula primordial no es el verbo, sino el adjetivo monosilábico. El sustantivo se forma por acumulación de adjetivos. No se dice luna: se dice *aéreo-claro sobre oscuro-redondo* o *anaranjado-tenue-del cielo* o cualquier otra agregación. [...] Hay objetos compuestos de dos términos, uno de carácter visual y otro auditivo: el color del naciente y el remoto grito de un pájaro. Los hay de muchos: el sol y el agua contra el pecho del nadador, el vago rosa trémulo que se ve con los ojos cerrados, la sensación de quien se deja llevar por un río y también por el sueño. Esos objetos de segundo grado pueden combinarse con otros; el proceso, mediante ciertas abreviaturas, es prácticamente infinito. Hay poemas famosos compuestos de una sola enorme palabra." ([Bor44])

This entails that the most natural notion of a "family" of $t$-structures is a *equivariant* $J$-family of $t$-structures, namely an equivariant map $J \to \text{TS}(\mathbf{D})$ from another $\mathbb{Z}$-poset $J$.

The formalism of equivariant families allows us to unify several constructions in the classical theory of $t$-structures: in particular

> The *semiorthogonal decompositions* of [BO95, Kuz11] are described as precisely those $J$-families $\mathfrak{t}\colon J \to \text{TS}(\mathbf{D})$ taking values on fixed points of the $\mathbb{Z}$-action; these are equivalently characterized as
>
> - the *stable* $t$-structures, where the torsion and torsionfree classes are themselves stable $\infty$-categories;
> - the equivariant $J$-families where $J$ has the trivial action.

And again

> The datum of a single $t$-structure $\mathfrak{t}\colon \{*\} \to \text{TS}(\mathbf{D})$ is equivalent to the datum of a whole $\mathbb{Z}$-*orbit* of $t$-structures, namely an equivariant map $\mathbb{Z} \to \text{TS}(\mathbf{D})$.

In light of these remarks, given a recollement $(i, q)\colon \mathbf{D}^0 \leftrightarrows \mathbf{D} \leftrightarrows \mathbf{D}^1$ it is natural to define the gluing of two $J$-families

$$\text{TS}(\mathbf{D}^0) \xleftarrow{\mathfrak{t}_0} J \xrightarrow{\mathfrak{t}_1} \text{TS}(\mathbf{D}^1) \tag{5.51}$$

to be the $J$-family $\mathfrak{t}_0 \uplus \mathfrak{t}_1 : J \to \text{TS}(\mathbf{D})\colon j \mapsto \mathfrak{t}_0(j) \uplus \mathfrak{t}_1(j)$.

It is now quite natural to ask how does the gluing operation interact with the two situations above: is the gluing of two $J$-families again a $J$-family? As we are going to show, the answer to this question is: yes. Indeed, it's easy to see that the gluing operation is an equivariant map, by recalling that $(\mathcal{E}_0 \uplus \mathcal{E}_1)[1] = \{f \in \hom(\mathbf{D}) \mid f[-1] \in \mathcal{E}_0 \uplus \mathcal{E}_1\}$, and that all of the functors $q, i_L, i_R$ preserves the pullouts (and so commute with the shift). We have

$$\begin{aligned}
(\mathcal{E}_0 \uplus \mathcal{E}_1)[1] &= \{f \in \hom(\mathbf{D}) \mid \{q, i_L\}(f[-1]) \in \mathcal{E}\} \\
&= \{f \in \hom(\mathbf{D}) \mid q(f[-1]) \in \mathcal{E}_1,\ i_L(f[-1]) \in \mathcal{E}_0\} \\
&= \{f \in \hom(\mathbf{D}) \mid q(f)[-1] \in \mathcal{E}_1,\ i_L(f)[-1] \in \mathcal{E}_0\} \\
&= \{f \in \hom(\mathbf{D}) \mid q(f) \in \mathcal{E}_1[1],\ i_L(f) \in \mathcal{E}_0[1]\} \\
&= \mathcal{E}_0[1] \uplus \mathcal{E}_1[1].
\end{aligned}$$

Given this, it is obvious that given two semiorthogonal decompositions $\mathfrak{t}_i\colon J \to \text{TS}(\mathbf{D}_i)$ on $\mathbf{D}^0, \mathbf{D}^1$, the $J$-family $\mathfrak{t}_0 \uplus \mathfrak{t}_1$ is again a semiorthogonal decomposition on $\mathbf{D}$ (the trivial action on $J$ remains the same; it is also possible to prove directly that if $\mathcal{E}_0, \mathcal{E}_1$ are left parts of two exact normal torsion theories $\mathbb{F}_0, \mathbb{F}_1$ on $\mathbf{D}^0, \mathbf{D}^1$, then the gluing $\mathcal{E}_0 \uplus \mathcal{E}_1$ is the left part of the exact normal torsion theory $\mathbb{F}_0 \uplus \mathbb{F}_1$ on $\mathbf{D}$). In some sense at the other side is the gluing of two $\mathbb{Z}$-orbits $\mathfrak{t}_0, \mathfrak{t}_1 \colon \mathbb{Z} \to \textsc{ts}(\mathbf{C})$ on $\mathbf{D}^0$ and $\mathbf{D}^1$. Namely, the glued $t$-structure $\mathfrak{t}_0 \uplus \mathfrak{t}_1$ on $\mathbf{D}$ is the $\mathbb{Z}$-orbit $(\mathfrak{t}_0 \uplus \mathfrak{t}_1)[k] = \mathfrak{t}_0[k] \uplus \mathfrak{t}_1[k]$.

The important point here is that this construction can be framed in the more general context of *perversity data* associated to a recollement, which we now discuss in the attempt to generalize at least part of the classical theory of "perverse sheaves" to the abstract, $\infty$-categorical and torsio-centric setting.

Definition 5.4.9. (Perversity datum): Let $p \colon \{0,1\} \to \mathbb{Z}$ be any function, called a *perversity datum*; suppose that a recollement

$$(i,q) \colon \mathbf{D}^0 \leftrightarrows \mathbf{D} \leftrightarrows \mathbf{D}^1$$

is given, and that $\mathfrak{t}_0, \mathfrak{t}_1$ are $t$-structures on $\mathbf{D}^0, \mathbf{D}^1$ respectively. We define the *($p$-)perverted $t$-structures* on $\mathbf{D}^0, \mathbf{D}^1$ as

$$\begin{aligned} {}^p\mathfrak{t}_0 &= \mathfrak{t}_0[p(0)] = (\mathbf{D}^0_{\geq p(0)}, \mathbf{D}^0_{<p(0)}) \\ {}^p\mathfrak{t}_1 &= \mathfrak{t}_1[p(1)] = (\mathbf{D}^1_{\geq p(1)}, \mathbf{D}^1_{<p(1)}) \end{aligned}$$

Definition 5.4.10. (Perverse objects): Let $p$ be a perversity datum, in the notation above; the *($p$-)glued $t$-structure* is the $t$-structure ${}^p(\mathfrak{t}_0 \uplus \mathfrak{t}_1) = {}^p\mathfrak{t}_0 \uplus {}^p\mathfrak{t}_1$. The heart of the $p$-perverted $t$-structure on $\mathbf{D}$ is called the ($\infty$-)category of *($p$-)perverse objects* of $\mathbf{D}$.

Notice that saying "the category of $p$-perverse objects of $\mathbf{D}$" is an abuse of notation: this category indeed does not depend only on $\mathbf{D}$ and $p$, but on all of the recollement data and on the $t$-structures $\mathfrak{t}_0$ and $\mathfrak{t}_1$. Also notice how for a constant perversity datum $p(0) = p(1) = k$, the $p$-perverted $t$-structure is nothing but the $t$-structure $\mathfrak{t}_0 \uplus \mathfrak{t}_1$ shifted by $k$.

We can extend the former discussion to the gluing of a whole $n$-tuple of $t$-structures, using a Urizen compass:

Remark 5.4.11. In the case of a Urizen compass of dimension $n$ (diagram **5.48**), whose leaves are the categories $\{\mathbf{D}^0, \dots, \mathbf{D}^n\}$, each endowed with a $t$-structure $\mathfrak{t}_i$; a perversity function $p \colon \{0, \dots, n\} \to \mathbb{Z}$ defines a perverted $t$-structure

$$ {}^p(\mathfrak{t}_0 \uplus \cdots \uplus \mathfrak{t}_n) = \mathfrak{t}_0[p(0)] \uplus \mathfrak{t}_1[p(1)] \uplus \cdots \uplus \mathfrak{t}_n[p(n)] \tag{5.52}$$

which is well-defined in any parenthesization thanks to the structure defining the Urizen compass. This result immediately generalizes to the case of a Urizen compass of $J$-families of $t$-structures, $\mathfrak{t}_i \colon J \to \textsc{ts}(\mathbf{D}_i)$, with $i = 0, \dots, n$. Indeed perversity data act on $J$-equivariant families of $t$-structures by

$$ {}^p\mathfrak{t}_i(j) = \mathfrak{t}_i(j)[p(i)] = (\mathbf{D}^i_{\geq j+p(i)}, \mathbf{D}^i_{<j+p(i)}), \tag{5.53}$$

where on the right we have adopted Notation **4.1.21**. This way, a $J$-perversity datum $p\colon \{0,\dots,n\} \to \mathbb{Z}$ induces a $p$-perverted $t$-structure

$$^{p}(\mathfrak{t}_0 \uplus \cdots \uplus \mathfrak{t}_n) = {}^{p}\mathfrak{t}_0 \uplus \cdots \uplus {}^{p}\mathfrak{t}_n : J \to \mathrm{TS}(\mathbf{D}^{[\![0,n]\!]}) \tag{5.54}$$

on $\mathbf{D}^{[\![0,n]\!]}$.

Remark 5.4.12. (Gluing of slicings.): Recall that a *slicing* on a stable $\infty$-category $\mathbf{D}$ consists of a $\mathbb{R}$-family of $t$-structures $\mathfrak{t}\colon \mathbb{R} \to \mathrm{TS}(\mathbf{D})$, where $\mathbb{R}$ is endowed with the usual total order. This means that we are given $t$-structures $\mathfrak{t}_\lambda = (\mathbf{D}_{\geq\lambda}, \mathbf{D}_{<\lambda})$, one for each $\lambda \in \mathbb{R}$, such that $\mathfrak{t}_{\lambda+1} = \mathfrak{t}_\lambda[1]$. Slicings on $\mathbf{D}$ are part of the abstract definition of a $t$-*stability* on a triangulated (or stable) category $\mathbf{D}$, see [Bri07, GKR04].

Grouping together all the above remarks, we obtain that the gluing of two slicings $\mathfrak{t}_i\colon \mathbb{R} \to \mathrm{TS}(\mathbf{D}^i)$ gives a slicing on $\mathbf{D}$ every time $\mathbf{D}^0 \leftrightarrows \mathbf{D} \leftrightarrows \mathbf{D}^1$ is a recollement on $\mathbf{D}$. Moreover, if $p\colon \{0,1\} \to \mathbb{Z}$ is a perversity datum, we have a corresponding notion of *$p$-perverted slicing* on $\mathbf{D}$. More generally one has a notion of $p$-perverted slicing on $\mathbf{D}^{[\![0,n]\!]}$ induced by a pervesity datum $p$ and by a Urizen compass of slicings $\mathsf{G}_n$.

# Chapter 6

# Operations on $t$-structures

In this chapter we collect several examples of operations on the set of $t$-structures on a fixed $\infty$-category $\mathbf{C}$, and functions between classes of $t$-structures on different categories.

From a formal point of view, this amounts to a study of closure properties of the $\infty$-category whose objects are categories with $t$-structure, $(\mathbf{C}, \mathfrak{t}_{\mathbf{C}})$. These objects will be called *$t$-structured categories* and they will be collected in the $\infty$-category $\mathbf{Cat}_{\mathfrak{t}}$. There is an obvious forgetful functor $U\colon \mathbf{Cat}_{\mathfrak{t}} \to \mathbf{Cat}_{\mathrm{st}}$, which has adjoints on both sides; again from a formal point of view, some of the results below amount to a(n easy) series of verification that $U$ has some nice properties.

As it has been observed in Ch. **4** the poset $\mathrm{TS}(\mathbf{C})$ of $t$-structures on $\mathbf{C}$ has a fairly rich structure: it is a partially ordered set, with a canonical action of the group of integers given by the shift functor (see **A.3.11**); it is often the case that nice properties on $\mathbf{C}$ turn $\mathrm{TS}(\mathbf{C})$ into a nicer poset: as we will see in §**6.3**, the class $\mathrm{TS}(\mathbf{Sp})$ of $t$-structures on the $\infty$-category of *spectra* becomes a monoid under the operation described there.

Among the most natural operations on categories there is their product: it is easy to show that given stable $\infty$-categories $\mathbf{C}, \mathbf{D}$ the product $\mathbf{C} \times \mathbf{D}$ is again stable;[1] from a torsio-centric perspective, it is then natural to give the following

Definition 6.0.1. (product $t$-structure): Given stable $\infty$-categories $\mathbf{C}, \mathbf{D}$ the *product* $t$-structure $\mathfrak{t}_{\mathbf{C}} \times \mathfrak{t}_{\mathbf{D}}$ on the product category $\mathbf{C} \times \mathbf{D}$ is defined to be the product (defined in **1.5.7**) of the factorization systems $\mathbb{F}(\mathfrak{t}_{\mathbf{C}}) \times \mathbb{F}(\mathfrak{t}_{\mathbf{D}})$ (the notation is the same of Thm. **3.1.1**).

[1] There are at least two ways to prove this; directly, or appealing to [Lur17, **1.1.4.2**].

## 6.1 Basic constructions.

DEFINITION 6.1.1. (INDUCED $t$-STRUCTURE): Let $\mathbf{B} \subseteq \mathbf{C}$ be a stable sub-$\infty$-category; let $0/\mathcal{E}|_{\mathbf{B}} = 0/\mathcal{E} \cap \mathbf{B}, \mathcal{M}/0|_{\mathbf{B}} = \mathcal{M}/0 \cap \mathbf{B}$ considered as full subcategories; then if the truncation and cotruncation functors $S_{\mathbf{C}}, R_{\mathbf{C}}$ restrict to functors $\mathbf{B} \to 0/\mathcal{E}|_{\mathbf{B}}, \mathcal{M}/0|_{\mathbf{B}}$ the category $\mathbf{B}$ inherits a $t$-structure called the *restricted* or *induced* $t$-structure $\mathfrak{t}|_{\mathbf{B}}$.

The proof is straightforward as restricting the co/truncation functors is a sufficient condition to ensure the existence of $\mathfrak{t}|_{\mathbf{B}}$. Note that in principle this is a weaker condition than having an induced normal torsion theory on $\mathbf{B}$, but that the latter stronger condition is the most natural to expect in concrete situations, as the following example shows.

EXAMPLE 6.1.2. Let $F\colon (\mathbf{C}, \mathfrak{t}_{\mathbf{C}}) \to (\mathbf{D}, \mathfrak{t}_{\mathbf{D}})$ be a $t$-exact functor (Def. **A.3.10**); the *fiber* of $F$ is defined by the pullback square (taken in the $(\infty, 2)$-category of stable $\infty$-categories)

$$\begin{array}{ccc} \mathrm{fib}(F) & \longrightarrow & \mathbf{C} \\ \big\downarrow & \lrcorner & \big\downarrow{\scriptstyle F} \\ \mathbf{0} & \longrightarrow & \mathbf{D}. \end{array} \tag{6.1}$$

In other words, $\mathrm{fib}(F)$ is the full subcategory on all those $X \in \mathbf{C}$ such that $FX \cong 0$. The fiber of $F$ inherits the induced $t$-structure, given that the factor (Def. **1.2.21**) of the $\mathbb{F}(\mathfrak{t}_{\mathbf{C}})$-factorization of a morphism in $\mathrm{fib}(F)$ lies again in $\mathrm{fib}(F)$.

DEFINITION 6.1.3. (CO-INDUCED $t$-STRUCTURE): The *Verdier quotient* $\mathbf{C}/\mathbf{B}$ of $\mathbf{C}$ by a (non necessarily thick) sub-$\infty$-category $\mathbf{B}$ is defined to be the universal functor out of $\mathbf{C}$ sending objects $\mathbf{B}$ to zero (or, equivalently, formally inverting those morphisms whose cofiber is in $\mathbf{B}$).

In the stable setting, the quotient $\mathbf{C}/\mathbf{B}$ is again a stable $\infty$-category.[2]

Under suitable assumptions, the quotient $\mathbf{C}/\mathbf{B}$ acquires a $t$-structure defined by [Lur17, Prop. **1.4.4.11**]: in the same notation as above, suppose $\mathbf{B}, \mathbf{C}$ are presentable and $i\colon \mathbf{B} \to \mathbf{C}$ is a fully faithful inclusion. Suppose $\mathfrak{t} \in \mathrm{TS}(\mathbf{C})$ is a presentable $t$-structure.

The quotient functor $q\colon \mathbf{C} \to \mathbf{C}/\mathbf{B}$ generates the left class of a $t$-structure $(q(\mathbf{C}_{\geq 0}), q(\mathbf{C}_{\geq 0})^{\perp})$ on $\mathbf{C}/\mathbf{B}$ by [Lur17, Prop. **1.4.4.11**].

REMARK 6.1.4. Here we offer a counterexample [Ant] showing that there are cases where this procedure can't induce a $t$-structure on the quotient: suppose that there is an exact and fully faithful inclusion $\mathbf{B} \to \mathbf{C}$ of categories with $t$-structure.

[2] It is the main aim of [] to show that this quotient operation enjoys the universal property of the cofiber $\varinjlim(0 \leftarrow \mathbf{B} \to \mathbf{C})$ in the $(\infty$-)category of stable $\infty$-categories.

From this we deduce an exact functor between hearts $\mathbf{B}^\heartsuit \to \mathbf{C}^\heartsuit$ on hearts; notice that $\mathbf{B}^\heartsuit$ is a so-called *weak Serre* subcategory of $\mathbf{C}^\heartsuit$ (it is closed under extensions, kernels, and cokernels). In order for there to be a $t$-structure on $\mathbf{C}/\mathbf{B}$ such that $q\colon \mathbf{C} \to \mathbf{C}/\mathbf{B}$ is exact, $\mathbf{B}^\heartsuit$ must be Serre, i.e. also closed under subobjects.

A concrete case where this doesn't happen is as follows. Suppose that $R$ is a coherent commutative ring. This means that every finitely generated ideal of $R$ is finitely presented, and it has the consequence that the category of finitely presented $R$-modules is abelian. Let's call this category $\textsc{Coh}(R)$ and view it as a full subcategory of $\mathbf{Mod}(R)$. We can consider $D^b_{\textsc{Coh}(R)}(\mathbf{Mod}(R)) \subseteq D^b(\mathbf{Mod}(R))$, where the category consists of bounded complexes of $R$-modules with homology modules in $\textsc{Coh}(R)$. There are clearly bounded $t$-structures on these. However, if $R$ is not noetherian, then the heart of the first, namely $\textsc{Coh}(R)$, will not be Serre inside the heart of the second, namely $\mathbf{Mod}(R)$.

## 6.2 The poset of $t$-structures.

Studying the order-theoretic properties of the set $\textsc{ts}(\mathbf{C})$ should be a natural step towards the classification of $t$-structures on $\mathbf{C}$.

It is natural, then, to ask whether $\textsc{ts}(\mathbf{C})$ admits finite joins and meets: a natural way to define these operations on $\mathfrak{t}_1 = (\mathbf{C}^{(1)}_{\geq 0}, \mathbf{C}^{(1)}_{<0})$ and $\mathfrak{t}_2 = (\mathbf{C}^{(2)}_{\geq 0}, \mathbf{C}^{(2)}_{<0})$ intersects respectively the aisle and the coaisle, setting

$$\mathfrak{t}_1 \cap \mathfrak{t}_2 = \left(\mathbf{C}^{(1)}_{\geq 0} \cap \mathbf{C}^{(2)}_{\geq 0}, perp.\right) \tag{6.2}$$

$$\mathfrak{t}_1 \cup \mathfrak{t}_2 = \left(perp.', \mathbf{C}^{(1)}_{<0} \cap \mathbf{C}^{(2)}_{<0}\right) \tag{6.3}$$

where *perp.* is a shorthand for $\left(\mathbf{C}^{(1)}_{\geq 0} \cap \mathbf{C}^{(2)}_{\geq 0}\right)^\perp$, and *perp.′* is a shorthand for ${}^\perp\left(\mathbf{C}^{(1)}_{<0} \cap \mathbf{C}^{(2)}_{<0}\right)$. These are called the *naïve join* and *naïve meet* respectively.

It is often the case, however, that the naïve join and meet operations in $\textsc{ts}(\mathbf{C})$ do not coincide with the "abstract" operations on the same set, definable via universal properties: [Bon13, **§1.2**] gives an example where the intersection $\mathbf{C}^{(1)}_{\geq 0} \cap \mathbf{C}^{(2)}_{\geq 0}$, seen as a subcategory of $\mathbf{C}$, can't be coreflective.

Because of this, $\textsc{ts}(\mathbf{C})$ seems to be rather poorly-behaved from the order-theoretic point of view. In fact, it is also possible to show that binary meet and join, when defined, do not distribute over each other.

It is however possible to give conditions on $\mathfrak{t}_1, \mathfrak{t}_2$ ensuring that the expected operations exist and behave nicely: this is the main aim of [Bon13], which we now follow closely: the final aim is to show that $\textsc{ts}(\mathbf{C})$ is a *set with consistencies* (or a *conset* for short), i.e. a partially ordered set where the domain of definition for joins and meets is determined by "consistency conditions" on the arguments, and where "partial distributivity laws" ([Bon13, §2.1]) hold.

Obviously, we state the consistency conditions for two $t$-structures on $\mathbf{C}$ in terms of the corresponding normal torsion theories.

Definition 6.2.1. (Upper and lower consistency): Let $\mathbb{F}_1, \mathbb{F}_2$ be two normal torsion theories on the stable $\infty$-category $\mathbf{C}$ ($\mathbb{F}_i = (\mathcal{E}_i, \mathcal{M}_i)$), and let $(S_i, R_i)$ be the pair coreflection/reflection of $\mathbb{F}_i$. Then $\mathbb{F}_1, \mathbb{F}_2$ are *lower consistent* (resp., *upper consistent*) if $S_1(0/\mathcal{E}_2) \subseteq 0/\mathcal{E}_2$ (resp., $R_2(\mathcal{M}_1/0) \subseteq \mathcal{M}_1/0$).

"Being upper/lower consistent" are symmetric binary relations on the set $\mathrm{TS}(\mathbf{C})$ denoted respectively $\dot{\curlyvee}$ and $\underset{\cdot}{\curlywedge}$. Two normal torsion theories $\mathbb{F}_1, \mathbb{F}_2$ which are both lower and upper consistent are simply called *consistent* and this relation is denoted $\mathbb{F}_1 \between \mathbb{F}_2$.

Proposition 6.2.2. Let $\mathfrak{t}_1, \mathfrak{t}_2$ be a lower consistent pair of $t$-structures on $\mathbf{C}$. Then the naive intersection (**6.2**) is the meet $\mathfrak{t}_1 \wedge \mathfrak{t}_2$. Dually, let $\mathfrak{t}_1, \mathfrak{t}_2$ be upper consistent; then the naive union (**6.3**) is the join $\mathfrak{t}_1 \vee \mathfrak{t}_2$.

Remark 6.2.3. It is easy to show that if $\mathfrak{t}_0 \preceq \mathfrak{t}_1$ then $\mathfrak{t}_0 \between \mathfrak{t}_1$, and the join/meet of the two is $\mathfrak{t}_1/\mathfrak{t}_0$.

Remark 6.2.4. [Bon13, Prop. **5**, **6**] prove that lower or upper consistency is a sufficient condition ensuring that the intersection or union of $t$-structures exists: given an $n$-tuple $\{\mathfrak{t}_1, \dots, \mathfrak{t}_n\}$ of $t$-structures

- if $\mathfrak{t}_i \underset{\cdot}{\curlywedge} \mathfrak{t}_j$ for each $i < j$, then the naïve intersection $\mathfrak{t}_1 \cap \dots \cap \mathfrak{t}_n$ is well-defined (it coincides with the abstract one) and associative;
- if $\mathfrak{t}_i \dot{\curlyvee} \mathfrak{t}_j$ for each $i < j$, then the naïve union $\mathfrak{t}_1 \cup \dots \cup \mathfrak{t}_n$ is well-defined (it coincides with the abstract one) and associative.

Consistency conditions on $t$-structures also ensure that the meet and join distribute over each other:

- if $\mathfrak{t}_1 \dot{\curlyvee} \mathfrak{t}_2$, $\mathfrak{t}_1 \dot{\curlyvee} \mathfrak{t}_3$, and $\mathfrak{t}_2 \underset{\cdot}{\curlywedge} \mathfrak{t}_3$ then $\mathfrak{t}_1 \dot{\curlyvee} (\mathfrak{t}_2 \cap \mathfrak{t}_3)$;
- if $\mathfrak{t}_1 \underset{\cdot}{\curlywedge} \mathfrak{t}_3$, $\mathfrak{t}_2 \underset{\cdot}{\curlywedge} \mathfrak{t}_3$, and $\mathfrak{t}_1 \dot{\curlyvee} \mathfrak{t}_2$ then $(\mathfrak{t}_1 \cup \mathfrak{t}_2) \underset{\cdot}{\curlywedge} \mathfrak{t}_3$.

The structure so determined is a *set with consistencies*, defined in [Bon13, §**2.1**].

## 6.3 Tensor product of $t$-structures.

> All God's children are not beautiful. Most of God's children are, in fact, barely presentable.
>
> F. Leibowitz

Let $\mathbf{C}, \mathbf{D}$ be two presentable $\infty$-categories ([Lur09, Ch. **5**]); then, for each presentable $\infty$-category $\mathbf{A}$ we consider the category $\mathrm{Bil}(\mathbf{C}, \mathbf{D}; \mathbf{A})$ of functors $F \colon \mathbf{C} \times \mathbf{D} \to \mathbf{A}$ such that each restriction $F(-, D)$ and $F(C, -)$ is cocontinuous. These functors are called *bilinear*.

It turns out ([Gro10, Lur17, Lur16]) that the functor $\mathbf{A} \mapsto \mathrm{Bil}(\mathbf{C}, \mathbf{D}; \mathbf{A})$ functor is representable for each pair of categories $\mathbf{C}, \mathbf{D}$ and represented by an object $\mathbf{C} \otimes \mathbf{D}$ called the *tensor product* of the two presentable $\infty$-categories $\mathbf{C}, \mathbf{D}$ (the analogy with the tensor product of vector spaces is evident).

Although we are only interested in the case where the categories involved are stable (and then their tensor product is again stable), this condition plays no rôle in the proof of

LEMMA 6.3.1. The $\infty$-category $\mathbf{C} \otimes \mathbf{D}$ such that

$$\mathrm{Bil}(\mathbf{C}, \mathbf{D}; \mathbf{A}) \cong \mathbf{QCat}(\mathbf{C} \otimes \mathbf{D}, \mathbf{A}) \tag{6.4}$$

is equivalent to $\mathbf{QCat}(\mathbf{C}^{\mathrm{op}}, \mathbf{D})_R$, the sub-$\infty$-category of functors $F \colon \mathbf{C}^{\mathrm{op}} \to \mathbf{D}$ that commute with limits.

*Proof.* It is a long and formal argument based on universal properties. □

PROPOSITION 6.3.2. If $\mathbf{C}, \mathbf{D}$ are *stable* and presentable, the category $\mathbf{C} \otimes \mathbf{D}$ is stable and presentable as well, and we have that $\mathbf{C} \otimes \mathbf{D} \cong \mathbf{QCat}(\mathbf{C}^{\mathrm{op}}, \mathbf{D})_R$.

*Proof.* A slick proof that $\mathbf{QCat}(\mathbf{C}^{\mathrm{op}}, \mathbf{D})_R$ is presentable is in [Gro10, p. 68]; showing that this category is also stable is a matter of unwinding definitions, or follows from our Lemma **4.4.1**. □

A natural question arises: does the tensor operation lifts from $\mathbf{Cat}_{\mathrm{st}}$ to $\mathbf{Cat}_{\mathfrak{t}}$? In other words, given $t$-structures $\mathfrak{t}_{\mathbf{C}}$ and $\mathfrak{t}_{\mathbf{D}}$ on categories $\mathbf{C}, \mathbf{D}$, how can we endow $\mathbf{C} \otimes \mathbf{D}$ with a $t$-structure $\mathfrak{t}_{\mathbf{C}} \otimes \mathfrak{t}_{\mathbf{D}}$ such that the following reasonable properties are satisfied?

- The operation $(\mathfrak{t}_{\mathbf{C}}, \mathfrak{t}_{\mathbf{D}}) \mapsto \mathfrak{t}_{\mathbf{C}} \otimes \mathfrak{t}_{\mathbf{D}}$ is "associative", namely the two $t$-structures $(\mathfrak{t}_0 \otimes \mathfrak{t}_1) \otimes \mathfrak{t}_2$ and $\mathfrak{t}_0 \otimes (\mathfrak{t}_1 \otimes \mathfrak{t}_2)$ correspond to each other via the equivalence $\mathbf{C}_0 \otimes (\mathbf{C}_1 \otimes \mathbf{C}_2) \cong (\mathbf{C}_0 \otimes \mathbf{C}_1) \otimes \mathbf{C}_2$, and "commutative", namely the $t$-structures $\mathfrak{t}_0 \otimes \mathfrak{t}_1$ and $\mathfrak{t}_1 \otimes \mathfrak{t}_0$ correspond to each other via the equivalence $\mathbf{C}_0 \otimes \mathbf{C}_1 \cong \mathbf{C}_1 \otimes \mathbf{C}_0$; moreover, $\otimes \colon \mathrm{TS}(\mathbf{C}_0) \otimes \mathrm{TS}(\mathbf{C}_1) \to \mathrm{TS}(\mathbf{C}_0 \otimes \mathbf{C}_1)$ is compatible with shifts, in the sense that

$$\mathfrak{t}_0[n] \otimes \mathfrak{t}_0[m] = (\mathfrak{t}_0 \otimes \mathfrak{t}_1)[n+m] \tag{6.5}$$

  for each $n, m \in \mathbb{Z}$.
- The canonical $t$-structure $\mathfrak{t}_{\mathbf{Sp}}$ on the category $\mathbf{Sp}$ of spectra is the unit for this monoidal composition, namely $\mathfrak{t}_{\mathbf{Sp}} \otimes \mathfrak{t} \in \mathrm{TS}(\mathbf{Sp} \otimes \mathbf{C})$ and $\mathfrak{t} \otimes \mathfrak{t}_{\mathbf{Sp}} \in \mathrm{TS}(\mathbf{C} \otimes \mathbf{Sp})$ both correspond to $\mathfrak{t} \in \mathrm{TS}(\mathbf{C})$ under the equivalence $\mathbf{Sp} \otimes \mathbf{C} \cong \mathbf{C} \cong \mathbf{C} \otimes \mathbf{Sp}$;
- Tensoring with a fixed $\mathfrak{t} \in \mathrm{TS}(\mathbf{Sp})$ gives an endofunction $\mathfrak{t} \otimes \mathfrak{t}_{\mathbf{C}}$ of $\mathrm{TS}(\mathbf{C})$, when composed with the equivalence $\mathrm{TS}(\mathbf{Sp} \otimes \mathbf{C}) \cong \mathrm{TS}(\mathbf{C})$; more precisely, there is an action $\mathrm{TS}(\mathbf{Sp}) \times \mathrm{TS}(\mathbf{C}) \to \mathrm{TS}(\mathbf{C})$, that becomes a monoid operation when $\mathbf{C} = \mathbf{Sp}$.

It turns out that this problem has a natural reformulation in terms of normal torsion theories, as it is rather easy to use the presentability of the categories involved to invoke the small object argument and produce a factorization system on $\mathbf{C} \otimes \mathbf{D}$. In particular, we rephrase the question in the following form:

> Given two normal torsion theories $\mathbb{F}_{\mathbf{C}}$ and $\mathbb{F}_{\mathbf{D}}$ on two presentable, stable $\infty$-categories $\mathbf{C}, \mathbf{D}$, define a new normal torsion theory $\mathbb{F}_{\mathbf{C}} \otimes \mathbb{F}_{\mathbf{D}}$ on $\mathbf{C} \otimes \mathbf{D}$, having the "nice" properties above.

To solve this problem, consider first of all the explicit formula in Lemma **6.3.1** which defines $\mathbf{C} \otimes \mathbf{D}$: as always, solving this universal problem gives a unique map $\mathbf{C} \times \mathbf{D} \to \mathbf{C} \otimes \mathbf{D}$, which is the unique bilinear functor corresponding to $1_{\mathbf{C}\otimes\mathbf{D}}$; this will be called the *canonical tensor*

$$\otimes \colon \mathbf{C} \times \mathbf{D} \to \mathbf{C} \otimes \mathbf{D} \tag{6.6}$$

This problem can be divided into several steps: first of all we want to induce a factorization system on $\mathbf{C} \otimes \mathbf{D}$ starting from factorization systems on the factors $\mathbf{C}, \mathbf{D}$. Next we want to see that this induced factorization system still has all the good features enjoyed by $\mathbb{F}_{\mathbf{C}}$ and $\mathbb{F}_{\mathbf{D}}$; the canonical tensor, together with its domain and codomain, plays an essential rôle here (notice that, again, stability is not a necessary condition, but only the most important case of interest for the present discussion):

- Consider the product of categories $\mathbf{C} \times \mathbf{D}$ and the product $t$-structure $\mathbb{F}_{\mathbf{C}} \times \mathbb{F}_{\mathbf{D}}$ (Def. **6.0.1**) on this category; by **1.3.9** the left class $\mathcal{E} \times \mathcal{E}'$ of this factorization system uniquely determines $\mathcal{M} \times \mathcal{M}'$ as its right orthogonal, so we will consider only $\mathcal{E} \times \mathcal{E}'$ in the following to simplify the discussion.
- Consider the image $\mathcal{E} \otimes \mathcal{E}'$ of $\mathcal{E} \times \mathcal{E}'$ via the canonical tensor $\otimes$; this, as a class of morphisms in $\mathbf{C} \otimes \mathbf{D}$ has a right orthogonal $(\mathcal{E} \otimes \mathcal{E}')^{\perp}$;
- The small object argument ([Joy] or [DS95]) applied to the class $\mathcal{E} \otimes \mathcal{E}'$ entails that the pair

$$(\mathcal{L}, \mathcal{R}) = \big(\mathsf{s}(\mathcal{E} \otimes \mathcal{E}'), (\mathcal{E} \otimes \mathcal{E}')^{\perp}\big) \tag{6.7}$$

  where $\mathsf{s}(-)$ is the *saturation* operator defined in **1.4.1**, forms a factorization system on $\mathbf{C} \otimes \mathbf{D}$.

## 6.4 Tilting of $t$-structures.

> Quando si vuole uccidere un uomo bisogna colpirlo al cuore, e un Winchester è l'arma più adatta.
>
> Ramón Rojo

Let us first recall (Def. **4.3.7**) that an abelian $\infty$-category is an $\infty$-category with biproducts, kernels and cokernels, and image-factorization which is in addition *homotopically discrete.*

It is not surprising that the language of abelian $\infty$-categories is rich enough to interpret the notion of normal torsion theory:[(3)] to be more precise, we can define a (normal) torsion theory on an abelian $\infty$-category $\mathbf{A}$ in a similar fashion of Def. **2.3.9**, paying attention to the fact that the stable setting endows the definition with several useful autodualities (like **2.3.16**) false in the abelian setting.

Start with the following example: let $\mathbf{C} = \mathbf{D}(\mathbf{A})$ be the derived $\infty$-category of an abelian category $\mathbf{A}$; it is interesting to ask which (factorization functors of) normal torsion theories $\mathfrak{t}_{\mathbf{C}}$ on $\mathbf{D}(\mathbf{A})$ factor through $\mathbf{A} = \mathbf{D}(\mathbf{A})^{\heartsuit} \subset \mathbf{D}(\mathbf{A})$; this means that

(1) we have (mild) co/completeness conditions on $\mathbf{A}$ (i.e. the existence in $\mathbf{A}$ of the co/limits involved in **2.3.17**);
(2) the reflection $R$ and coreflection $S$ of $\mathfrak{t}_{\mathbf{C}}$ factor as follows:

$$\begin{array}{ccc} \mathbf{C}^{\heartsuit} & \xrightarrow{\;S|_{\mathbf{C}^{\heartsuit}}\;} & 0/\mathcal{E} \\ & {\scriptstyle\hat{S}}\searrow \qquad \nearrow & \\ & 0/\mathcal{E} \cap \mathbf{C}^{\heartsuit} & \end{array} \qquad \begin{array}{ccc} \mathbf{C}^{\heartsuit} & \xrightarrow{\;R|_{\mathbf{C}^{\heartsuit}}\;} & \mathcal{M}/0 \\ & {\scriptstyle\hat{R}}\searrow \qquad \nearrow & \\ & \mathcal{M}/0 \cap \mathbf{C}^{\heartsuit} & \end{array} \tag{6.8}$$

This informal definition is needed to cope with a torsio-centric reformulation of *tilting theory.* Our aim here is not to delve into the details of such an intricate and vast topic, but only to skim the surface of it: indeed, at the level of generality we are interested in, *tilting* of a $t$-structure $\mathfrak{t}$ is a device to produce another $t$-structure out of $\mathfrak{t}$ and a (normal) torsion theory on the heart $\mathbf{C}^{\heartsuit,\mathfrak{t}}$; the $t$-structure on $\mathbf{C}$ are acted by (normal) torsion theories on their hearts.

Definition 6.4.1. Let $\mathbf{C}$ be a stable $\infty$-category, $\mathbb{F} = (\mathcal{E}, \mathcal{M})$ a normal torsion theory on $\mathbf{C}$ and $\mathbb{T} = (\mathcal{X}, \mathcal{Y})$ a (normal?) torsion theory on the heart $\mathbf{C}^{\heartsuit,\mathfrak{t}}$. We define the two classes

$$\begin{aligned} \mathcal{E} \mathbin{\text{🝭}} \mathcal{X} &= \{f \in \hom(\mathbf{C}) \mid f \in \mathcal{E}[-1],\ h_{\mathfrak{t}}(f) \in \mathcal{X}\} \\ \mathcal{M} \mathbin{\text{🝭}} \mathcal{Y} &= \{g \in \hom(\mathbf{C}) \mid g \in \mathcal{M}[1],\ h_{\mathfrak{t}}(g) \in \mathcal{Y}\},^{(4)} \end{aligned} \tag{6.9}$$

where $h_{\mathfrak{t}} \colon \mathbf{C} \to \mathbf{C}^{\heartsuit,\mathfrak{t}}$ is the canonical functor of projection to the heart. These two classes define a new normal torsion theory $\mathfrak{t} \mathbin{\text{🝭}} \mathbb{T}$ on $\mathbf{C}$, called the *tilting* of $\mathfrak{t}$ by $\mathbb{T}$.

---

(3) This is the context where historically torsion theories were introduced [Dic66]; in some sense, stable categories are ontologically more primitive since "all" abelian categories arise as hearts of suitable $t$-structures.

(4) The symbol 🝭 (pron. *retort*) recalls the alchemical token for an alembic; here the term hints at the double meaning of the word retort.

REMARK 6.4.2. The idea behind the definition of tilting is to have a way to factor morphisms "until the upper half-plane", and below the horizontal line $Y = 1$;[(5)] specifying a (normal) torsion theory on $\mathbf{C}^{\heartsuit,\mathfrak{t}}$ amounts to specifying a factorization on the objects of the strip $[0,1)$.

PROPOSITION 6.4.3. Let $\mathbb{T}$ be a (normal) torsion theory on $\mathbf{C}^{\heartsuit,\mathfrak{t}}$, and let $\mathbb{S}$ another (normal) torsion theory on $\mathbf{C}^{\heartsuit,\mathfrak{t}\frown\mathbb{T}}$. Now, the tilting operation "behaves like an action", namely

- $(\mathfrak{t} \frown \mathbb{T}) \frown \mathbb{S} = \mathfrak{t} \frown (\mathbb{T} \star \mathbb{S})$, for an operation $\star$ between (normal) torsion theories on the heart;
- $\mathfrak{t} \frown \mathbb{T}_\mathfrak{t} = \mathfrak{t}$, if $\mathbb{T}_\mathfrak{t}$ is the factorization system induced by $\mathfrak{t}$ on its heart.

DEFINITION 6.4.4. (COMPATIBLE $t$-STRUCTURES): Let $\mathbb{F}, \mathbb{F}'$ be two $t$-structures on the stable $\infty$-category $\mathbf{C}$; then $\mathbb{F}'$ is *compatible with* $\mathbb{F}$ (or $\mathbb{F}$-*compatible*) if the $\mathbb{F}'$-factorization of every object $X \in \mathbf{C}^{\heartsuit,\mathfrak{t}}$ belongs again to $\mathbf{C}^{\heartsuit,\mathfrak{t}}$.

In view of the definition of the heart functor $h_\mathfrak{t} \colon \mathbf{C} \to \mathbf{C}^{\heartsuit,\mathfrak{t}}$ as $X \mapsto R_1 S_0 X$, and since an object $A \in \mathbf{C}$ lies in $\mathbf{C}^{\heartsuit,\mathfrak{t}}$ if and only if $h_\mathfrak{t} A \cong A$, we have that $\mathbb{F}'$ is $\mathbb{F}$-compatible if and only if its coreflection/reflection pair $(S', R')$ is such that

$$R_1 S_0 S' = S'; \qquad R_1 S_0 R' = R'. \tag{6.10}$$

REMARK 6.4.5. Let $J$ be a $\mathbb{Z}$-poset and $\mathfrak{t} \colon J \to \text{TS}(\mathbf{C})$ a $J$-slicing on $\mathbf{C}$; let $\bar{\jmath}$ a specified element of $J$ and $\mathfrak{t}_{\bar{\jmath}} = (\mathcal{E}_{\bar{\jmath}}, \mathcal{M}_{\bar{\jmath}})$ its image under $\mathfrak{t}$; then, every $\mathfrak{t}_j$ such that $\mathfrak{t}_{\bar{\jmath}+1} \preceq \mathfrak{t}_j \preceq \mathfrak{t}_{\bar{\jmath}}$ is $\mathfrak{t}_{\bar{\jmath}}$-compatible.

PROPOSITION 6.4.6. Given $\mathbb{F}, \mathbb{F}'$ compatible $t$-structures on $\mathbf{C}$, $\mathbb{F}'$ induces a normal torsion theory on the heart $\mathbf{C}^{\heartsuit,\mathfrak{t}}$, denoted $\mathbb{F}'|_\mathfrak{t}$, and $\mathbb{F} \frown (\mathbb{F}'|_\mathfrak{t}) = \mathbb{F}'$.

PROPOSITION 6.4.7. There is a bijective correspondence between tiltings of $\mathbb{F}$ by (normal) torsion theories on $\mathbf{C}^{\heartsuit,\mathfrak{t}}$ and $\mathbb{F}$-compatible normal torsion theories on $\mathbf{C}$.

The situation is best depicted in the following picture giving the factorization rule; in view of Remark **1.4.8**, this also yields orthogonality of the two classes so determined.

---

[(5)] The $Y$ axis is oriented downwards: see Figure **6.1** below.

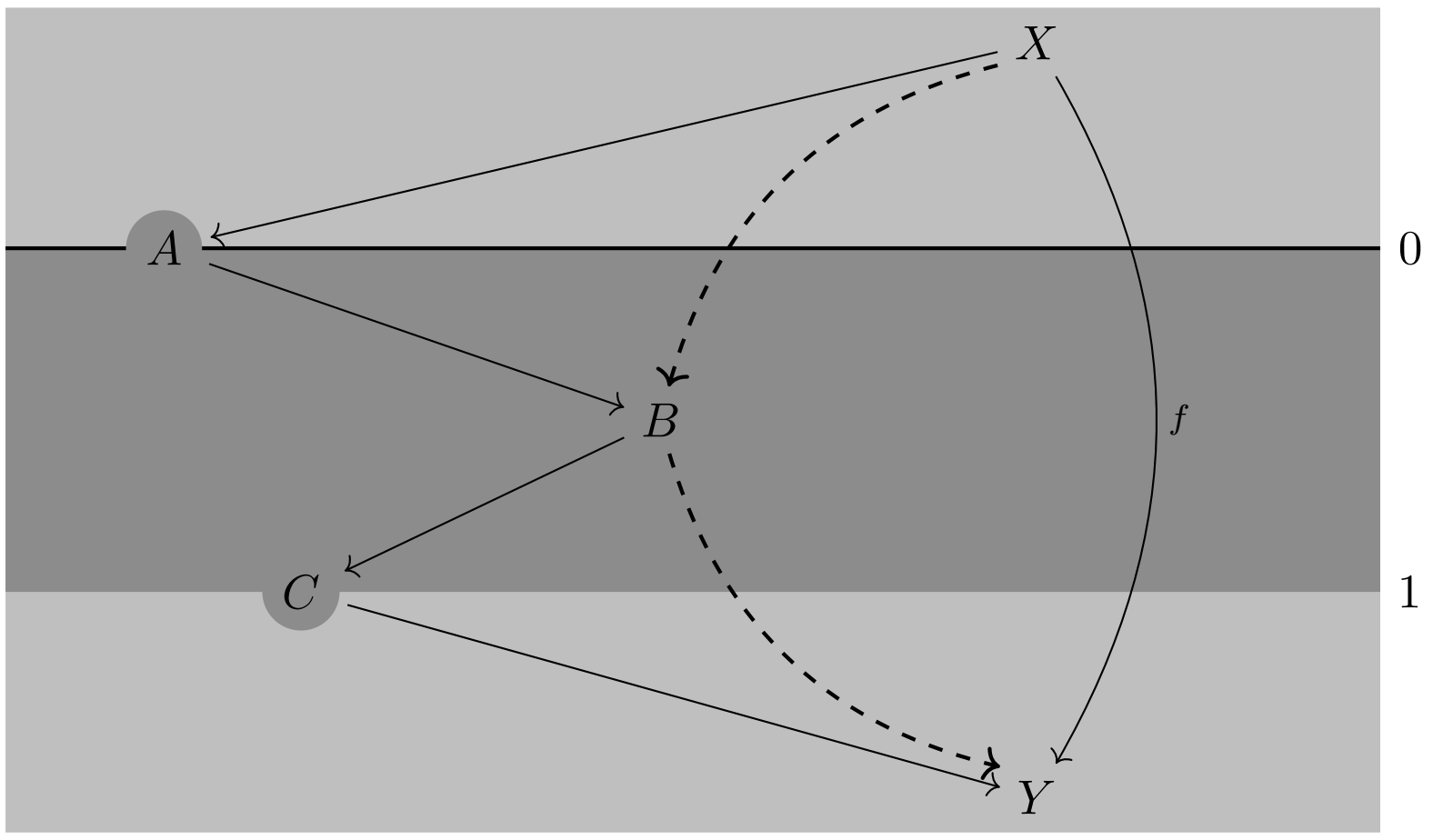

Figure 6.1: Tilting factorization of $f$.

Even if the result can also be obtained from a direct argument, as a consequence of the following general fact about ternary FS:

LEMMA 6.4.8. (TILTING OF FACTORIZATIONS): Let $\mathfrak{t}\colon \mathbb{Z} \to \mathrm{TS}(\mathbf{C})$ be a $\mathbb{Z}$-family of normal torsion theories on a stable $\infty$-category $\mathbf{C}$, having values $\mathfrak{t}_i = (\mathcal{E}_i, \mathcal{M}_i)$ for $i \in \mathbb{Z}$ (here, we will only consider the values $\mathfrak{t}_0, \mathfrak{t}_1 = \mathfrak{t}_0[1]$); let $(\mathcal{L}, \mathcal{R})$ be a torsion theory on the heart $\mathbf{C}^{\heartsuit} = \mathbf{C}_{[0,1)}$ such that $\mathcal{E}_1 \subseteq \mathcal{L} \subseteq \mathcal{E}_0$ (equivalently, $\mathcal{M}_0 \subseteq \mathcal{R} \subseteq \mathcal{M}_1$). Then, the factorization $(e_{\frown}, m_{\frown})$

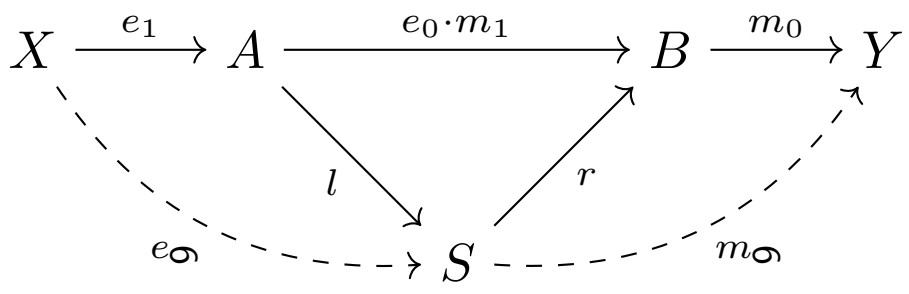

(6.11)

of a morphism $f\colon X \to Y$ in $\mathbf{C}$, obtained from the synergy of the ternary factorization induced by $\mathfrak{t}_0 \preceq \mathfrak{t}_1$ (see **1.5.11**), plus the $(\mathcal{L}, \mathcal{R})$-factorization of its middle part $A \to B \in \mathcal{E}_0 \cap \mathcal{M}_1$, defines a factorization system on $\mathbf{C}$, called the *tilting* of $\mathfrak{t}$ (confused with its 0-value $\mathfrak{t}_0$, in view of Remark **4.1.17**) by $(\mathcal{L}, \mathcal{R})$, and denoted $\mathfrak{t} \frown (\mathcal{L}, \mathcal{R})$.

*Proof.* We have to show that the rule outlined above constitutes a factorization system; the stretagy is to summon [KT93, Thm. **A**] again (see the proof of the "Rosetta stone" **3.1.1**); there is, however, also a direct proof of this fact, appealing the "$\fatsemi$" notation of **1.3.5**: it's easy to see that (in the notation of the statement) $\mathcal{E}_1 \perp \mathcal{R} \fatsemi \mathcal{M}_0$ and $\mathcal{L} \perp \mathcal{R} \fatsemi \mathcal{M}_0$. Now, this allows us to conclude since given a lifting problem

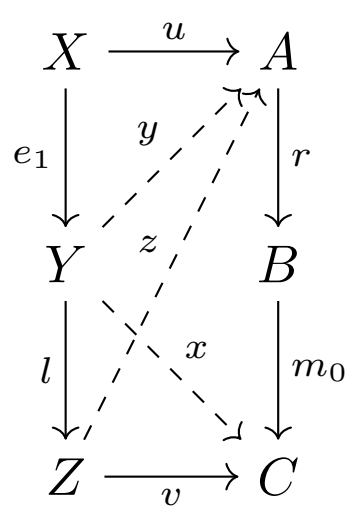

(6.12)

the arrows $x, y, z$ obtained respectively as composition $v \circ l$, and as solutions to suitable lifting problems, give the desired orthogonality. □

# 6.5 Algebras for a monad

In the present section we sketch a general method to induce a $t$-structure on the stable $\infty$-category of algebras for a monad $T$ on a stable $\mathbf{C}$. Apart from some locally-defined new conventions, the same notation as in the rest of the text applies here.

LEMMA 6.5.1. Let $\mathbb{F} = (\mathcal{E}, \mathcal{M})$ be a factorization system on $\mathbf{C}$ and $T$ a monad on $\mathbf{C}$ that preserves the marking $\mathcal{E}$ of $\mathbb{F}$, i.e. such that $T\mathcal{E} \subset \mathcal{E}$; then there is a factorization system $(\mathcal{E}', \mathcal{M}') = U^{\leftarrow}(\mathbb{F})$ on $\mathbf{C}^T$ (the EM category of algebras for the monad $T$) defined by $\mathcal{E}' = U^{\leftarrow}(\mathcal{E}), \mathcal{M}' = U^{\leftarrow}(\mathcal{M})$.

PROPOSITION 6.5.2. If $\mathbf{C}$ is a stable $\infty$-category and $T \colon \mathbf{C} \to \mathbf{C}$ a monad on $\mathbf{C}$ which preserves finite colimits, then the category $\mathbf{C}^T$ of $T$-algebras is again stable.

PROPOSITION 6.5.3. ($t$-STRUCTURE ON $T$-ALGEBRAS): Given a $t$-structure $\mathfrak{t}$ on $\mathbf{C}$, whose normal torsion theory is $\mathbb{F} = (\mathcal{E}, \mathcal{M})$, the procedure above defines a $t$-structure on the category of $T$-algebras $\mathbf{C}^T$, for $T$ a $\mathcal{E}$-preserving finitely cocontinuous monad on $\mathbf{C}$.

*Proof.* To show that $U^{\leftarrow}(\mathbb{F})$ is a normal torsion theory we have to show that

(1) $\mathbb{F}'$ is bireflective, i.e. both $\mathcal{E}', \mathcal{M}'$ are 3-for-2 classes;
(2) $\mathbb{F}'$ is normal i.e. one of the equivalent conditions blabla is satisfied.

The preimage of a 3-for-2 class under any functor is again 3-for-2. This proves the first item. Normality follows from the assumptions in the following form:

the arrow $(KX, k) \to 0$ lies in $\mathcal{E}'$, for each $(X, x) \in \mathbf{C}^T$, if we take $(KX, k)$ to be the fiber

$$\begin{array}{ccc} KX & \longrightarrow & (X, x) \\ \downarrow & & \downarrow \\ 0 & \longrightarrow & (RX, r) \end{array} \tag{6.13}$$

of $(X, x) \to (RX, r)$ (the reflection associated to $\mathbb{F}$). □

Application: let **C** be stable and monoidal; any internal monoid $M$ in **C** induces a monad $- \otimes M$; unders suitable assumptions, the category of $M$-objects (algebras for $- \otimes M$) inherits a $t$-structure.

# Chapter 7

# Stability Conditions

## 7.1 Introduction

The cleanest cut is the one you withhold

Anonymous?

The notion of *Bridgeland stability*[(1)] comes from theoretical Physics, and was proposed by T. Bridgeland in order to better understand a construction in String Theory, the so-called $\Pi$-*stability* of [Dou02, Dou01]; Bridgeland showed that this notion has a natural interpretation in the language of triangulated categories (the idea of identifying objects of the derived category of sheaves on a space with physical D-branes dates back to the work of Moore and Harvey [HM98]).

The main result outlined in [Bri07, Bri09] is that the set of all stability conditions on a given triangulated category $\mathbf{T}$ can be naturally endowed with a topology, induced by a generalized metric. This allows one to define interesting geometric structures out from a triangulated category.

Up to now, a great effort has been put (sometimes, unfortunately, to no avail) into explicitly describing the spaces of stability conditions attached to derived categories of certain algebraic varieties, and to study some of their geometric properties; at the moment of writing, a general theory of these spaces is missing[(2)]

---

(1) There is an unavoidable clash of notation between stability conditions as described here, and the abstract notion of stability for a category. We underline here that this analogy does not exist.

(2) See [Bri09], where the author says:

> there is some yet-to-be discovered construction that will allow one to define interesting geometric structures on these spaces. [...] the agreement between spaces of stability conditions and moduli spaces of conformal field theories is impressive enough to suggest that stability conditions do indeed capture some part of the mathematics of string theory. My own feeling is that at some point

The main aim of the present chapter is to re-enact the classical theory of [Bri07] in the framework of stable $\infty$-categories. In this respect, this is one of the important chapters of the present thesis, as it constitutes one of the main applications of the language initiated by the "Rosetta stone" theorem **3.1.1**. Nevertheless, we only concentrate on a single piece of the rather vast theory of stability conditions on categories, limiting uourselves to showing that given two "close enough" stability functions $Z$ and $W$ and a slicing $J$ compatible with $Z$ then there exists a slicing compatible with $W$, close enough to $J$. A more detailed recovering of other major results about the space of stability conditions will hopefully be discussed in a forthcoming article [FL16a]. Although our proof will closely follow the original argument by Bridgeland, there are a few points where the use of the language developed in the previous chapters of this thesis allow us to give a somehow neater treatment.

Bridgeland's theory exploits some notable[3] properties of increasing families of $t$-structures on a triangulated category $\mathbf{C}$, indexed by the set of real numbers, i.e. monotonic $\mathbb{Z}$-equivariant functions $\mathbb{R} \to \text{TS}(\mathbf{C})$; we paved the way for this definition in our Ch. **4**.

These collections are called ($\mathbb{R}$-)*slicings* in the stable setting; an extremely remarkable result, hidden in Bridgeland's original formulation and made clear by the torsio-entric perspective, is the following:

> simple topological properties of $\mathbb{R}$ (completeness as a metric space, properties of the standard euclidean topology and of the topology of lower convergence generated by the base $\{[a, b) \mid a, b \in \mathbb{Q}\}$...) reflect into categorical properties of slicings

A deeper analysis of this phenomenon occupies §**7.2**.

NOTATION 7.1.1. We make a number of blanket assumptions throughout the chapter: $\mathbf{C}$ is, as always, a stable $\infty$-category, and $\mathfrak{t}$ is a $t$-structure on $\mathbf{C}$; we often demand that $\mathbf{C}$ is cocomplete, and $\mathfrak{t}$ is left, right or two-sided complete. If $J : \mathbb{R} \to \text{TS}(\mathbf{C})$ is a slicing, we define $\mathbf{H}_t = \mathbf{C}_{[t,t+1)}$; the collection $\{\mathbf{H}_t\}$ is called the *heart* of the slicing. The set of slicings $J : \mathbb{R} \to \text{FS}(\mathbf{C})$ is denoted 切$_{\mathbb{R}}(\mathbf{C})$[4]. The real line has to be thought as a time-axis, in such a way that a slicing consists of "a collection of cuttings at prescribed time"; the value of the slicing $J$ at time $\lambda$, $J(\lambda) = (\mathbf{C}_{\geq\lambda}, \mathbf{C}_{<\lambda})$, will often be called the *slice at time* $\lambda$, or the $\lambda$-*slice* of $J$. One has the

---

> in the near future the notion of a stability condition will be subsumed into some more satisfactory framework.

The present chapter is a –clumsy or not, the reader will decide– first step towards this more satisfactory framework.

[3] Peculiar to the standard topological structure of the set of real numbers, but not fully essential: see [GKR04] for an enlightening "formal theory of stability functions", which has been a constant source of inspiration for the present chapter.

[4] The Japanese verb 切る ("kiru", *to cut*) contains the radical 切, the same of *katana*.

inclusion $\mathbf{C}_{\geq\lambda_0} \subseteq \bigcap_{\lambda<\lambda_0} \mathbf{C}_{\geq\lambda}$. A slicing will be called *continuous* at $\lambda_0$ if

$$\bigcap_{\lambda>\lambda_0} \mathbf{C}_{\geq\lambda} = \mathbf{C}_{\geq\lambda_0}. \tag{7.1}$$

It will be called continuous if it is continuous at $\lambda_0$ for every $\lambda_0 \in \mathbb{R}$. We also set

$$\mathbf{C}_{\leq\lambda_0} = \bigcap_{\lambda>\lambda_0} \mathbf{C}_{<\lambda}. \tag{7.2}$$

Notice that if $\lambda_0 < \lambda_1$, then $\mathbf{C}_{\leq\lambda_0} \cap \mathbf{C}_{\geq\lambda_1} = \{0\}$ since, by definition of $\mathbf{C}_{\leq\lambda_0}$, we have $\mathbf{C}_{\leq\lambda_0} \subseteq \mathbf{C}_{<\lambda_1}$. Finally, for $\lambda_0 \leq \lambda_1$ we set

$$\mathbf{C}_{[\lambda_0,\lambda_1]} = \mathbf{C}_{\geq\lambda_0} \cap \mathbf{C}_{\leq\lambda_1}. \tag{7.3}$$

Also, as a shorthand notation, we write $\mathbf{C}_\lambda = \mathbf{C}_{[\lambda,\lambda]}$ for any $\lambda \in \mathbb{R}$.

Definition 7.1.2. Let $\mathbf{C}_0$ be a full subcategory of an stable $\infty$-category $\mathbf{C}$, and let $X$ be an object of $\mathbf{C}_0$. If we have a pullout diagram

$$\begin{array}{ccc} X_s & \longrightarrow & X \\ \downarrow & & \downarrow \\ 0 & \longrightarrow & X_q \end{array} \tag{7.4}$$

in $\mathbf{C}$ with $X_s$ and $X_q$ in $\mathbf{C}_0$, then we say that $X_s$ is a *subobject* of $X$ and that $X_q$ is a quotient of $X$ (relative to $\mathbf{C}_0$).

Definition 7.1.3. A full subcategory $\mathbf{C}_0$ of a stable $\infty$-category $\mathbf{C}$ is called *of finite length* (or simply *finite*) if for each object $A \in \mathbf{C}_0$ there is no infinite ascending chain of subobjects of $A$ (equivalently, there is no infinite descending chain of quotients of $A$).

## 7.2 Slicings

疾く斬るって...剣はそんな小さなものかね

Kamīzumi Nobutsuna

Let $J\colon \mathbb{R} \to \mathrm{FS}(\mathbf{C})$ be a continuous slicing.

Definition 7.2.1. (Suprema and infima): For any object $A$ of $\mathbf{C}$ we set

$$\sup(A) = \inf\{\lambda \in \mathbb{R} : A \in \mathbf{C}_{<\lambda}\};$$
$$\inf(A) = \sup\{\lambda \in \mathbb{R} : A \in \mathbf{C}_{\geq\lambda}\}$$

with the convention $\sup(0) = -\infty$ and $\inf(0) = +\infty$ (if $\mathbf{C}$ is *left/right complete*, [Lur17, Def. **1.2.1.19**], the zero object is the only object whose sup and inf are not finite).

REMARK 7.2.2. It follows directly from the definition that $A \in \mathbf{C}_{\geq\mu}$ implies $\inf(A) \geq \mu$ and $A \in \mathbf{C}_{<\mu}$ implies $\sup(A) \leq \mu$. In particular, if $A \in \mathbf{C}_{[a,b)} = \mathbf{C}_{\geq a} \cap \mathbf{C}_{<b}$ then $a \leq \inf(A)$ and $\sup(A) \leq b$.

DEFINITION 7.2.3. A continuous slicing $J$ will be called *regular* if for any nonzero object $A$ in $\mathbf{C}_{[a,b)}$ one has $\sup(A) < b$.

Unless otherwise stated, all slicings considered in the following will be regular.

LEMMA 7.2.4. If $\inf(A) > \mu$ then $A \in \mathbf{C}_{\geq\mu}$ and if $\sup(A) < \mu$ then $A \in \mathbf{C}_{<\mu}$. In particular, it follows that for any $A \neq 0$ one has $\inf(A) \leq \sup(A)$ and $A \in \mathbf{C}_{[\inf(A),\sup(A)]}$.

*Proof.* If $\inf(A) > \mu$ there exists $\lambda_\mu > \mu$ such that $A \in \mathbf{C}_{\geq\lambda_\mu}$. Since $\lambda_\mu > \mu$, one immediately gets $A \in \mathbf{C}_{\geq\mu}$. The proof for $\sup(A)$ is analogous. It follows from this that $A \in \bigcap_{\mu<\inf(A)} \mathbf{C}_{\geq\mu} = \mathbf{C}_{\geq\inf(A)}$ and $A \in \bigcap_{\mu>\sup(A)} \mathbf{C}_{<\mu} = \mathbf{C}_{\leq\sup(A)}$. If $A \neq 0$ this gives $\mathbf{C}_{\geq\inf(A)} \cap \mathbf{C}_{\leq\sup(A)} \neq \{0\}$ and so $\inf(A) \leq \sup(A)$ and $A \in \mathbf{C}_{[\inf(A),\sup(A)]}$. □

This proves that for a regular slicing the two inequalities $a \leq \inf(A)$ and $\sup(A) < b$ for a nonzero object $A$ in $\mathbf{C}_{[a,b)}$ form a chain:

COROLLARY 7.2.5. Let $J$ be a regular slicing and let $A \in \mathbf{C}_{[a,b)}$ be a nonzero object. Then $a \leq \inf(A) \leq \sup(A) < b$.

An important result links together the contractibility of mapping spaces $\mathbf{C}(X,Y)$ and suitable inequalities between infima and suprema of co/domains of these maps.

LEMMA 7.2.6. If $\inf(X) > \sup(Y)$ then $\mathbf{C}(X,Y) = \{0\}$.

*Proof.* Let $t \in \mathbb{R}$ be such that $\sup(Y) < t < \inf(X)$. Then $X \in \mathbf{C}_{\geq t}$ and $Y \in \mathbf{C}_{<t}$, by Lemma **7.2.4**; the object-orthogonality of classes in the slice at time $t$ allows us to conclude. □

The situation is depicted as follows: there is a "natural direction" in which nonzero morphisms of $(\mathbf{C}, J)$ go: if $\inf(X)$ is greater than $\sup(Y)$, then $Y$ only receives zero morphisms from $X$.

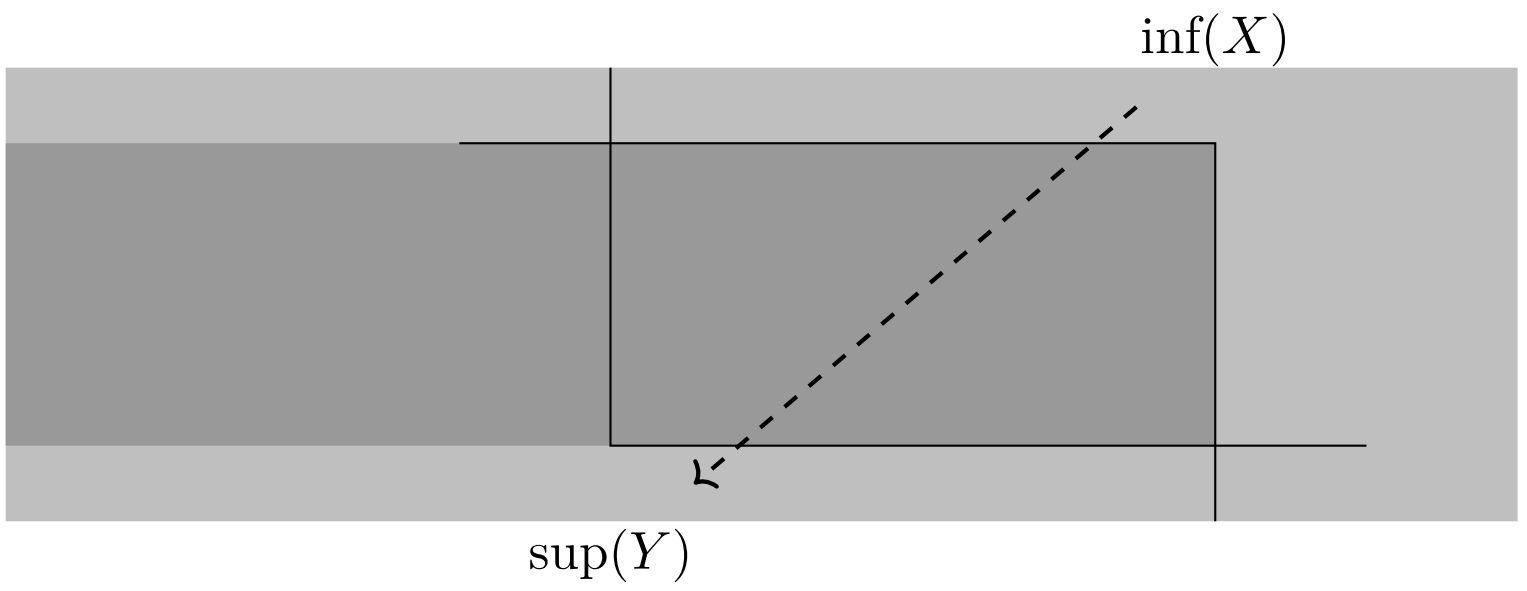

Figure 7.1: Lemma **7.2.6**.

Taking the contrapositive, the above Lemma gives the following.

Lemma 7.2.7. Let $f\colon X \to Y$ be a nonzero morphism in $\mathbf{C}$. Then $\inf(X) \le \sup(Y)$.

Remark 7.2.8. Lemma **7.2.7** provides an additional proof of the fact that for a nonzero object $A$ in $\mathbf{C}$ one has $\inf(A) \le \sup(A)$. Indeed, if $A$ is nonzero, then $\mathrm{id}_A\colon A \to A$ is a nonzero morphism.

Lemma 7.2.9. Let $A$ be an object in $\mathbf{C}$ and let $\mu < \sup(A)$. Then there exists a nonzero morphism $f\colon A_\mu \to A$ with $\inf(A_\mu) \ge \mu$. Similarly, if $\mu > \inf(A)$ then there exists a nonzero morphism $f\colon A \to A_\mu$ with $\sup(A_\mu) \le \mu$.

*Proof.* Since $\mu < \sup(A)$, we have $\mu \notin \{\lambda \in \mathbb{R} : A \in \mathbf{C}_{<\lambda}\}$ and so $A \notin \mathbf{C}_{<\mu} = \mathbf{C}_{\ge\mu}^{\perp}$, and so there exists $A_\mu \in \mathbf{C}_{\ge\mu}$ and a nonzero morphism $f\colon A_\mu \to A$. Since $A \in \mathbf{C}_{\ge\mu}$, we have $\inf(A) \ge \mu$ by Remark **7.2.2**. The proof of the second part of the statement is analogous. □

Definition 7.2.10. (Thin subcategory): The subcategories $\mathbf{C}_{[a,b)}$ of a slicing show an extremely notable behaviour when $[a, b)$ is a "sufficiently small" interval: we call every such $\mathbf{C}_{[a,b)}$ a *thin* subcategory, having in mind [Bri07, Def. **7.2**]; alternatively, we will call $\mathbf{C}_{[a,b)}$ the $[a, b)$-endocardium of the slicing $J$(the reason for this quaint choice of notation is explained in §**7.4**).

Lemma 7.2.11. Let

$$\begin{array}{ccc} A & \longrightarrow & B \\ \downarrow & & \downarrow \\ 0 & \longrightarrow & C \end{array} \tag{7.5}$$

be a fiber sequence in $\mathbf{C}$ with $A, B$ and $C$ in $\mathbf{C}_{[a,b)}$ with $b - a \leq 1$. Then $\sup(A) \leq \sup(B)$ and $\inf(B) \leq \inf(C)$.

*Proof.* We only prove $\sup(A) \leq \sup(B)$ , the other proof being dual. Assume $\sup(A) > \sup(B)$. Then there exists $\mu$ with $\sup(A) > \mu > \sup(B)$ and so by Lemma **7.2.9** there exists a nonzero morphism $f\colon A_\mu \to A$, with $\inf(A_\mu) \geq \mu > \sup(B)$. By Lemma **7.2.6**, the composition $A_\mu \xrightarrow{f} A \to B$ is zero, and so (by the universal property and the 3-for-2 property of pullbacks) the morphism $f$ factors through $C[-1]$. Since the composition $f\colon A_\mu \to C[-1] \to A$ is nonzero, so is the morphism $A_\mu \to C[-1]$, and so by Lemma **7.2.7** $\sup(B) < \inf(A_\mu) \leq \sup(C[-1]) = \sup(C) - 1$. This gives $|\sup(B) - \sup(C)| > 1$. On the other hand, since $B, C \in \mathbf{C}_{[a,b)}$, by Corollary **7.2.5**, both $\sup(B)$ and $\sup(C)$ lie in the interval $[a, b]$ and so $|\sup(B) - \sup(C)| \leq |a - b| \leq 1$. ☐

REMARK 7.2.12. Lemma **7.2.11** in particular implies that, if $b - a \leq 1$ and

$$\begin{array}{ccc} A & \longrightarrow & B \\ \big\downarrow & \lrcorner \qquad \ulcorner & \big\downarrow \\ 0 & \longrightarrow & C \end{array} \tag{7.6}$$

is a fiber sequence with vertices in $\mathbf{C}_{[a,b)}$ and with $B \in \mathbf{C}_{[\tilde{a},\tilde{b})}$, for some $a \leq \tilde{a} < \tilde{b} \leq b$, then $A \in \mathbf{C}_{[a,\tilde{b})}$ and $C \in \mathbf{C}_{[\tilde{a},b)}$.

We also record a direct proof of this fact, independent from **7.2.11**. Since $(\mathbf{C}_{<\tilde{b}}, \mathbf{C}_{\geq\tilde{b}})$ is a $t$-structure on $\mathbf{C}$, we have a pullout diagram

$$\begin{array}{ccccc} 0 & \xrightarrow{e_{\tilde{b}}} & A_{\geq\tilde{b}} & \xrightarrow{m_{\tilde{b}}} & A \\ & & \big\downarrow & & \big\downarrow \\ & & 0 & \longrightarrow & A_{<\tilde{b}} \\ & & & & \big\downarrow{\scriptstyle m_{\tilde{b}}} \\ & & & & 0 \end{array} \tag{7.7}$$

with $A_{\geq\tilde{b}}$ in $\mathbf{C}_{\geq\tilde{b}}$ and $A_{<\tilde{b}}$ in $\mathbf{C}_{<\tilde{b}}$. Since $a \leq \tilde{b}$, we have $\mathcal{E}_{\tilde{b}} \subseteq \mathcal{E}_a$ and so $A \to A_{<\tilde{b}}$ is in $\mathcal{E}_a$. Since $A \in \mathbf{C}_{[a,b)} \subseteq \mathbf{C}_{\geq a}$, the terminal morphism $A \to 0$ is in $\mathcal{E}_a$. So, by the 3-for-2 property of $\mathcal{E}_a$ also $A_{<\tilde{b}} \to 0$ is in $\mathcal{E}_a$, i.e., $A_{<\tilde{b}} \in \mathbf{C}_{\geq a}$. Therefore $A_{<\tilde{b}} \in \mathbf{C}_{[a,\tilde{b})}$; we will write $A_{<\tilde{b}} = A_{[a,\tilde{b})}$ to emphasize this fact. Similarly we have $A_{\geq\tilde{b}} \in \mathbf{C}_{[\tilde{b},b)}$ and we write $A_{\geq\tilde{b}} = A_{[\tilde{b},b)}$.

Consider now the pasting of pullout diagrams

$$\begin{array}{ccccc} & & C[-1] & \longrightarrow & 0 \\ & & \downarrow & & \downarrow \\ A_{[\tilde{b},b)} & \longrightarrow & A & \longrightarrow & B \\ \downarrow & & \downarrow & & \downarrow \\ 0 & \longrightarrow & A_{[a,\tilde{b})} & \longrightarrow & K \\ & & \downarrow & & \downarrow \\ & & 0 & \longrightarrow & C \end{array} \tag{7.8}$$

Since $A_{[\tilde{b},b)}$ is in $\mathbf{C}_{[\tilde{b},b)}$ and $B \in \mathbf{C}_{[\tilde{a},\tilde{b})}$, the morphism $A_{[\tilde{b},b)} \to A \to B$ is the zero morphism and so $A_{[\tilde{b},b)} \to A$ factors through $C[-1]$. But $C[-1] \in \mathbf{C}_{[a-1,b-1)}$ and $b-1 \leq a < \tilde{b}$, so that $\mathbf{C}(A_{[\tilde{b},b)}, C[-1]) = 0$. This implies that $A_{[\tilde{b},b)} \to A$ is the zero morphism, and so $A_{[\tilde{b},b)} = 0$ and $A = A_{[a,\tilde{b})}$. The proof for $C$ is dual.

Remark 7.2.13. If $X \in \mathbf{C}_{[0,1)} \smallsetminus \mathbf{C}_{\{0\}}$, then there is a nonzero morphism $Y_\varepsilon \to X$ for some $\varepsilon > 0$ and $\bar{Y} \in \mathbf{C}_{[\varepsilon,1)}$. Indeed, it is immediate to notice that if $X \in \mathbf{C}_{[0,1)} \smallsetminus \mathbf{C}_0$, then there exists an $0 < \varepsilon < 1$ such that $X \notin \mathbf{C}_{[0,\varepsilon)}$, so $X \notin \mathbf{C}^{\perp}_{[\varepsilon,+\infty)} = \mathbf{C}_{<\varepsilon}$, hence it receives a nonzero morphism $Y_\varepsilon \to X$ from an object $Y_\varepsilon \in \mathbf{C}_{[\varepsilon,+\infty)}$; now $\mathbb{F}_1$-factor this morphism:

$$Y_\varepsilon \xrightarrow{e_1} \bar{Y} \xrightarrow{m_1} X; \tag{7.9}$$

the object $\bar{Y}$ now lies in $\mathbf{C}_{[\varepsilon,1)}$.

### 7.2.1 A topology on slicings

In [Bri07] the author defines a generalized metric (and hence a topology) on the set Stab($\mathbf{D}$) of stability conditions on the triangulated category $\mathbf{D}$; now, we show that this definition corresponds, in the torsio-centric approach, to a generalized metric (and hence a topology) on the set of slicings.

#### 7.2.1.1 A metric on 切$_{\mathbb{R}}(\mathbf{C})$

Definition 7.2.14. Let $I$ and $J$ two slicings on $\mathbf{C}$ and denote by $(\mathbf{C}^I_{<t}, \mathbf{C}^I_{\geq t})$ and $(\mathbf{C}^J_{<t}, \mathbf{C}^J_{\geq t})$ the corresponding families of $t$-structures. We set

$$d(I,J) = \inf\{\varepsilon > 0 \mid \mathbf{C}^I_{<t} \subseteq \mathbf{C}^J_{<t+\varepsilon} \text{ and } \mathbf{C}^I_{\geq t} \subseteq \mathbf{C}^J_{\geq t-\varepsilon} \text{ any for } t \in \mathbb{R}\}. \tag{7.10}$$

This defines a function

$$d\colon 切_{\mathbb{R}}(\mathbf{C}) \times 切_{\mathbb{R}}(\mathbf{C}) \to [0,+\infty] \tag{7.11}$$

REMARK 7.2.15. One can equivalently define $d$ as

$$d(I, J) = \inf\{\varepsilon > 0 \mid \mathbf{C}^J_{\geq t} \subseteq \mathbf{C}^I_{\geq t-\varepsilon} \text{ and } \mathbf{C}^I_{\geq t} \subseteq \mathbf{C}^J_{\geq t-\varepsilon} \text{ any for } t \in \mathbb{R}\}. \tag{7.12}$$

Namely, the condition $\mathbf{C}^I_{<t} \subseteq \mathbf{C}^J_{<t+\varepsilon}$ is equivalent to $\mathbf{C}^{I,\perp}_{\geq t} \subseteq \mathbf{C}^{J,\perp}_{\geq t+\varepsilon}$ and so to $\mathbf{C}^J_{\geq t+\varepsilon} \subseteq \mathbf{C}^I_{\geq t}$. Since this has to hold for every $t$, this is equivalent to $\mathbf{C}^J_{\geq t} \subseteq \mathbf{C}^I_{\geq t-\varepsilon}$.

We split the proof that the function $d$ is a metric on 切$_{\mathbb{R}}(\mathbf{C})$ in lemmas **7.2.16**, **7.2.18**, **7.2.19** below.

LEMMA 7.2.16. The function $d$ is symmetric.

*Proof.* Manifest from the expression for $d$ given in Remark **7.2.15**. □

LEMMA 7.2.17. If $d(I, J)$ is finite, then $\mathbf{C}^I_{\geq t} \subseteq \mathbf{C}^J_{\geq t-d(I,J)}$ and $\mathbf{C}^J_{\geq t} \subseteq \mathbf{C}^I_{\geq t-d(I,J)}$, for any $t \in \mathbb{R}$.

*Proof.* Let $t_0 \in \mathbb{R}$. By definition of $d$ and by Remark **7.2.15**, for any $\varepsilon > 0$ there exists $\delta$ with $d(I, J) \leq \delta < d(I, J) + \varepsilon$ such that $\mathbf{C}^J_{\geq t} \subseteq \mathbf{C}^I_{\geq t-\delta}$ and $\mathbf{C}^I_{\geq t} \subseteq \mathbf{C}^J_{\geq t-\delta}$ for any $t \in \mathbb{R}$. In particular this implies $\mathbf{C}^J_{\geq t_0} \subseteq \mathbf{C}^I_{\geq t_0-d(I,J)-\varepsilon}$ and $\mathbf{C}^I_{\geq t_0} \subseteq \mathbf{C}^J_{\geq t_0-d(I,J)-\varepsilon}$. Since this holds for any $\varepsilon > 0$, we get $\mathbf{C}^J_{\geq t_0} \subseteq \mathbf{C}^I_{\geq t_0-d(I,J)}$ and $\mathbf{C}^I_{\geq t_0} \subseteq \mathbf{C}^J_{\geq t_0-d(I,J)}$. Since $t_0$ was arbitrary, this concludes the proof. □

LEMMA 7.2.18. One has $d(I, J) = 0$ if and only if $I = J$.

*Proof.* Clearly, if $I = J$ then $d(I, J) = 0$. Conversely, assume $d(I, J) = 0$. Then, by Lemma **7.2.17**, we get $\mathbf{C}^J_{\geq t} \subseteq \mathbf{C}^I_{\geq t}$ and $\mathbf{C}^I_{\geq t} \subseteq \mathbf{C}^J_{\geq t}$, i.e., $\mathbf{C}^I_{\geq t} = \mathbf{C}^J_{\geq t}$, for any $t \in \mathbb{R}$. □

LEMMA 7.2.19. The function $d$ satisfies the triangle inequality, i.e. for any three slicings $I, J, K$ one has

$$d(I, K) \leq d(I, J) + d(J, K) \tag{7.13}$$

*Proof.* If either $d(I, J)$ or $d(J, K)$ are infinite then there is nothing to prove. Assume then that both $d(I, J)$ and $d(J, K)$ are finite. By Lemma **7.2.17**, for any $t \in \mathbb{R}$ we have $\mathbf{C}^I_{\geq t} \subseteq \mathbf{C}^J_{\geq t-d(I,J)} \subseteq \mathbf{C}^K_{\geq t-d(I,J)-d(J,K)}$ and $\mathbf{C}^K_{\geq t} \subseteq \mathbf{C}^J_{\geq t-d(J,K)} \subseteq \mathbf{C}^I_{\geq t-d(J,K)-d(I,J)}$. □

DEFINITION 7.2.20. Let $\varepsilon > 0$ be a real number and $J\colon \mathbb{R} \to \mathrm{FS}(\mathbf{C})$ be a slicing; we define

$$U_\varepsilon(J) = \left\{J^\star\colon \mathbb{R} \to \mathrm{FS}(\mathbf{C}) \mid \exists\delta > 0 : (\forall t \in \mathbb{R})\ \mathcal{E}_{t+\varepsilon} \subseteq \mathcal{E}^\star_{t+\delta} \subseteq \mathcal{E}^\star_{t-\delta} \subseteq \mathcal{E}_{t-\varepsilon}\right\} \tag{7.14}$$

where $\mathcal{E}_\lambda$ is the left class of $J(\lambda)$, and similarly $\mathcal{E}^\star_\lambda$ is the left class of $J^\star(\lambda)$ for each $\lambda \in \mathbb{R}$.

Proposition 7.2.21. The set $\mathcal{U} = \{U_\varepsilon(J) \mid \varepsilon > 0,\ J \in 切_{\mathbb{R}}(\mathbf{C})\}$ forms a basis for a topology $\tau_{\mathcal{U}}$ on $切_{\mathbb{R}}(\mathbf{C})$.

*Proof.* As always, we have to prove that

(1) The family $\mathcal{U}$ forms a covering of $切_{\mathbb{R}}(\mathbf{C})$;
(2) every nonempty intersection $U_\alpha(J_1) \cap U_\beta(J_2)$ containing $J$ contains also a basis element containing $J$.

The first point is obvious, as every $J \in 切_{\mathbb{R}}(\mathbf{C})$ lies in $U_\varepsilon(J)$ for $\varepsilon > 0$.

Now, if $J \in U_\alpha(J_1) \cap U_\beta(J_2)$ for $J_1, J_2 \in 切_{\mathbb{R}}(\mathbf{C})$ and $\alpha, \beta > 0$, then we have inequalities

$$\begin{aligned}\mathcal{E}^1_{t-\alpha} \supseteq \mathcal{E}_{t-\delta_1} \supseteq \mathcal{E}_{t+\delta_1} \supseteq \mathcal{E}^1_{t+\alpha}\\ \mathcal{E}^2_{t-\beta} \supseteq \mathcal{E}_{t-\delta_2} \supseteq \mathcal{E}_{t+\delta_2} \supseteq \mathcal{E}^2_{t+\beta}\end{aligned}$$

for suitable $\delta_1, \delta_2 > 0$; it is enough to choose $\gamma > 0$ such that the inequalities

$$\mathcal{E}_{t+\delta_1} \subseteq \mathcal{E}_{t+\gamma} \subseteq \mathcal{E}_{t-\delta_1} \quad \text{and} \quad \mathcal{E}_{t+\delta_2} \subseteq \mathcal{E}_{t+\gamma} \subseteq \mathcal{E}_{t-\delta_2} \tag{7.15}$$

both hold: once this choice has been made, every other $J^\star \in U_\gamma(J)$ satisfies $\mathcal{E}_{t+\gamma} \subseteq \mathcal{E}^\star_{t+\delta} \subseteq \mathcal{E}^\star_{t-\delta} \subseteq \mathcal{E}_{t-\gamma}$ and hence belongs to $U_\alpha(J_1) \cap U_\beta(J_2)$. Any $\gamma < \min\{\delta_1, \delta_2\}$ does the job. □

Proposition 7.2.22. The topology $\tau_{\mathcal{U}}$ on $切_{\mathbb{R}}(\mathbf{C})$ is induced by the metric of Definition **7.2.14**.

*Proof.* A slicing $J^\star$ belongs to the radius $\epsilon/2$ open ball centered at $J$ if and only if

$$\inf\{\delta > 0 \mid \mathbf{C}^{J^\star}_{\geq t} \subseteq \mathbf{C}^{J}_{\geq t-\delta} \text{ and } \mathbf{C}^{J}_{\geq t} \subseteq \mathbf{C}^{J^\star}_{\geq t-\delta} \text{ for any } t \in \mathbb{R}\} < \epsilon/2, \tag{7.16}$$

i.e., if and only if

$$\inf\{\delta > 0 \mid \mathcal{E}^\star_t \subseteq \mathcal{E}_{t-\delta} \text{ and } \mathcal{E}_t \subseteq \mathcal{E}^\star_{t-\delta} \text{ for any } t \in \mathbb{R}\} < \epsilon/2, \tag{7.17}$$

and so if and only if there exists a $\delta > 0$ with $\delta < \varepsilon/2$ such that $\mathcal{E}^\star_t \subseteq \mathcal{E}_{t-\delta}$ and $\mathcal{E}_t \subseteq \mathcal{E}^\star_{t-\delta}$ for all $t \in \mathbb{R}$. Since $t$ is arbitrary, this is equivalent to $\mathcal{E}^\star_{t-\delta} \subseteq \mathcal{E}_{t-2\delta}$ and $\mathcal{E}_{t+2\delta} \subseteq \mathcal{E}^\star_{t+\delta}$ for all $t \in \mathbb{R}$. Since $2\delta < \varepsilon$ we have

$$\mathcal{E}_{t+\varepsilon} \subseteq \mathcal{E}_{t+2\delta} \subseteq \mathcal{E}^\star_{t+\delta} \subseteq \mathcal{E}^\star_{t-\delta} \subseteq \mathcal{E}_{t-2\delta} \subseteq \mathcal{E}_{t-\varepsilon}, \tag{7.18}$$

so $J^\star \in U_\varepsilon(J)$. In other words, $B_{\varepsilon/2}(J) \subseteq U_\varepsilon(J)$. Vice versa, if $J^\star \in U_\varepsilon(J)$ then there exists $\delta_0 > 0$ such that $\mathcal{E}_{t+\varepsilon} \subseteq \mathcal{E}^\star_{t+\delta_0} \subseteq \mathcal{E}^\star_{t-\delta_0} \subseteq \mathcal{E}_{t-\varepsilon}$. Since for every $\delta$ with $0 < \delta < \delta_0$ we have $\mathcal{E}_{t+\delta_0} \subseteq \mathcal{E}^\star_{t+\delta} \subseteq \mathcal{E}^\star_{t-\delta} \subseteq \mathcal{E}_{t-\delta_0}$, we see that for every $0 < \delta < \delta_0$ we have $\mathcal{E}_{t+\varepsilon} \subseteq \mathcal{E}^\star_{t+\delta} \subseteq \mathcal{E}^\star_{t-\delta} \subseteq \mathcal{E}_{t-\varepsilon}$. This gives $\mathcal{E}^\star_t \subseteq \mathcal{E}_{t-\varepsilon+\delta}$ and $\mathcal{E}_t \subseteq \mathcal{E}^\star_{t-\varepsilon+\delta}$ and so $d(J, J^\star) \leq \varepsilon + \delta$ for every $0 < \delta < \delta_0$. In particular, we have $d(J, J^\star) < 2\epsilon$, i.e., $U_\varepsilon(J) \subseteq B_{2\varepsilon}(J)$. □

## 7.3 Stability conditions

> It is essential to think that anything you are doing has to become the occasion for slashing.
>
> M. Musashi

Notation 7.3.1. (Cones and half planes): We adopt the following shorthand to denote certain subsets of the complex plane which we will extensively use from now on:

- Given an interval $[a, b) \subseteq \mathbb{R}$, we denote by $K_{[a,b)} \subseteq \mathbb{C}$ the *cone*

$$K_{[a,b)} = \{z \in \mathbb{C} \mid z = \rho e^{\pi i \theta} \text{ with } \rho \in \mathbb{R}_{\geq 0} \text{ and } \theta \in [a, b)\}. \tag{7.19}$$

  We also adopt all variants like $K_{[a,b]}, K_{(a,b]}, K_{(a,b)}$, all with the obvious meaning.
- Whenever $b - a = \pi$, we call the cone $\mathcal{H}_a := K_{[a,a+\pi)}$ a *half plane* of slope $a$; the half plane of slope 0 will be called the *standard half plane* and denoted $\mathcal{H}$.

Definition 7.3.2. Let $\mathbf{C}$ be a stable $\infty$-category. A *stability condition* on $\mathbf{C}$ is a pair $\sigma = (Z, J)$, where:

sc1) $J \colon \mathbb{R} \to \mathrm{ts}(\mathbf{C})$ is a *slicing* on $\mathbf{C}$;

sc2) $Z \colon \mathbf{C} \to \mathbb{C}$ is a functor[5] which factors through the Grothendieck group of $\mathbf{C}$, i.e., for every fiber sequence

$$\begin{array}{ccc} A & \longrightarrow & B \\ \downarrow & & \downarrow \\ 0 & \longrightarrow & C \end{array} \tag{7.20}$$

in $\mathbf{C}$, one has $Z(A) - Z(B) + Z(C) = 0$ (this property is called *additivity* for $Z$);

sc3) $Z$ is *compatible with the slicing*, i.e. for any $a < b$ in $\mathbb{R}$ one has $Z(\mathbf{C}_{[a,b)}) \subseteq K_{[a,b)}$.

sc4) $Z$ is *nondegenerate on the hearts*, i.e. for any $t \in \mathbb{R}$ one has $Z(X) \neq 0$ for any nonzero object of $\mathbf{H}_t$.

The functor $Z$ is often regarded as a mere function, and called a *stability function* on $\mathbf{C}$.

It follows immediately from sc2 and sc3 that

[5] The set $\mathbb{C}$ is seen as the small groupoid having the set $\mathbb{C}$ as set of objects and exactly one arrow between any two objects; because of this, $Z$ is determined by a function $\mathbf{C}_0 \to \mathbb{C}$.

Remark 7.3.3. The complex number $Z(A)$ only depends on the equivalence class of $A$. Moreover, $Z(A[\pm 1]) = -Z(A)$, so that $Z(A[2]) = Z(A)$ for any object $A$ in $\mathbf{C}$.

Remark 7.3.4. One has $Z(\mathbf{C}_a) \subseteq K_{\{a\}}$, i.e. $Z(X) = \rho(X)e^{i\pi a}$ for any nonzero object $X$ in $\mathbf{C}_a$.

Definition 7.3.5. (Bridgeland subcategory): Let $\mathbf{C}$ be a stable $\infty$-category and $Z\colon \mathbf{C} \to \mathbb{C}$ a stability function on $\mathbf{C}$; a *Bridgeland subcategory* is a full extension closed subcategory $\mathbf{B} \subseteq \mathbf{C}$ such that

(1) The image of the stability function $Z|_{\mathbf{B}}\colon \mathbf{B} \to \mathbb{C}$ is contained in a half-plane $\mathcal{H}_a$;
(2) Each morphism $f\colon X \to Y$ in $\mathbf{B}$ admits a factorization $X \to Z \to Y$ such that there are pullout diagrams

$$\begin{array}{ccc} Z' & \longrightarrow & X \\ \downarrow & & \downarrow \\ 0 & \longrightarrow & Z \end{array} \qquad \begin{array}{ccc} Z & \longrightarrow & Y \\ \downarrow & & \downarrow \\ 0 & \longrightarrow & Z'' \end{array} \tag{7.21}$$

namely, the object $Z$ is at the same time a subobject of $Y$ and a quotient of $X$.

Definition 7.3.6. (Bridgeland cover): A *Bridgeland cover* $\mathfrak{B} = \{\mathbf{B}_\lambda\}_{\lambda\in\Lambda}$ is a family of Bridgeland subcategories $\mathbf{B}_\lambda$ whose extension closure $\langle \bigcup \mathbf{B}_\lambda \rangle$ equals the whole $\mathbf{C}$.

Definition 7.3.7. A stable $\infty$-category $\mathbf{C}$ is said to be *locally finite* with respect to a Bridgeland cover $\mathfrak{B}$ if the subcategories $\mathbf{B}_\lambda$ are of finite length in the sense of Definition **7.1.3**.

We are going to show in the following section that, given a stability condition $(J, Z)$ on the stable $\infty$-category $\mathbf{C}$, the collection $\{\mathbf{C}_{[a,b)}\}_{a\le b\le a+1}$ is a Bridgeland cover of $\mathbf{C}$.

## 7.4 Hearts and endocardia

> If your heart is large enough to envelop your adversaries, you can see right through them and avoid their attacks. And once you envelop them, you will be able to guide them along the path indicated to you by heaven and earth.
>
> M. ‘O Sensei’ Ueshiba

By Thm. **4.3.9**, every category $\mathbf{H}_\lambda = \mathbf{C}_{[\lambda,\lambda+1)}$ is abelian; this subcategory is the $\lambda$-*heart*, i.e. the heart of the slice at time $\lambda$.

We now want to extend the validity of this result to thin (Def. **7.2.10**) subcategories $\mathbf{C}_{[a,b)}$, by showing (in Thm. **7.4.1** below) that all these $\mathbf{C}_{[a,b)}$ are abelian $\infty$-categories in the sense of Def. **4.3.7**.

Theorem 7.4.1. (Abelianity of endocardia): Every $[a,b)$-endocardium is an abelian $\infty$-category; it is, in particular, a category with kernel and cokernel functors, respectively $\ker_{[a,b)}$ and $\mathrm{coker}_{[a,b)}$, and these kernels and cokernels fit into pullback and pushout diagrams

$$\begin{array}{ccc} \ker_{[a,b)} & \longrightarrow & X \\ \downarrow & \lrcorner & \downarrow{\scriptstyle f} \\ 0 & \longrightarrow & Y \end{array} \qquad \begin{array}{ccc} X & \longrightarrow & 0 \\ {\scriptstyle f}\downarrow & \ulcorner & \downarrow \\ Y & \longrightarrow & \mathrm{coker}_{[a,b)} \end{array} \tag{7.22}$$

for each $f\colon X \to Y$. There is, moreover, a canonical isomorphism

$$\mathrm{coker}_{[a,b)}\left[\begin{smallmatrix} \ker_{[a,b)}(f) \\ \downarrow \\ X \end{smallmatrix}\right] \xrightarrow{\ \simeq\ } \ker_{[a,b)}\left[\begin{smallmatrix} Y \\ \downarrow \\ \mathrm{coker}_{[a,b)}(f) \end{smallmatrix}\right] \tag{7.23}$$

whose domain and codomain are called the *coimage* and *image* of $f$ respectively.

*Proof.* We re-draw the diagram constructed in **4.3.9**, and reproduce the argument therein: refer to (**7.25**) below, where $f\colon X \to Z \to Y$, $F = \mathrm{fib}(f), C = \mathrm{cofib}(f)$ are the fiber and cofiber of $f$, and we refer to $\ker_{[a,b)}(f) = S_aF$ and $\mathrm{coker}_{[a,b)}(f) = R_bC$ as the objects emerging from the ternary factorizations

$$\begin{array}{ccccccc} 0 & \longrightarrow & S_bF & \longrightarrow & S_aF & \longrightarrow & F \\ C & \longrightarrow & R_bC & \longrightarrow & R_aC & \longrightarrow & 0 \end{array} \tag{7.24}$$

obtained from the normal torsion theories $\mathbb{F}_a \preceq \mathbb{F}_b$. Now, notice that by definition $\ker_{[a,b)}(f) \in \mathcal{E}_a$ and $\mathrm{coker}_{[a,b)}(f) \in \mathcal{M}_b$, hence the two objects belong to the $[a,b)$-endocardium if and only if $\ker_{[a,b)}(f) \in \mathcal{M}_b$ and $\mathrm{coker}_{[a,b)}(f) \in \mathcal{E}_a$. But this easily follows from the fact that $X, Y \in \mathbf{C}_{[a,b)} = \mathcal{E}_a \cap \mathcal{M}_b$ and from the closure properties of each $\mathcal{E}_\lambda, \mathcal{M}_\lambda$: we are in the following situation,

$$\begin{array}{ccccccccc} \ker_{[a,b)}(f) & \xrightarrow{m_a} & F & \xrightarrow{m_b} & X & \xrightarrow{e_a m_b} & 0 & & \\ {\scriptstyle e_a}\downarrow & & {\scriptstyle e_a}\downarrow & \boxed{\zeta} & {\scriptstyle e_a}\downarrow & & \downarrow{\scriptstyle e_a[1]} & & \\ 0 & \xrightarrow{m_b[-1]} & R_bC[-1] & \xrightarrow{m_b} & Z & \xrightarrow{e_a} & S_aF[1] & \xrightarrow{e_a[1]} & 0 \\ & & {\scriptstyle m_b[-1]}\downarrow & & {\scriptstyle m_b}\downarrow & & \downarrow{\scriptstyle m_b} & & \downarrow{\scriptstyle m_b} \\ & & 0 & \xrightarrow[e_a m_b]{} & Y & \xrightarrow[e_a]{} & C & \xrightarrow[e_a[1]]{} & \mathrm{coker}_{[a,b)}(f) \end{array} \tag{7.25}$$

and hence, by the closure properties of $\mathcal{M}_b$ and $\mathcal{E}_a$, we conclude. □

To conclude the proof we must show that $\ker_{[a,b)}(f)$ and $\mathrm{coker}_{[a,b)}(f)$ indeed have the desired universal properties of kernel and cokernel, namely that in each endocardium the diagrams

$$\begin{array}{ccccc} \ker_{[a,b)}(f) & \longrightarrow & X & \longrightarrow & 0 \\ \downarrow & \lrcorner & \downarrow{\scriptstyle f} & \ulcorner & \downarrow \\ 0 & \longrightarrow & Y & \longrightarrow & \mathrm{coker}_{[a,b)}(f) \end{array} \tag{7.26}$$

are, respectively a pullback and a pushout. This, together with the fact that in every $[a,b)$-endocardium there is a canonical isomorphism $\mathrm{coim}(f) \to \mathrm{im}(f)$, follows from a slight modification of the argument given in Lemma **4.3.16**, **4.3.18** and Prop. **4.3.19** in Ch. **4**.

Again, it remains to show that in every $[a,b)$-endocardium there is a canonical isomorphism $\mathrm{coim}(f) \to \mathrm{im}(f)$; again, this follows by adapting the proof of **4.3.9** in a similar way.

Proposition 7.4.2. Let $J \in \text{切}_{\mathbb{R}}(\mathbf{C})$ be a slicing on a stable $\infty$-category $\mathbf{C}$. Then any $[a,b)$-endocardium of $J$ is a Bridgeland subcategory.

*Proof.* Conditions (1) and (2) of Def. **7.3.5** are rather immediate:

(1) It is obvious since $Z(\mathbf{C}_{[a,b)})$ is contained in an half-plane by property (sc3) of Def. **7.3.2**.

(2) It is a direct consequence of Thm. **7.4.1**, since the object $Z$ in the pullout $(\zeta)$ does the job.

So we are left to prove that each $\mathbf{C}_{[a,b)}$ is an extension closed subcategory; in fact, more is true, since each $[a,b)$-endocardium is also closed under subobjects and quotients.

To see this, consider the diagram

$$\begin{array}{ccccccc} 0 & \xrightarrow{e_a m_b} & A & \xrightarrow{m_b} & B & & \\ & & {\scriptstyle e_a m_b}\downarrow & & \downarrow{\scriptstyle e_a} & & \\ & & 0 & \xrightarrow[e_a m_b]{} & C & \xrightarrow{e_a m_b} & 0 \end{array} \tag{7.27}$$

the assumption that $A, C \in \mathbf{C}_{[a,b)}$, together with stability of $\mathcal{E}_a, \mathcal{M}_b$ under composition and pushout/pullback entails that also $B \in \mathbf{C}_{[a,b)}$. If now $B, C \in \mathbf{C}_{[a,b)}$ in the same diagram, the fact that $A \notin \mathbf{C}_{[a,b)}$ would contradict Lemma **7.2.11**, since $A \in \mathbf{C}_{[a,a+1)} \smallsetminus \mathbf{C}_{[a,b)} = \mathbf{C}_{[b,a+1)}$ entails $b \leq \sup(A) < a+1$, but we must have $\sup(A) \leq \sup(B) < b$.

In a similar way, if $A, B \in \mathbf{C}_{[a,b)}$, if $C \in \mathbf{C}_{[b-1,a)} = \mathbf{C}_{[b-1,b)} \smallsetminus \mathbf{C}_{[a,b)}$, then $\inf(C) < a$, whereas Lemma **7.2.11** entails that $a \leq \inf(B) \leq \inf(C)$. □

REMARK 7.4.3. The $[a,b)$-endocardia clearly cover the whole of $\mathbf{C}$. This is true both with $a, b$ ranging among all pairs of real numbers with $a < b$, and with the constraint $a < b \leq a+1$, or even with a narrower constraint like $a < b \leq a + \varepsilon$ for some $\varepsilon > 0$. In other words $[a,b)$-endocardia are a Bridgeland cover in the sense of definition **7.3.6**.

DEFINITION 7.4.4. A slicing $J$ is called *locally finite* if $\mathbf{C}$ can be covered by finite length endocardia. Equivalently, this means that $\mathbf{C}$ is locally finite with respect to a suitable Bridgeland cover of endocardia in the sense of Definition **7.3.7**.

In what follows, unless otherwise stated, we will assume that the slicings are locally finite.

## 7.5 Deformation of stability conditions

Let $\sigma = (Z, J)$ be a stability condition on $\mathbf{C}$.

NOTATION 7.5.1. We call an object $E \in \mathbf{C} \smallsetminus \{0\}$ *J-thin* (or simply *thin*) if it is contained in some $[a,b)$-endocardium $\mathbf{C}_{[a,b)}$. We denote by $\mathbf{C}^{\asymp}$ the full subcategory of $\mathbf{C}$ on $J$-thin objects.

NOTATION 7.5.2. Let $\| - \|_\sigma$ be the norm on additive functions $\mathbf{C} \to \mathbb{C}$, defined by

$$\|U\|_\sigma = \sup_{t\in\mathbb{R}} \Big( \sup_{E\in\mathbf{C}_t \smallsetminus\{0\}} \left\{ \frac{|U(E)|}{|Z(E)|} \right\} \Big) \tag{7.28}$$

LEMMA 7.5.3. There exists a unique collection of functions

$$\varphi_{[a,b)} \colon \mathbf{C}_{[a,b)} \smallsetminus \{0\} \to [a,b), \tag{7.29}$$

with $(a,b)$ ranging over the set of all pairs of real numbers with $a < b < a+1$, such that

- $Z(E) = \rho(E)\, e^{i\pi\varphi_{[a,b)}(E)}$ for every $E$ in $\mathbf{C}_{[a,b)} \smallsetminus \{0\}$;
- if $a \leq a' < b' \leq b$, then $\varphi_{[a,b)}\big|_{\mathbf{C}_{[a',b')}\smallsetminus\{0\}} = \varphi_{[a',b')}$;
- if $t \in [a,b)$, then $\varphi_{[a,b)}\big|_{\mathbf{C}_t\smallsetminus\{0\}} \equiv t$.
- $\varphi_{[a+1,b+1)}(E[1]) = \varphi_{[a,b)}(E) + 1$.

*Proof.* Since $Z(\mathbf{C}_{[a,b)}) \subseteq K_{[a,b)}$ and $Z(E) \neq 0$ for $E \neq 0$, for every $E$ in $\mathbf{C}_{[a,b)} \smallsetminus \{0\}$ there is a well defined argument of $Z(E)$ with $\arg(Z(E)) \in [a,b)$. Defining $\varphi_{[a,b)}$ as $\arg(Z(E))$ satisfies all the conditions in the statement of the Lemma. Uniqueness is obvious. □

COROLLARY 7.5.4. There exists a unique function

$$\varphi \colon \mathbf{C}^{\asymp} \to \mathbb{R} \tag{7.30}$$

such that $\varphi\big|_{\mathbf{C}_{[a,b)}\smallsetminus\{0\}} = \varphi_{[a,b)}$ for every $a < b < a+1$.

*Proof.* By uniqueness of the argument of $Z(E)$ in a given interval, if $X$ is an object both in $\mathbf{C}_{[a,b)}$ and in $\mathbf{C}_{[a',b')}$ (and so also an object in $\mathbf{C}_{[\max\{a,a'\},\min\{b,b'\})}$) we have

$$\varphi_{[a,b)}(X) = \varphi_{[\max\{a,a'\},\min\{b,b'\})}(X) = \varphi_{[a',b')}(X). \tag{7.31}$$

So the "local" functions $\varphi_{[a,b)}$ glue together into a sigle "global" function $\varphi$. □

Definition 7.5.5. The function $\varphi$ whose existence and uniqueness has been shown in the previous corollary will be called the *$Z$-phase* of $J$-thin objects of $\mathbf{C}$.

Notice that the $Z$-phase $\varphi$ satisfies $\varphi(E[1]) = \varphi(E)+1$ for every nonzero object $E$ of $\mathbf{C}$.

We now come to the main aim of the present section, which is to show that every additive function $W$ in a suitably small neighborhood of a fixed $Z$, is in fact another stability function linked to a "slightly modified" slicing and forming a deformed stability condition $(J^{(W)}, W)$ (Def. **7.3.2**).

Definition/Proposition 7.5.6. (prestability functions preserve cones): Let $W\colon \mathbf{C} \to \mathbb{C}$ an additive function (Def. **7.3.2**) such that $\|Z - W\|_\sigma \le \sin\varepsilon$, where $0 \le \varepsilon \ll \pi/2$: these functions are the *prestability functions* around $Z$. Then, for every $a < b$ in $\mathbb{R}$ one has

$$W(\mathbf{C}_{[a,b)} \smallsetminus \{0\}) \subseteq K_{[a-\varepsilon,b+\varepsilon)} \tag{7.32}$$

*Proof.* This follows immediately from the inequality

$$|W(E) - Z(E)| \le \sin\varepsilon\, |Z(E)| \tag{7.33}$$

for every nonzero $E$ with $E \in \mathbf{C}_t$ with $t \in [a,b)$. □

Remark 7.5.7. The main result outlined in this section can be summarized as "every sufficiently near prestability function around $Z$ is in fact a stability function and it is part of a single stability condition $(J^{(W)}, W)$ around $(J, Z)$" (the "fundamental deformation theorem" **7.5.25**).

As an immediate consequence we get:

Corollary 7.5.8. Let $W\colon \mathbf{C} \to \mathbb{C}$ an additive function (Def. **7.3.2**) such that $\|Z - W\|_\sigma \le \sin\varepsilon$, with $0 \le \varepsilon < 1/4$. Then there exists a unique collection of functions

$$\psi_{[a,b)}\colon \mathbf{C}_{[a,b)} \smallsetminus \{0\} \to \mathbb{R} \tag{7.34}$$

with $(a,b)$ ranging in all pairs of real numbers with $a < b < a+1$, such that

- $W(E) = \rho(E)\, e^{i\pi\psi_{[a,b)}(E)}$ per ogni $E$ in $\mathbf{C}_{[a,b)} \smallsetminus \{0\}$;
- se $a < a' < b' < b$, allora $\psi_{[a,b)}\big|_{\mathbf{C}_{[a',b')}\smallsetminus\{0\}} = \psi_{[a',b')}$;

- per ogni $E$ in $\mathbf{C}_{[a,b)} \smallsetminus \{0\}$ vale $|\psi(E) - \varphi(E)| < \varepsilon$.
- $\psi_{[a+1,b+1)}(E[1]) = \psi_{[a,b)}(E) + 1$.

Moreover, there exists a unique function

$$\psi\colon \mathbf{C}^{\asymp} \to \mathbb{R} \tag{7.35}$$

such that

$$\psi\big|_{\mathbf{C}_{[a,b)}\smallsetminus\{0\}} = \psi_{[a,b)} \tag{7.36}$$

for every $a < b < a+1$.

Notation 7.5.9. Throughout what follows, $\mathbf{C}$ will be a stable $\infty$-category, $J\colon \mathbb{R} \to \mathrm{TS}(\mathbf{C})$ will be a fixed slicing on $\mathbf{C}$, and $\varepsilon$ will be a suitably small real number; the general assumption is that $\varepsilon \ll 1$, but in some special cases we will be able to give sharp estimates. We will also denote by $Z, W$ two stability functions, and by $\phi$ and $\psi$ the phase functions of $Z$ and $W$, respectively. We will also assume that $|\varphi - \psi| < \varepsilon/2$.

Definition 7.5.10. ($\varepsilon$-envelops): Let $[a, b)$ be an interval with $a < b < a + \varepsilon \ll a + 1$, and let $[\alpha, \beta)$ be an thin interval (see Def. **7.5.1**) containing $[a, b)$. We say that $[\alpha, \beta)$ *$\varepsilon$-envelops* $[a, b)$ if $[a - \varepsilon, b + \varepsilon) \subseteq [\alpha, \beta)$; notice that $[\alpha, \beta)$ $\varepsilon$-envelops $[a, b)$ if, and only if, $[a, b) \subseteq [\alpha + \varepsilon, \beta - \varepsilon)$. We denote this situation by "$[a, b) \subseteq_\varepsilon [\alpha, \beta)$".

Remark 7.5.11. Notice that the set $\varepsilon\text{-ENV}(a, b)$ of all intervals $\varepsilon$-enveloping a fixed $[a, b)$ is a directed poset: if $[\alpha, \beta), [\alpha', \beta') \in \varepsilon\text{-ENV}(a, b)$ then $\big[\max\{\alpha, \alpha'\}, \min\{\beta, \beta'\}\big) \in \varepsilon\text{-ENV}(a, b)$ is contained in both intervals.

Definition 7.5.12. ($W$-semistable objects in a thin interval): Under the assumptions of Notation **7.5.9** we say that a nonzero object $E \in \mathbf{C}$ is *$W$-semistable in* $\mathbf{C}_{[a,b)}$ if

- $E \in \mathbf{C}_{[a,b)}$;
- For each nontrivial fiber sequence $A \to E \to B$ in $\mathbf{C}_{[\alpha,\beta)}$, where $[\alpha, \beta) \in \mathrm{ENV}(a, b)$, we have

$$\psi(A) \le \psi(E). \tag{7.37}$$

Definition 7.5.13. If an object $E$ in $\mathbf{C}_{[a,b)}$ is *not* $W$-semistable, then there exists a nontrivial fiber sequence $A \to E \to B$ in $\mathbf{C}_{[\alpha,\beta)}$ with $[\alpha, \beta) \in \mathrm{ENV}(a, b)$ and $\psi(A) > \psi(E)$. Such a fiber sequence will be called a *destabilizing sequence* for $E$ on $[\alpha, \beta)$.

Proposition 7.5.14. Assume $\mathbf{C}_{[\alpha,\beta)}$ is of finite length. Then an object $E$ in $\mathbf{C}_{[a,b)}$ which is not $W$-semistable has a $W$-semistable quotient in $\mathbf{C}_{[\alpha,\beta)}$. Moreover we can choose the $W$-phase of the semistable quotient to be minimal.

*Proof.* If $E$ is not $W$-semistable then there is a quotient $E_1$ of $E$ with $\psi(E_1) < \psi(E)$. This inequality in particular implies that the quotient map is nontrivial. If $E_1$ is $W$-semistable then we are done. Otherwise we have a quotient $E_2$ of $E_1$ with $\psi(E_2) < \psi(E_1)$. Clearly $E_2$ is also a quotient of $E$, so if $E_2$ is $W$-semistable we are done. Proceeding this way, we either end up with a $W$-semistable quotient or we build an infinite sequence of nontrivial quotients. But since the endocardium we are working in is of finite length there can not be infinite sequences of nontrivial quotients. To see that the $W$-phase of the $W$-semistable quotient can be chosen to be minimal, assume that $E \to E_\lambda$ and $E \to E_\mu$ are two semistable quotients of $E$, with $\psi(E_\lambda) = \lambda < \mu = \psi(E_\mu)$. Let $F_\lambda$ be the fiber of $E \to E_\lambda$. Since $\psi(F_\lambda) \geq \psi(E) > \psi(E_\lambda)$ we have $\psi(F_\lambda) > \psi(E_\mu)$. Since $E_\mu$ is $W$-semistable this implies that there are no nontrivial morphisms from $F_\lambda$ to $E_\mu$, so that the diagram

$$\begin{array}{ccc} F_\lambda & \longrightarrow & E \\ \downarrow & & \downarrow \\ 0 & \longrightarrow & E_\mu \end{array} \tag{7.38}$$

commutes. But then, by the universal property of pullouts it factors as

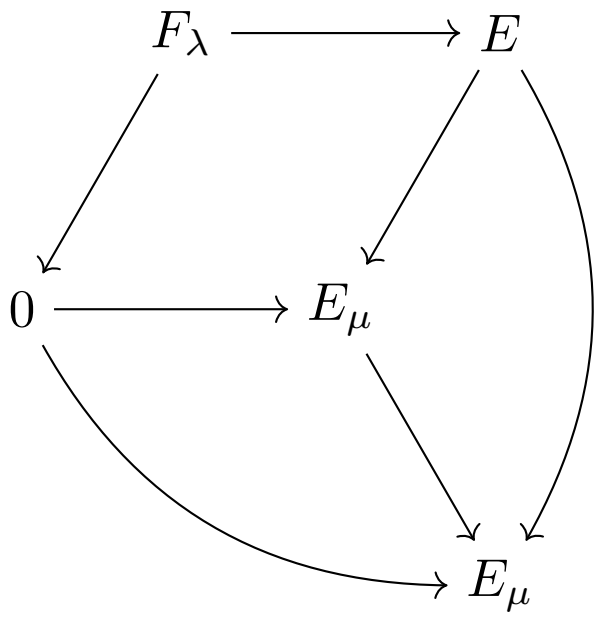

(7.39)

So we see that if there were not a minimal phase $W$-semistable quotient we could build an infinite sequence of nontrivial quotients. $\square$

Proposition 7.5.15. If $[a, b)$ is a sufficiently small interval, then every object in $\mathbf{C}_{[a,b)}$ has a Postnikov tower whose weaves are $W$-semistable objects with decreasing $W$-phases.

*Proof.* If $E$ is $W$-semistable there is nothing to prove. If $E$ is not, consider the fiber sequence

$$\begin{array}{ccc} F_\lambda & \longrightarrow & E \\ \downarrow & & \downarrow \\ 0 & \longrightarrow & E_\lambda \end{array} \tag{7.40}$$

where $\lambda$ is minimal. Then repeat the reasoning on $F_\lambda$. If $F_\lambda$ is $W$-semistable we are done, since $\psi(F_\lambda) > \psi(E_\lambda)$. Indeed, if we set $A_1 = F_\lambda$ and $A_2 = E_\lambda$ we see we have a Postnikov tower with $W$-semistable weaves $A_1, A_2$,

$$\begin{array}{ccccc} 0 & & & & \\ \downarrow & & & & \\ A_1 & \longrightarrow & 0 & & \\ \downarrow & & \downarrow & & \\ E & \longrightarrow & A_2 & \longrightarrow & 0 \end{array} \tag{7.41}$$

with $\psi(A_2) < \psi(A_1)$. If $F_\lambda$ is not $W$-semistable, then we can consider the fiber sequence

$$\begin{array}{ccc} F_{\lambda'} & \longrightarrow & F_\lambda \\ \downarrow & & \downarrow \\ 0 & \longrightarrow & E_{\lambda'} \end{array} \tag{7.42}$$

where $\lambda'$ is minimal. The composite of pullout diagrams

$$\begin{array}{ccccc} F_{\lambda'} & \longrightarrow & F_\lambda & \longrightarrow & E \\ \downarrow & & \downarrow & & \downarrow \\ 0 & \longrightarrow & E_{\lambda'} & \longrightarrow & K \\ & & \downarrow & & \downarrow \\ & & 0 & \longrightarrow & E_\lambda \end{array} \tag{7.43}$$

shows that $\lambda' > \lambda$. Indeed, $\psi(K) > \lambda$ by the minimality assumption on $\lambda$ and so $\lambda' > \psi(K)$; this gives in particular $\lambda' > \lambda$. Similarly, $\psi(F_{\lambda'}) > \lambda'$. So if $F_{\lambda'}$ is $W$-semistable we are done: write $E_1 = F_\lambda$, $A_1 = F_{\lambda'}$, $A_2 = E_{\lambda'}$ and $A_3 = E_\lambda$ to see that we have a Postnikov tower with $W$-semistable weaves $A_1, A_2, A_3$

$$
\begin{array}{ccccccc}
0 & & & & & & \\
\downarrow & & & & & & \\
A_1 & \longrightarrow & 0 & & & & \\
\downarrow & & \downarrow & & & & \\
E_1 & \longrightarrow & A_2 & \longrightarrow & 0 & & \\
\downarrow & & \downarrow & & \downarrow & & \\
E & \longrightarrow & K & \longrightarrow & A_3 & \longrightarrow & 0
\end{array}
\tag{7.44}
$$

with $\psi(A_3) < \psi(A_2) < \psi(A_1)$. If $F_{\lambda'}$ is not $W$-semistable, we iterate the process. This will eventually end due to the finite length assumption on $\mathbf{C}_{[a,b)}$. $\square$

Notation 7.5.16. If $[\alpha,\beta) \in \mathrm{ENV}(a,b)$, we denote by $\mathbf{C}^{\measuredangle(W,J)}_{[a,b)\subseteq_\varepsilon[\alpha,\beta)}$ the full subcategory of $\mathbf{C}_{[a,b)}$ on $W$-semistable objects in $\mathbf{C}_{[a,b)}$.

Proposition 7.5.17. If $[\alpha',\beta')$ is another interval $\varepsilon$-enveloping $[a,b)$, one has

$$
\mathbf{C}^{\measuredangle(W,J)}_{[a,b)\subseteq_\varepsilon[\alpha,\beta)} = \mathbf{C}^{\measuredangle(W,J)}_{[a,b)\subseteq_\varepsilon[\alpha',\beta')}. \tag{7.45}
$$

*Proof.* Without loss of generality we can assume that $[\alpha',\beta') \subseteq [\alpha,\beta)$, and also that $\alpha = \alpha'$ or $\beta = \beta'$ (this follows directly from Remark **7.5.11**). In other words we want to prove that

$$
\begin{aligned}
\mathbf{C}^{\measuredangle(W,J)}_{[a,b)\subseteq_\varepsilon[\alpha,\beta)} &= \mathbf{C}^{\measuredangle(W,J)}_{[a,b)\subseteq_\varepsilon[\alpha,\beta')} \\
\mathbf{C}^{\measuredangle(W,J)}_{[a,b)\subseteq_\varepsilon[\alpha,\beta)} &= \mathbf{C}^{\measuredangle(W,J)}_{[a,b)\subseteq_\varepsilon[\alpha',\beta)}.
\end{aligned}
$$

Suppose $\alpha = \alpha' < \beta' \leq \beta$: a similar argument proves the result when $\alpha \leq \alpha' < \beta = \beta'$; we start noticing that $\mathbf{C}^{\measuredangle(W,J)}_{[a,b)\subseteq_\varepsilon[\alpha,\beta)} \subseteq \mathbf{C}^{\measuredangle(W,J)}_{[a,b)\subseteq_\varepsilon[\alpha,\beta')}$ is immediate, and we show the other inclusion arguing by contradiction. Let $E \in \mathbf{C}^{\measuredangle(W,J)}_{[a,b)\subseteq_\varepsilon[\alpha,\beta')}$, and let's show that a nontrivial fiber sequence $A \to E \to B$ destabilizing $E$ on $[\alpha,\beta)$ induces another fiber sequence destabilizing $E$ on $[\alpha,\beta')$. Given such a fiber sequence $A \to E \to B$, we have $\psi A > \psi E > \psi B$; in particular $\psi E > \psi B$; then, we can build the diagram

$$
\begin{array}{ccccc}
A_{[\alpha,\beta')} & \longrightarrow & K_{[\alpha,\beta')} & \longrightarrow & E_{[\alpha,\beta')} \\
\downarrow & & \downarrow & & \downarrow \\
0 & \longrightarrow & (S_{\beta'}B)_{[\beta',\beta)} & \longrightarrow & B_{[\alpha,\beta)} \\
 & & \downarrow & & \downarrow \\
 & & 0 & \longrightarrow & (R_{\beta'}B)_{[\alpha,\beta')}
\end{array}
\tag{7.46}
$$

where $S_{\beta'}B \to B \to R_{\beta'}B$ is the fiber sequence induced by the normal torsion theory $J(\beta')$, and we use a subscript on objects to denote in which endocardium $\mathbf{C}_{[a,b)}$ they lie.

Now we can find the desired contradiction. If we show that $\psi E < \psi B_1$, we have $\psi B < \psi B_1$, hence $\psi B_2 < \psi B$. But now, $\psi B_2 (< \psi B) < \psi E$ entails $\psi E < \psi K$, and this is a thin fiber sequence in $[\alpha, \beta')$ destabilizing $E$ on such an interval. So, we are left with the proof that $\psi E < \psi B_1$: to this end, we estimate the $W$-phase using the $Z$-phase $\phi$ and its proximity with the $W$-phase $\psi$. We have $\phi(E) \in [\alpha+\varepsilon, \beta'-\varepsilon)$ and $|\phi - \psi| < \varepsilon/2$, so if $\phi E < \beta' - \varepsilon$ one has $\psi E < \beta' - \varepsilon/2$, and if $\phi B_1 \geq \beta'$ one has $\psi B_1 \geq \beta' - \varepsilon/2 > \psi E$. This concludes the proof. □

DEFINITION 7.5.18. ($W$-SEMISTABLE THIN OBJECTS): As a consequence of Proposition **7.5.17**, given $a < b \leq a + \varepsilon$, we can define

$$\mathbf{C}^{\measuredangle(W,J)}_{[a,b)} = \mathbf{C}^{\measuredangle(W,J)}_{[a,b)\subseteq_\varepsilon [\alpha,\beta)} \tag{7.47}$$

for *any* $[\alpha, \beta) \in \text{ENV}_\varepsilon(a, b)$. We call $\mathbf{C}^{\measuredangle(W,J)}_{[a,b)}$ the subcategory of $W$-semistable objects in $\mathbf{C}_{[a,b)}$.

PROPOSITION 7.5.19. If $a \leq a' < b' \leq b < a + \varepsilon$ then

$$\mathbf{C}^{\measuredangle(W,J)}_{[a',b')} = \mathbf{C}^{\measuredangle(W,J)}_{[a,b)} \cap \mathbf{C}_{[a',b')}; \tag{7.48}$$

as a consequence, we have the following equalities

$$\mathbf{C}^{\measuredangle(W,J)}_{[a,b)\cap[a',b')} = \mathbf{C}^{\measuredangle(W,J)}_{[a,b)} \cap \mathbf{C}_{[a',b')} = \mathbf{C}^{\measuredangle(W,J)}_{[a',b')} \cap \mathbf{C}_{[a,b)} \tag{7.49}$$

for *any* (not only those contained one into the other) pair of intervals $[a, b), [a', b')$ (with the consistent convention that $\mathbf{C}_\varnothing := \mathbf{0}$).

*Proof.* Prop. **7.5.17** ensures that both sides consist of objects $E$ which are $W$-semistable on the same interval $[\alpha^*, \beta^*)$ in $\varepsilon$-ENV$([a, b) \cap [a', b'))$ (it suffices to choose the biggest among two, one $\varepsilon$-enveloping $[a, b)$ and the other $\varepsilon$-enveloping $[a', b')$). □

Prop. **7.5.17** ensures that the following definition is sound:

DEFINITION 7.5.20. ($W$-SEMISTABLE OBJECTS): Let $W$ be a prestability function with $\|Z - W\|_\sigma \leq \sin \varepsilon$. We define the subcategory $\mathbf{C}^{\measuredangle(W,J)}$ of ($J$-thin and) *$W$-semistable objects* to be the full subcategory of $\mathbf{C}$ having objects $E \in \mathbf{C}$ such that there is a suitably small interval $[a, b) \subseteq [a, a+\varepsilon/2)$ for which $E \in \mathbf{C}_{[a,b)}$, and it is $W$-semistable. In other words,

$$\mathbf{C}^{\measuredangle(W,J)} = \bigcup_{a<b<a+\varepsilon/2} \mathbf{C}^{\measuredangle(W,J)}_{[a,b)}. \tag{7.50}$$

LEMMA 7.5.21. One has $\mathbf{C}^{\measuredangle(W,J)}_{[a,b)} = \mathbf{C}^{\measuredangle(W,J)} \cap \mathbf{C}_{[a,b)}$.

*Proof.* It is a direct consequence of Prop. **7.5.17**:

$$\left(\bigcup_{c<d<c+\varepsilon/2} \mathbf{C}^{\measuredangle(W,J)}_{[c,d)}\right) \cap \mathbf{C}_{[a,b)} = \bigcup_{c<d<c+\varepsilon/2} \mathbf{C}^{\measuredangle(W,J)}_{[c,d)} \cap \mathbf{C}_{[a,b)} = \bigcup_{\substack{c<d<c+\varepsilon/2 \\ [c,d)\cap[a,b)\neq\varnothing}} \mathbf{C}^{\measuredangle(W,J)}_{[c,d)\cap[a,b)}$$

This is, by construction, contained in $\mathbf{C}^{\measuredangle(W,J)}_{[a,b)}$, and it obviously contains it as one of the summands. □

Lemma 7.5.22. If $\|Z - W\|_\sigma < \sin(\varepsilon/2)$, then $\mathbf{C}^{\measuredangle(W,J)}_t \subseteq \mathbf{C}_{[t-\varepsilon,t+\varepsilon)}$.

*Proof.* Let $E \neq 0$ be an object in $\mathbf{C}^{\measuredangle(W,J)}_t$. Then $E \in \mathbf{C}_{[a,b)}$ for real numbers $a < b < a+\varepsilon/2$. Therefore, $\varphi(E) \in [a,b)$ and then $t = \psi(E) \in [a-\varepsilon/2, b+\varepsilon/2)$. It follows that $[a,b) \subseteq [t-\varepsilon, t+\varepsilon)$ and then $E \in \mathbf{C}_{[t-\varepsilon,t+\varepsilon)}$. □

Proposition 7.5.23. Let $E_1 \in \mathbf{C}^{\measuredangle(W,J)}_{t_1}$ and $E_2 \in \mathbf{C}^{\measuredangle(W,J)}_{t_2}$ with $t_1 > t_2$. Then $\mathbf{C}(E_1, E_2) = 0$.

*Proof.* Start by assuming $t_1 - t_2 < 2\varepsilon$. Then, Lemma **7.5.22** entails that $E_1 \in \mathbf{C}_{[t_1-\varepsilon,t_1+\varepsilon)}$ and $E_2 \in \mathbf{C}_{[t_2-\varepsilon,t_2+\varepsilon)}$, hence both lie in $\mathbf{C}_{[t_1-\varepsilon,t_2+\varepsilon)}$. The interval $[t_1-\varepsilon, t_2+\varepsilon)$ is $\varepsilon$-enveloped by $[t_1-2\varepsilon, t_2+2\varepsilon)$. So, $E_1, E_2$ are both $W$-semistable in $[t_1-\varepsilon, t_2+\varepsilon)$ with $\psi(E_1) > \psi(E_2)$. Then there are no nontrivial morphisms from $E_1$ to $E_2$, i.e., $\mathbf{C}(E_1,E_2) = 0$. Now when $t_1 - t_2 \geq 2\varepsilon$, we have $[t_1-\varepsilon, t_1+\varepsilon) \cap [t_2-\varepsilon, t_2+\varepsilon) = \emptyset$, so the subcategory $\mathbf{C}_{[t_1-\varepsilon,t_1+\varepsilon)}$ containing $E_1$ is left-orthogonal to the subcategory $\mathbf{C}_{[t_2-\varepsilon,t_2+\varepsilon)}$ containing $E_2$. □

Now, let $\mathbf{S} \subseteq \mathbf{C}$ be a full subcategory, seen as a set of objects, and denote $\langle \mathbf{S} \rangle$ (or $\langle \mathbf{S} \rangle_{\mathbf{C}}$ when the context does not uniquely specify the embedding) the extension-closure of $\mathbf{S}$, full in $\mathbf{C}$. Define

$$\mathbf{C}^{(W)}_{<t} = \langle \{\mathbf{C}^{\measuredangle(W,J)}_{\psi=s}\}_{s<t} \rangle \qquad \mathbf{C}^{(W)}_{\geq t} = \langle \{\mathbf{C}^{\measuredangle(W,J)}_{\psi=s}\}_{s\geq t} \rangle. \tag{7.51}$$

Lemma 7.5.24. The subcategories $\mathbf{C}^{(W)}_{\geq t}, \mathbf{C}^{(W)}_{<t}$ form an orthogonal pair.

*Proof.* Generally, if $\mathbf{S}_1$ and $\mathbf{S}_2$ are two subcategories of $\mathbf{C}$ with $\mathbf{S}_1 \perp \mathbf{S}_2$, then $\langle \mathbf{S}_1 \rangle \perp \langle \mathbf{S}_2 \rangle$. This can be easily proved by double induction on the "length" of extensions: one first shows by induction on the length of the iterated extension by objects in $\mathbf{S}_2$ leading to an object $Y$ in $\langle \mathbf{S}_2 \rangle$ that $S_1 \perp Y$ for every $S_1 \in \mathbf{S}_1$. Next, one shows by induction on the length of the iterated extension by objects in $\mathbf{S}_1$ leading to an object $X$ in $\langle \mathbf{S}_1 \rangle$ that $X \perp Y$. □

THEOREM 7.5.25. (FUNDAMENTAL DEFORMATION THEOREM): The pair of subcategories $(\mathbf{C}^{(W)}_{\geq t}, \mathbf{C}^{(W)}_{<t})$ defined in (**7.51**) determines a normal torsion theory on $\mathbf{C}$, for each $t \in \mathbb{R}$. Moreover, the function $J^{(W)}\colon t \mapsto (\mathbf{C}^{(W)}_{\geq t}, \mathbf{C}^{(W)}_{<t})$ is monotone and $\mathbb{Z}$-equivariant, and hence defines a slicing on $\mathbf{C}$, which is $\varepsilon$-near to the initial slicing $J$, in the metric of Def. **7.2.14**.

This slicing $J^{(W)}$ is part of the stability condition $(W, J^{(W)})$ (in particular, $W$ is compatible with the slicing $J^{(W)}$ in the sense of Def. **7.3.2**.SC3), and this is called the *standard deformation* of $\sigma = (Z, J)$.

To simplify the exposition of a somewhat involved argument, we split the proof into several preliminary results; the final argument will follow almost directly from all the preceding considerations.

As a preparatory remark, we prove the following two results, proving that each object in a sufficiently thin interval falls into a $W$-fiber sequence:

LEMMA 7.5.26. Let $X$ be an object in a thin endocardium; then $X$ falls into a fiber sequence

$$X^{(W)}_{\geq t} \to X \to X^{(W)}_{<t} \tag{7.52}$$

where $X^{(W)}_{\geq t} \in \mathbf{C}^{(W)}_{\geq t}$ and $X^{(W)}_{<t} \in \mathbf{C}^{(W)}_{<t}$.

*Proof.* By Prop. **7.5.15** $X$ has a (finite) Postnikov tower with $W$-semistable weaves $\{A_1, \ldots, A_n\}$ of decreasing $W$-phases $\lambda_1 > \lambda_2 > \cdots > \lambda_n$. If we consider this tower, i.e. the diagram

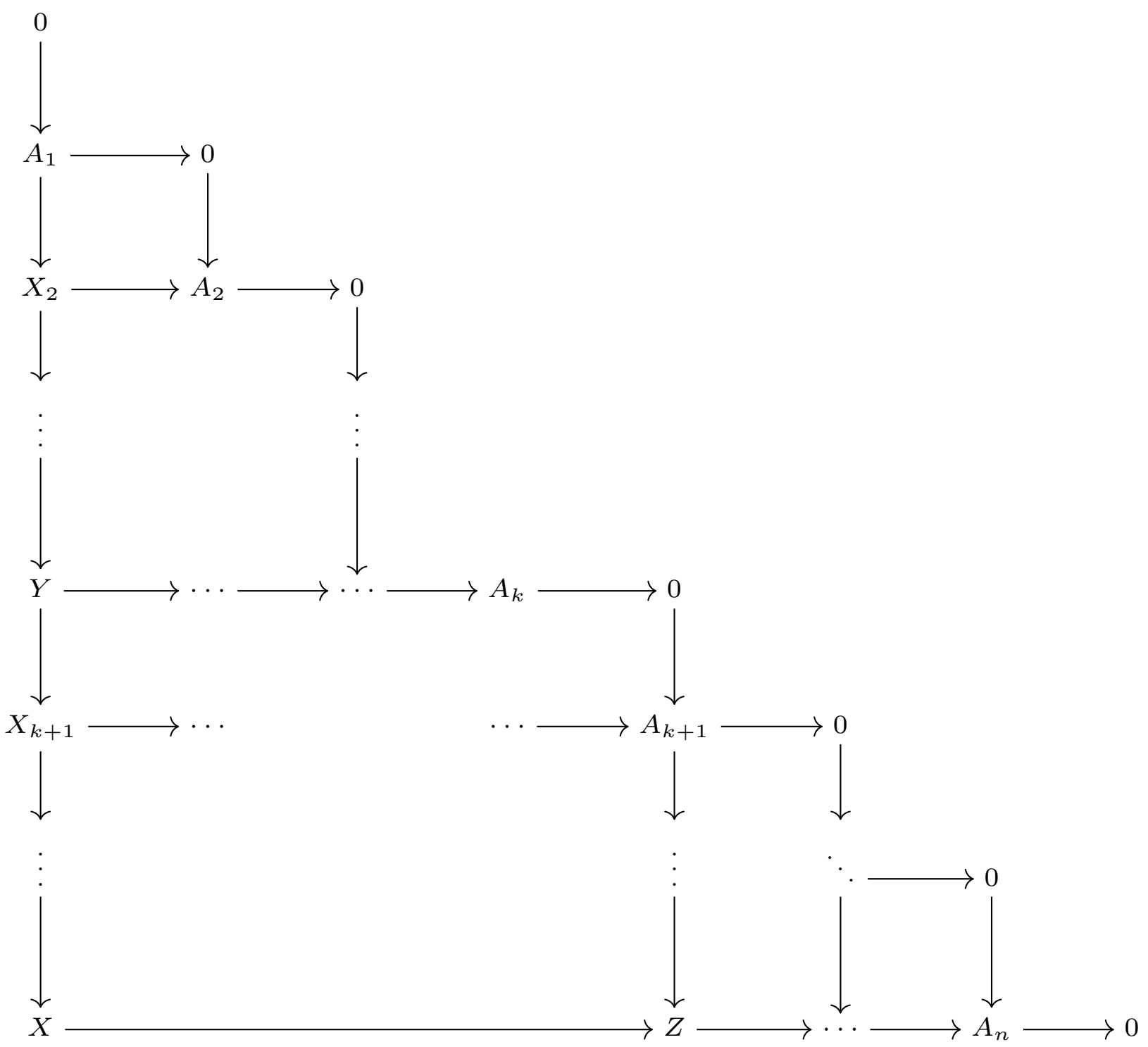

(7.53)

we can extract the pullout subdiagram

$$\begin{array}{ccc} Y & \longrightarrow & 0 \\ \downarrow & & \downarrow \\ X & \longrightarrow & Z \end{array} \tag{7.54}$$

and the two subdiagrams defining the fiber sequence:

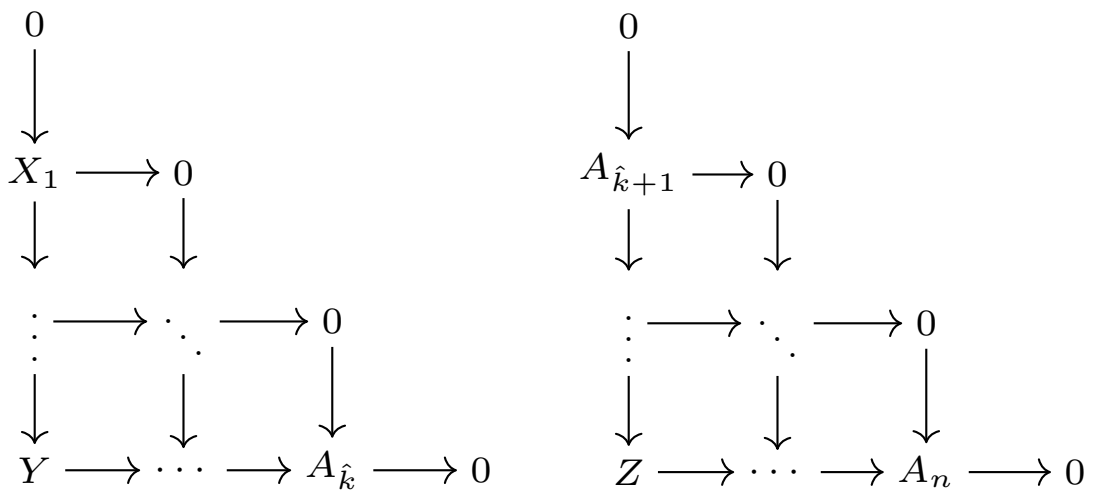

(7.55)

The first diagram says that $Y \in \mathbf{C}^{(W)}_{\geq t}$, and the second that $Z \in \mathbf{C}^{(W)}_{<t}$, since both classes $\mathbf{C}^{(W)}_{\geq t}$, $\mathbf{C}^{(W)}_{<t}$ are extension closed and $Y, Z$ result from iterated extensions done in these classes. □

Lemma 7.5.27. $\mathbf{C}_{<0} = \langle \mathbf{C}_{[a,b)} \mid [a,b) \text{ thin}, b < 0 \rangle$ (dually, $\mathbf{C}_{\geq 0} = \langle \mathbf{C}_{[a,b)} \mid [a,b) \text{ thin}, a \geq 0 \rangle$).

Corollary 7.5.28. $\mathbf{C}_{\geq 0} \subseteq \mathbf{C}^{(W)}_{\geq -\varepsilon}$ (dually, $\mathbf{C}_{<0} \subseteq \mathbf{C}^{(W)}_{<\varepsilon}$).

We now prove the final result as a consequence of Prop. **7.5.15** and Lemma **7.5.26**: each object $X \in \mathbf{C}$ fits into a fiber sequence $X^{(W)}_{\geq t} \to X \to X^{(W)}_{<t}$.

Let us consider the factorization of the initial morphism $0 \to X$ with respect to the slicing $J$: we obtain the diagram

$$\begin{array}{ccccccc}
0 & & & & & & \\
\downarrow & & & & & & \\
S_{t+\varepsilon}X & \longrightarrow & 0 & & & & \\
\downarrow & & \downarrow & & & & \\
S_{t-\varepsilon}X & \longrightarrow & A & \longrightarrow & 0 & & \\
\downarrow & & \downarrow & & \downarrow & & \\
X & \longrightarrow & R_{t+\varepsilon}X & \longrightarrow & R_{t-\varepsilon}(X) & \longrightarrow & 0
\end{array} \tag{7.56}$$

where $X_{t+\varepsilon} \in \mathbf{C}_{\geq t+\varepsilon}$, $C \in \mathbf{C}_{<t-\varepsilon}$ and $A \in \mathbf{C}_{[t-\varepsilon,t+\varepsilon)}$. Now we have $X_{t+\varepsilon} \in \mathbf{C}^{(W)}_{\geq t}$ and $X_{t-\varepsilon} \in \mathbf{C}^{(W)}_{<t}$. Moreover, the interval $[t-\varepsilon, t+\varepsilon)$ is thin, so we have a fiber sequence

$$\begin{array}{ccc}
A^{(W)}_{\geq t} & \longrightarrow & 0 \\
\downarrow & & \downarrow \\
A & \longrightarrow & A^{(W)}_{<t}
\end{array} \tag{7.57}$$

which leads to a refinement of the starting factorization as

$$\begin{array}{ccccccccc}
0 & & & & & & & & \\
\downarrow & & & & & & & & \\
S_{t+\varepsilon}X & \longrightarrow & 0 & & & & & & \\
\downarrow & & \downarrow & & & & & & \\
Y & \longrightarrow & A^{(W)}_{\geq t} & \longrightarrow & 0 & & & & \\
\downarrow & & \downarrow & & \downarrow & & & & \\
S_{t-\varepsilon}X & \longrightarrow & A & \longrightarrow & A^{(W)}_{<t} & \longrightarrow & 0 & & \\
\downarrow & & \downarrow & & \downarrow & & \downarrow & & \\
X & \longrightarrow & R_{t+\varepsilon}X & \longrightarrow & Z & \longrightarrow & R_{t-\varepsilon}X & \longrightarrow & 0
\end{array} \tag{7.58}$$

Notably, the objects $Y, Z$ factor the arrows $S_{t+\varepsilon}X \to S_{t-\varepsilon}X$ and $R_{t+\varepsilon}X \to R_{t-\varepsilon}X$ which "approximate" in some sense the desired reflections; the idea behind this proof is to show that $Y$ plays the rôle of the desired coreflection $S^{(W)}_t X$ and $Z$ plays the rôle of the reflection $R^{(W)}_t X$ at level $t$, cutting at time $t$ and falling in the desired classes: since both $\mathbf{C}^{(W)}_{\geq t}$ and $\mathbf{C}^{(W)}_{<t}$ are extension closed, we have

- $Y \in \mathbf{C}^{(W)}_{\geq t}$, since in the fiber sequence $S_{t+\varepsilon}X \to Y \to A^{(W)}_{\geq 0}$ the extremal objects both lie in $\mathbf{C}^{(W)}_{\geq t}$ (by construction, and Cor. **7.5.28**),
- $Z \in \mathbf{C}^{(W)}_{<t}$, since in the fiber sequence $A^{(W)}_{<0} \to Z \to R_{t-\varepsilon}X$ the extremal objects both lie in $\mathbf{C}^{(W)}_{<t}$ (again by construction, and invoking the dual of Cor. **7.5.28**: $R_{t-\varepsilon}X \in \mathbf{C}_{<-\varepsilon} \subseteq \mathbf{C}^{(W)}_{<0}$).

This, in particular, shows that the process of building the fiber sequence

$$Y \to X \to Z \tag{7.59}$$

is entirely canonical, and since [RT07, **3.1**] holds in the stable setting the objects $Y, Z$ are the coreflection and reflection of a normal torsion theory on $\mathbf{C}$; we denote these functors $(S^{(W)}, R^{(W)})$. Now, the characterization of torsion/free classes from the pair coreflection/reflection entails that the normal torsion theory is completely determined from the pair $(S^{(W)}, R^{(W)})$ via the relations

$$\begin{aligned}
\mathcal{E}^{(W)} &= \{f \in \hom(\mathbf{C}) \mid R^{(W)}f \text{ iso}\} \\
\mathcal{M}^{(W)} &= \{g \in \hom(\mathbf{C}) \mid S^{(W)}g \text{ iso}\}.
\end{aligned}$$

To show that this defines a slicing, we have to prove $\mathbb{Z}$-equivariancy and monotonicity. Since $\psi(E[1]) = \psi(E) + 1$, we have $\mathbf{C}^{\measuredangle(W,J)}_{\psi=s}[1] = \mathbf{C}^{\measuredangle(W,J)}_{\psi=s+1}$,

where we exploited the fact that shifts preserve semistable objects and the fact that $\psi(E[1]) = \psi(E) + 1$. Since shift commutes with the $\langle - \rangle$ operation on classes, we have $\mathbf{C}^{(W)}_{<t}[1] = \langle \{\mathbf{C}^{\measuredangle(W,J)}_{\psi=s+1}\}_{s\geq t}\rangle = \langle \{\mathbf{C}^{\measuredangle(W,J)}_{\psi=s}\}_{s\geq t+1} = \mathbf{C}^{(W)}_{<t+1}$. Similarly, $\mathbf{C}^{(W)}_{\geq t}[1] = \mathbf{C}^{(W)}_{\geq t+1}$.

To conclude the proof we show that $J^{(W)} \in B_\varepsilon(J)$: following Def. **7.2.20** and Prop. **7.2.21**, we must show that there exists $\delta > 0$ such that

$$(\forall t \in \mathbb{R})\ \mathcal{E}_{t+\varepsilon} \subseteq \mathcal{E}^{(W)}_{t+\delta} \subseteq \mathcal{E}^{(W)}_{t-\delta} \subseteq \mathcal{E}_{t-\varepsilon} \tag{7.60}$$

This is immediate, as a consequence of Lemma **7.5.27** and Cor. **7.5.28**.

# Appendix A

# Stable $\infty$-categories

The present chapter serves as a reference for the rest of the thesis, outlining the fundamentals of stable $\infty$-category theory. Apart from classical literature on triangulated categories ([HJ10, Nee01]) we follow the only available source on stable $\infty$-categories [Lur17], deviating a little from the presentation given there, to add some new considerations and complete, explicit proofs of certain useful classical constructions (like an extensive proof, alternative to that in [Lur17] of the validity of triangulated category axioms in the homotopy category of a stable $\infty$-category).

We start by trying to outline a bit of history of homological algebra to motivate the quest for a higher-categorical formulation of its basic principles. For this account (which makes no claim of originality or completeness), the survey [Wei99] has been an essential source of inspiration.

## A.1 Triangulated higher categories.

> Otra escuela declara [...] que nuestra vida es apenas el recuerdo o reflejo crepuscular, y sin duda falseado y mutilado, de un proceso irrecuperable.
>
> [Bor44], *Tlön, Uqbar, Orbis Tertius*

The notion of triangulated category is deeply linked to homotopy theory. The native language in which Def. **A.1.1** below was originally formulated was *stable homotopy theory*, where suitable sequences of arrows

$$X \to Y \to Z \to \Sigma X \tag{A.1}$$

played an essential rôle in the definition of the *stable homotopy category* of topological spectra and the endofunctor $\Sigma$ acts as the (reduced) suspension, i.e. as the homotopy pushout

$$\begin{array}{ccc} X & \longrightarrow & CX \\ \downarrow & & \downarrow \\ CX & \longrightarrow & \Sigma X \end{array}$$

The invertibility of $\Sigma$ is an essential feature of stable homotopy theory, and the construction giving the universal category where $\Sigma$ becomes an equivalence is part of the so-called *Spanier-Whitehead stabilization* $\mathrm{SW}(\mathbf{Spc})$ of the category of CW-complexes **Spc**. We briefly investigate the construction of $\mathrm{SW}(\mathbf{C})$ in §**A.4**.

A first axiomatization for the phenomena giving rise to these structures dates back to A. Dold and D. Puppe's [DP61]; subsequently, motivated by this result, Grothendieck and Verdier recognized a similar structure on the homotopy category of $\mathbf{Ch}(\mathbf{A})$ (chain complexes on the abelian category $\mathbf{A}$), and encoded this procedure of modding out null-homotopic maps to construct the *derived category* $\mathbf{D}(\mathbf{A})$ of $\mathbf{Ch}(\mathbf{A})$.

Verdier outlined in his [Ver96] a (ingenious but rather cumbersome) set of axioms, aimed at capturing the behaviour of these notable classes of the *distinguished triangles* (**A.1**), acting like exact sequences and involving an additive autoequivalence $\Sigma\colon \mathbf{C} \to \mathbf{C}$, generalizing the reduced suspension $\Sigma$.

Subsequently, D. Quillen axiomatized the notion of *abstract homotopy theory* [Bau89] with his definition of a *model category*; this in some sense unified the language of homotopy and homology theory, giving a more profound intuition of the latter being an additive manifestation of the former, and in particular conveying the idea that homotopies behave the same way also outside the category of spaces (and exist, for example, between maps of chain complexes, or maps of simplicial sets).

Even at this point however the systematization of the theory of triangulated categories was far from being satisfactory, since the origin of the axioms was obscure and really far from being canonical. This "bad behaviour" shows up in several practical situations, populating the dense literature on the subject: after having given the definition of a triangulated category, we embark on a deep analysis of their meaning. Convenient shorthands to denote a distinguished triangle in a triangulated category $\mathbf{C}$ are the following;

$$X \to Y \to Z \to^{+},\; X \to Y \to Z \to,\; X \to Y \to Z \to X[1]$$

(see **A.2.10**) or even $X \to Y \to Z$, when no ambiguity can arise from this compactness.

Definition A.1.1. (Triangulated category): A category $\mathbf{C}$ is called *suspended* if it is endowed with an endofunctor $\Sigma$; an additive category with suspension $(\mathbf{C}, \Sigma)$ is said to be *triangulated* if the following axioms are satisfied:

pt 1) The suspension endofunctor is an equivalence of categories;

pt 2) There exists a class of diagrams in $\mathbf{C}$, called *distinguished triangles* of the form $X \to Y \to Z \to \Sigma X$ (often denoted $X \to Y \to Z \to^+$ for short) which is closed under isomorphism and contains every sequence of the form $X \xrightarrow{\mathrm{id}_X} X \to 0 \to \Sigma X$;

pt 3) Any arrow $f\colon \Delta[1] \to \mathbf{C}$ fits into at least one distinguished triangle $X \xrightarrow{f} Y \to Z \to \Sigma X$;

pt 4) (rotation) The diagram $X \xrightarrow{u} Y \xrightarrow{v} Z \xrightarrow{w} \Sigma X$ is distinguished if and only if the "rotated diagram" $Y \xrightarrow{-v} Z \xrightarrow{-w} \Sigma X \xrightarrow{-\Sigma u} \Sigma Y$ is distinguished;

pt 5) (completion) In any diagram of the form

$$\begin{array}{ccccccc} X & \longrightarrow & Y & \longrightarrow & Z & \longrightarrow & \Sigma X \\ {\scriptstyle f}\downarrow & & {\scriptstyle g}\downarrow & & & & \downarrow{\scriptstyle \Sigma f} \\ X' & \longrightarrow & Y' & \longrightarrow & Z' & \longrightarrow & \Sigma X' \end{array} \tag{A.2}$$

where the rows are distinguished triangles, there exists a morphism $h\colon Z \to Z'$ making the whole diagram a morphism of triangles (which, once triangles are regarded as suitable functors $J \to \mathbf{C}$ are simply natural transformations between two such functors).

tr) Given *three* distinguished triangles

$$X \xrightarrow{f} Y \to Y/X \qquad Y \xrightarrow{g} Z \to Z/Y \qquad X \xrightarrow{gf} Z \to Z/X \tag{A.3}$$

(where the cone of each arrow is temporarily represented as a quotient to suggest the meaning of the cone construction) arranged in a "braid" diagram

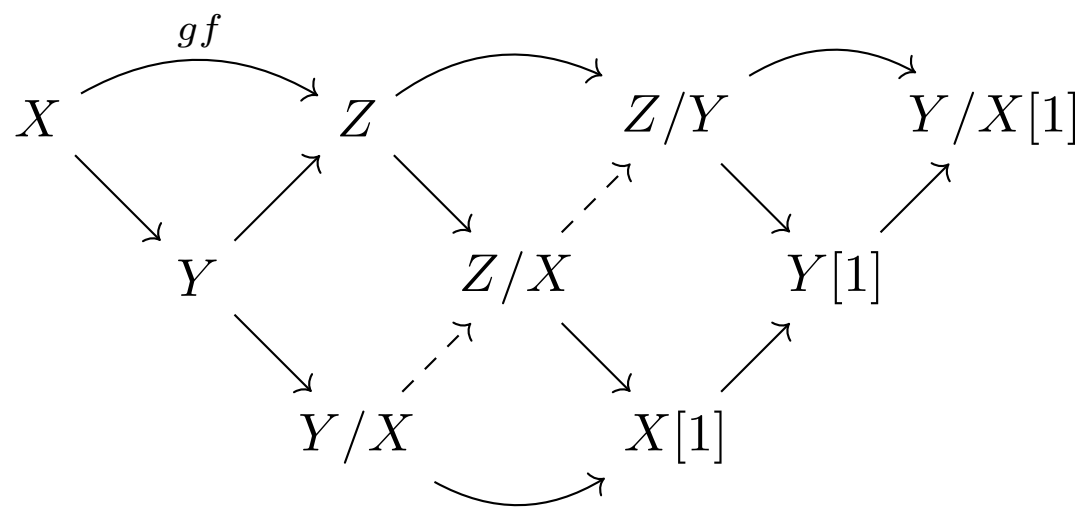

then there is a (non-unique) way to complete it with the arrows $s, t$ indicated.

Now, a deeper analysis of the design behind these axioms shows several drawbacks:

- axiom pt **3)** embeds a map $f\colon X \to Y$ in a distinguished triangle $X \xrightarrow{f} Y \to Z \to \Sigma X$, with a procedure which is not canonical, and yet all the most important examples of triangulated category show this property by means of “weakly canonical” constructions (the object $Z$ in the axiom is determine “up to a contractible space of choices” as the *homotopy colimits* or *mapping cone* of $f$, in some flavour of higher category theory).

- On the same lines, property PT **5)**, which asserts that each “morphism of triangles”

$$\begin{array}{ccccccc} A & \longrightarrow & B & \longrightarrow & C & \longrightarrow & \Sigma A \\ \downarrow f & & \downarrow g & & \downarrow h & & \\ A' & \longrightarrow & B' & \longrightarrow & C' & \longrightarrow & \Sigma A' \end{array}$$

is determined by only two elements, is not canonical: there is no unique choice of a third element, the only hope being that there is a choice which is well-suited for “some” other purpose, since again in the most important cases like $\mathbf{D}(\mathbf{A})$ or $\mathrm{Ho}(\mathbf{Sp})$ the completion axiom holds as a consequence of a universal property (of the homotopy co/-limits involved).

- (This is a more conceptual, but important drawback.) As it is noted in [MK07], the derived category of an abelian category $\mathbf{A}$, taken as a triangulated category alone, has no universal property;

From a modern perspective, is is easy to see that this situation reflects some deep features of homotopy theory: the category of chain complexes $\mathbf{Ch}(\mathbf{A})$ has a fairly natural choice of a model structure; this entails that $\mathbf{Ch}(\mathbf{A})$ is a fairly rich environment; the localization procedure outlined by Verdier does not retain these additional pieces of information encoded in the homotopy co/limits in $\mathbf{Ch}(\mathbf{A})$, because they are hidden in a higher categorical structure that the localization procedure is not able to preserve.

It must be said, however, that despite this highly unsatisfactory situation, a great deal of refined mathematics stemmed from the theory of triangulated categories:

- One of Verdier-Grothendieck’s primary tasks (to shed a light on the construction of derived functors) is easily achieved (the language of model categories clarifies best the meaning and construction of derived functors);
- In a suitable sense the derived category of sheaves on a good space contains enough information to rebuild the space from scratch (this is a result in reconstruction theory, mainly worked out in [BO01]);
- Several properties of an abelian category $\mathbf{A}$ can be deduced from the study of a notable kind of subcategories of $\mathbf{D}(\mathbf{A})$ (the adjacent classes of a “$t$-structure”) on $\mathbf{D}(\mathbf{A})$) and of a generic triangulated category $\mathbf{D}$ (this is by far the most important application for the purposes of the present thesis).

In light of this, one could argue that Def. **A.1.1** behaves like the definition of topological spaces to a certain extent: the definition is not modeled to be user-friendly, but to be pervasive, and concrete examples of the definition often enjoy additional properties making them more wieldy.

However, with the passing of time, understanding the deep meaning of the axioms in Def. **A.1.1** became more and more a priority. It became evident that triangulated category where the "false and mutilated memory of an irrecoverable process", behaving like 1-dimensional shadows of a higher dimensional notion:

> triangulated categories arise as *decategorification* of some structure taking place in the $\infty$-categorical world, and the axioms defining them are designed to keep track of the 1-categorical trace of this more refined notion.

Shadows of objects retain no information about their colours; in the same spirit, triangulated categories retain little or no information about the higher structure generating them.[(1)]

Because of these reasons, it would be desirable to have at our disposal a more intrinsic notion of triangulated category, satisfying some reasonable requests of universality: whenever a higher category $\mathbf{C}$ enjoys a property which we will call "stability", then

SC1) its homotopy category $\text{Ho}(\mathbf{C})$ carries a triangulated structure in the sense of Def. **A.1.1**;

SC2) the axioms characterizing a triangulated structure are "easily verified and well-motivated consequences of evident universal arguments"([Lur17, Remark **1.1.2.16**]);

SC3) classical derived categories arising in Homological Algebra can be regarded as homotopy categories of stable $\infty$-categories functorially associated to an abelian category $\mathbf{A}$ (see [Lur17, §**1.3.1**]).

The most common examples show that finding a triangulated structure on $\text{Ho}(\mathbf{C})$ is often sufficient for most practical purposes where one only needs information that survive the homotopy identification process. However, as soon as one needs to take into account additional information about homotopy co/limits that existed in $\mathbf{C}$, its stability comes into play.

(1) Albeit seldom spelled out explicitly, we can trace in this remark a fundamental tenet of the theory outlined in [Lur17]:

> In the same way every shadow comes from an object, produced once the sun sheds a light on it, every "non-pathological" triangulated category is the 1-dimensional shadow (i.e. the homotopy category) of an higher-dimensional object.

No effort is made here to hide that this fruitful metaphor is borrowed from [Car10], even if with a different meaning and in a different context.

## A.2 Building stable categories.

> A stable mind is *fudoshin*, a mind not disturbed or upset by verbal mistreatment
>
> M. Hatusmi

Pathological examples aside (see [MSS07], from which the following distinction is taken verbatim), there are essentially two procedures to build “nice” triangulated categories:

- In Algebra they often arise as the stable category of a Frobenius category ([Hel68, **4.4**], [GM96, **IV.3** Exercise **8**]).
- In algebraic topology they usually appear as a full triangulated subcategory of the homotopy category of a Quillen stable model category [Hov99, **7.1**].

The (closure under equivalence of) these two classes contain respectively the so-called *algebraic* and *topological* triangulated categories described in [Sch10].

Classical triangulated categories can also be seen as Spanier-Whitehead stabilizations of the homotopy category $\mathrm{Ho}(\mathbf{M})$ of a pointed model category $\mathbf{M}$ (see [Del04] for an exhaustive treatment of this construction, which we sketch in §**A.4** below).

So, several different models for triangulated higher categories arose as a reaction to different needs in abstract homological algebra (where derived categories of rings play a central rôle), algebraic geometry (where one is led to study derived categories of –modules of– sheaves of rings) or in a fairly non-additive setting as algebraic topology (where the main example of such a structure is the homotopy category of spectra $\mathrm{Ho}(\mathbf{Sp})$); there’s no doubt that allowing a certain play among different models may be more successful in describing a particular phenomenon (or a wider range of phenomena), whereas being forced to a particular one may turn out to be insufficient.

Now, according to the “principle of equivalence” between models of higher category theory there must be a similar notion in the language of $\infty$-categories, i.e. some property of an $\infty$-category $\mathbf{C}$ ensuring that the “requests” sc1)—sc3) above are satisfied.

Building this theory is precisely the aim of [Lur17, Ch. **1.1**]. As this is the most interesting and well-developed model at the moment of writing, and the one we constantly had in mind, we now give a rapid account of the main lines of stable $\infty$-category theory.

We invite the reader to take [Lur17] as a permanent reference for this section, hoping to convince those already acquainted with the theory of triangulated categories that stable $\infty$-categories are in fact a simpler and more manageable reformulation of the basic principles they already know how to manipulate.

### A.2.1 Stable $\infty$-categories.

Let $\square = \Delta[1] \times \Delta[1]$ be the "prototype of a square",

$$\begin{array}{ccc} (0,0) & \longrightarrow & (0,1) \\ \Big\downarrow & & \Big\downarrow \\ (1,0) & \longrightarrow & (1,1) \end{array}$$

such that the category of functors $\square \to \mathbf{C}$ consists of commutative squares in $\mathbf{C}$. With this identification in mind, we can give the following

Definition A.2.1. ((Co)cartesian square): A diagram $F \colon \square \to \mathbf{C}$ in a (finitely bicomplete) $\infty$-category is said to be *cocartesian* (resp., *cartesian*) if the square

$$\begin{array}{ccc} F(0,0) & \longrightarrow & F(0,1) \\ \Big\downarrow & & \Big\downarrow \\ F(1,0) & \longrightarrow & F(1,1) \end{array}$$

is a homotopy pushout (resp., a homotopy pullback).

Alternatively, one can characterize the category $\square$ as $\Delta[1] \times \Delta[1] = (\Lambda^2_2)^{\lhd} = (\Lambda^2_0)^{\rhd}$ (see the diagrams below, and [Lur09] for the notation;

$$\begin{array}{ccc} (0,0) & \longrightarrow & (1,0) \\ \Big\downarrow & {\scriptstyle \Lambda^2_0} & \\ (0,1) & & \end{array} \qquad\qquad \begin{array}{ccc} & & (1,0) \\ & {\scriptstyle \Lambda^2_2} & \Big\downarrow \\ (0,1) & \longrightarrow & (1,1) \end{array} \tag{A.4}$$

each of these descriptions will turn out to be useful). In the same way, we denote pictorially the two horn-inclusions

$$\begin{aligned} i_{\ulcorner} &\colon \ulcorner \to \square \quad \big( = \Lambda^2_0 \to (\Lambda^2_0)^{\rhd} \big) \\ i_{\lrcorner} &\colon \lrcorner \to \square \quad \big( = \Lambda^2_2 \to (\Lambda^2_2)^{\lhd} \big) \end{aligned}$$

(see [Lur09, Notation **1.2.8.4**]) and the induced maps

$$i^*_{\ulcorner} \colon \mathrm{Map}(\square, \mathbf{C}) \to \mathrm{Map}(\ulcorner, \mathbf{C}) \tag{A.5}$$

$$i^*_{\lrcorner} \colon \mathrm{Map}(\square, \mathbf{C}) \to \mathrm{Map}(\lrcorner, \mathbf{C}) \tag{A.6}$$

from the category of commutative squares in $\mathbf{C}$, "restricting" a given diagram to its top or bottom part, respectively. These functors are part of a

string of adjoints

$$(i_\ulcorner)_! \dashv \boxed{i_\ulcorner^* \dashv (i_\ulcorner)_*} : \mathrm{Map}(\Box, \mathbf{C}) \leftrightarrows \mathrm{Map}(\ulcorner, \mathbf{C}) \tag{A.7}$$

$$\boxed{(i_\lrcorner)_! \dashv i_\lrcorner^*} \dashv (i_\lrcorner)_* : \mathrm{Map}(\Box, \mathbf{C}) \leftrightarrows \mathrm{Map}(\lrcorner, \mathbf{C}) \tag{A.8}$$

where $(i_\ulcorner)_!$ and $(i_\lrcorner)_*$ are easily seen to be evaluations at the initial and terminal object of $\ulcorner$ and $\lrcorner$, respectively.

Given $F \in \mathrm{Map}(\Box, \mathbf{C})$ the canonical morphisms obtained from the boxed adjunctions,

$$\eta_{\ulcorner,F} : F \to (i_\ulcorner)_* i_\ulcorner^* F$$
$$\epsilon_{\lrcorner,F} : (i_\lrcorner)_! i_\lrcorner^* F \to F$$

give the canonical "comparison" arrow $F(1,1) \to \varinjlim i_\ulcorner^* F$ and $\varprojlim i_\lrcorner^* F \to F(0,0)$.

With these notations we can say that

- $F \in \mathrm{Map}(\Box, \mathbf{C})$ is *cartesian* if $\eta_{\ulcorner,F}$ is invertible;
- $F \in \mathrm{Map}(\Box, \mathbf{C})$ is *cocartesian* if $\epsilon_{\lrcorner,F}$ is invertible.

DEFINITION A.2.2. (STABLE $\infty$-CATEGORY): A $\infty$-category $\mathbf{C}$ is called *stable* if

(1) it has any finite (homotopy) limit and colimit;
(2) A square $F : \Box \to \mathbf{C}$ is cartesian if and only if it is cocartesian.

NOTATION A.2.3. Squares which are both pullback and pushout are called *pulation squares* or *bicartesian squares* (see [AHS90, Def. **11.32**]) in the literature. We choose to call them *pullout squares* and we refer to axiom **2** above as the *pullout axiom*: in such terms, a stable $\infty$-category is a finitely bicomplete $\infty$-category satisfying the pullout axiom.

REMARK A.2.4. The pullout axiom is by far the most characteristic feature of stable $\infty$-categories; it is the most ubiquitously applied property of diagrams in such a setting, to the point that in some sense the rest of the present section is devoted to a better understanding of the consequences of this statement alone.

We being with the simplest remark: most of the arguments in the following discussion are a consequence of the following

REMARK A.2.5. (A 3-FOR-2 PROPERTY FOR PULLOUTS): The pullout axiom implies that the class $\mathcal{P}$ of pullout squares in a category $\mathbf{C}$ satisfies a 3-for-2 property: in fact, it is a classical, easy result (see [AHS90, Prop. **11.10**] and its dual) that pullback squares, regarded as morphisms in $\mathbf{C}^{\Delta[1]}$, form a R32 class and dually, pushout squares form a L32 class (these are called *pasting laws* for pullback and pushout squares) in the sense of our Definition **1.4.6**.

Notation A.2.6. It is a common practice to denote diagrammatically a (co)cartesian square by “enhancing” the corner where the universal object sits (this well-established convention has been used with no further mention throughout our discussion): as a “graphical” representation of the auto-duality of the pullout axiom, we choose to denote a pullout square by enhancing *both* corners:

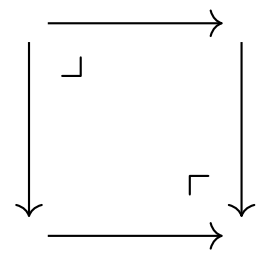

Remark A.2.7. Any 1-category **C** satisfying the pullout axiom with respect to 1-dimensional pullbacks and pushouts is equivalent to the terminal category.

*Proof.* First of all notice that in a stable $\infty$-category the functors $\Sigma \dashv \Omega$ form an equivalence of $\infty$-categories; this follows from the fact that in the diagram the object $X$ has the universal property of both $\Omega\Sigma X$ and $\Sigma\Omega X$. Now, from the fact that $\Sigma X$ is the pushout of $0 \leftarrow X \to 0$, we deduce that $\Sigma X \cong 0$. □

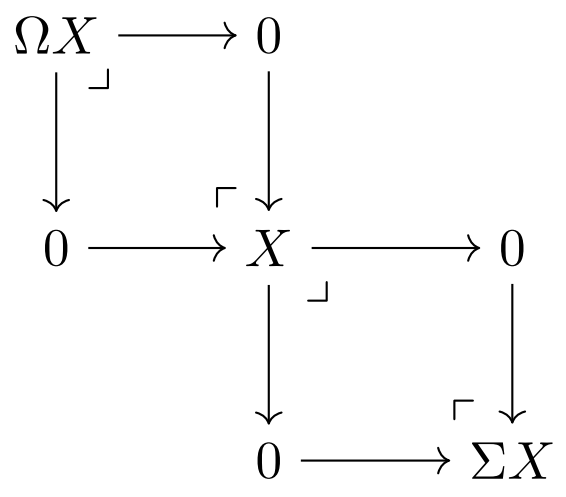

Among the most essential features of stability, there is the fact that all stable categories are naturally enriched over abelian groups (or, rather, over a “homotopy meaningful” notion of abelian group): Remark **A.2.7** above showed that the pullout axiom is a really strong assumption on a category, so strong that it can only live in the “weakly-universal” world of $(\infty, 1)$-categories. Now, we learn that the pullout axiom characterizes almost completely the structure of a stable $\infty$-category.

Remark A.2.8. (The pullout axiom induces an enrichment.): A stable $\infty$-category **C**

- has a zero object, i.e. there exists an arrow $1 \to \varnothing$ (which is forced to be an isomorphism);
- **C** has biproducts, i.e. $X \times Y \simeq X \amalg Y$ for any two $X, Y \in \mathbf{C}$, naturally in both $X$ and $Y$.

We skip the proof of this statement; the interested reader can take it as an exercise and test the power of the pullout axioms.

REMARK A.2.9. The proof of the above statement heavily relies on a result of Freyd's [Fre64]; the *biproduct* of objects $X, Y$ in $\mathbf{C}$ can be characterized as an object $S = S_{X,Y}$ such that

- There are arrows $Y \leftrightarrows S \leftrightarrows X$;
- The arrow $Y \to S \to Y$ compose to the identity of $Y$, and the arrows $X \to S \to X$ compose to the identity of $X$;
- There are "exact sequences" (in the sense of a pointed, finitely bicomplete category) $0 \to Y \to S \to X \to 0$ and $0 \to X \to S \to Y \to 0$.

The biproduct of $X, Y$ is denoted $X \oplus Y \cong X \times Y \cong X \amalg Y$. A pleasant consequence of Freyd's characterization is that in any additive category the enrichment over the category of abelian groups is *canonical*; in fact, exploiting the isomorphism $Y \times Y \cong Y \amalg Y$ one is able to define the *sum* of $f, g \colon X \rightrightarrows Y$ as

$$f + g \colon X \to X \times X \xrightarrow{(f,g)} Y \times Y \cong Y \amalg Y \to Y \tag{A.9}$$

In fact, this result can be retrieved in the setting of stable $\infty$-categories (see [Lur17, Lemma **1.1.2.9**]); we do not want to reproduce the whole argument: instead we want to investigate the construction of the *loop* and *suspension* functors in a pointed category.

The *suspension* $\Sigma X$ of an object $X$ in a finitely cocomplete, pointed $\infty$-category $\mathbf{C}$ can be defined as the (homotopy) colimit of the diagram $0 \leftarrow X \to 0$; dually, the *looping* (or *loop object*) $\Omega X$ of an object $X$ in such a $\mathbf{C}$ is defined as the (homotopy) limit of $0 \to X \leftarrow 0$.

This notation is natural with a topological intuition in mind, where these operations amount to the *reduced suspension* (see (**A.1**)) and *loop space* of $X$ (thought of as the fiber of the fibration $PX \to X$, where $PX$ is the *path space* of $X$); evaluating a square $F \colon \square \to \mathbf{C}$ at its right-bottom vertex gives an endofunctor $\Sigma \colon \mathbf{C} \to \mathbf{C}$, and the looping $\Omega$ is the right adjoint of this functor $\Sigma$. We depict the objects $\Sigma X, \Omega X$ as vertices of the diagrams

$$\begin{array}{ccc} X & \longrightarrow & * \\ \downarrow & & \downarrow \\ * & \longrightarrow & \Sigma X \end{array} \qquad \begin{array}{ccc} \Omega Y & \longrightarrow & * \\ \downarrow & & \downarrow \\ * & \longrightarrow & Y. \end{array}$$

The pullout axiom defining a stable $\infty$-category implies that these two correspondences (which in general are adjoint functors between $\infty$-categories: see [Lur17, Remark **1.1.2.8**]) are a pair of mutually inverse equivalences ([Gro10, Prop. **5.8**]).

NOTATION A.2.10. In a stable setting, we will often denote the image of $X$ under the suspension $\Sigma$ as $X[1]$, and by extension $X[n]$ will denote, for any $n \geq 2$ the object $\Sigma^n X$ (obviously, $X[0] := X$). Dually, $X[-n] := \Omega^n X$ for any $n \geq 1$.

This notation is in line with the long tradition of denoting by $X[1]$ the *shift* of an object $X$ in a triangulated category; this notation adds to the already existing ones like $X \to Y \to Z \to^+$ and will be used together with the others with no further mention.

Remark A.2.11. (Stable $\infty$-categories are nice): Due to the non-canonical behaviour of axioms PT **1)**–PT **5)**, during the years there have been several attempts to produce a better-behaved axiomatics, more canonical but still general enough to encompass the interesting examples. One of these was the notion of a *Neeman triangulated category*: we address the reader to [Nee91] to get acquainted with the definition.

Here, we show that *every* (homotopy category of a) stable $\infty$-category is Neeman-triangulated.

Proposition A.2.12. ([Nee91]): Let $\mathbf{C}$ be a stable $\infty$-category, and $A \to A' \to A''$, $B \to B' \to B''$ two fiber sequences on $\mathbf{C}$.

Then the commutative square

$$\begin{array}{ccc} A' \oplus B & \to & A'' \oplus B' \\ \downarrow & & \downarrow \\ 0 & \longrightarrow & A[1] \oplus B'' \end{array}$$

is a fiber sequence.

Example A.2.13. (A complete proof of the octahedral axiom): Among all axioms stated in Def. **A.1.1**, the octahedral axiom PT **5)** is the most difficult to motivate. At first sight, it seems like a god-given condition ensuring that some fairly unnatural things happen. On a second thought, however, there are at least two ways to motivate it:

- the axiom is motivated by the desire to see the *freshman algebraist's theorem* hold in triangulated categories: using the same notation as in PT **5)**, the axiom asserts that $\frac{Z/X}{Y/X} \cong Z/Y$;
- the axiom is motivated by the fact that, in the category of spaces, the classical geometric definition of mapping cone of $f\colon X \to Y$, fitting in a sequence $X \to Y \to C(f)$ ensures the presence of a canonical morphism $C(f) \to C(g \circ f)$, and the cofiber of this map is homotopy equivalent to $C(g)$.

In a stable $\infty$-category $\mathbf{C}$ we are in the following situation:

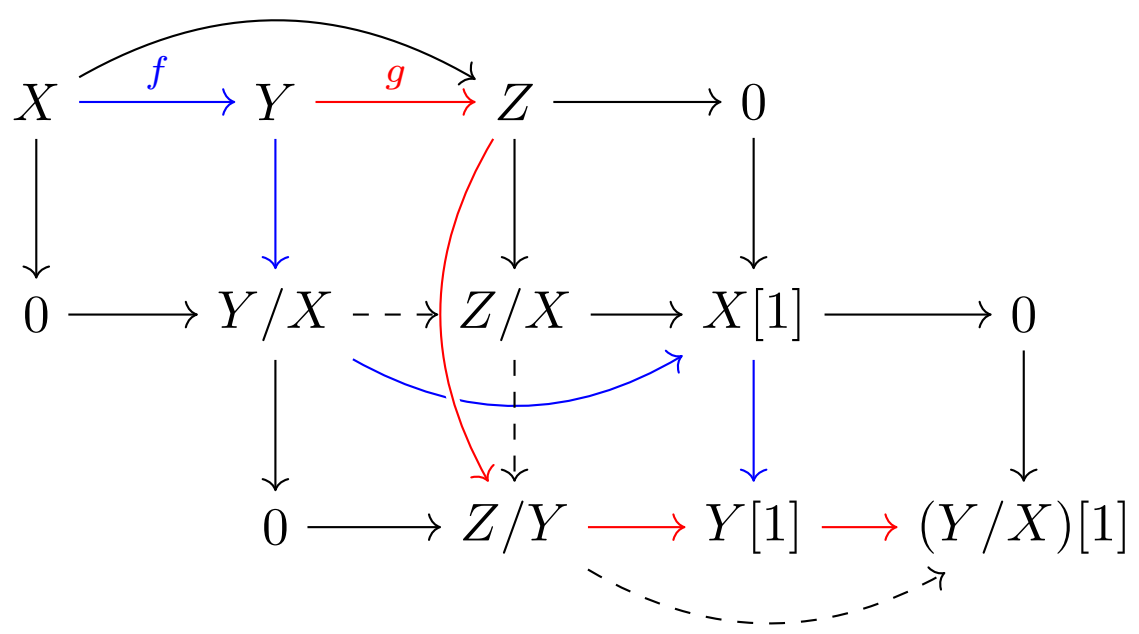

where different colours denote different fiber sequences (i.e., triangles in the homotopy category). Axiom PT **5)** says that we can find arrows $Y/X \to Z/X \to Z/Y$ such that the triangle $Y/X \to Z/X \to Z/Y \to (Y/X)[1]$ is distinguished.

Here is a sketch of a direct, elementary proof for the octahedral axiom.

First of all one must notice that all the preceding axioms PT **1)**–PT **5)** hold almost immediately thanks to the universal properties of the homotopy co/limits involved: in particular, the completion axiom is a consequence of the universal property of a pullback/pushout square, and it implies that the diagram

$$\begin{array}{ccccccc} X & \xrightarrow{f} & Y & \longrightarrow & Y/Z & \longrightarrow & X[1] \\ {\scriptstyle f}\downarrow & & \downarrow{\scriptstyle g} & & & & \downarrow{\scriptstyle f[1]} \\ Y & \xrightarrow[g]{} & Z & \longrightarrow & Z/X & \longrightarrow & Y[1] \end{array} \tag{A.10}$$

can be completed with an arrow $Y/X \xrightarrow{\phi} Z/X$, fitting in the square

$$\begin{array}{ccc} Y & \longrightarrow & Z \\ \downarrow & & \downarrow \\ Y/X & \xrightarrow[\varphi]{} & Z/X \end{array} \tag{A.11}$$

Now consider the objects $V, W$ respectively obtained as pushouts of $Y/X \leftarrow Y \xrightarrow{g} Z$ and $0 \leftarrow Y/X \xrightarrow{\phi} Z/X$; these data fit in a diagram

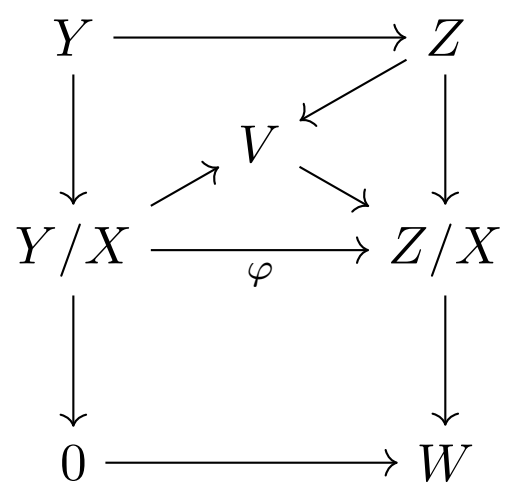

(A.12)

and the 3-for-2 property for pullout squares **A.2.5** now implies that the outer rectangle is a pushout, hence $W \cong Z/Y$. It remains to prove that $V \cong Z/X$; this follows again from the 3-for-2 property applied to

$$\begin{array}{ccccc} X & \longrightarrow & Y & \longrightarrow & Z \\ \downarrow & & \downarrow & & \downarrow \\ 0 & \longrightarrow & Y/X & \longrightarrow & V. \end{array} \tag{A.13}$$

## A.3 *t*-structures.

どのように急須奇妙な
同時に表すことができます
孤独の快適さ
そして、会社の喜び。

Zen *haiku*

The notion of $t$-structure appears in [BBD82] to try to axiomatize the following situation:

Definition A.3.1. (The canonical $t$-structure in $\mathbf{D}(R)$): Let $R$ be a ring, and $\mathbf{D}(R)$ the derived category of modules over $R$; in $\mathbf{D}(R)$ we can find two full subcategories

$$\mathbf{D}_{\geq 0}(R) = \{A_* \in \mathbf{D}(R) \mid H^n(A_*) = 0;\ n \leq 0\}$$
$$\mathbf{D}_{\leq 0}(R) = \{B_* \in \mathbf{D}(R) \mid H^n(B_*) = 0;\ n \geq 0\}$$

such that

- (**orthogonality**): $\hom(A_*[1], B_*) = 0$;
- (**closure under shifts**) $\mathbf{D}_{\geq 0}(R)[1] \subseteq \mathbf{D}_{\geq 0}(R)$ and $\mathbf{D}_{\leq 0}(R)[-1] \subseteq \mathbf{D}_{\leq 0}(R)$;
- (**factorization**) every object $X_* \in \mathbf{D}(R)$ fits into a distinguished triangle
$$X_{\geq 0} \longrightarrow X \longrightarrow X_{\leq 0} \to X_{\geq 0}[1] \tag{A.14}$$

These classes naturally determine an *abelian* subcategory of $\mathbf{D}(R)$, the *heart* $\mathbf{D}(R)^\heartsuit$ of the $t$-structure.

In the following section we briefly sketch some of the basic classical definitions taken from [KS] and the classical [BBD82]; the $\infty$-categorical analogue of the theory has been defined in [Lur17, §**1.2.1**]. Here we merely recall a couple of definitions for the ease of the reader: from [Lur17, Def. **1.2.1.1** and **1.2.1.4**] one obtains the following translation of the definition of $t$-structure.

DEFINITION A.3.2. Let $\mathbf{C}$ be a stable $\infty$-category. A *$t$-structure* on $\mathbf{C}$ consists of a pair $\mathfrak{t} = (\mathbf{C}_{\geq 0}, \mathbf{C}_{<0})$ of full sub-$\infty$-categories satisfying the following properties:

(i) orthogonality: $\mathbf{C}(X, Y)$ is a contractible simplicial set for each $X \in \mathbf{C}_{\geq 0}$, $Y \in \mathbf{C}_{<0}$;
(ii) Setting $\mathbf{C}_{\geq 1} = \mathbf{C}_{\geq 0}[1]$ and $\mathbf{C}_{<-1} = \mathbf{C}_{<0}[-1]$ one has $\mathbf{C}_{\geq 1} \subseteq \mathbf{C}_{\geq 0}$ and $\mathbf{C}_{<-1} \subseteq \mathbf{C}_{<0}$;
(iii) Any object $X \in \mathbf{C}$ fits into a (homotopy) fiber sequence $X_{\geq 0} \to X \to X_{<0}$, with $X_{\geq 0}$ in $\mathbf{C}_{\geq 0}$ and $X_{<0}$ in $\mathbf{C}_{<0}$.

The subcategories $\mathbf{C}_{\geq 0}, \mathbf{C}_{<0}$ are called respectively the *coaisle* and the *aisle* of the $t$-structure (see [KV88]).

REMARK A.3.3. The definition as it is stated is a slight reformulation of the classical one given in [BBD82]; it is rather curious that the authors of the book do not give any reasonable rationale to explain what does the "t" stand for. A natural explanation is that it is a truncation (!) of the word "**t**runcation" (see [Hum] and the discussion therein).

REMARK A.3.4. The assignments $X \mapsto X_{\geq 0}$ and $X \mapsto X_{<0}$ define two functors $\tau_{\geq 0}$ and $\tau_{<0}$ which are, respectively, a right adjoint to the inclusion functor $\mathbf{C}_{\geq 0} \hookrightarrow \mathbf{C}$ and a left adjoint to the inclusion functor $\mathbf{C}_{<0} \hookrightarrow \mathbf{C}$. In other words, $\mathbf{C}_{\geq 0}, \mathbf{C}_{<0} \subseteq \mathbf{C}$ are respectively [Lur17, **1.2.1.5-8**] a coreflective and a reflective subcategory of $\mathbf{C}$.

This in particular implies that

- the full subcategories $\mathbf{C}_{\geq n} = \mathbf{C}_{\geq}[n]$, are coreflective via a coreflection $\tau_{\geq n}$; dually $\mathbf{C}_{<n} = \mathbf{C}_{<0}[n]$ are reflective via a reflection $\tau_{<n}$,
- $\mathbf{C}_{<n}$ is stable under all limits which exist in $\mathbf{C}$, and colimits are computed by applying the reflector $\tau_{<n}$ to the colimit computed in $\mathbf{C}$; dually, $\mathbf{C}_{\geq n}$ is stable under all colimits, and limits are $\mathbf{C}$-limits coreflected via $\tau_{\geq n}$; from the last of these remarks we deduce a useful corollary:

COROLLARY A.3.5. The functor $\tau_{<n}$ maps a pullout in $\mathbf{C}$ to a pushout in $\mathbf{C}_{<n}$ while $\tau_{\geq n}$ maps a pullout in $\mathbf{C}$ to a pullback in $\mathbf{C}_{\geq n}$.

NOTATION A.3.6. This is an important notational remark: the subcategory that we here denote $\mathbf{C}_{<0}$ is the subcategory which would be denoted $\mathbf{C}_{\leq 0}[-1]$ in [Lur17].

REMARK A.3.7. It's easy to see that Definition **A.3.2** is modeled on the classical definition of a $t$-structure ([KS], [BBD82]). In fact a $t$-structure $\mathfrak{t}$ on $\mathbf{C}$, following [Lur17], can also be characterized as a $t$-structure (in the classical sense) on the homotopy category of $\mathbf{C}$ ([Lur17, Def. **1.2.1.4**]), once $\mathbf{C}_{\geq 0}, \mathbf{C}_{<0}$ are identified with the subcategories of the homotopy category of $\mathbf{C}$ spanned by those objects which belong to the (classical) $t$-structure $\mathfrak{t}$ on the homotopy category.

Remark A.3.8. The datum of a $t$-structure via both classes $(\mathbf{C}_{\geq 0}, \mathbf{C}_{<0})$ is a bit redundant: in fact, each of the two classes uniquely determines the other via the object-orthogonality relation **1.2.14**.

Remark A.3.9. The notation $\mathbf{C}_{\geq 1}$ for $\mathbf{C}_{\geq 0}[1]$ is powerful but potentially misleading: one is led to view $\mathbf{C}_{\geq 0}$ as the seminfinite interval $[0, +\infty)$ in the real line and $\mathbf{C}_{\geq 1}$ as the seminfinite interval $[1, +\infty)$. This is indeed a very useful analogy (see Remark **4.3.3**) but one should always keep in mind that as a particular case of the inclusion condition $\mathbf{C}_{\geq 1} \subseteq \mathbf{C}_{\geq 0}$ also the extreme case $\mathbf{C}_{\geq 1} = \mathbf{C}_{\geq 0}$ is possible, in blatant contradiction of the real line half-intervals mental picture.

Definition A.3.10. ($t$-exact functor): Let $\mathbf{C}, \mathbf{D}$ be two stable $\infty$-categories, endowed with $t$-structures $\mathfrak{t}_{\mathbf{C}}, \mathfrak{t}_{\mathbf{D}}$; a functor $F\colon \mathbf{C} \to \mathbf{D}$ is *left $t$-exact* if it is exact and $F(\mathbf{C}_{\geq 0}) \subseteq \mathbf{D}_{\geq 0}$. It is called *right $t$-exact* if it is exact and $F(\mathbf{C}_{<0}) \subseteq \mathbf{D}_{<0}$.

Remark A.3.11. The collection $\text{TS}(\mathbf{C})$ of all $t$-structures on $\mathbf{C}$ has a natural partial order defined by $\mathfrak{t} \preceq \mathfrak{t}'$ iff $\mathbf{C}_{<0} \subseteq \mathbf{C}'_{<0}$. The ordered group $\mathbb{Z}$ acts **4.1.8** on $\text{TS}(\mathbf{C})$ with the generator **4.1.9** $+1$ mapping a $t$-structure $\mathfrak{t} = (\mathbf{C}_{\geq 0}, \mathbf{C}_{<0})$ to the $t$-structure $\mathfrak{t}[1] = (\mathbf{C}_{\geq 1}, \mathbf{C}_{<1})$. Since by **A.3.2**(ii) $\mathfrak{t} \preceq \mathfrak{t}[1]$, one sees that $\text{TS}(\mathbf{C})$ is naturally a $\mathbb{Z}$-poset.

In light of this remark, it is natural to consider *families* of $t$-structures with values in a generic $\mathbb{Z}$-poset $J$; this is discussed in our Ch. **4**.

Remark A.3.12. ($t$-structures are localizations): An alternative description for a $t$-structure is given in [Lur17, Prop. **1.2.1.16**] via a *$t$-localization* $L$, i.e. a reflection functor $L$ satisfying one of the following equivalent properties:

- The class of $L$-local morphisms[(2)] is generated (as a quasisaturated marking) by a family of initial arrows $\{0 \to X\}$;
- The class of $L$-local morphisms is generated (as a quasisaturated marking) by the class of initial arrows $\{0 \to X \mid LX \simeq 0\}$;
- The essential image $L\mathbf{C} \subset \mathbf{C}$ is an extension-closed class.

The $t$-structure $\mathfrak{t}(L)$ determined by the $t$-localization $L\colon \mathbf{C} \to \mathbf{C}$ is given by the pair of subcategories

$$\mathbf{C}_{\geq 0}(L) := \{A \mid LA \simeq 0\}, \qquad \mathbf{C}_{<0}(L) := \{B \mid LB \simeq B\}. \tag{A.15}$$

It is no surprise that the obvious example of $t$-localization is the truncation $\tau_{<0} : \mathbf{C} \to \mathbf{C}_{<0}$ associated with a $t$-structure $(\mathbf{C}_{\geq 0}, \mathbf{C}_{<0})$, and that one has $\mathbf{C}_{\geq 0}(\tau_{<0}) = \mathbf{C}_{\geq 0}$ and $\mathbf{C}_{<0}(\tau_{<0}) = \mathbf{C}_{<0}$.

This connection is precisely what motivated us to exploit the theory of factorization systems to give an alternative description of the data contained

(2) An arrow $f$ in $\mathbf{C}$ is called *$L$-local* if it is inverted by $L$; it's easy to see that $L$-local objects form a quasisaturated class in the sense of [Lur17, Def. **1.2.1.14**].

in a $t$-structure: the synergy between orthogonality encoded in **A.3.2.(i)** and reflectivity of the subcategories generated by $\mathfrak{t}$, suggests taking the "torsio-centric" approach.

# A.4 Spanier-Whitehead stabilization.

Let $\mathbf{A}$ be any category, endowed with an endofunctor $\Sigma\colon \mathbf{A} \to \mathbf{A}$. The problem adressed by the Spanier-Whitehead construction is the following: how to produce a category with endofunctor $(\mathrm{SW}(\mathbf{A}), \hat{\Sigma})$ such that

(1) there is an embedding $\mathbf{A} \hookrightarrow \mathrm{SW}(\mathbf{A})$;
(2) $\hat{\Sigma}|_{\mathbf{A}} = \Sigma$;
(3) $\hat{\Sigma}$ is an equivalence of categories

and such that the pair $(\mathrm{SW}(\mathbf{A}), \hat{\Sigma})$ is initial with these properties?

There are two ways to formalize the problem. We analyze them both, borrowing equally from Tierney's [Tie69] and [Del04]. The treatment of SW-stabilization given here motivates very well the meaning of [Lur17, **1.4.1**, **1.4.2**].

## A.4.1 Construction via monads.

Let $\mathbb{N}$ be the monoid of natural numbers, considered as a category: it has a monoidal product given by the sum operation, such that the unit object is zero. Since $\mathbb{N}$ is a monoid in ($\mathbf{Set} \subset$)$\mathbf{Cat}$, the functor $T_{\mathbb{N}} = (-) \times \mathbb{N}$ is a monad, and the category of $T_{\mathbb{N}}$-algebras can be described as the category whose objects are pairs $(\mathbf{A}, \Sigma\colon \mathbf{A} \to \mathbf{A})$; more explicitly, a $T_{\mathbb{N}}$-algebra is a pair $(\mathbf{A}, \Sigma)$ where $\mathbf{A}$ is a category, and $\Sigma\colon \mathbf{A} \to \mathbf{A}$ is a functor such that the diagrams

$$
\begin{array}{ccc}
\mathbf{A} \times \mathbf{1} & \xrightarrow{\mathbf{A}\times\eta} & \mathbf{A} \times \mathbb{N} \\
& \overset{\sim}{\searrow} & \downarrow{\scriptstyle \tilde{\Sigma}} \\
& & \mathbf{A}
\end{array}
\qquad
\begin{array}{ccc}
\mathbf{A} \times \mathbb{N} \times \mathbb{N} & \xrightarrow{\tilde{\Sigma}\times\mathbb{N}} & \mathbf{A} \times \mathbb{N} \\
{\scriptstyle \mathbf{A}\times\mu}\downarrow & & \downarrow{\scriptstyle \tilde{\Sigma}} \\
\mathbf{A} \times \mathbb{N} & \xrightarrow{\tilde{\Sigma}} & \mathbf{A}
\end{array}
\tag{A.16}
$$

(where $\tilde{\Sigma}(A, n) = \Sigma^n A$ and $\eta, \mu$ are the monoid maps of $\mathbb{N}$) all commute.

NOTATION A.4.1. In the following, $T_{\mathbb{N}}$-algebras will be called *categories with endomorphism.*

Let now $\mathbb{N} \hookrightarrow \mathbb{Z}$ the obvious inclusion. When regarded as a category, the group of integers is a groupoid, so $S_{\mathbb{Z}} = (-) \times \mathbb{Z}$ is again a monad on $\mathbf{Cat}$.

The category of $S_{\mathbb{Z}}$-algebras consists of pairs $(\mathbf{A}, \Sigma)$ where $\Sigma\colon \mathbf{A} \to \mathbf{A}$ is an *auto*morphism (so in particular every $S_{\mathbb{Z}}$-algebra is a $T_{\mathbb{N}}$-algebra). Similar diagrams are requested to commute, so that if we consider the restriction $\Sigma_{(n)} = \Sigma|_{\mathbf{A}\times\{n\}}$ for any $n \in \mathbb{Z}$, and we identify $\mathbf{A} \times \{n\} \cong \mathbf{A}$, then we have that $\Sigma_{(1)} = \Sigma$, $\Sigma_{(-1)} = \Sigma^{-1}$ and so on.

Remark A.4.2. The homomorphism $\iota\colon \mathbb{N} \hookrightarrow \mathbb{Z}$ induces a morphism of monads $T \to S$, which we call again $\iota$; this in turns induces a "forgetful" functor

$$U\colon \mathbf{Cat}^{\mathbb{Z}} \hookrightarrow \mathbf{Cat}^{\mathbb{N}} \tag{A.17}$$

(the forgetful action of $U$ is clear when its action is explicited: it simply forgets that an automorphism $\Sigma$ of $\mathbf{A}$ has an inverse.

We want to give a left adjoint $F\colon \mathbf{Cat}^{\mathbb{N}} \to \mathbf{Cat}^{\mathbb{Z}}$ to the functor $U$, obtaining a precise description of its action on objects of $\mathbf{Cat}$. To this end, given $(A, \Sigma) \in \mathbf{Cat}^{\mathbb{N}}$ let us consider the coequalizer diagram in $\mathbf{Cat}$:

$$\mathbf{A} \times \mathbb{N} \times \mathbb{Z} \underset{(\mathbf{A}\times\mu)\circ(\iota\times\mathbb{Z})}{\overset{\Sigma\times\mathbb{Z}}{\rightrightarrows}} \mathbf{A} \times \mathbb{Z} \longrightarrow F(\mathbf{A}, \Sigma) \tag{A.18}$$

Now, all monads like $S_{\mathbb{Z}}$, i.e. all monads of the form $(-)\times M$ for $M$ a monoid in a monoidal(ly cocomplete) category $(\mathbf{A}, \times)$ preserve colimits, hence there is a unique $S_{\mathbb{Z}}$-algebra structure on $F(\mathbf{A}, \Sigma)$ such that

$$\hat{\Sigma}\colon F(\mathbf{A}, \Sigma) \times \mathbb{Z} \to F(\mathbf{A}, \Sigma) \tag{A.19}$$

is an automorphism of $F(\mathbf{A}, \Sigma)$ and the correspondence $\tilde{F}\colon (\mathbf{A}, \Sigma) \mapsto (F(\mathbf{A}, \Sigma), \hat{\Sigma})$ is the desired left adjoint. The category $F(\mathbf{A}, \Sigma)$ can be considered the *free category with automorphism* on the category with endomorphism $(\mathbf{A}, \Sigma)$.

This category satisfies the desired universal property: there exists a functor

$$\alpha\colon (\mathbf{A}, \Sigma) \to \tilde{F}(\mathbf{A}, \Sigma) \tag{A.20}$$

(the unit of the adjunction we built) such that for any $S_{\mathbb{Z}}$-algebra morphism $H\colon (\mathbf{A}, \Sigma) \to (\mathbf{B}, \Theta)$ where $(\mathbf{B}, \Theta)$ is a $T$-algebra, there is a unique $T_{\mathbb{N}}$-algebra morphism $\bar{H}\colon \tilde{F}(\mathbf{A}, \Sigma) \to (\mathbf{B}, \Theta)$ such that the following diagram commutes:

$$\begin{array}{ccc} (\mathbf{A}, \Sigma) & \longrightarrow & \tilde{F}(\mathbf{A}, \Sigma) \\ & \searrow \quad \nearrow & \\ & (\mathbf{B}, \Theta). & \end{array} \tag{A.21}$$

The category of *topological spectra* consists of the Spanier-Whitehead stabilization of the category of CW-complexes, as well as the category of chain complexes of abelian groups (or modules over a ring $R$); these examples are discussed in [Tie69].

## A.5 Stability in different models.

Whirl in circles
Around a stable center.

M. Ueshiba

### A.5.1 Stable Model categories.

Every pointed model category [Hov99, Ch. **7**] $\mathbf{M}$ carries an adjunction between endofunctors

$$\Sigma \dashv \Omega \colon \mathbf{M} \leftrightarrows \mathbf{M} \tag{A.22}$$

defined respectively as the homotopy pushout and homotopy pullback below:

$$\begin{array}{ccc} X & \longrightarrow & * \\ \downarrow & & \downarrow \\ * & \longrightarrow & \Sigma X \end{array} \qquad \begin{array}{ccc} \Omega Y & \longrightarrow & * \\ \downarrow & & \downarrow \\ * & \longrightarrow & Y. \end{array}$$

It is a matter of unraveling definition to show that these two functors are mutually adjoint. A pointed model category is said to be *stable* if the above adjunction is a Quillen equivalence.

#### A.5.1.1 $k$-linear DG-categories.

A *$k$-linear* DG*-category* is a category enriched ([Kel82, Gen15]) over the category of chain complexes of vector spaces over the field $k$; a DG-category $\mathbb{D}$ is called *pretriangulated* if the following two axioms hold:

- For every object $X \in \mathbb{D}$ the shifted representable DG-module $\mathbb{D}(-,X)[k] \in \widehat{\mathbb{D}}$ is homotopic to a representable $\mathbb{D}(-,X\langle k\rangle)$;
- For every $f\colon \mathbb{D}(-,X) \to \mathbb{D}(-,Y)$ a morphism of representable DG-modules in $\widehat{\mathbb{D}}$, the DG-module

$$\underline{C}\big(\mathbb{D}(-,f)\big)\colon \underline{C}(\mathbb{D}(-,X)) \to \underline{C}(\mathbb{D}(-,Y)) \tag{A.23}$$

  is homotopic to a representable $\mathbb{D}(-,c(f))$.

The homotopy category of a DG-category is defined by taking the $H^0$ of each hom-space $\mathbb{D}(X,Y)$ (or, more formally, the image of $\mathbb{D}$ under the 2-functor $H_{0,*}\colon$ DG-**Cat** $\to$ **Cat**). The homotopy category of a pretriangulated DG-category is triangulated, in the sense of definition **A.1.1**. We define an *enhancement* for a triangulated category $\mathbf{D}$ to be a pretriangulated DG-category $\mathbb{D}$ such that there is an equivalence $[\mathbb{D}] \cong \mathbf{D}$.

Quoting [Lur17]:

The theory of differential graded categories is closely related to the theory of stable $\infty$-categories. More precisely, one can show that the data of a (pretriangulated) differential graded category over a field $k$ is equivalent to the data of a stable $\infty$-category $\mathbf{C}$ equipped with an enrichment over the monoidal $\infty$-category of $k$-module spectra. The theory of differential graded categories provides a convenient language for working with stable $\infty$-categories of algebraic origin (for example, those which arise from chain complexes of coherent sheaves on algebraic varieties), but is inadequate for treating examples which arise in stable homotopy theory.

#### A.5.1.2 Stable $\infty$-categories.

Stable $\infty$-categories are extensively described in [Lur17], throughout the present chapter, and throughout the present thesis; here, we outline how, in the setting of $\infty$-categories, the lack of universality for the construction of $\mathbf{D}(\mathbf{A})$ is completely solved: first of all, recall that the *Dold-Kan correspondence* [Kan58, GJ99] establishes an equivalence of categories between the category $\mathrm{Ch}^+(\mathbf{Ab})$ (chain complexes of abelian groups, concentrated in positive degree) and $\mathbf{sAb}$ (simplicial sets whose sets of $n$-simplices all are abelian groups).

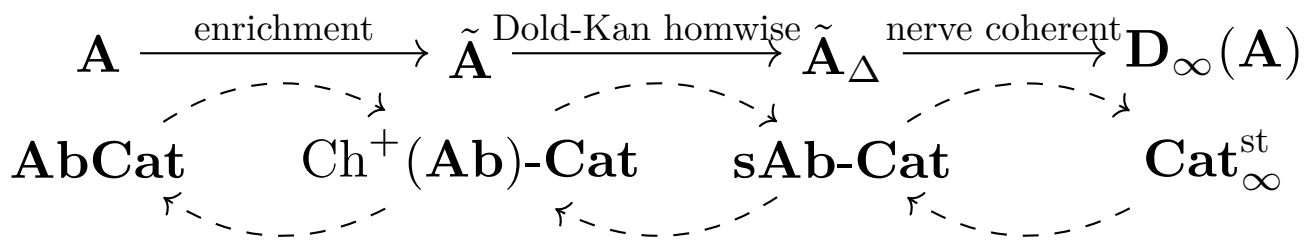

Figure A.1: Construction of the derived $\infty$-category of $\mathbf{A}$.

# Index